\documentclass[11pt]{amsart}
\usepackage{amsmath,amssymb,amscd}
\usepackage[mathscr]{eucal}
\usepackage{xypic}
\input xy
\xyoption{all}

\theoremstyle{plain}
\newtheorem{thm}{Theorem}[section]
\newtheorem{prop}[thm]{Proposition}
\newtheorem{lem}[thm]{Lemma}
\newtheorem{cor}[thm]{Corollary}
\newtheorem{defn}[thm]{Definition}

\newcommand{\A}{\mathbb A}
\newcommand{\C}{\mathbb C}
\renewcommand{\P}{\mathbb P}
\newcommand{\Q}{\mathbb Q}
\newcommand{\R}{\mathbb R}
\newcommand{\Z}{\mathbb Z}

\newcommand{\scC}{\mathscr C}
\newcommand{\scD}{\mathscr D}
\newcommand{\scH}{\mathscr H}

\newcommand{\scZ}{\mathscr Z}

\newcommand{\caC}{\mathcal C}
\newcommand{\caD}{\mathcal D}
\newcommand{\caE}{\mathcal E}
\newcommand{\caF}{\mathcal F}
\newcommand{\caG}{\mathcal G}
\newcommand{\caH}{\mathcal H}
\newcommand{\caO}{\mathcal O}
\newcommand{\caP}{\mathcal P}

\newcommand{\frA}{\mathfrak A}
\newcommand{\frB}{\mathfrak B}

\newcommand{\frP}{\mathfrak P}
\newcommand{\frS}{\mathfrak S}

\newcommand{\IIm}{\operatorname{Im}}
\newcommand{\CH}{\operatorname{ch}}

\newcommand{\IId}{\operatorname{Id}}
\newcommand{\KER}{\operatorname{Ker}}

\newcommand{\CPT}{\operatorname{cpt}}

\newcommand{\TOT}{\operatorname{Tot}}
\newcommand{\TTR}{\operatorname{tr}}

\newcommand{\PRJ}{\operatorname{pr}}
\newcommand{\sgn}{\operatorname{sgn}}
\newcommand{\alt}{\operatorname{alt}}
\newcommand{\Alt}{\operatorname{Alt}}

\newcommand{\sbinom}[2]{{\textstyle \binom{#1}{#2}}}
\newcommand{\ssum}{\textstyle \sum}
\newcommand{\PR}[1]{#1^{\prime }}
\newcommand{\PPRR}[1]{#1^{\prime \prime }}

\newcommand{\TRI}{\vartriangle }
\newcommand{\PRA}{\bullet_{\mathbb A}}
\newcommand{\PRP}{\bullet_{\mathbb P}}
\newcommand{\PRAP}{\bullet_{\mathbb A,\mathbb P}}
\newcommand{\EPRAP}{\ast_{\mathbb A,\mathbb P}}
\newcommand{\EPRA}{\ast_{\mathbb A}}
\newcommand{\ahdif}{\overline{\partial }}
\newcommand{\PP}[1]{\square^{#1}}
\newcommand{\Sps}[1]{\left<#1\right>}
\newcommand{\CST}[2]{\frac{#1}{(2\pi i)^{#2}}}
\newcommand{\OV}[1]{\overline{#1}}
\newcommand{\QQ}[1]{#1_{\mathbb Q}}

\newcommand{\mmD}{\tau \mathscr D}

\newcommand{\WHmmD}{\tau \widehat{\mathscr D}}

\newcommand{\WH}[1]{\widehat{#1}}
\newcommand{\WT}[1]{\widetilde{#1}}

\newcommand{\TRNS}[2]{\operatorname{tr}_{#1}(\lambda #2)}
\newcommand{\ALT}[2]{C^{\Alt}_{#1}(#2)}
\newcommand{\WALT}[2]{\WH{C}^{\Alt}_{#1}(#2)}

\title{Higher arithmetic Chern character}
\author{Yuichiro Takeda}
\address{Faculty of Mathematics, Kyushu University, Motooka 744, 
Nishi-ku, Fukuoka 819-0395, Japan}
\email{yutakeda@math.kyushu-u.ac.jp}

\thanks{JSPS Grant-in-Aid for Scientific Research (No.21540019)}

\begin{document}
\maketitle

\begin{abstract}
A map from the higher arithmetic $K$-group defined by the author 
to the higher arithmetic Chow group constructed by Burgos and Feliu 
is given.
It is a higher extension of the arithmetic Chern character 
established by Gillet and Soul\'{e}, and it can also be regarded as 
an analogue of Beilinson's regulator in Arakelov geometry.
It is shown that this map is compatible with the pull-back 
map and the product structure.
\end{abstract}

\vspace{0.5cm}

\setcounter{tocdepth}{1}
\tableofcontents

\section{Introduction}
\vskip 1pc

In a celebrated paper \cite{gilletsoule2}, Gillet and 
Soul\'{e} established a theory of arithmetic Chern character 
associated with a hermitian vector bundle in Arakelov geometry.
This theory is very useful in generating elements of arithmetic 
Chow groups, and a lot of interesting equalities of numbers 
have been deduced by calculating the arithmetic intersection 
of these elements.

The arithmetic Chern characters of hermitian vector bundles on 
an arithmetic variety $X$ provide a map from the arithmetic 
$K_0$-group to the arithmetic Chow group 
\begin{equation}
\WH{\CH}_0^p:\WH{K}_0(X)\to \WH{CH}^p(X). \label{sec1:map1}
\end{equation}
Recently, the both sides of this map were extended to higher degree.
For instance, in \cite{takeda} the author gave a definition 
of higher arithmetic $K$-groups $\WH{K}_r(X)$.
On the other hand, in \cite{goncharov1} and \cite{BF} 
two distinct definitions of higher arithmetic Chow groups 
$\WH{CH}^p(X, r)$ were proposed when $X$ is defined 
over an arithmetic field.
The both definitions are based on descriptions of the regulator 
map 
$$
CH^p(X, r)\to H_{\scD}^{2p-r}(X, \R(p))
$$
by means of integral on algebraic cycles, and 
in \cite{BFT} it was shown that these two definitions 
are essentially the same.

The aim of this paper is to construct a map from 
$\WH{K}_r(X)$ to $\WH{CH}^p(X, r)$ called 
the {\it higher arithmetic Chern character}.
It can be seen as an extension of the map \eqref{sec1:map1} 
to higher degree, and also as an analogue in Arakelov geometry of 
the higher Chern character $\CH_r^p:K_r(X)\to CH^p(X, r)$.

In order to define the higher Chern character map we usually 
make use of homotopy theory.
However, it seems difficult to redefine the higher arithmetic 
Chow groups $\WH{CH}^p(X, r)$ in 
the context of homotopy theory.
Hence we adopt an alternative approach.

In \cite{levine} Levine introduced a scheme called 
{\it iterated double} in order to prove several properties of 
the higher Chow group, such as localization sequence and 
contravariant functoriality, without relying 
on the moving lemma.
In doing this, he showed that the higher $K$-group of a scheme 
can be embedded into the $K_0$-group of the associated 
iterated double.
Let us explain this in detail.

Let $X$ be a regular scheme, and $\square =\P^1-\{1\}$ 
the affine line with a subscheme $\partial \square =\{0, \infty \}$.
Let $\PP{r}$ be the cartesian product of $r$ copies of $\square $.
Then the subscheme $\partial \square $ on each component 
gives a normal crossing divisor $\partial \PP{r}\subset \PP{r}$.
As for the precise definition of $\partial \PP{r}$, see \S 3.6.
It is well-known that the $r$-th $K$-group $K_r(X)$ is isomorphic 
to the multi-relative $K_0$-group of the scheme $X\times \PP{r}$ 
with the normal crossing divisor $X\times \partial \PP{r}$.
Denote by $T=D(X\times \PP{r}; X\times \partial \PP{r})$ 
the associated iterated double.
Then we obtain closed subschemes $T_1, \ldots , T_r\subset T$ 
with a closed embedding 
$$
i_{\emptyset }:(X\times \PP{r}; X\times \partial \PP{r})\hookrightarrow 
(T; T_1, \ldots , T_r).
$$
Levine showed in \cite{levine} that the pull-back map of 
multi-relative $K$-groups 
$$
i_{\emptyset }^*:K_*(T; T_1, \ldots , T_r)\to 
K_*(X\times \PP{r}; X\times \partial \PP{r})
$$
is bijective.
Hence we have the following sequence of maps:
\begin{equation}
K_r(X)\simeq K_0(X\times \PP{r}; X\times \partial \PP{r})
\overset{i_{\emptyset }^*}{\overset{\sim }{\leftarrow }}
K_0(T; T_1, \ldots , T_r)\subset K_0(T). \label{sec1:map2}
\end{equation}
This means that any element of $K_r(X)$ can be represented by 
a virtual vector bundle on the iterated double $T$.

In this paper we pursue an analogue of \eqref{sec1:map2} 
in Arakelov geometry.
To accomplish this purpose, we define arithmetic analogues of 
$K$-groups appearing in \eqref{sec1:map2}.
Then we obtain the following sequence:
\begin{equation}
\QQ{\WH{K}_r(X)}\simeq 
\QQ{\WH{K}_0(X\times \PP{r}; X\times \partial \PP{r})}
\overset{\WH{i}_{\emptyset }^*}{\twoheadleftarrow }
\QQ{\WH{K}_0(T; T_1, \ldots , T_r)}
\subset \QQ{\WH{K}_0(T)},  \label{sec1:map3}
\end{equation}
where $\QQ{A}$ stands for $A\otimes \Q$.
Note that in the arithmetic case this sequence holds only in rational 
coefficients, and the middle arrow $\WH{i}_{\emptyset }^*$ 
is shown to be surjective.
We do not know whether $\WH{i}_{\emptyset }^*$ 
is bijective or not, nevertheless we know how any element 
of $\KER \WH{i}_{\emptyset }^*$ is described.
We next establish a theory of arithmetic Chern character 
of a hermitian vector bundle on the iterated double $T$.
This provides a map 
\begin{equation}
\WH{\CH}^p_{T,0}:\WH{K}_0(T)\to \WH{CH}^p(X, r). \label{sec1:map4}
\end{equation}
Combining \eqref{sec1:map3} with \eqref{sec1:map4}, we can 
obtain the desired map 
\begin{equation}
\WH{\CH}_r^p:\QQ{\WH{K}_r(X)}\to \WH{CH}^p(X, r). \label{sec1:map5}
\end{equation}

If we ignore the metric structure in the definition of 
\eqref{sec1:map5}, we can define 
$$
\CH_r^p:\QQ{K_r(X)}\to H_{\scD}^{2p-r-1}(X, \R(p))
$$
which make the following diagram commutative:
$$
\begin{CD}
\QQ{\WH{K}_r(X)} @>>> \QQ{K_r(X)} \\
@VV{\WH{\CH}_r^p}V @VV{\CH_r^p}V \\
\WH{CH}^p(X, r) @>>> CH^p(X, r).
\end{CD}
$$
However, in this paper we could not prove that 
the map $\CH_r^p$ agrees with the higher Chern character.
Nevertheless we can show that the composite 
with the regulator map 
$$
\QQ{K_r(X)}\overset{\CH_r^p}{\to }CH^p(X, r)\to H^{2p-r}(X, \R(p))
$$
agrees with Beilinson's regulator.
This is an evidence that the map $\CH_r^p$ is 
nothing but the higher Chern character.

We describe the content of each section.

In \S 2, we first recall a chain complex of exact cubes 
which computes the rational $K$-theory \cite{mccarthy}.
Then we introduce a subcomplex 
on which the symmetric group acts alternatingly, and show that 
these two chain complexes are quasi-isomorphic.
We next introduce $\scC$-complex \cite{hanamura}, and we employ 
it to define a chain complex called {\it multi-relative complex 
of exact cubes}, which computes the multi-relative rational $K$-theory.
We examine functorial properties of the multi-relative complex 
of exact cubes.

In \S 3, we introduce several complexes of differential 
forms and of currents which compute the Deligne cohomology.
We then introduce Wang's forms and 
the higher Bott-Chern forms \cite{BW}.
We also give a map from the multi-relative complex of exact cubes 
introduced in \S 3 to a complex of differential forms, which can 
be seen as a generalization of the higher Bott-Chern forms 
to the relative case.

In \S 4, we introduce iterated double, 
and establish a theory of Chern forms of hermitian vector bundles 
on it.
In \S 5, we define relative arithmetic $K$-theory, and 
arithmetic $K_0$-group of an iterated double.
We then deduce the sequence \eqref{sec1:map3}.

In \S 6, after recalling the definition of the higher 
arithmetic Chow groups $\WH{CH}^p(X, r)$ given by Burgos nad Feliu 
\cite{BF}, we give a definition of arithmetic Chern character 
of hermitian vector bundles on an iterated double.
Then we show that it induces the map \eqref{sec1:map4}.
In \S 7, we deduce the desired map \eqref{sec1:map5}, 
and in \S 8 we show the compatibility 
of this map with the pull-back map.

The last three sections are devoted to showing the compatibility 
of \eqref{sec1:map5} with product with 
the arithmetic $K_0$-group and the arithmetic Chow group.
In \S 9, we put a tensor product structure on 
the multi-relative complexes of exact cubes.
In \S 10, we show that all the arithmetic $K$-groups defined in \S 5 
admit products with $\WH{K}_0(X)$, and the sequence \eqref{sec1:map3} 
respects the $\WH{K}_0(X)$-module structures.
Finally in \S 11, we show that the arithmetic Chern character 
of a hermitian vector bundle on an iterated double has a 
multiplicative property with the tensor product.
Collecting these results, we can show that \eqref{sec1:map5} 
is compatible with the product structures.

{\it Acknowledgment}: 
During the preparation of this paper, the author visited 
Centre de Recerca Matem\'{a}tica (CRM) in Bellaterra from March 
to April of 2010, 
and Universitat de Barcelona several times.
The author would like to thank all the staff there for their 
hospitality and support.
In particular, he would like to thank sincerely 
Jose Ignacio Burgos Gil for his hospitality and fruitful discussions.

\vskip 1pc
\noindent
{\bf Notations and conventions.} \ 
Throughout this paper, we fix an universe and assume that 
all the sets we will consider are contained in this fixed universe.
Hence the category of vector bundles on a scheme is a small 
exact category.

Any small exact category is assumed to have a distinguished 
zero object denoted  by $0$, and any exact functor between 
exact categories is assumed to preserve the distinguished 
zero object.
In particular, we fix a distinguished zero object of 
the category of abelian groups.

For a scheme $X$, we denote by $\frP(X)$ the category of 
locally free $\caO_X$-modules of finite rank on $X$.
We identify locally free $\caO_X$-modules of finite rank with 
vector bundles of finite rank on $X$ in the usual way.
Then $\frP(X)$ is a small exact category, and a distinguished 
zero object of $\frP(X)$ is given as follows: For any open set 
$U\subset X$, $0(U)$ is the distinguished zero object 
of the category of abelian groups.

For any morphism $f:X\to Y$ of schemes, there is a pull-back 
functor $f^*:\frP(Y)\to \frP(X)$, which is an exact functor 
preserving the distinguished zero object.
When $f$ is the identity $\IId_X:X\to X$, we identify 
$\IId_X^*$ with the identity functor of $\frP(X)$.

We fix some conventions on complexes.
With a complex of abelian groups $(A^*, d_A)$, we can associate 
a chain complex $(A_*, \partial_A)$ in the way that $A_*=A^{-*}$ 
and $\partial_A=d_A$.
In this paper we always identify a complex and the associated chain 
complex in this way.
For a chain complex of abelian groups $A_*=(A_n, \partial_A)$, 
let $A[r]_*$ be a chain complex such that $A[r]_n=A_{n-r}$ with 
$\partial_{A[r]}=(-1)^r\partial_A:A_{n-r}\to A_{n-r-1}$.
Similarly, for a complex $A^*=(A^*, d_A)$, let $A[r]^*$ be a complex 
such that $A[r]^n=A^{n+r}$ with $d_{A[r]}=(-1)^rd_A$.
We define the {\it simple complex} $s(f)$ of a map $f:A_*\to B_*$ 
of chain complexes to be 
$$
s(f)_n=A_n\oplus B_{n+1}
$$
with the boundary map 
$$
\partial :s(f)_n\to s(f)_{n-1}, \quad 
\partial (a, b)=(\partial a, f(a)-\partial b).
$$
Then we have a long exact sequence on homology 
$$
\cdots \to H_{n+1}(B_*)\to H_n(s(f)_*)\to H_n(A_*)\to H_n(B_*)
\to \cdots . 
$$
Throughout this paper, we denote by $[a]\in H_n(A_*)$ 
the homology class represented by $a\in \KER \partial $.
Moreover, we denote by $\WT{a}\in A_n/\IIm \partial $ 
the element represented by $a\in A_n$.

Finally, we define the {\it truncated relative cohomology groups} 
of a map of complexes $A^*\to B^*$ to be 
$$
\WH{H}^n(A^*, B^*)=\{(a, \WT{b})\in 
A^n\oplus B^{n-1}/\IIm d; \ da=0, f(a)=db\}
$$
\cite[Def.4.2]{burgos2}.
We are going to use these groups to define 
Green objects associated with an algebraic cycle in \S 6.

\setcounter{equation}{0}
\vskip 2pc
\section{Multi-relative complexes of exact cubes}
\vskip 1pc

\subsection{Complexes of exact cubes}
Let $\frA$ be a small exact category.
We see the ordered set $\{-1, 0, 1\}$ as a category.
For $n\geq 0$, an $n$-{\it cube} of $\frA$ is a covariant functor 
from the product $\{-1, 0, 1\}^{\times n}$ to $\frA$.
Note that a $0$-cube of $\frA$ is an object of $\frA$.
For an $n$-cube $\caF$ of $\frA$, denote by 
$\caF_{\alpha_1, \ldots , \alpha_n}$ the image of the 
object $(\alpha_1, \ldots , \alpha_n)\in \{-1, 0, 1\}^{\times n}$ 
by $\caF$.
For $1\leq j\leq n$ and $-1\leq i\leq 1$, a {\it face} 
of $\caF$ is an $(n-1)$-cube $\partial_j^i\caF$ defined by 
$$
(\partial_j^i\caF)_{\alpha_1, \ldots , \alpha_{n-1}}=
\caF_{\alpha_1, \ldots , \alpha_{j-1}, i, \alpha_j, 
\ldots , \alpha_{n-1}}.
$$
For $\alpha =(\alpha_1, \ldots , \alpha_{n-1})\in 
\{-1, 0, 1\}^{n-1}$ and for an integer $1\leq j\leq n$, 
an {\it edge} of $\caF $ is an $1$-cube 
$\partial_{j^c}^{\alpha }\caF $ described as 
$$
\caF_{\alpha_1, \ldots , \alpha_{j-1}, -1, \alpha_j, \ldots , 
\alpha_{n-1}}\to \caF_{\alpha_1, \ldots , \alpha_{j-1}, 0, \alpha_j, 
\ldots , \alpha_{n-1}}\to \caF_{\alpha_1, \ldots , \alpha_{j-1}, 1, 
\alpha_j, \ldots , \alpha_{n-1}}.
$$
An $n$-cube $\caF $ of $\frA $ is said to be {\it exact} if all 
the edges of $\caF $ are short exact sequences.

For an exact $(n-1)$-cube $\caF$ of $\frA$ and for an integer 
$1\leq j\leq n$, let $s_j^1\caF$ be an exact $n$-cube 
such that the edge $\partial_{j^c}^{\alpha }(s_j^1\caF)$ is 
$\caF_{\alpha }\overset{\IId }{\to }\caF_{\alpha }\to 0$
for any $\alpha \in \{-1, 0, 1\}^{n-1}$.
Similarly, let $s_j^{-1}\caF$ be an exact $n$-cube such that 
$\partial_{j^c}^{\alpha }(s_j^{-1}\caF)$ is 
$0\to \caF_{\alpha }\overset{\IId }{\to }\caF_{\alpha }$
for any $\alpha \in \{-1, 0, 1\}^{n-1}$.
These exact $n$-cubes are said to be {\it degenerate}.

Denote by ${\Q}C_n(\frA)$ the $\Q$-vector space generated 
by all exact $n$-cubes of $\frA$.
The alternating sum of faces gives a map of $\Q$-vector spaces 
$$
\partial =\sum_{j=1}^r\sum_{i=-1}^1(-1)^{i+j}\partial_j^i:
{\Q}C_r(\frA)\to {\Q}C_{r-1}(\frA), 
$$
by which $({\Q}C_*(\frA), \partial )$ is a chain complex.
Denote by $D_n$ the subgroup of $\Q C_n(\frA)$ generated by 
all degenerate $n$-cubes.
Then $D_*$ is a subcomplex of $\Q C_*(\frA)$.

\vskip 1pc
\begin{thm}\cite{mccarthy}
The homology group of the quotient complex  
$$
\WT{\Q}C_*(\frA)=\Q C_*(\frA)/D_*
$$
is canonically isomorphic to the rational algebraic $K$-theory 
of $\frA$:
$$
H_n(\WT{\Q}C_*(\frA), \partial )\simeq \QQ{K_n(\frA)}.
$$
\end{thm}

\vskip 1pc
\subsection{The subcomplex $\WT{\Q}\ALT{*}{\frA}$}
Let $\frA$ be a small exact category, and $\frS_n$ 
the $n$-th symmetric group.
For an exact $n$-cube $\caF$ of $\frA$ and for 
$\sigma \in \frS_n$, let $\sigma \caF$ be an 
exact $n$-cube such that 
$(\sigma \caF)_{\alpha_1, \ldots , \alpha_n}=
\caF_{\alpha_{\sigma(1)}, \ldots , 
\alpha_{\sigma(n)}}$.
Then $\frS_n$ acts on $\Q C_n(\frA)$, and 
we can obtain the alternating part 
$$
\Q\ALT{n}{\frA}=\{x\in \Q C_n(\frA); \ 
\sigma (x)=(\sgn \sigma )x  \ \text{for any} \ 
\sigma \in \frS_n\}
$$
with the section $\Alt_n:\Q C_n(\frA)\to \Q\ALT{n}{\frA}$ defined by 
\begin{equation}
\Alt_n(x)=\frac{1}{n!}\sum_{\sigma \in \frS_n}
(\sgn \sigma )\sigma (x). \label{sec2:map1}
\end{equation}
We can show in the same way as \cite[Lem.1.1]{bloch2} 
that $\Alt_{n-1}\partial =\partial \Alt_n$, which implies that 
$\Q\ALT{*}{\frA}$ is a subcomplex of $\Q C_*(\frA)$ and that 
$\Alt_*$ is a map of complexes.

Set 
$$
\WT{\Q}\ALT{n}{\frA}=\Q\ALT{n}{\frA}/(\Q\ALT{n}{\frA}\cap D_n).
$$
Then $\WT{\Q}\ALT{*}{\frA}$ is a subcomplex of $\WT{\Q}C_*(\frA)$.
Since $\sigma \caF $ is degenerate if so is $\caF $, 
the map \eqref{sec2:map1} gives rise to a map of complexes 
$$
\Alt_*:\WT{\Q}C_*(\frA)\to \WT{\Q}\ALT{*}{\frA}, 
$$
which is also a section of the inclusion 
$\WT{\Q}\ALT{*}{\frA}\subset \WT{\Q}C_*(\frA )$.

\vskip 1pc
\begin{thm}
The inclusion 
$\WT{\Q}\ALT{*}{\frA}\hookrightarrow \WT{\Q}C_*(\frA)$ 
is a quasi-isomorphism.
Hence we have a canonical isomorphism 
$$
H_n(\WT{\Q}\ALT{*}{\frA})\simeq 
\QQ{K_n(\frA)}.
$$
\end{thm}

\vskip 1pc
To prove this theorem, we need some preparation.
For an exact $1$-cube 
$\caF:A\overset{f}{\to }B\overset{g}{\to }C$ of $\frA $, 
consider the following exact $2$-cube:
$$
\rho (\caF)=\qquad \begin{CD}
A @>{\IId }>> A  @.  \\
@V{\IId }VV @V{f}VV @.  \\
A @>{f}>> B @>{g}>> C  \\
@. @V{g}VV @V{\IId }VV  \\
@. C @>{\IId }>> C.
\end{CD}
$$
Suppose $n\geq 1$.
For an exact $n$-cube $\caF$ of $\frA$ and for $1\leq j\leq n$, 
let $\rho_{n,j}(\caF)$ be an exact $(n+1)$-cube such that 
$$
\partial_1^{\alpha_1}\cdots \partial_{j-1}^{\alpha_{j-1}}
\partial_{j+2}^{\alpha_{j+1}}\cdots \partial_{n+1}^{\alpha_n}
\rho_{n,j}(\caF)=\rho \left(\partial_1^{\alpha_1}\cdots 
\partial_{j-1}^{\alpha_{j-1}}\partial_{j+1}^{\alpha_{j+1}}
\cdots \partial_n^{\alpha_n}\caF\right).
$$
It is obvious that $\rho_{n,j}(\caF)$ is degenerate 
if so is $\caF$.
The lemma below follows immediately from the definition:

\vskip 1pc
\begin{lem}
\begin{align*}
\partial_j^0\rho_{n,j}(\caF)&=\caF=\partial_{j+1}^0
\rho_{n,j}(\caF), \\
\partial_j^{-1}\rho_{n,j}(\caF)&=s_j^1\partial_j^{-1}
\caF=\partial_{j+1}^{-1}\rho_{n,j}(\caF), \\
\partial_j^1\rho_{n,j}(\caF)&=s_j^{-1}\partial_j^1\caF=
\partial_{j+1}^1\rho_{n,j}(\caF), 
\end{align*}
and 
$$
\partial_k^i\rho_{n,j}(\caF)=\begin{cases}
\rho_{n-1,j-1}(\partial_k^i\caF), \ & k<j,  \\
\rho_{n-1,j}(\partial_{k-1}^i\caF), \ & k>j+1.
\end{cases}
$$
Moreover, 
$$
\rho_{n+1,j+1}(\rho_{n,j}(\caF))=\rho_{n+1,j}(\rho_{n,j}(\caF)).
$$
\end{lem}

\vskip 1pc
Let $\WT{\Q}C_{n,m}(\frA)=\WT{\Q}C_{n+m}(\frA)$ with 
boundary maps 
\begin{align*}
\PR{\partial }=\sum_{j=1}^n\sum_{i=-1}^1(-1)^{i+j}
\partial_j^i&:\WT{\Q}C_{n,m}(\frA)
\to \WT{\Q}C_{n-1,m}(\frA), \\
\PPRR{\partial }=\sum_{j=1}^m\sum_{i=-1}^1(-1)^{i+j}
\partial_{n+j}^i&:\WT{\Q}C_{n,m}(\frA)\to 
\WT{\Q}C_{n,m-1}(\frA).
\end{align*}
Then $(\WT{\Q}C_{n,m}(\frA), \PR{\partial }, \PPRR{\partial })$ 
is a double chain complex.

\vskip 1pc
\begin{lem}
If $m\geq 1$, then 
$(\WT{\Q}C_{*,m}(\frA), \PR{\partial })$ is acyclic, and 
if $n\geq 1$, then 
$(\WT{\Q}C_{n,*}(\frA), \PPRR{\partial })$ is acyclic.
\end{lem}

{\it Proof}: 
For $m\geq 1$, define a map 
$$
\varPhi_{n,m}:\WT{\Q}C_{n,m}(\frA)\to \WT{\Q}C_{n+1,m}(\frA)
$$
by $\varPhi_{n,m}(\caF)=(-1)^{n+1}\rho_{n+m,n+1}(\caF)$ 
for an exact $(n+m)$-cube $\caF$ of $\frA$.
Then $(\varPhi_{*,m})$ turns out to be a homotopy from the zero map 
to the identity of the chain complex 
$(\WT{\Q}C_{*,m}(\frA), \PR{\partial })$.
Similarly, if we define 
$$
\varPsi_{n,m}:\WT{\Q}C_{n,m}(\frA)\to \WT{\Q}C_{n,m+1}(\frA)
$$
for $n\geq 1$ by $\varPsi_{n,m}(\caF)=-\rho_{n+m,n}(\caF)$, then 
$(\varPsi_{n,*})$ turns out to be a homotopy from the zero map 
to the identity of the chain complex 
$(\WT{\Q}C_{n,*}(\frA), \PPRR{\partial })$.
\qed

\vskip 1pc
{\it Proof of Thm.2.2}:
We will prove that $\Alt_*$ induces an isomorphism on homology.
We see $\frS_n$ as a subgroup of $\frS_{n+m}$ in the way that 
$$
\frS_n=\left\{\sigma \in \frS_{n+m}; \ \sigma (j)=j \ \ 
\text{if} \ \ n+1\leq j\leq n+m\right\}.
$$
Set 
$$
\Q C^0_{n,m}(\frA)=\left\{x\in \Q C_{n+m}(\frA); 
\sigma (x)=(\sgn \sigma )x \ \text{for any} \ \sigma 
\in \frS_n\right\}, 
$$
and 
$$
\WT{\Q}C^0_{n,m}(\frA)=\Q C^0_{n,m}(\frA)/\Q C^0_{n,m}(\frA)\cap 
D_{n+m}.
$$
Then $(\WT{\Q}C^0_{n,m}(\frA), \PR{\partial }, \PPRR{\partial })$ 
is a subcomplex of 
$(\WT{\Q}C_{n,m}(\frA), \PR{\partial }, \PPRR{\partial })$.
Define a map 
$$
\Alt_{n,m}:\WT{\Q}C_{n,m}(\frA)\to \WT{\Q}C_{n,m}^0(\frA)
$$
by 
$$
\Alt_{n,m}(\caF)=\frac{1}{n!}\sum_{\sigma \in \frS_n}(\sgn \sigma )
\sigma (\caF)
$$
for any exact $(n+m)$-cube $\caF$ of $\frA$.
Since the action of $\frS_n$ commutes with $\rho_{n+m,n+1}$, 
it holds that 
\begin{equation}
\Alt_{n,m+1}\varPhi_{n,m}=\varPhi_{n,m}\Alt_{n,m}. \label{sec2:eq1}
\end{equation}

If $m\geq 1$, then the chain complex 
$(\WT{\Q}C^0_{*,m}(\frA), \PR{\partial })$ is proved to be acyclic.
In fact, the map 
$$
\varPhi^{\alt}_{n,m}=\Alt_{n+1,m}\varPhi_{n,m}:
\WT{\Q}C^0_{n,m}(\frA)\to \WT{\Q}C^0_{n+1,m}(\frA)
$$
turns out to be a homotopy from the zero map to the identity.
On the other hand, since $(\WT{\Q}C^0_{n,*}(\frA), \PPRR{\partial })$ 
is a subcomplex of $(\WT{\Q}C_{n,*}(\frA), \PPRR{\partial })$ 
with the splitting map $\Alt_{n,*}$, Lem.2.4 implies that 
$(\WT{\Q}C_{n,*}(\frA), \PPRR{\partial })$ is also acyclic for $n\geq 1$.
Consequently, the inclusions 
\begin{align*}
\varepsilon_1:& \ \WT{\Q}C_*(\frA)=\WT{\Q}C^0_{0,*}(\frA)\to 
\TOT (\WT{\Q}C^0_{*,*}(\frA)),  \\
\varepsilon_2:& \ \WT{\Q}\ALT{*}{\frA}=\WT{\Q}C^0_{*,0}(\frA)
\to \TOT (\WT{\Q}C^0_{*,*}(\frA))
\end{align*}
are quasi-isomorphisms, therefore the composite 
$$
(\varepsilon_2)_*^{-1}(\varepsilon_1)_*:H_*(\WT{\Q}C_*(\frA))
\overset{\sim }{\to }H_*(\WT{\Q}\ALT{*}{\frA})
$$
is an isomorphism.

Suppose $m\geq 1$ and let 
$x\in \KER \PPRR{\partial }\subset \WT{\Q}C^0_{0,m}(\frA)$.
Then $x$ is homologous in \linebreak 
$\TOT (\WT{\Q}C^0_{*,*}(\frA))$ to 
$$
x-(\PR{\partial }-\PPRR{\partial })\varPhi^{\alt}_{0,m}(x)=
\PPRR{\partial }\varPhi^{\alt}_{0,m}(x)\in \WT{\Q}C^0_{1,m-1}(\frA).
$$
It is also homologous to 
\begin{align*}
\PPRR{\partial }\varPhi^{\alt}_{0,m}(x)-(\PR{\partial }+
\PPRR{\partial })\varPhi^{\alt}_{1,m-1}\PPRR{\partial }
\varPhi^{\alt}_{0,m}(x)
&=\varPhi^{\alt}_{0,m-1}\PR{\partial }\PPRR{\partial }
\varPhi^{\alt}_{0,m}(x)-\PPRR{\partial }\varPhi^{\alt}_{1,m-1}
\PPRR{\partial }\varPhi^{\alt}_{0,m}(x)  \\
&=\varPhi^{\alt}_{0,m-1}\PPRR{\partial }(x)-\PPRR{\partial }
\varPhi^{\alt}_{1,m-1}\PPRR{\partial }\varPhi^{\alt}_{0,m}(x)  \\
&=-\PPRR{\partial }\varPhi^{\alt}_{1,m-1}\PPRR{\partial }
\varPhi^{\alt}_{0,m}(x)\in \WT{\Q}C^0_{2,m-2}(\frA).
\end{align*}
Repeating this procedure we can say that $x$ is homologous in 
$\TOT (\WT{\Q}C^0_{*,*}(\frA))$ to 
$$
(-1)^{\frac{1}{2}m(m-1)}\PPRR{\partial }\varPhi^{\alt}_{m-1,1}
\PPRR{\partial }\varPhi^{\alt}_{m-2,2}\cdots \PPRR{\partial }
\varPhi^{\alt}_{0,m}(x)\in \WT{\Q}C^0_{m,0}(\frA).
$$
This means that the isomorphism 
$(\varepsilon_2)_*^{-1}(\varepsilon_1)_*$ is given by 
$$
[x]\mapsto (-1)^{\frac{1}{2}m(m-1)}[\PPRR{\partial }
\varPhi^{\alt}_{m-1,1}\PPRR{\partial }\varPhi^{\alt}_{m-2,2}\cdots 
\PPRR{\partial }\varPhi^{\alt}_{0,m}(x)].
$$

Suppose $m>1$ and $k\leq m-2$.
For an exact $m$-cube $\caF$ of $\frA$, 
\begin{align*}
\PPRR{\partial }\varPhi_{k,m-k}(\caF)=&
\sum_{j=1}^{m-k}\sum_{i=-1}^1(-1)^{i+j+k+1}\partial_{k+1+j}^i
\rho_{m,k+1}(\caF)  \\
=&\, (-1)^k\caF +\sum_{j=2}^{m-k}\sum_{i=-1}^1(-1)^{i+j+k+1}
\rho_{m-1,k+1}(\partial_{k+j}^i\caF).
\end{align*}
This implies that for $y\in \WT{\Q}C^0_{k,m-k}(\frA)$ we have 
\begin{equation}
\PPRR{\partial }\varPhi^{\alt}_{k,m-k}(y)=(-1)^k\Alt_{k+1,m-k-1}(y)-
\varPhi^{\alt}_{k,m-k-1}(\PPRR{\partial }_{k+1}y), \label{sec2:eq2}
\end{equation}
where 
$$
\PPRR{\partial }_{k+1}=\sum_{j=1}^{m-k-1}\sum_{i=-1}^1
(-1)^{i+j}\partial_{k+1+j}^i.
$$
For an exact $(m-1)$-cube $\caG$ of $\frA$, 
$$
\Alt_{k+2,m-k-1}\varPhi_{k+1,m-k-1}\varPhi_{k,m-k-1}(\caG)
=-\Alt_{k+2,m-k-1}\rho_{m,k+2}\rho_{m-1,k+1}(\caG), 
$$
and it is zero since 
$\rho_{m,k+2}\rho_{m-1,k+1}(\caG)=\rho_{m,k+1}\rho_{m-1,k+1}(\caG)$ 
is invariant by the transposition $(k+1, k+2)\in \frS_{k+2}$.
Hence 
$$
\Alt_{k+2,m-k-1}\varPhi_{k+1,m-k-1}\varPhi_{k,m-k-1}(z)=0
$$
for $z\in \WT{\Q}C^0_{k,m-k-1}(\frA)$, and it follows from 
\eqref{sec2:eq1} that 
\begin{align*}
\varPhi^{\alt}_{k+1,m-k-1}\varPhi^{\alt}_{k,m-k-1}(z)=&
\Alt_{k+2,m-k-1}\varPhi_{k+1,m-k-1}\Alt_{k+1,m-k-1}
\varPhi_{k,m-k-1}(z) \\
=&\Alt_{k+2,m-k-1}\varPhi_{k+1,m-k-1}\varPhi_{k,m-k-1}(z)=0.
\end{align*}
Taking the image by $\varPhi^{\alt}_{k+1,m-k-1}$ 
of the both sides of \eqref{sec2:eq2} we have 
\begin{align*}
\varPhi^{\alt}_{k+1,m-k-1}\PPRR{\partial }\varPhi^{\alt}_{k,m-k}(y)
=&\, (-1)^k\varPhi^{\alt}_{k+1,m-k-1}\Alt_{k+1,m-k-1}(y)  \\
=&\, (-1)^k\Alt_{k+2,m-k-1}\varPhi_{k+1,m-k-1}\Alt_{k+1,m-k-1}(y)
\end{align*}
for $y\in \WT{\Q}C^0_{k,m-k}(\frA )$, and by \eqref{sec2:eq1} 
it is equal to
$$
(-1)^k\Alt_{k+2,m-k-1}\varPhi_{k+1,m-k-1}(y)
=(-1)^k\varPhi^{\alt}_{k+1,m-k-1}(y).
$$
Hence for $m\geq 1$ and for $x\in \WT{\Q}C_{0,m}(\frA)$, 
\begin{align*}
(-1)^{\frac{1}{2}m(m-1)}&
\PPRR{\partial }\varPhi^{\alt}_{m-1,1}\PPRR{\partial }
\varPhi^{\alt}_{m-2,2}\cdots \PPRR{\partial }
\varPhi^{\alt}_{0,m}(x)  \\
=&(-1)^{\frac{1}{2}m(m-1)}\PPRR{\partial }
\varPhi^{\alt}_{m-1,1}\PPRR{\partial }\varPhi^{\alt}_{m-2,2}
\cdots \PPRR{\partial }\varPhi^{\alt}_{1,m-1}(x)  \\
=&-(-1)^{\frac{1}{2}m(m-1)}\PPRR{\partial }
\varPhi^{\alt}_{m-1,1}\PPRR{\partial }\varPhi^{\alt}_{m-2,2}
\cdots \PPRR{\partial }\varPhi^{\alt}_{2,m-2}(x)  \\
=&\cdots  \\
=&(-1)^{m-1}\PPRR{\partial }\varPhi^{\alt}_{m-1,1}(x).
\end{align*}
Since 
\begin{align*}
\PPRR{\partial }\varPhi^{\alt}_{m-1,1}(\caG)
=&(-1)^m\PPRR{\partial }\Alt_{m,1}\rho_{m,m}(\caG)  \\
=&(-1)^m\Alt_m\PPRR{\partial }\rho_{m,m}(\caG)=(-1)^{m-1}\Alt_m(\caG)
\end{align*}
for any exact $m$-cube $\caG$, we conclude that 
$$
(-1)^{\frac{1}{2}m(m-1)}
\PPRR{\partial }\varPhi^{\alt}_{m-1,1}\PPRR{\partial }
\varPhi^{\alt}_{m-2,2}\cdots \PPRR{\partial }
\varPhi^{\alt}_{0,m}(x)=\Alt_m(x).
$$
This means that the isomorphism 
$$
(\varepsilon_2)_*^{-1}(\varepsilon_1)_*:H_m(\WT{\Q}C_*(\frA))
\to H_m(\WT{\Q}\ALT{*}{\frA})
$$
agrees with the map induced by $\Alt_*:\WT{\Q}C_*(\frA)\to 
\WT{\Q}\ALT{*}{\frA}$, 
which completes the proof.
\qed

\vskip 1pc
\subsection{$\scC$-complexes}
In this subsection we will introduce $\scC$-complexes 
\cite{hanamura}.
A {\it $\scC$-complex} $A=(A_*^m, F_A^{m,n})$ is a family of 
chain complexes $(A^m_*, \partial_{A^m})$ for 
$m\in \Z$ such that $A_*^m=0$ except for finitely many $m$, 
with maps of abelian groups 
$$
F_A^{m,n}:A^m_*\to A^n_{*+n-m-1}
$$
for $m<n$ satisfying the relation 
$$
(-1)^mF_A^{m,n}\partial_{A^m}+(-1)^n\partial_{A^n}F_A^{m,n}
+\sum_{m<l<n}F_A^{l,n}F_A^{m,l}=0.
$$
Define the {\it total chain complex} $\TOT (A)$ of the 
$\scC$-complex $A$ by 
$$
\TOT (A)_p=\underset{m}{\oplus }A^m_{m+p}
$$
with the boundary map 
$$
\partial (x)^m=(-1)^m\partial_{A^m}(x^m)+
\sum_{l<m}F_A^{l,m}(x^l)
$$
for $x=(x^m)\in \TOT (A)_p=\underset{m}{\oplus }A^m_{m+p}$.

Let $A=(A^m_*, F_A^{m,n})$ and $B=(B^m_*, F_B^{m,n})$ be 
$\scC$-complexes.
A {\it map of $\scC$-complexes} from $A$ to $B$ is a family of 
maps of abelian groups 
$$
f^{m,n}:A^m_*\to B^n_{*+n-m}
$$
for $m\leq n$ satisfying the relation 
$$
(-1)^n\partial_{B^n}f^{m,n}+\sum_{l}F_B^{l,n}f^{m,l}=
(-1)^mf^{m,n}\partial_{A^m}+\sum_{l}f^{l,n}F_A^{m,l}.
$$
This condition is equivalent to that 
$$
\TOT (f)=\oplus f^{m,n}:\TOT (A)\to \TOT (B)
$$
is a map of chain complexes.
For two maps $f=(f^{m,n}):A\to B$ and $g=(g^{m,n}):B\to C$ of 
$\scC$-complexes, define the {\it composite} $gf$ by 
$$
(gf)^{m,n}=\sum_lg^{l,n}f^{m,l}:A_*^m\to C_{*+n-m}^n.
$$
It is obvious that $gf$ is also a map of $\scC$-complexes.

Let $f, g:A\to B$ be maps of $\scC$-complexes.
A {\it homotopy} $\varPhi =(\varPhi^{m,n})$ from $f$ to $g$ is 
a family of maps of abelian groups
$$
\varPhi^{m.n}:A^m_*\to B^n_{*+n-m+1}
$$
for any $m\leq n+1$ satisfying the relation 
$$
(-1)^m\varPhi^{m,n}\partial_{A^m}+\sum_l\varPhi^{l,n}F_A^{m,l}+
(-1)^n\partial_{B^n}\varPhi^{m,n}+\sum_lF_B^{l,n}\varPhi^{m,l}=
g^{m,n}-f^{m,n}.
$$
This condition is equivalent to that 
$$
\oplus \, \varPhi^{m,n}:\TOT (A)_*\to \TOT (B)_{*+1}
$$
is a chain homotopy from $\TOT (f)$ to $\TOT (g)$.

For a $\scC$-complex $A=(A_*^m, F_A^{m,n})$ and for 
$r\in \Z$, let $A[r]=(A[r]^m, F_{A[r]}^{m,n})$ be 
the $\scC$-complex given as follows:
$$
A[r]_*^m=A_*^{m+r}, \quad F_{A[r]}^{m,n}=(-1)^rF_A^{m+r,n+r}:
A_*^{m+r}\to A_*^{n+r}.
$$
Then there is a natural isomorphism of chain complexes 
$\TOT (A[r])\simeq \TOT (A)[r]$.

Let $f:A\to B$ be a map of $\scC$-complexes.
Define the {\it simple complex} $s(f)=(C^m_*, F_C^{m,n})$ of $f$ 
as follows: 
\begin{align*}
C^m_*&=A^m_*\oplus B^{m-1}_*,  \\
\partial_C(a, b)&=(\partial_A(a), \partial_B(b)),  \\
F_C^{m,n}(a, b)&=(F_A^{m,n}(a), f^{m,n-1}(a)-F_B^{m-1,n-1}(b)).
\end{align*}
It is easy to see that $s(f)$ is a $\scC$-complex such that 
its total chain complex is canonically isomorphic 
to the simple complex of the map $\TOT (f):\TOT (A)\to \TOT (B)$.
Let $p^{m,n}:s(f)_*^m=A_*^m\oplus B_*^{m-1}\to A_*^n$ be the map 
defined by 
$$
p^{m,n}(a, b)=\begin{cases}a, &\ m=n,  \\
0, &\ m\not= n, \end{cases}
$$
and 
$i^{m,n}:B[-1]_*^m=B_*^{m-1}\to s(f)_*^n=A_*^n\oplus B_*^{n-1}$ 
the map defined by 
$$
i^{m,n}(b)=\begin{cases}(0, b), &\ m=n,  \\
(0, 0), &\ m\not= n. \end{cases}
$$
Then $p=(p^{m,n}):s(f)\to A$ and $i=(i^{m,n}):B[-1]\to s(f)$ 
are maps of $\scC$-complexes.
Taking the total chain complexes of a sequence of maps of 
$\scC$-complexes 
$B[-1]\overset{i}{\to }s(f)\overset{p}{\to }A\overset{f}{\to }B$ 
yields a distinguished triangle of chain complexes.

The proposition below can be easily verified.

\vskip 1pc
\begin{prop}
Let 
$$
\begin{CD}
A @>{f}>>  B  \\
@V{\varphi_A}VV  @VV{\varphi_B}V  \\
\PR{A} @>{\PR{f}}>>  \PR{B}
\end{CD}
$$
be a diagram of maps of $\scC$-complexes which is commutative 
up to homotopy.
Let \linebreak 
$\varPhi_f:A\to \PR{B}$ be a homotopy from $\varphi_Bf$ to 
$\PR{f}\varphi_A$.
Then the family of maps $\varphi_s^{m,n}:s(f)^m\to s(\PR{f})^n$ 
defined by 
$$
\varphi_s^{m,n}(a, b)=
(\varphi_A^{m,n}(a), \varphi_B^{m-1,n-1}(b)+\varPhi_f^{m,n-1}(a))
$$
forms a map of $\scC$-complexes 
$\varphi_s=(\varphi_s^{m,n}):s(f)\to s(\PR{f})$.
Moreover, the diagram 
$$
\begin{CD}
B[-1] @>{i}>> s(f) @>{p}>> A \\
@V{\varphi_B[-1]}VV @VV{\varphi_s}V @VV{\varphi_A}V  \\
\PR{B}[-1] @>{\PR{i}}>> s(\PR{f}) @>{\PR{p}}>> \PR{A}
\end{CD}
$$
is strictly commutative.
\end{prop}

\vskip 1pc
\begin{prop}
Let $f:A\to B$ be a map of $\scC$-complexes, and $p:s(f)\to A$ 
the map of \linebreak 
$\scC$-complexes defined above.
Assume that $f$ has a right inverse map up to homotopy, namely, 
there is a map $g:B\to A$ such that $fg$ is homotopy equivalent to 
the identity of $B$.
Let $\varPsi $ be a homotopy from the identity of $B$ to $fg$.
Then $t=(t^{m,n}):A\to s(f)$ defined by 
$$
t^{m,n}(a)=\left(\delta^{m,n}(a)-\sum_lg^{l,n}f^{m,l}(a), 
-\sum_l\varPsi^{l,n-1}f^{m,l}(a)\right), 
$$
where 
$$
\delta^{m,n}(a)=\begin{cases}a, &\ m=n,  \\
0, &\ m\not= n,   \end{cases}
$$
is a map of $\scC$-complexes such that $pt=\IId_A-gf$, and 
it is a left inverse map of $p$ up to homotopy, namely, 
$tp$ is homotopy equivalent to the identity of $s(f)$.
Moreover, the composite $tg:B\to s(f)$ is homotopy equivalent 
to the zero map.
\end{prop}

{\it Proof}: 
It is easy to see that $t$ is a map of $\scC$-complexes such that 
$pt=\IId_A-gf$.
Let $\varPsi_1^{m,n}:s(f)^m_*\to s(f)^n_{*+n-m+1}$ be the map 
defined by 
$$
\varPsi_1^{m,n}(a, b)=\left(-g^{m-1,n}(b), 
-\varPsi^{m-1, n-1}(b)\right).
$$
Then $\varPsi_1=(\varPsi_1^{m,n})$ turns out to be a homotopy 
from the identity of $s(f)$ to $tp$.
Let $\varPsi_2^{m,n}:B_*^m\to s(f)_{*+n-m+1}^n$ 
be the map defined by 
$$
\varPsi_2^{m,n}(b)=\left(-\sum_lg^{l,n}\varPsi^{m,l}(b), 
-\sum_l\varPsi^{l,n-1}\varPsi^{m,l}(b)\right), 
$$
then $\varPsi_2=(\varPsi_2^{m,n})$ turns out to be a homotopy 
from the zero map to $tg$.
\qed

\vskip 1pc
Let 
$$
\begin{CD}
A @>{f}>>  B  \\
@V{\varphi_A}VV  @VV{\varphi_B}V  \\
\PR{A} @>{\PR{f}}>>  \PR{B}
\end{CD}
$$
be a diagram of maps of $\scC$-complexes which is commutative 
up to homotopy.
Suppose that $f$ and $\PR{f}$ have right inverse maps 
$g:B\to A$ and $\PR{g}:\PR{B}\to \PR{A}$ up to homotopy 
such that $\varphi_Ag$ is homtopy equivalent to $\PR{g}\varphi_B$.
Let 
\begin{align*}
\varPhi_f:&\ A\to \PR{B}, \\
\varPhi_g:&\ B\to \PR{A}, \\
\varPsi:&\ B\to B, \\
\PR{\varPsi }:&\ \PR{B}\to \PR{B}
\end{align*}
be homotopies from $\varphi_Bf$ to $\PR{f}\varphi_A$, from 
$\varphi_Ag$ to $\PR{g}\varphi_B$, from $\IId_B$ to $fg$, 
and from $\IId_{\PR{B}}$ to $\PR{f}\PR{g}$ respectively.
Then we have two homotopies 
$$
\varPhi_fg+\PR{f}\varPhi_g+\varphi_B\varPsi, \ 
\PR{\varPsi }\varphi_B:B\to \PR{B}
$$
from $\varphi_B$ to $\PR{f}\PR{g}\varphi_B$.

\vskip 1pc
\begin{defn}
With the above notations, a second homotopy from 
$\varPhi_fg+\PR{f}\varPhi_g+\varphi_B\varPsi$ to 
$\PR{\varPsi }\varphi_B$ is a family of maps 
$$
\varTheta^{m,n}:B^m_*\to {\PR{B}}^n_{*+n-m+2}
$$
for $m\leq n+2$ satisfying the relation 
\begin{align*}
(-1)^n&\partial \varTheta^{m,n}+\sum_lF_{\PR{B}}^{l,n}\varTheta^{m,l}
-(-1)^m\varTheta^{m,n}\partial -\sum_l\varTheta^{l,n}F_B^{m,l} \\
=&\sum_l{\PR{\varPsi }}^{l,n}\varphi_B^{m,l}
-\sum_l\varPhi_f^{l,n}g^{m,l}-\sum_l{\PR{f}}^{l,n}\varPhi_g^{m,l}
-\sum_l\varphi_B^{l,n}\varPsi^{m,l}.
\end{align*}
\end{defn}

\vskip 1pc
The proposition below is verified by a direct calculation:

\vskip 1pc
\begin{prop}
With the above notations, let $t:A\to s(f)$ and 
$\PR{t}:\PR{A}\to s(\PR{f})$ be the maps given in Prop.2.6.
Then $\varphi_st$ is homotopy equivalent to $\PR{t}\varphi_A$, 
and the homotopy $\varPi :A\to s(\PR{f})$ from $\varphi_st$ to 
$\PR{t}\varphi_A$ is given by 
$$
\varPi^{m,n}(a)=\left(-\sum_l\varPhi_g^{l,n}f^{m,l}(a)-
\sum_l{\PR{g}}^{l,n}\varPhi_f^{m,l}(a), 
-\sum_l{\PR{\varPsi }}^{l,n-1}\varPhi_f^{m,l}(a)
+\sum_l\varTheta^{l,n-1}f^{m,l}(a)\right).
$$
\end{prop}

\vskip 1pc
\subsection{Multi-relative complexes of exact cubes}
Let $X$ be a scheme.
As we mentioned in the introduction, the category $\frP(X)$ 
of vector bundles of finite rank on $X$ is a small exact category.
Hence we can obtain the chain complex of exact cubes of vector bundles 
$\WT{\Q}C_*(X)$ and its alternating part 
$\WT{\Q}\ALT{*}{X}$ as in \S 2.2.
Then Thm.2.1 and Thm.2.2 imply the following:

\vskip 1pc
\begin{thm}
The homology groups of the chain complexes $\WT{\Q}C_*(X)$ and 
$\WT{\Q}\ALT{*}{X}$ are canonically isomorphic to the rational $K$-theory of $X$:
$$
H_n(\WT{\Q}\ALT{*}{X})\simeq H_n(\WT{\Q}C_*(X))\simeq \QQ{K_n(X)}.
$$
\end{thm}

\vskip 1pc
Consider a sequence of morphisms of schemes 
$$
X_0\overset{f_1}{\to }X_1\overset{f_2}{\to }
\cdots \overset{f_r}{\to }X_r.
$$
For a vector bundle $\caF$ on $X_r$, let us define an exact 
$(r-1)$-cube $(f_1, \ldots , f_r)^*\caF $ on $X_0$ as follows: 
In the case of $r\geq 2$, 
\begin{multline*}
\left((f_1, \ldots , f_r)^*\caF\right)_{\alpha_1, \ldots , 
\alpha_{r-1}} \\
=\begin{cases}
0, & \ \exists j\ \text{such that} \ 
\alpha _j=1,   \\
(f_{j_1}\cdots f_1)^*(f_{j_2}\cdots f_{j_1+1})^*\cdots 
(f_{j_l}\cdots f_{j_{l-1}+1})^*(f_r\cdots f_{j_l+1})^*\caF, 
& \ \text{otherwise}, \end{cases}
\end{multline*}
where $\alpha ^{-1}(-1)=\{j_1, \ldots , j_l\}$ with 
$j_1<\cdots <j_l$.
The maps in $(f_1, \ldots , f_r)^*\caF $ are the zero maps 
or the natural isomorphisms.
When $r=1$, define $(f_1)^*\caF =f_1^*\caF$.
We can generalize this procedure to exact cubes in the 
following way: 
For an exact $n$-cube $\caF$ on $X_r$, 
$(f_1, \ldots , f_r)^*\caF $ is an exact $(n+r-1)$-cube 
on $X_0$ so that 
$$
\partial_r^{\alpha_r}\cdots \partial_{n+r-1}^{\alpha_{n+r-1}}
((f_1, \ldots , f_r)^*\caF )=
(f_1, \ldots , f_r)^*(\partial_1^{\alpha_r}\cdots 
\partial_n^{\alpha_{n+r-1}}\caF ).
$$

\vskip 1pc
\begin{prop}
For an exact $n$-cube $\caF$ on $X_r$, the faces 
of the exact cube $(f_1, \ldots , f_r)^*\caF $ are as follows:
For $1\leq j\leq r-1$, 
$$
\partial_j^i((f_1, \ldots , f_r)^*\caF )=\begin{cases}
(f_1, \ldots , f_j)^*\left((f_{j+1}, \ldots , f_r)^*\caF \right), 
& \ i=-1, \\
(f_1, \ldots , f_{j+1}f_j, \ldots , f_r)^*\caF, & \ i=0,  \\
0, & \ i=1
\end{cases}
$$
and for $r\leq j$, 
$$
\partial_j^i((f_1, \ldots , f_r)^*\caF )=
(f_1, \ldots , f_r)^*(\partial_{j-r+1}^i\caF ).
$$
\end{prop}

\vskip 1pc
Let $Y_1, \ldots , Y_r$ be closed subschemes of a scheme $X$.
Let $Y_J=\underset{j\in J}{\cap }Y_j$ for 
$J\subset \{1, \ldots , r\}$, and 
$\iota_k:Y_{J\cup \{k\}}\hookrightarrow Y_J$ 
the embedding for $k\notin J$.
For a division $J=K\coprod I$ such that 
$K=\{k_1, \ldots , k_{n-m}\}$ with $k_1<\cdots <k_{n-m}$ 
and for $\sigma \in \frS_{n-m}$, we have a sequence of embeddings 
$$
Y_J\overset{\iota_{k_{\sigma (1)}}}{\hookrightarrow }
Y_{J-\{k_{\sigma (1)}\}}
\overset{\iota_{k_{\sigma (2)}}}{\hookrightarrow }
Y_{J-\{k_{\sigma (1)}, k_{\sigma (2)}\}}
\overset{\iota_{k_{\sigma (3)}}}{\hookrightarrow }\cdots 
\overset{\iota_{k_{\sigma (n-m)}}}{\hookrightarrow }Y_I.
$$
Note that $(\iota_{k_{\sigma (1)}}, \ldots , 
\iota_{k_{\sigma (n-m)}})^*\caF $ is degenerate if so is $\caF $.
Hence with the above sequence we can associate a map 
$$
\varXi_K=\sum_{\sigma \in \frS_{n-m}}
(\sgn \sigma )(\iota_{k_{\sigma (1)}}, \ldots , 
\iota_{k_{\sigma (n-m)}})^*:\WT{\Q}C_*(Y_I)\to 
\WT{\Q}C_{*+n-m-1}(Y_J).
$$

Let us consider the composite of $\varXi_K$ with the boundary map.
To do this, we need to introduce signature of a division of subsets of 
$\{1, \ldots , r\}$.
For a division $K\coprod I=J$ such that 
\begin{align*}
K=\{k_1, \ldots , k_{n-m}\}, &\ \ k_1<\cdots <k_{n-m},  \\
I=\{i_1, \ldots , i_m\}, &\ \ i_1<\cdots <i_m,  \\
J=\{j_1, \ldots , j_n\}, &\ \ j_1<\cdots <j_n, 
\end{align*}
define the signature as follows: 
$$
\sgn \binom{K \ I}{J}=\sgn 
\begin{pmatrix}
k_1 & \ldots & k_{n-m} & i_1 & \ldots & i_m \\
j_1 & \hdotsfor{4} & j_n
\end{pmatrix}.
$$
The following lemma is easily verified.

\vskip 1pc
\begin{lem}
Let $J=L\coprod \PR{L}\coprod I$ be a division and write 
$K=L\coprod \PR{L}$ and $P=\PR{L}\coprod I$.
Then 
$$
\sgn \binom{K \ I}{J}\sgn \binom{L \ \PR{L}}{K}
=\sgn \binom{L \ P}{J}\sgn \binom{\PR{L} \ I}{P}.
$$
\end{lem}

\vskip 1pc
\begin{prop}
For $x\in \WT{\Q}C_*(Y_I)$, 
$$
\partial \varXi_K(x)+(-1)^{n-m}\varXi_K(\partial x)=
\sum_{a=1}^{n-m-1}(-1)^{a+1}\underset{|L|=a}{\sum_{L\coprod \PR{L}=K}}
\sgn \sbinom{L \ \PR{L}}{K}\varXi_L\varXi_{\PR{L}}(x).
$$
\end{prop}

{\it Proof}: 
Prop.2.10 implies that 
\begin{align}
\partial \varXi_K(x)=&\sum_{\sigma \in \frS_{n-m}}
(\sgn \sigma )\sum_{a=1}^{n-m-1}
(-1)^{a+1}(\iota_{k_{\sigma (1)}}, \ldots , \iota_{k_{\sigma (a)}})^*
\left((\iota_{k_{\sigma (a+1)}}, \ldots , 
\iota_{k_{\sigma (n-m)}})^*(x)\right) \notag \\
&+\sum_{\sigma \in \frS_{n-m}}(\sgn \sigma )\sum_{a=1}^{n-m-1}
(-1)^a(\iota_{k_{\sigma (1)}}, \ldots , 
\iota_{k_{\sigma (a+1)}}\iota_{k_{\sigma (a)}}, \ldots , 
\iota_{k_{\sigma (n-m)}})^*(x) \label{sec2:eq3} \\
&+(-1)^{n-m+1}\sum_{\sigma \in \frS_{n-m}}(\sgn \sigma )
(\iota_{k_{\sigma (1)}}, \ldots , \iota_{k_{\sigma (n-m)}})^*
(\partial x) \label{sec2:eq4}.
\end{align}
In the above, \eqref{sec2:eq3} is obviously zero, and 
\eqref{sec2:eq4} is equal to $(-1)^{n-m+1}\varXi_K(\partial x)$.
Hence 
\begin{align*}
\partial &\varXi_K(x)+(-1)^{n-m}\varXi_K(\partial x) \\
=&\sum_{\sigma \in \frS_{n-m}}(\sgn \sigma )\sum_{a=1}^{n-m-1}
(-1)^{a+1}(\iota_{k_{\sigma (1)}}, \ldots , \iota_{k_{\sigma (a)}})^*
\left((\iota_{k_{\sigma (a+1)}}, \ldots , \iota_{k_{\sigma (n-m)}})^*
(x)\right). 
\end{align*}
If we denote 
\begin{align*}
L=\{k_{\sigma (1)}, \ldots , k_{\sigma (a)}\}=&\ 
\{l_1, \ldots , l_a\}, \\
\PR{L}=\{k_{\sigma (a+1)}, \ldots , k_{\sigma (n-m)}\}
=&\ \{\PR{l}_1, \ldots , \PR{l}_{n-m-a}\}
\end{align*}
with $l_1<\cdots <l_a$ and $\PR{l}_1<\cdots <\PR{l}_{n-m-a}$, then 
\begin{align*}
\partial \varXi_K(x)&\ +(-1)^{n-m}\varXi_K(\partial x) \\
=&\sum_{a=1}^{n-m-1}(-1)^{a+1}\underset{|L|=a}{\sum_{L\coprod \PR{L}=K}}
\sgn \sbinom{L \ \PR{L}}{K}\times \\
&\hskip 2pc \sum_{\tau \in \frS_a}(\sgn \tau )(\iota_{\tau (l_1)}, 
\ldots , \iota_{\tau (l_a)})^*\left(\sum_{\eta \in \frS_{n-m-a}}
(\sgn \eta )(\iota_{\eta (\PR{l}_1)}, \ldots , 
\iota_{\eta (\PR{l}_{n-m-a})})^*(x)\right) \\
=&\sum_{a=1}^{n-m-1}(-1)^{a+1}\underset{|L|=a}{\sum_{L\coprod \PR{L}=K}}
\sgn \sbinom{L \ \PR{L}}{K}\varXi_L\varXi_{\PR{L}}(x), 
\end{align*}
which completes the proof.
\qed

\vskip 1pc
Define a map 
$$
F^{m,n}:\underset{|I|=m}{\oplus }\WT{\Q}C_*(Y_I)\to 
\underset{|J|=n}{\oplus }\WT{\Q}C_{*+n-m-1}(Y_J)
$$
by 
$$
F^{m,n}(x)_J=(-1)^n\underset{|I|=m}{\sum_{K\coprod I=J}}
\sgn \sbinom{K \ I}{J}\varXi_K(x_I)
$$
for $x=(x_I)\in \underset{|I|=m}{\oplus }\WT{\Q}C_*(Y_I)$.

\vskip 1pc
\begin{cor}
$$
\WT{\Q}C_*(X; Y_1, \ldots , Y_r)=
\left(\underset{|I|=m}{\oplus }\WT{\Q}C_*(Y_I), F^{m,n}\right)
$$
is a $\scC$-complex.
\end{cor}

{\it Proof}: 
Let $x=(x_I)\in \underset{|I|=m}{\oplus }\WT{\Q}C_*(Y_I)$.
Then it follows from Prop.2.12 and Lem.2.11 that 
\begin{align*}
(-1)^n&\partial F^{m,n}(x)_J+(-1)^mF^{m,n}(\partial x)_J \\
=&\underset{|I|=m}{\sum_{K\coprod I=J}}\sgn \sbinom{K \ I}{J}
\sum_{a=1}^{n-m-1}(-1)^{a+1}\underset{|L|=a}{\sum_{L\coprod \PR{L}=K}}
\sgn \sbinom{L \ \PR{L}}{K}\varXi_L\varXi_{\PR{L}}(x) \\
=&\sum_{a=1}^{n-m-1}(-1)^{a+1}\underset{|L|=a}{\sum_{L\coprod P=J}}
\sgn \sbinom{L \ P}{J}\underset{|I|=m}{\sum_{\PR{L}\coprod I=P}}
\sgn \sbinom{\PR{L} \ I}{P}\varXi_L\varXi_{\PR{L}}(x_I) \\
=&\sum_{a=1}^{n-m-1}(-1)^{n+1}\underset{|L|=a}{\sum_{L\coprod P=J}}
\sgn \sbinom{L \ P}{J}\varXi_L\left(F^{m,n-a}(x)_P\right) \\
=&-\sum_{a=1}^{n-m-1}F^{n-a,n}F^{m,n-a}(x)_J, 
\end{align*}
which completes the proof.
\qed

\vskip 1pc
We next consider pull-back map associated with a morphism 
of schemes.
Let $T$ be another scheme with closed subschemes 
$D_1, \ldots , D_r$, and $f:X\to T$ a morphism 
such that $f(Y_j)\subset D_j$ for $1\leq j\leq r$.
In what follows, we denote such a morphism by 
$$
f:(X; Y_1, \ldots , Y_r)\to (T; D_1, \ldots , D_r).
$$
Denote the restriction of $f$ to $Y_I\to D_I$ by 
the same symbol $f$ to simplify the notations.
For a division $K\coprod I=J$ of subsets of $\{1, \ldots , r\}$ 
such that $K=\{k_1, \ldots , k_{n-m}\}$ with $k_1<\cdots <k_{n-m}$, 
define a map 
$\varXi_{K,f}: \WT{\Q}C_*(D_I)\to \WT{\Q}C_{*+n-m}(Y_J)$ by 
$$
\varXi_{K,f}=\sum_{p=0}^{n-m}(-1)^p
\sum_{\sigma \in \frS_{n-m}}(\sgn \sigma )
(\iota_{k_{\sigma (1)}}, \ldots , \iota_{k_{\sigma(p)}}, f, 
\iota_{k_{\sigma(p+1)}}, \ldots , \iota_{k_{\sigma (n-m)}})^*.
$$
In particular, in the case that $K$ is the emptyset, 
$\varXi_{\emptyset ,f}=f^*:\WT{\Q}C_*(D_I)\to \WT{\Q}C_*(Y_I)$.

\vskip 1pc
\begin{prop}
For $x\in \WT{\Q}C_*(D_I)$, 
\begin{multline*}
\partial \varXi_{K,f}(x)+(-1)^{n-m+1}\varXi_{K,f}(\partial x) \\
=-\sum_{a=1}^{n-m}\underset{|L|=a}{\sum_{L\coprod \PR{L}=K}}
\sgn \sbinom{L \ \PR{L}}{K}\varXi_L\varXi_{\PR{L},f}(x)
+\sum_{a=0}^{n-m-1}(-1)^a\underset{|L|=a}{\sum_{L\coprod \PR{L}=K}}
\sgn \sbinom{L \ \PR{L}}{K}\varXi_{L,f}\varXi_{\PR{L}}(x).
\end{multline*}
\end{prop}

{\it Proof}: 
Prop.2.10 implies that 
$$
\partial \varXi_{K,f}(x)+(-1)^{n-m+1}\varXi_{K,f}(\partial x)
=\sum_{a=1}^{n-m}(-1)^{a+1}\partial_a^{-1}\varXi_{K,f}(x)
+\sum_{a=1}^{n-m}(-1)^a\partial_a^0\varXi_{K,f}(x), 
$$
and 
\begin{align*}
\sum_{a=1}^{n-m}(-1)^a\partial_a^0&\varXi_{K,f}(x)  \\
=\sum_{\sigma \in \frS_{n-m}}(\sgn \sigma )&\left\{
\sum_{p=0}^{n-m}\sum_{a=1}^{p-1}(-1)^{a+p}(\ldots , 
\iota_{k_{\sigma (a+1)}}\iota_{k_{\sigma (a)}}, \ldots , 
\iota_{k_{\sigma (p)}}, f, \iota_{k_{\sigma (p+1)}}, \ldots )^*(x)
\right. \\
&\hskip 1pc +\sum_{p=0}^{n-m}\sum_{a=p+1}^{n-m-1}(-1)^{a+p-1}
(\ldots , \iota_{k_{\sigma (p)}}, f, \iota_{k_{\sigma (p+1)}}, 
\ldots , \iota_{k_{\sigma (a+1)}}\iota_{k_{\sigma (a)}}, 
\ldots )^*(x) \\
&\hskip 1pc +\left.\sum_{p=1}^{n-m}(\ldots , 
f\iota_{k_{\sigma (p)}}, \ldots )^*(x)-\sum_{p=0}^{n-m-1}(\ldots , 
\iota_{k_{\sigma (p+1)}}f, \ldots )^*(x)\right\}, 
\end{align*}
which is obviously zero.
On the other hand, 
\begin{align*}
\sum_{a=1}^{n-m}&(-1)^{a+1}\partial_a^{-1}\varXi_{K,f}(x) \\
=&\sum_{p=1}^{n-m}\sum_{a=1}^p(-1)^{a+p+1}
\sum_{\sigma \in \frS_{n-m}}(\sgn \sigma )\times \\
&\hskip 3pc (\iota_{k_{\sigma (1)}}, \ldots , 
\iota_{k_{\sigma (a)}})^*\left((\iota_{k_{\sigma (a+1)}}, \ldots , 
\iota_{k_{\sigma (p)}}, f, \iota_{k_{\sigma (p+1)}}, \ldots , 
\iota_{k_{\sigma (n-m)}})^*(x)\right)  \\
&+\sum_{p=0}^{n-m-1}\sum_{a=p}^{n-m-1}(-1)^{a+p}
\sum_{\sigma \in \frS_{n-m}}(\sgn \sigma )\times   \\
&\hskip 2pc (\iota_{k_{\sigma (1)}}, \ldots , \iota_{k_{\sigma (p)}}, 
f, \iota_{k_{\sigma (p+1)}}, \ldots , \iota_{k_{\sigma (a)}})^*
\left((\iota_{k_{\sigma (a+1)}}, \ldots , 
\iota_{k_{\sigma (n-m)}})^*(x)\right)  \\
=&-\sum_{a=1}^{n-m}\underset{|L|=a}{\sum_{L\coprod \PR{L}=K}}
\sgn \sbinom{L \ \PR{L}}{K}\varXi_L\varXi_{\PR{L},f}(x)
+\sum_{a=0}^{n-m-1}\underset{|L|=a}{\sum_{L\coprod \PR{L}=K}}
(-1)^a\sgn \sbinom{L \ \PR{L}}{K}\varXi_{L,f}\varXi_{\PR{L}}(x), 
\end{align*}
which completes the proof.
\qed

\vskip 1pc
Define a map 
$$
(f^*)^{m,n}:\underset{|I|=m}{\oplus }\WT{\Q}C_*(D_I)\to 
\underset{|J|=n}{\oplus }\WT{\Q}C_{*+n-m}(Y_J)
$$
by 
$$
(f^*)^{m,n}(x)_J=\underset{|I|=m}{\sum_{K\coprod I=J}}
\sgn \sbinom{K \ I}{J}\varXi_{K,f}(x_I)
$$
for $x=(x_I)\in \underset{|I|=m}{\oplus }\WT{\Q}C_*(D_I)$.

\vskip 1pc
\begin{cor}
$$
f^*={(f^*)^{m,n}}:\WT{\Q}C_*(T; D_1, \ldots , D_r)\to 
\WT{\Q}C_*(X; Y_1, \ldots , Y_r)
$$
is a map of $\scC$-complexes.
\end{cor}

{\it Proof}: 
It follows from Prop.2.14 and Lem.2.11 that 
\begin{align*}
(-1)^n&\partial (f^*)^{m,n}(x)_J-(-1)^m(f^*)^{m,n}(\partial x)_J \\
=&(-1)^n\underset{|I|=m}{\sum_{K\coprod I=J}}
\sgn \sbinom{K \ I}{J}\left\{
-\sum_{a=1}^{n-m}\underset{|L|=a}{\sum_{L\coprod \PR{L}=K}}
\sgn \sbinom{L \ \PR{L}}{K}\varXi_L\varXi_{\PR{L},f}(x_I)\right. \\
&\hskip 10pc +\left.\sum_{a=0}^{n-m-1}(-1)^a
\underset{|L|=a}{\sum_{L\coprod \PR{L}=K}}
\sgn \sbinom{L \ \PR{L}}{K}\varXi_{L,f}
\varXi_{\PR{L}}(x_I)\right\} \\
=&(-1)^{n+1}\sum_{a=1}^{n-m}\underset{|L|=a}{\sum_{L\coprod P=J}}
\sgn \sbinom{L \ P}{J}\underset{|I|=m}{\sum_{\PR{L}\coprod I=P}}
\sgn \sbinom{\PR{L} \ I}{P}\varXi_L\varXi_{\PR{L},f}(x_I) \\
&+\sum_{a=0}^{n-m-1}(-1)^{n+a}
\underset{|L|=a}{\sum_{L\coprod P=J}}\sgn \sbinom{L \ P}{J}
\underset{|I|=m}{\sum_{\PR{L}\coprod I=P}}
\sgn \sbinom{\PR{L} \ I}{P}\varXi_{L,f}\varXi_{\PR{L}}(x_I)  \\
=&-\sum_{a=1}^{n-m}F^{n-a,n}(f^*)^{m,n-a}(x)_J+
\sum_{a=0}^{n-m-1}(f^*)^{n-a,n}F^{m,n-a}(x)_J, 
\end{align*}
which completes the proof.
\qed

\vskip 1pc
Comparing the definition of the simple complex of the pull-back map 
$\iota_r^*$ with the definition of $\WT{\Q}C_*(X; Y_1, \ldots , Y_r)$, 
we obtain the following corollary:

\vskip 1pc
\begin{cor}
The $\scC$-complex $\WT{\Q}C_*(X; Y_1, \ldots , Y_r)$ is 
canonically isomorphic to the simple complex of the map 
$$
\iota_r^*:\WT{\Q}C_*(X; Y_1, \ldots , Y_{r-1})\to 
\WT{\Q}C_*(Y_r; Y_1\cap Y_r, \ldots , Y_{r-1}\cap Y_r)
$$
induced by the embedding $\iota_r:Y_r\hookrightarrow X$.
Hence the homology groups of 
$\WT{\Q}C_*(X; Y_1, \ldots , Y_r)$ are isomorphic to the rational 
multi-relative $K$-theory of $(X; Y_1, \ldots , Y_r)$: 
$$
H_n(\WT{\Q}C_*(X; Y_1, \ldots , Y_r))\simeq 
\QQ{K_n(X; Y_1, \ldots , Y_r)}.
$$
\end{cor}

\vskip 1pc
Let 
$$
(X; Y_1, \ldots , Y_r)\overset{f}{\to }
(T; D_1, \ldots , D_r)\overset{g}{\to }
(S; E_1, \ldots , E_r)
$$
be morphisms of schemes with closed subschemes.
For a division $K\coprod I=J$ of subsets of $\{1, \ldots , r\}$ 
such that $K=\{k_1, \ldots , k_{n-m}\}$ with 
$k_1<\cdots <k_{n-m}$, define a map 
$\varXi_{K,f,g}: \WT{\Q}C_*(E_I)\to \WT{\Q}C_{*+n-m+1}(Y_J)$ by 
$$
\varXi_{K,f,g}=\sum_{0\leq p\leq q\leq n-m}(-1)^{p+q}
\sum_{\sigma \in \frS_{n-m}}(\sgn \sigma )
(\ldots , \iota_{k_{\sigma (p)}}, f, \iota_{k_{\sigma (p+1)}}, 
\ldots , \iota_{k_{\sigma (q)}}, g, \iota_{k_{\sigma (q+1)}}, 
\ldots )^*.
$$

\vskip 1pc
\begin{prop}
If we write $h=gf$, then for $x\in \WT{\Q}C_*(E_I)$, 
\begin{align*}
\partial \varXi_{K,f,g}&(x)+(-1)^{n-m}\varXi_{K,f,g}(\partial x) \\
=&\sum_{a=1}^{n-m}\underset{|L|=a}{\sum_{L\coprod \PR{L}=K}}
(-1)^{a+1}\sgn \sbinom{L \ \PR{L}}{K}\varXi_L\varXi_{\PR{L},f,g}(x)
+\sum_{a=0}^{n-m}\underset{|L|=a}{\sum_{L\coprod \PR{L}=K}}
\sgn \sbinom{L \ \PR{L}}{K}\varXi_{L,f}\varXi_{\PR{L},g}(x)  \\
&+\sum_{a=0}^{n-m-1}\underset{|L|=a}{\sum_{L\coprod \PR{L}=K}}
(-1)^{a+1}\sgn \sbinom{L \ \PR{L}}{K}\varXi_{L,f,g}\varXi_{\PR{L}}(x)
-\varXi_{K,h}(x).
\end{align*}
\end{prop}

{\it Proof}:
Prop.2.10 implies that 
\begin{align*}
\partial \varXi_{K,f,g}(x)+&\, (-1)^{n-m}\varXi_{K,f,g}(\partial x) \\
=&\sum_{a=1}^{n-m+1}(-1)^{a+1}\partial_a^{-1}\varXi_{K,f,g}(x) 
+\sum_{a=1}^{n-m+1}(-1)^a\partial_a^0\varXi_{K,f,g}(x).
\end{align*}
In the same way as in the proof of Prop.2.12 and 
Prop.2.14 we can show that 
\begin{align*}
\sum_{a=1}^{n-m+1}&(-1)^a\partial_a^0\varXi_{K,f,g}(x) \\
=&\sum_{p=0}^{n-m}\sum_{\sigma \in \frS_{n-m}}(-1)^{p+1}
(\sgn \sigma )(\iota_{k_{\sigma (1)}}, \ldots , 
\iota_{k_{\sigma (p)}}, gf, \iota_{k_{\sigma (p+1)}}, \ldots , 
\iota_{k_{\sigma (n-m)}})^*(x)  \\
=&-\varXi_{K,h}(x).
\end{align*}
On the other hand, 
\begin{align*}
\sum_{a=1}^{n-m+1}&\, (-1)^{a+1}\partial_a^{-1}\varXi_{K,f,g}(x)  \\
=&\, \sum_{\sigma \in \frS_{n-m}}(\sgn \sigma )
\sum_{0\leq p\leq q\leq n-m}(-1)^{p+q}\times \\
&\left\{\sum_{a=1}^p(-1)^{a+1}(\iota_{k_{\sigma (1)}}, 
\ldots , \iota_{k_{\sigma (a)}})^*\left((\iota_{k_{\sigma (a+1)}}, 
\ldots , \iota_{k_{\sigma (p)}}, f, \ldots , 
\iota_{k_{\sigma (q)}}, g, \ldots )^*(x)\right)\right. \\
&\hskip 1pc +\sum_{a=p}^q(-1)^a(\ldots, 
\iota_{k_{\sigma (p)}}, f, \ldots , \iota_{k_{\sigma (a)}})^*
\left((\iota_{k_{\sigma (a+1)}}, \ldots , 
\iota_{k_{\sigma (q)}}, g, \ldots )^*(x)\right) \\
&\hskip 1pc +\left.\sum_{a=q}^{n-m-1}(-1)^{a+1}(\ldots, 
\iota_{k_{\sigma (p)}}, f, \ldots , \iota_{k_{\sigma (q)}}, g, 
\ldots , \iota_{k_{\sigma (a)}})^*\left((\iota_{k_{\sigma (a+1)}}, 
\ldots , \iota_{k_{\sigma (n-m)}})^*(x)\right)\right\} \\ 
=&\sum_{a=1}^{n-m}\underset{|L|=a}{\sum_{L\coprod \PR{L}=K}}
(-1)^{a+1}\sgn \sbinom{L \ \PR{L}}{K}\varXi_L\varXi_{\PR{L},f,g}(x)
+\sum_{a=1}^{n-m}\underset{|L|=a}{\sum_{L\coprod \PR{L}=K}}
\sgn \sbinom{L \ \PR{L}}{K}\varXi_{L,f}\varXi_{\PR{L},g}(x) \\
&+\sum_{a=1}^{n-m}\underset{|L|=a}{\sum_{L\coprod \PR{L}=K}}
(-1)^{a+1}\sgn \sbinom{L \ \PR{L}}{K}\varXi_{L,f,g}\varXi_{\PR{L}}(x), 
\end{align*}
which completes the proof.
\qed

\vskip 1pc
Define a map 
$$
\varPhi^{m,n}:\underset{|I|=m}{\oplus }\WT{\Q}C_*(E_I)\to 
\underset{|J|=n}{\oplus }\WT{\Q}C_{*+n-m+1}(Y_J)
$$
by 
$$
\varPhi^{m,n}(x)_J=(-1)^n\underset{|I|=m}{\sum_{K\coprod I=J}}
\sgn \sbinom{K \ I}{J}\varXi_{K,f,g}(x_I)
$$
for $x=(x_I)\in \underset{|I|=m}{\oplus }\WT{\Q}C_*(E_I)$.

\vskip 1pc
\begin{cor}
$$
\varPhi ={\varPhi^{m,n}}:\WT{\Q}C_*(S; E_1, \ldots , E_r)\to 
\WT{\Q}C_{*+1}(X; Y_1, \ldots , Y_r)
$$
is a homotopy from $h^*$ to $f^*g^*$.
\end{cor}

{\it Proof}: 
It follows from Prop.2.17 and Lem.2.11 that 
\begin{align*}
(-1)^n&\partial \varPhi^{m,n}(x)_J+(-1)^m
\varPhi^{m,n}(\partial x)_J \\
=&\sum_{a=1}^{n-m}(-1)^{a+1}
\underset{|L|=a}{\sum_{L\coprod P=J}}\sgn \sbinom{L \ P}{J}
\underset{|I|=m}{\sum_{\PR{L}\coprod I=P}}
\sgn \sbinom{\PR{L} \ I}{P}\varXi_L\varXi_{\PR{L},f,g}(x_I) \\
&+\sum_{a=0}^{n-m}\underset{|L|=a}{\sum_{L\coprod P=J}}
\sgn \sbinom{L \ P}{J}\underset{|I|=m}{\sum_{\PR{L}\coprod I=P}}
\sgn \sbinom{\PR{L} \ I}{P}\varXi_{L,f}\varXi_{\PR{L},g}(x_I)  \\
&+\sum_{a=0}^{n-m-1}(-1)^{a+1}
\underset{|L|=a}{\sum_{L\coprod P=J}}\sgn \sbinom{L \ P}{J}
\underset{|I|=m}{\sum_{\PR{L}\coprod I=P}}
\sgn \sbinom{\PR{L} \ I}{P}\varXi_{L,f,g}\varXi_{\PR{L}}(x) \\
&-\underset{|I|=m}{\sum_{K\coprod I=J}}
\sgn \sbinom{K \ I}{J}\varXi_{K,h}(x)  \\
=&-\sum_{a=1}^{n-m}F^{n-a,n}\varPhi^{m,n-a}(x_I)
+\sum_{a=0}^{n-m}(f^*)^{n-a,n}(g^*)^{m,n-a}(x_I) \\
&-\sum_{a=0}^{n-m-1}\varPhi^{n-a,n}F^{m,n-a}(x_I)-(h^*)^{m,n}(x_I), 
\end{align*}
which completes the proof.
\qed

\vskip 1pc
\subsection{The $\scC$-complex $\WT{\Q}\ALT{*}{X; Y_1, \ldots , Y_r}$}
Consider the sequence of embeddings 
$$
Y_1\cap Y_2\overset{\iota_1}{\to }Y_2\overset{\IId }{\to }Y_2
\overset{\iota_2}{\to }X
$$
and take a vector bundle $\caF$ on $X$.
Then the $2$-cube $(\iota_1, \IId , \iota_2)^*\caF $ on 
$Y_1\cap Y_2$ is described as 
$$
\begin{CD}
\iota_1^*\iota_2^*\caF @>>> \iota_1^*\iota_2^*\caF  \\
@VVV  @VVV  \\
\iota_1^*\iota_2^*\caF @>>> (\iota_2\iota_1)^*\caF , 
\end{CD}
$$
which is not degenerate.
This means that $(\IId^*)^{0,2}$ is not the zero map.
More generally, the pull-back map by the identity morphism 
$$
(\IId_X)^*:\WT{\Q}C_*(X; Y_1, \ldots , Y_r)\to 
\WT{\Q}C_*(X; Y_1, \ldots , Y_r).
$$
is not the identity, hence we can not apply Prop.2.6 and Prop.2.8 
to the multi-relative complexes of exact cubes.
However, since the above cube is symmetric, we can overcome 
this drawback by using the alternating part $\WT{\Q}\ALT{*}{X}$.

Let 
$$
X_0\overset{f_1}{\to }X_1\overset{f_2}{\to }\cdots \overset{f_r}{\to }X_r
$$
be a sequence of morphisms of schemes, and $\caF$ an exact 
$n$-cube on $X_r$.
For any $\sigma \in \frS_n$, define $\PR{\sigma }\in \frS_{n+r-1}$ by 
$$
\PR{\sigma }=\begin{pmatrix}
1 & \cdots & r-1 & r & \cdots & n+r-1 \\
1 & \cdots & r-1 & \sigma (1)+r-1 & \cdots & \sigma (n)+r-1 
\end{pmatrix}.
$$
Then it holds that 
$$
(f_1, \ldots f_r)^*(\sigma (\caF))=
\PR{\sigma }((f_1, \ldots , f_r)^*(\caF)).
$$
This implies that for a division 
$K\coprod I=J$ of $\{1, \ldots , r\}$ we have 
$$
\Alt_*\varXi_K\Alt_*=\Alt_*\varXi_K:
\WT{\Q}C_*(Y_I)\to \WT{\Q}C_{*+n-m-1}(Y_J), 
$$
which leads to the following proposition:

\vskip 1pc
\begin{prop}
If we put 
$$
\varXi^{\alt}_K=\Alt_*\varXi_K:\WT{\Q}\ALT{*}{Y_I}\to 
\WT{\Q}\ALT{*+n-m-1}{Y_J}
$$
and 
$$
F^{m,n}=(-1)^n\underset{|I|=m}{\sum_{K\coprod I=J}}
\sgn \sbinom{K \ I}{J}\varXi^{\alt}_K:
\underset{|I|=m}{\oplus }\WT{\Q}\ALT{*}{Y_I}\to 
\underset{|J|=n}{\oplus }\WT{\Q}\ALT{*+n-m-1}{Y_J}, 
$$
then the family 
$$
\WT{\Q}\ALT{*}{X; Y_1, \ldots , Y_r}=
\left(\underset{|I|=m}{\oplus }
\WT{\Q}\ALT{*}{Y_I}, F^{m,n}\right)
$$
is a $\scC$-complex.
Similarly, for morphisms 
$$
(X; Y_1, \ldots , Y_r)\overset{f}{\to }(T; D_1, \ldots , D_r)
\overset{g}{\to }(S; E_1, \ldots , E_r), 
$$
write $\varXi_{K,f}^{\alt}=\Alt_*\varXi_{K,f}$, 
$\varXi_{K,f,g}^{\alt}=\Alt_*\varXi_{K,f,g}$ and 
\begin{align*}
(f^*)^{m,n}=\underset{|I|=m}{\sum_{K\coprod I=J}}
\sgn \sbinom{K \ I}{J}\varXi_{K,f}^{\alt}:&\ 
\underset{|I|=m}{\oplus }\WT{\Q}\ALT{*}{D_I}\to 
\underset{|J|=n}{\oplus }\WT{\Q}\ALT{*+n-m}{Y_J},  \\
\varPhi^{m,n}=(-1)^n\underset{|I|=m}{\sum_{K\coprod I=J}}
\sgn \sbinom{K \ I}{J}\varXi_{K,f,g}^{\alt}:&\ 
\underset{|I|=m}{\oplus }\WT{\Q}\ALT{*}{E_I}\to 
\underset{|J|=n}{\oplus }\WT{\Q}\ALT{*+n-m+1}{Y_J}.
\end{align*}
Then $f^*=((f^*)^{m,n})$ is a map of $\scC$-complexes and 
$\varPhi=(\varPhi^{m,n})$ is a homotopy from $h^*$ to $f^*g^*$.
\end{prop}

\vskip 1pc
\begin{prop}
For the identity map $\IId_X:X\to X$, 
$$
\IId_X^*: \WT{\Q}\ALT{*}{X; Y_1, \ldots , Y_r}\to 
\WT{\Q}\ALT{*}{X; Y_1, \ldots , Y_r}
$$
is the identity map of the $\scC$-complex.
Hence if we assume that \linebreak 
$(S; E_1, \ldots , E_r)=(X; Y_1, \ldots , Y_r)$ 
and $gf=\IId_X$ in the above notations, then $f^*g^*$ is homotopy 
equivalent to the identity, and a homotopy from the identity to 
$g^*f^*$ is given by 
$$
\varPhi^{m,n}=(-1)^n\underset{|I|=m}{\sum_{K\coprod I=J}}
\sgn \sbinom{K \ I}{J}\varXi_{K,f,g}^{\alt}:
\underset{|I|=m}{\oplus }\WT{\Q}\ALT{*}{Y_I}\to 
\underset{|J|=n}{\oplus }\WT{\Q}\ALT{*+n-m+1}{Y_J}.
$$
\end{prop}

{\it Proof}: 
Suppose $n-m\geq 1$.
It is easy to see that for any exact cube $\caF$, 
$$
(\iota_{k_{\sigma (1)}}, \ldots , \iota_{k_{\sigma(p)}}, \IId_X, 
\iota_{k_{\sigma(p+1)}}, \ldots , \iota_{k_{\sigma (n-m)}})^*\caF 
$$
is invariant under the action of the transposition $(p, p+1)$ if 
$1\leq p\leq n-m-1$, and it is degenerate if $p=0$ or $n-m$.
Hence 
$$
\Alt_*(\iota_{k_{\sigma (1)}}, \ldots , \iota_{k_{\sigma(p)}}, \IId_X, 
\iota_{k_{\sigma(p+1)}}, \ldots , \iota_{k_{\sigma (n-m)}})^*(\caF)=0, 
$$
which completes the proof.
\qed

\setcounter{equation}{0}
\vskip 2pc
\section{The multi-relative $K$-theory and the higher Bott-Chern forms}

\vskip 1pc
\subsection{The cocubical schemes $(\P^1)^*$ and $\PP{*}$}
A {\it cubical complex} is a family of abelian groups 
$\{C_n\}_{n=0,1,\ldots }$ with maps 
\begin{align*}
\partial_i^0, \partial_i^{\infty }:&\ C_n\to C_{n-1}, 
\ 1\leq i\leq n, \\
s_i:&\ C_{n-1}\to C_n, \ 1\leq i\leq n,
\end{align*}
satisfying that 
\begin{align*}
\partial_i^k\partial_j^l=&\partial_{j+1}^k\partial_i^l, 
\qquad \ i<j, \\
\partial_i^ks_j=&\begin{cases}s_{j-1}\partial_i^k, &\ i<j, \\
\IId, &\ i=j,  \\
s_j\partial_{i-1}^k, &\ i>j,  \end{cases}  \\
s_is_j=&s_{j+1}s_i, \qquad \ \ i\leq j
\end{align*}
for $k, l=0$ or $\infty $.
Then the same family of groups $\{C_n\}$ with 
$$
\partial =\sum_{i=1}^n(-1)^i(\partial_i^0-\partial_i^{\infty }):
C_n\to C_{n-1}
$$
forms a chain complex.
Set 
$$
D_n=\sum_{i=1}^n\IIm (s_i:C_{n-1}\to C_n)\subset C_n.
$$
Then $D_*$ is a subcomplex of $C_*$.
An element of $D_*$ is said to be {\it degenerate}.
The subcomplex defined by 
$$
(NC)_n=\underset{j=1}{\overset{n}{\cap }}\KER \partial_j^{\infty }
$$
is called the {\it normalized subcomplex} of $C_*$.
Then there is a direct sum decomposition 
$$
C_*=D_*\oplus NC_*.
$$
Hence, if we write $\WT{C}_*=C_*/D_*$, then there is a canonical 
isomorphism of chain complexes 
$$
NC_*\simeq \WT{C}_*.
$$
When an element $x\in C_*$ is decomposed as $x=x_0+x_d$ 
such that $x_0\in NC_*$ and $x_d\in D_*$, $x_0$ is called 
the {\it normalized component} of $x$, 
and $x$ is said to be {\it normalized} if $x_d=0$.

Let $\P^1$ be the projective line.
Let $z$ be the canonical coordinate of $\P^1$, and $z_j=\pi_j^*z$, 
where $\pi_j:(\P^1)^r\to \P$ is the  $i$-th projection.
Define coface and codegeneracy maps 
\begin{align*}
\delta_j^0, \delta_j^{\infty }&:(\P^1)^r\to (\P^1)^{r+1}, \ 
1\leq j\leq r+1,  \\
\sigma_j&:(\P^1)^r\to (\P^1)^{r-1}, \ 1\leq j\leq r
\end{align*}
as follows:
\begin{align*}
\delta_j^0(z_1, \ldots , z_r)&=
(z_1, \ldots , z_{j-1}, 0, z_j, \ldots , z_r),  \\
\delta_j^{\infty }(z_1, \ldots , z_r)&=
(z_1, \ldots , z_{j-1}, \infty , z_j, \ldots , z_r),  \\
\sigma_j(z_1, \ldots , z_r)&=
(z_1, \ldots , z_{j-1}, z_{j+1}, \ldots , z_r).
\end{align*}
Then $((\P^1)^*, \delta_j^i, \sigma_j)$ is a cocubical scheme.
For any subset $J=\{j_1, \ldots , j_n\}\subset \{1, \ldots , r\}$ 
with $j_1<\cdots <j_n$ and 
for any map $\iota :J\to \{0, \infty \}$, the closed subscheme 
$$
D_{J, \iota }=\IIm \left(\delta_{j_n}^{\iota (j_n)}\cdots 
\delta_{j_1}^{\iota (j_1)}:(\P^1)^{r-n}\to (\P^1)^r\right)
$$
is called a {\it face} of $(\P^1)^r$.
In other words, 
$$
D_{J, \iota }=\{(z_1, \ldots , z_r)\in (\P^1)^r; 
z_{j_k}=\iota (j_k), \ k=1, \ldots , n\}.
$$
Let 
$D_j=\{z_j=0 \ \text{or} \ \infty \}\subset (\P^1)^r$ 
and $D_J=\underset{j\in J}{\cap }D_j$.
Then $D_J$ is a disjoint union of faces of $(\P^1)^r$:
$$
D_J=\underset{\iota:J\to \{0,\infty \}}{\coprod }D_{J,\iota }
\subset (\P^1)^r. 
$$

Let $\square =\P^1-\{1\}$.
Then $\PP{*}\subset (\P^1)^*$ with the same coface and 
codegeneracy maps is also a cocubical scheme.
Denote the faces of $\PP{r}$ by the same symbol:
\begin{align*}
D_{J, \iota }&=\IIm \left(\delta_{j_1}^{\iota (j_1)}\cdots 
\delta_{j_n}^{\iota (j_n)}:\PP{r-n}\to \PP{r}\right), \\
D_J&=\underset{\iota:J\to \{0,\infty \}}{\coprod }D_{J,\iota }
\subset \PP{r}.
\end{align*}

\vskip 1pc
\subsection{Complexes of differential forms}
In this subsection we will introduce several complexes of 
differential forms which we will use throughout the paper.
In what follows, by {\it complex algebraic manifold} 
we mean the complex manifold associated with a smooth algebraic 
variety defined over the complex number field $\C$.
Let $X$ be a complex algebraic manifold.
Denote by $(E_{\log,\R}^*(X),d)$ the complex of real smooth 
differential forms on $X$ with logarithmic singularities along 
infinity, which is defined in \cite{burgos1}.
Then $(E_{\log,\R}^*(X), d)$ with the natural bigrading 
$$
E_{\log,\R}^n(X)\otimes \C=\underset{p+q=n}{\oplus }
E_{\log}^{p,q}(X)
$$
forms a Dolbeault complex.
It is proved in \cite{burgos1} that the complex 
$E_{\log,\R}^*(X)$ with the above bigrading computes the 
cohomology groups of $X$ with the usual Hodge structure.
Denote by $(\scD_{\log}^*(X, p), d_{\scD})$ the associated 
Deligne complex \cite[\S 2]{burgos2}, and by 
$\mmD_{\log}^*(X, p)=\tau_{\leq 2p}\scD_{\log}^*(X, p)$ 
the subcomplex of $\scD_{\log}^*(X, p)$ truncated 
at the degree $2p$.
That is to say, 
$$
\mmD_{\log}^n(X, p)=\begin{cases}
E^{n-1}_{\log,\R}(X)(p-1)\cap 
\underset{p^{\prime }<p, q^{\prime }<p}
{\underset{p^{\prime }+q^{\prime }=n-1}{\oplus }}
E_{\log}^{p^{\prime }, q^{\prime }}(X),  &  n<2p,  \\
E^{2p}_{\log,\R}(X)(p)\cap E_{\log}^{p, p}(X)\cap \KER d, 
&  n=2p, \\
0, & n>2p, 
\end{cases}
$$
where $\R(p)=(2\pi i)^p\R\subset \C$ and $E_{\log,\R}^*(X)(p)$ 
is the space of differential forms on $X$ with 
values in $\R(p)$.
The differential 
$d_{\scD}:\mmD_{\log}^n(X, p)\to \mmD_{\log}^{n+1}(X, p)$ 
is given by 
$$
d_{\scD}(\omega )=\begin{cases}
-\pi (d\omega ),  &  n<2p-1,  \\
-2\partial \overline{\partial }\omega , & n=2p-1,  \\
0, & n>2p-1,  \end{cases}
$$
where $\pi :E_{\log}^n(X)\otimes \C\to \mmD_{\log}^{n+1}(X, p)$ 
is the canonical projection.
Then it is shown in \cite[Cor.2.7]{burgos2} that if $n\leq 2p$, then 
the cohomology groups of $(\mmD_{\log}^*(X, p), d_{\scD})$ are canonically 
\linebreak 
isomorphic to the Deligne cohomology of $X$:
$$
H^n(\mmD_{\log}^*(X, p), d_{\scD})\simeq H_{\scD}^n(X, \R(p)).
$$

Let us recall the multiplicative structure of $\mmD_{\log}^*(X, p)$ 
given in \cite[\S 3]{burgos2}.
Define a product 
$$
\bullet :\mmD_{\log}^n(X, p)\otimes \mmD_{\log}^m(X, q)
\to \mmD_{\log}^{n+m}(X, p+q)
$$
as follows: 
If $n<2p$ and $m<2q$, then 
$$
\omega \bullet \eta =(-1)^n(\partial \omega^{p-1,n-p}-
\ahdif \omega^{n-p,p-1}))\wedge \eta +\omega \wedge 
(\partial \eta^{q-1,m-q}-\ahdif \eta^{m-q,q-1}), 
$$
where $\omega^{p,q}$ is the $(p,q)$-part of the differential form 
$\omega $.
If $n=2p$ or $m=2q$, then define 
$\omega \bullet \eta =\omega \wedge \eta $.
It is shown in \cite[\S 3]{burgos2} that 
$$
d_{\scD}(\omega \bullet \eta )=d_{\scD}\omega \bullet \eta +
(-1)^n\omega \bullet d_{\scD}\eta 
$$
and the induce map 
on cohomology agrees with the product in the Deligne cohomology 
defined in \cite[\S 2]{EV}.

Set 
$$
\scD_{\A}^{n,-r}(X, p)=\mmD^n_{\log}(X\times \PP{r}, p)
$$
and define differentials by 
\begin{gather*}
d_{\scD}:\scD_{\A}^{n,-r}(X, p)\to \scD_{\A}^{n+1,-r}(X, p), \\
\delta_{\A}=\sum_{j=1}^r(-1)^j((\delta_j^0)^*-(\delta_j^{\infty })^*):
\scD_{\A}^{n,-r}(X, p)\to \scD_{\A}^{n,-r+1}(X, p).
\end{gather*}
Then $(\scD_{\A}^{n,-r}(X, p), d_{\scD}, \delta_{\A})$ is 
a double complex.
When we fix the first index $n$, $r\mapsto \scD_{\A}^{n,-r}(X, p)$ 
has a cubical structure.
Denote by $D_{\A}^{n,-*}(X, p)$ the subcomplex of 
degenerate elements with respect to this cubical structure, and set 
$$
\WT{\scD}_{\A}^{n,-*}(X, p)=\scD_{\A}^{n,-*}(X, p)/D_{\A}^{n,-*}(X, p).
$$
Denote by $(\WT{\scD}_{\A}^*(X, p), d_s)$ the single complex 
associated with $\WT{\scD}_{\A}^{*,*}(X, p)$.
Then it is shown in \cite[Prop.2.8]{BF} that the inclusion 
$$
\mmD^*_{\log}(X, p)=\WT{\scD}_{\A}^{*,0}(X, p)\to 
\WT{\scD}_{\A}^*(X, p)
$$
is a quasi-isomorphism.

We introduce another double complex.
Set 
$$
\scD_{\P}^{n,-r}(X, p)=\mmD^n_{\log}(X\times (\P^1)^r, p)
$$
and define differentials by 
\begin{gather*}
d_{\scD}:\scD_{\P}^{n,-r}(X, p)\to \scD_{\P}^{n+1,-r}(X, p), \\
\delta_{\P}=\sum_{j=1}^r(-1)^j((\delta_j^0)^*-(\delta_j^{\infty })^*):
\scD_{\P}^{n,-r}(X, p)\to \scD_{\P}^{n,-r+1}(X, p).
\end{gather*}
Then $(\scD_{\P}^{n,-r}(X, p), d_{\scD}, \delta_{\P})$ is also 
a double complex.
Let $z$ be the canonical coordinate of $\P^1$, and fix a K\"{a}hler 
form $\Omega =\partial \ahdif \log (1+|z|^2)\in \scD_{\P}^{2,-1}(X, 1)$ 
on $\P^1$.
Set 
$$
D_{\P}^{n,-r}(X, p)=\sum_{j=1}^r\left(
\sigma_j^*(\scD_{\P}^{n,-r+1}(X, p))+\pi_j^*\Omega \wedge 
\sigma_j^*(\scD_{\P}^{n-2,-r+1}(X, p-1))\right), 
$$
where $\pi_j:X\times (\P^1)^r\to \P^1$ is the projection to the $j$-th 
component of $(\P^1)^r$, and $\sigma_j$ is the codegeneracy map 
of the cocubical scheme $(\P^1)^*$.
Set 
$$
\WT{\scD}_{\P}^{n,-*}(X, p)=\scD_{\P}^{n,-*}(X, p)/D_{\P}^{n,-*}(X, p), 
$$
and let $(\WT{\scD}_{\P}^*(X, p), d_s)$ be the associated single complex.
Then it is shown in \cite[Prop.1.2]{BW} that the inclusion 
$$
\mmD^*_{\log}(X, p)=\WT{\scD}_{\P}^{*,0}(X, p)\to 
\WT{\scD}_{\P}^*(X, p)
$$
is a quasi-isomorphism.

Finally we introduce a triple complex mixing the above construction.
Set 
$$
\scD_{\A,\P}^{n,-r,-s}(X, p)=
\mmD_{\log}^n(X\times \PP{r}\times (\P^1)^s, p)
$$
with differentials $d_{\scD}, \delta_{\A}, \delta_{\P}$ defined 
in the same way as above.
Define a subcomplex of degenerate elements as follows:
$$
D_{\A,\P}^{n,-r,-s}(X, p)=D_{\A}^{n,-r}(X\times (\P^1)^s, p)+
D_{\P}^{n,-s}(X\times \PP{r}, p)\subset \scD_{\A,\P}^{n,-r,-s}(X, p). 
$$
Set 
$$
\WT{\scD}_{\A,\P}^{n,-r,-s}(X, p)=
\scD_{\A,\P}^{n,-r,-s}(X, p)/D_{\A,\P}^{n,-r,-s}(X, p), 
$$
and let $(\WT{\scD}_{\A,\P}^*(X, p), d_s)$ be the associated 
single complex.
Then the inclusion 
$$
\mmD^*_{\log}(X, p)=\WT{\scD}_{\A,\P}^{*,0,0}(X, p)\to 
\WT{\scD}_{\A,\P}^*(X, p)
$$
is also a quasi-isomorphism.
To sum up, we obtain the following commutative diagram of 
quasi-isomorphisms:
$$
\begin{CD}
\mmD_{\log}^*(X, p) @>>> \WT{\scD}_{\A}^*(X, p)  \\
@VVV  @VVV   \\
\WT{\scD}_{\P}^*(X, p) @>>> \WT{\scD}_{\A,\P}^*(X, p).
\end{CD}
$$

We next consider $\mmD^*_{\log}(X, p)$-module structures on 
the above mentioned complexes.
For $\omega \in \mmD^n_{\log}(X, p)$ and 
$\eta \in \scD_{\A}^{m,-r}(X, q)$, define a product by 
$$
\omega \PRA \eta =\pi_X^*\omega \bullet \eta 
\in \scD_{\A}^{n+m,-r}(X, p+q), 
$$
where $\pi_X:X\times \PP{r}\to X$ is the projection.
Then it induces a map 
$$
\PRA :\mmD^*_{\log}(X, p)\times \WT{\scD}_{\A}^*(X, q)
\to \WT{\scD}_{\A}^*(X, p+q)
$$
which satisfies the relation 
$$
d_s(\omega \PRA \eta )=d_{\scD}\omega \PRA \eta +
(-1)^{\deg \omega }\omega \PRA d_s\eta .
$$
Similarly, we can define products 
\begin{align*}
\PRP :&\ \mmD^*_{\log}(X, p)\times \WT{\scD}_{\P}^*(X, q)
\to \WT{\scD}_{\P}^*(X, p+q),   \\ 
\PRAP :&\ \mmD^*_{\log}(X, p)\times \WT{\scD}_{\A,\P}^*(X, q)
\to \WT{\scD}_{\A,\P}^*(X, p+q), 
\end{align*}
which satisfy the same relation.

\vskip 1pc
\subsection{Complexes of currents}
Let $X$ be a complex algebraic manifold of dimension $d_X$.
Denote by $E_{\R}^n(X)_{\CPT}$ the space of real smooth 
differential forms on $X$ of degree $n$ with compact support, 
and by $D_{\R}^n(X)$ the topological dual of 
$E_{\R}^{2d_X-n}(X)_{\CPT}(d_X)$.
Define $d:D_{\R}^n(X)\to D_{\R}^{n+1}(X)$ by 
$dT(\omega )=(-1)^nT(d\omega )$ for 
$\omega \in E_{\R}^{2d_X-n-1}(X)_{\CPT}(d_X)$.
Denote by $E^{p,q}(X)_{\CPT}$ the space of $(p,q)$-forms 
on $X$ with compact support, and by $D^{p,q}(X)$ the topological dual 
of $E^{d_X-p,d_X-q}(X)_{\CPT}$.
Then $(D_{\R}^*(X), d)$ with the bigrading 
$$
D_{\R}^n(X)\otimes \C=\underset{p+q=n}{\oplus }D^{p.q}(X)
$$
forms a Dolbeault complex as well.
Denote by $\mmD_D^*(X, p)$ the associated Deligne complex 
truncated at $2p$.

Let 
$$
[ \quad ]:E_{\R}^n(X)\to D_{\R}^n(X)
$$
be the map given by the integral 
$$
[\eta ](\omega )=\CST{1}{d_X}\int_X\omega \wedge \eta 
$$
for $\omega \in E_{\R}^{2d_X-n}(X)_{\CPT}$.
It satisfies $d[\eta ]=[d\eta ]$, and gives a quasi-isomorphism 
of the Dolbeault complexes.
Hence it induces a quasi-isomorphism of the associated 
Deligne complexes:
$$
[\quad ]:\mmD^*(X, p)\to \mmD_D^*(X, p).
$$

For an integral subvariety $V\subset X$ of codimension $p$, 
let $\delta_V\in \mmD_D^{2p}(X, p)$ be the current given by the integral 
$$
\delta_V(\omega )=\frac{1}{(2\pi i)^{d_X-p}}\int_{\WT{V}}\iota^*\omega 
$$
for $\omega \in E_{\R}^{2d_X-2p}(X)_{\CPT}(d_X-p)$, 
where $\iota :\WT{V}\to V\subset X$ is a disingularization of $V$.
We can extend this construction linearly to any cycle $z$ 
of codimension $p$ and we can define a current 
$\delta_z\in \mmD_D^{2p}(X, p)$.

Let $Y$ be a complex algebraic manifold of dimension $d_Y$, 
and $f:X\to Y$ a proper morphism.
We can define the {\it direct image map} 
$$
f_*:D_{\R}^*(X)\to D_{\R}^*(Y)(d_Y-d_X)[2d_Y-2d_X]
$$
by $f_*(T)(\omega )=T(f^*\omega )$ for $T\in D_{\R}^*(X)$ and 
$\omega \in E_{\R}^*(Y)_{\CPT}$.
It induces a map of the Deligne complexes 
$$
f_*:\mmD_D^*(X, p)\to \mmD_D^*(Y, p+d_Y-d_X)[2d_Y-2d_X].
$$

\vskip 1pc
\subsection{Wang's forms}
In this subsection we will introduce Wang's forms \cite{BW}.
Let $z$ be the canonical coordinate of the projective line $\P^1$ 
and $z_i=\pi_i^*z$, where $\pi_i:(\P^1)^r\to \P^1$ is the $i$-th 
projection.
For $r\leq 1$, define a differential form $W_r$ on $(\P^1)^r$ by 
$$
W_r=\frac{1}{2r!}\sum_{i=1}^r(-1)^iS_r^i, 
$$
where 
$$
S_r^i=\sum_{\sigma \in \frS_r}(\sgn \sigma )\log |z_{\sigma (1)}|^2
\frac{dz_{\sigma (2)}}{z_{\sigma (2)}}\wedge \cdots \wedge 
\frac{dz_{\sigma (i)}}{z_{\sigma (i)}}\wedge 
\frac{d\OV{z}_{\sigma (i+1)}}{\OV{z}_{\sigma (i+1)}}\wedge \cdots 
\wedge \frac{d\OV{z}_{\sigma (r)}}{\OV{z}_{\sigma (r)}}.
$$
The form $W_r$ has logarithmic singularities along 
the codimension one faces of $(\P^1)^r$, and it is locally 
integrable on $(\P^1)^r$ such that $\OV{W}_r=(-1)^{r-1}W_r$.
Hence as a current, $[W_r]\in \mmD_D^r((\P^1)^r, r)$.
When $r=0$, suppose that $W_0=1$.

\vskip 1pc
\begin{prop}\cite[Thm.6.7]{BFT}
When $r\geq 1$, the current $[W_r]$ satisfies the relation 
$$
d_{\scD}[W_r]=\sum_{j=1}^r(-1)^j\left((\delta_j^0)_*
[W_{r-1}]-(\delta_j^{\infty })_*[W_{r-1}]\right).
$$
\end{prop}

\vskip 1pc
We now construct several maps between the Deligne complexes 
given in \S 3.2 and \S 3.3 in the case that $X$ is compact.
Let 
\begin{align*}
\pi_{\P}:&\ X\times (\P^1)^r\to (\P^1)^r,  \\
\pi_X:&\ X\times (\P^1)^r\to X
\end{align*}
be the projections.
For any $\omega \in \scD_{\P}^{n,-r}(X, p)=\mmD^n(X\times (\P^1)^r, p)$, 
consider the integral along fibers 
$$
\frac{1}{(2\pi i)^r}\int_{(\P^1)^r}\omega \bullet \pi_{\P}^*W_r\in 
\mmD^{n-r}(X, p).
$$
We can show in the same way as \cite[Lem.6.8]{BW} that the above 
integral is zero if $\omega \in D_{\P}^{n,-r}(X, p)$.

\vskip 1pc
\begin{prop}
The map 
$$
\kappa_{\P}:\WT{\scD}_{\P}^*(X, p)\to \mmD^{n-r}(X, p)
$$
induced by the above integral is a map of complexes.
\end{prop}

{\it Proof}: 
For $\omega \in \scD_{\P}^{n,-r}(X, p)$, 
\begin{align*}
d_{\scD}\kappa_{\P}(\omega )=&\frac{1}{(2\pi i)^r}d_{\scD}\left(
\int_{(\P^1)^r}\omega \bullet \pi_{\P}^*W_r\right)  \\
=&\frac{1}{(2\pi i)^r}\int_{(\P^1)^r}(d_{\scD}-d_{\scD,\P})
(\omega \bullet \pi_{\P}^*W_r), 
\end{align*}
where $d_{\scD,\P}$ is the differential of 
$\mmD^*(X\times (\P^1)^r, p)$ on the component $(\P^1)^r$.
Since $d_{\scD}W_r=0$ as a differential form, 
$$
\frac{1}{(2\pi i)^r}\int_{(\P^1)^r}d_{\scD}(\omega \bullet \pi_{\P}^*W_r)
=\frac{1}{(2\pi i)^r}\int_{(\P^1)^r}d_{\scD}\omega \bullet \pi_{\P}^*W_r.
$$
On the other hand, it follows from Prop.3.1 that 
$$
\frac{1}{(2\pi i)^r}\int_{(\P^1)^r}
d_{\scD,\P}(\omega \bullet \pi_{\P}^*W_r)=
\frac{(-1)^{n-1}}{(2\pi i)^{r-1}}\int_{(\P^1)^{r-1}}
(\delta_{\P}\omega )\bullet \pi_{\P}^*W_{r-1}.
$$ 
Hence 
\begin{align*}
d_{\scD}\kappa_{\P}(\omega )=&\ \frac{1}{(2\pi i)^r}\int_{(\P^1)^r}
d_{\scD}\omega \bullet \pi_{\P}^*W_r+\frac{(-1)^n}{(2\pi i)^{r-1}}
\int_{(\P^1)^{r-1}}(\delta_{\P}\omega )\bullet \pi_{\P}^*W_{r-1}  \\
=&\ \kappa_{\P}(d_{\scD}\omega )+(-1)^n\kappa_{\P}(\delta_{\P}\omega ), 
\end{align*}
which completes the proof.
\qed

\vskip 1pc
It is easy to see that the map $\kappa_{\P}$ is a left inverse of 
the inclusion $\mmD^*(X, p)\to \WT{\scD}_{\P}^*(X, p)$.
In particular, $\kappa_{\P}$ is a quasi-isomorphism.

Let us next consider a similar map on the complex 
$\WT{\scD}_{\A}^*(X, p)$.
In this case, we can not take the integration along fibers.
However, it is shown in \cite[Prop.6.5]{BFT} that for any 
$\omega \in \scD_{\A}^{n,-r}(X, p)=\mmD_{\log}^n(X\times \PP{r}, p)$ 
the product $\omega \bullet \pi_{\P}^*W_r$ is locally 
integrable on the compactification $X\times (\P^1)^r$ of 
$X\times \PP{r}$, therefore we can define the current 
$$
{\pi_X}_*[\omega \bullet \pi_{\P}^*W_r]\in 
\mmD_D^{n-r}(X, p).
$$
Moreover, \cite[Prop.6.5]{BFT} says that 
$$
d_{\scD}{\pi_X}_*[\omega \bullet \pi_{\P}^*W_r]=
{\pi_X}_*[d_{\scD}\omega \bullet \pi_{\P}^*W_r]+
(-1)^n{\pi_X}_*[\delta_{\A}\omega \bullet \pi_{\P}^*W_{r-1}].
$$
Since $\pi_{X,*}[\omega \bullet \pi_{\P}^*W_r]=0$ if 
$\omega \in D_{\A}^{n,-r}(X, p)$, we conclude that 
the above current gives a map of complexes 
$$
\kappa_{\A}:\WT{\scD}_{\A}^*(X, p)\to \mmD_D^*(X, p), 
$$
which make the diagram 
$$
\xymatrix{
\mmD^*(X, p)\ar[rr] \ar[rd]_{\scriptscriptstyle{[ \ \ ]}} & & 
\WT{\scD}_{\A}^*(X, p)\ar[ld]^{\kappa_{\A}} \\
& \mmD_D^*(X, p)
}
$$
commutative.
In particular, $\kappa_{\A}$ is also a quasi-isomorphism.

Finally let us consider a map from the complex 
$\WT{\scD}_{\A,\P}^*(X, p)$.
In the same way as the previous case we can show that for any 
$\omega \in \scD_{\A}^{n,-r,-s}(X, p)=
\mmD_{\log}^n(X\times \PP{r}\times (\P^1)^s, p)$, 
the product $\omega \bullet \pi_{\P}^*W_{r+s}$ is locally 
integrable on the compactification $X\times (\P^1)^{r+s}$ of 
$X\times \PP{r}\times (\P^1)^s$, and that the current 
${\pi_X}_*[\omega \bullet \pi_{\P}^*W_{r+s}]\in 
\mmD_D^{n-r-s}(X, p)$ gives a quasi-isomorphism 
$$
\kappa_{\A,\P}:\WT{\scD}_{\A}^*(X, p)\to \mmD_D^*(X, p).
$$

Summing up the results in this subsection, 
we obtain the commutative diagram 
$$
\xymatrix{
\WT{\scD}_{\A}^*(X, p)\ar[r]\ar[d]_{\kappa_{\A}} & 
\WT{\scD}_{\A,\P}^*(X, p) \ar[dl]^{\kappa_{\A,\P}} & 
\WT{\scD}_{\P}^*(X, p) \ar[l]\ar[d]^{\kappa_{\P}} \\
\mmD_D^*(X, p) & & \mmD^*(X, p), 
\ar[ll]_{\scriptscriptstyle{[ \ \ ]}}
}
$$
all the maps in which are quasi-isomorphisms.

\vskip 1pc
\subsection{The higher Bott-Chern forms}
In this subsection we will recall the higher Bott-Chern forms 
\cite{BW}.
Let $X$ be a complex algebraic manifold and $\caF$ a vector 
bundle on $X$.
Given a smooth hermitian metric $h$ on $\caF$, there exists 
a unique connection on $\caF$ which is compatible both with 
the complex structure and with the metric.
Using the curvature form of this connection, we obtain 
a differential form 
$$
\CH_0(\caF, h)\in \underset{p}{\oplus }E_{\R}^{2p}(X)(p)\cap 
E^{p,p}(X)\cap \KER d
$$
which represents the Chern character of $\caF$.
This from is called the {\it Chern form} of $(\caF, h)$.

A smooth hermitian metric $h$ on $\caF$ is said to be 
{\it smooth at infinity} if there is a vector bundle $\PR{\caF}$ 
on a smooth compactification $\OV{X}$ of $X$ with a smooth 
hermitian metric $\PR{h}$ such that the restriction 
$(\PR{\caF}|_X, \PR{h}|_X)$ is isometric to $(\caF, h)$.
In this case the Chern form $\CH_0(\OV{\caF})$ can be extended to 
a smooth differential form on $\OV{X}$.
In particular, 
$$
\CH_0(\OV{\caF})\in \underset{p}{\oplus }\mmD_{\log}^{2p}(X, p).
$$
In what follows, any hermitian metric on a vector bundle is 
supposed to be smooth at infinity.

An {\it exact hermitian $n$-cube} on $X$ is an exact $n$-cube 
consisting of vector bundles on $X$ with smooth hermitian metrics.
Denote by ${\Q}\WH{C}_*(X)$ the chain complex of exact hermitian 
cubes on $X$ with smooth at infinity metrics, and 
by $\WT{\Q}\WH{C}_*(X)$ the quotient complex of $\WT{\Q}\WH{C}_*(X)$ 
by the subcomplex of degenerate cubes.
Denote by $\WT{\Q}\WALT{*}{X}$ the alternating part 
of $\WT{\Q}\WH{C}_*(X)$.
Since the space of smooth at infinity metrics on 
a vector bundle is convex, the map 
$\WT{\Q}\WH{C}_*(X)\to \WT{\Q}C_*(X)$ forgetting metrics is 
a quasi-isomorphism.
Hence Thm.2.9 implies the following:

\vskip 1pc
\begin{thm}
The homology groups of $\WT{\Q}\WH{C}_*(X)$ and $\WT{\Q}\WALT{*}{X}$ 
are canonically isomorphic to the rational algebraic $K$-theory 
of $X$:
$$
H_n(\WT{\Q}\WALT{*}{X})\simeq H_n(\WT{\Q}\WH{C}_*(X))
\simeq \QQ{K_n(X)}.
$$
\end{thm}

\vskip 1pc
An exact hermitian $n$-cube $\OV{\caF}$ on $X$ is {\it emi} 
if for any $1\leq j\leq n$, the metric on 
$\partial_j^{-1}\OV{\caF}$ is induced from the metric on 
$\partial_j^0\OV{\caF}$ by means of the inclusion 
$\partial_j^{-1}\OV{\caF}\hookrightarrow 
\partial_j^0\OV{\caF}$.
Let $\WT{\Q}\WH{C}^{emi}_*(X)$ denote the subcomplex of 
$\WT{\Q}\WH{C}_*(X)$ consisting of emi-cubes.
As seen in \cite[\S 3]{BW}, we can obtain a map of complexes 
$$
\lambda :\WT{\Q}\WH{C}_*(X)\to \WT{\Q}\WH{C}^{emi}_*(X)
$$
such that the composite of $\lambda $ with the inclusion 
$\WT{\Q}\WH{C}^{emi}_*(X)\hookrightarrow \WT{\Q}\WH{C}_*(X)$ 
is homotopy equivalent to the identity.

With any emi-$n$-cube $\OV{\caF}$ on $X$ we can associate 
a hermitian vector bundle $\TTR_n(\OV{\caF})$ on 
$X\times (\P^1)^n$ such that there are isometries 
\begin{equation}
\left\{
\begin{aligned}
\TTR_n(\OV{\caF})|_{\{z_j=0\}}&\simeq 
\TTR_{n-1}(\partial_j^0\OV{\caF}), \\
\TTR_n(\OV{\caF})|_{\{z_j=\infty \}}&\simeq 
\TTR_{n-1}(\partial_j^{-1}\OV{\caF})\oplus 
\TTR_{n-1}(\partial_j^1\OV{\caF})
\end{aligned}\right. \label{sec3:iso1}
\end{equation}
for $1\leq j\leq n$.
The hermitian vector bundle $\TTR_n(\OV{\caF})$ is called 
{\it transgression bundle} of $\OV{\caF}$ \cite[Def.3.8]{BW}.
Note that the metric on $\TTR_n(\OV{\caF})$ is smooth at infinity 
if each metric on $\OV{\caF}$ is smooth at infinity, as 
mentioned in \cite[Def.3.8]{BW}.
Define the {\it Bott-Chern form} of an exact hermitian $n$-cube 
$\OV{\caF}$ on $X$ to be the Chern form of 
$\TRNS{n}{\OV{\caF}}$: 
$$
\CH_n(\OV{\caF})_{\P}=\CH_0(\TRNS{n}{\OV{\caF}})\in 
\underset{p}{\oplus }\WT{\scD}_{\P}^{2p-n}(X, p).
$$
It follows from the isometries \eqref{sec3:iso1} that 
$\delta_{\P}\CH_n(\OV{\caF})_{\P}=\CH_{n-1}(\partial \OV{\caF})_{\P}$.
Since $\CH_n(\OV{\caF})_{\P}=0$ if $\OV{\caF}$ is degenerate 
\cite[Prop.3.11]{BW}, it induces a map of complexes 
$$
\CH_{*,\P}:\WT{\Q}\WH{C}_*(X)\to 
\underset{p}{\oplus }\WT{\scD}_{\P}^{2p-*}(X, p).
$$
In the case that $X$ is compact, we can obtain 
$$
\CH_n(\OV{\caF})=\kappa_{\P}(\CH_n(\OV{\caF})_{\P})\in 
\underset{p}{\oplus }\mmD^{2p-n}(X, p).
$$
Note that in the case that $n=1$, 
$\WT{\CH}_1(\OV{\caF})\in \underset{p}{\oplus }
\mmD^{2p-n}(X, p)/\IIm d_{\scD}$ 
agrees with the Bott-Chern secondary class of a short exact sequence 
$\OV{\caF}$ defined by Gillet and Soul\'{e} in \cite{gilletsoule2}.

\vskip 1pc
\begin{thm}\cite[Thm.5.2]{BW}
When $X$ is compact, the map on homology induced by 
the Bott-Chern forms 
$$
\CH_n=\underset{p}{\oplus }\CH_n^p\QQ{K_n(X)}\to 
\underset{p}{\oplus }H_{\scD}^{2p-n}(X, \R(p))
$$
agrees with Beilinson's regulator.
\end{thm}

\vskip 1pc
{\it Remark}: 
The target of Beilinson's regulator is the absolute Hodge cohomology, 
not the Deligne cohomology.
Hence in \cite{BW} the higher Bott-Chern forms sit in 
a complex which computes the absolute Hodge cohomology.
However, since the both cohomology theories are canonically 
isomorphic for compact manifolds, the above theorem follows.

\vskip 1pc
\subsection{A multi-relative complex of exact hermitian cubes and 
the higher Bott-Chern forms}
We begin by introducing a metrized version of the $\scC$-complex 
of exact cubes defined in \S 2.4.
Let $X$ be a complex algebraic manifold and $Y_1, \ldots , Y_r$ 
closed submanifolds of $X$.
Then the same construction as in \S 2.4 and \S 2.5 gives 
a $\scC$-complex 
$$
\WT{\Q}\WH{C}_*(X; Y_1, \ldots , Y_r)=\left(
\underset{|I|=m}{\oplus }\WT{\Q}\WH{C}_*(Y_I), F^{m,n}\right)
$$
and 
$$
\WT{\Q}\WALT{*}{X; Y_1, \ldots , Y_r}=\left(
\underset{|I|=m}{\oplus }\WT{\Q}\WALT{*}{Y_I}, F^{m,n}\right).
$$
It follows from Cor.2.16 and Thm.3.3 that there are canonical isomorphisms 
$$
H_n(\WT{\Q}\WALT{*}{X; Y_1, \ldots , Y_r})\simeq 
H_n(\WT{\Q}\WH{C}_*(X; Y_1, \ldots , Y_r))\simeq 
\QQ{K_n(X; Y_1, \ldots , Y_r)}.
$$

Let us recall the notations introduced in \S 3.1.
Consider the product $X\times \PP{r}$ with the normal 
crossing divisor 
$$
X\times \partial \PP{r}=X\times D_1+\cdots +X\times D_r, 
$$
where $D_j=\{z_j=0, \infty \}\subset \PP{r}$.
We identify $X$ with $X\times \{(\infty , \ldots , \infty )\}$, 
which is a connected component  of $X\times D_{\{1, \ldots , r\}}$.
This gives an embedding of chain complexes 
\begin{equation}
i_X:\WT{\Q}\WH{C}_*(X)\hookrightarrow 
\WT{\Q}\WH{C}_*(X\times \PP{r}; X\times \partial \PP{r})[r]. 
\label{sec3:map1}
\end{equation}
The homotopy invariant property of $K$-theory implies that 
the map \eqref{sec3:map1} induces an isomorphism of $K$-groups:
$$
\QQ{K_{n+r}(X)}\simeq 
\QQ{K_n(X\times \square^r; X\times \partial \square^r)}.
$$

Let $\OV{\caF}_I$ be an exact hermitian $n$-cube on $X\times D_I$.
In other words, $\OV{\caF}_I$ is a family 
$\{\OV{\caF}_{I,\iota }\}_{\iota:I\to \{0, \infty \}}$ 
such that $\OV{\caF}_{I,\iota }$ is an exact hermitian $n$-cube on 
$X\times D_{I,\iota }$ with a smooth at infinity metric.
We identify $D_{I,\iota }$ with $\PP{r-|I|}$, and let 
$$
\CH_n(\OV{\caF}_I)_{\P}
=\sum_{\iota :I\to \{0,\infty \}}
(-1)^{|\iota |}\CH_n(\OV{\caF}_{I,\iota })_{\P}\in 
\underset{p}{\oplus }\WT{\scD}_{\P}^{2p-n}(X\times \PP{r-|I|}, p), 
$$
where $|\iota |$ is the cardinarity of the set 
$\{i\in I; \iota (i)=\infty \}$.
Then it gives a map 
$$
\CH_{n,\P}:\WT{\Q}\WH{C}_n(X\times D_I)\to 
\underset{p}{\oplus }\WT{\scD}_{\A,\P}^{2p-r+|I|-n}(X, p).
$$
Since 
$\delta_{\P}\CH_n(\OV{\caF}_I)_{\P}
=\CH_{n-1}(\partial \OV{\caF}_I)_{\P}$, 
we have the following:

\vskip 1pc
\begin{prop}
For any $x_I\in \WT{\Q}\WH{C}_n(X\times D_I)$, we have 
$$
d_s\CH_n(x_I)_{\P}=\sum_{l=1}^{r-|I|}(-1)^l
\CH_n(x_I|_{X\times D_{I_l}})_{\P}+(-1)^{r-|I|} 
\CH_{n-1}(\partial x_I)_{\P}
$$
in $\WT{\scD}_{\A,\P}^*(X, p)$, where $\{i_1, \ldots , i_{r-|I|}\}$ 
is the complement of $I$ with $i_1<\cdots <i_{r-|I|}$ 
and $I_l=I\cup \{i_l\}$.
\end{prop}

\vskip 1pc
\begin{defn}
An element $x\in \WT{\Q}\WH{C}_n(X)$ is said to be 
{\it isometrically equivalent to a degenerate element} 
if there is a lift 
$$
\sum_ir_i[\OV{\caF}_i]\in \Q\WH{C}_n(X)
$$
of $x$ such that each $\OV{\caF}_i$ is isometric to 
a degenerate cube.
\end{defn}

\vskip 1pc
It is obvious from the definition of the map $\lambda $ 
in \cite[\S 3]{BW} that $\lambda \OV{\caF}$ is isometrically 
equivalent to a degenerate element if so is $\OV{\caF}$.
Furthermore, it is obvious from the definition of 
the transgression bundle that an isometry 
$\OV{\caF}\simeq \OV{\caG}$ of emi-$n$-cubes 
on $X$ induces an isometry $\TTR_n\OV{\caF}\simeq \TTR_n\OV{\caG}$ 
of hermitian vector bundles.
Hence if $x\in \WT{\Q}\WH{C}_n(X)$ is isometrically equivalent to 
a degenerate element, then $\CH_n(x)_{\P}=0$.

Consider a sequence of morphisms of complex algebraic manifolds: 
$$
X_0\overset{f_1}{\to }X_1\overset{f_2}{\to }\cdots 
\overset{f_r}{\to }X_r.
$$
Then for a hermitian vector bundle $\OV{\caF}$ on $X_r$, 
$(f_1, \ldots , f_r)^*\OV{\caF}$ is isometrically 
equivalent to a degenerate element if $r\geq 2$, since all 
the maps in $(f_1, \ldots , f_r)^*\OV{\caF}$ are isometries or 
the zero maps.
More generally, for an exact hermitian cube $\OV{\caF}$ on $X_r$, 
$(f_1, \ldots , f_r)^*\OV{\caF}$ is isometrically equivalent 
to a degenerate element if $r\geq 2$.
Hence we have the following:

\vskip 1pc
\begin{prop}
Consider the maps of complex algebraic manifolds 
with closed submanifolds 
$$
(X; Y_1, \ldots , Y_r)\overset{f}{\to }(T, D_1, \ldots , D_r)\overset{g}{\to }
(S; E_1, \ldots , E_r).
$$
Take $x=(x_I)\in \underset{|I|=m}{\oplus }\WT{\Q}\WH{C}_*(E_I)$ 
and $J\subset \{1, \ldots , r\}$ with $|J|=n$.
Then $F^{m,n}(x)_J$ for $n-m\geq 2$, $g^{m,n}(x)_J$ for 
$n-m\geq 1$, and $\varPhi^{m,n}(x)_J$ for any $m$ and $n$ 
are isometrically equivalent to degenerate elements.
Hence 
\begin{align*}
\CH_*(F^{m,n}(x)_J)_{\P}=0 &\ \ \text{if} \ n-m\geq 2,   \\
\CH_*(g^{m,n}(x)_J)_{\P}=0 &\ \ \text{if} \ n-m\geq 1,   \\
\CH_*(\varPhi^{m,n}(x)_J)_{\P}=0 &\ \ \text{for any $m$ and $n$}.
\end{align*}
\end{prop}

\vskip 1pc
\begin{prop}
For 
$$
x=(x_I)\in \WT{\Q}\WH{C}_n(X\times \PP{r}; 
X\times \partial \PP{r})=\underset{I}{\oplus }
\WT{\Q}\WH{C}_{n+|I|}(X\times D_I),
$$
let 
$$
\CH_{n}(x)=\sum_I
(-1)^{\frac{1}{2}|I|(|I|+1)+r|I|+\Sigma I}\CH_{n+|I|}(x_I)_{\P}
\in \underset{p}{\oplus }\WT{\scD}_{\A,\P}^{2p-r-n}(X, p), 
$$
where $|I|$ is the cardinarity of $I$ and $\Sigma I$ is 
the sum of elements of $I$.
Then 
$$
\CH_*:\WT{\Q}\WH{C}_*(X\times \PP{r}; X\times \partial \PP{r})[r]
\to \underset{p}{\oplus }\WT{\scD}_{\A,\P}^{2p-*}(X, p)
$$
is a map of chain complexes.
Hence it induces a map 
$$
\CH_n=\underset{p}{\oplus }\CH_n^p:
\QQ{K_n(X\times \PP{r}; X\times \partial \PP{r})}
\to \underset{p}{\oplus }H_{\scD}^{2p-n-r}(X, \R(p)).
$$
\end{prop}

{\it Proof}: 
Let $\{i_1, \ldots , i_{r-|I|}\}$ be the complement of $I$ with 
$i_1<\cdots <i_{r-|I|}$ and $I_l=I\cup \{i_l\}$.
Then Prop.3.5 implies that 
\begin{align*}
d_s\CH_n(x)=&
\sum_I\sum_{l=1}^{r-|I|}(-1)^{\frac{1}{2}|I|(|I|+1)+
r|I|+\Sigma I+l}\CH_{n+|I|}(x_I|_{X\times D_{I_l}})_{\P}  \\
&+\sum_I(-1)^{\frac{1}{2}|I|(|I|+1)+r|I|+\Sigma I+r-|I|}
\CH_{n+|I|-1}(\partial x_I)_{\P}.
\end{align*}
If we write $J=I_l=\{j_1, \ldots , j_{|J|}\}$ with 
$j_1<\ldots <j_{|J|}$ and $i_l=j_k$, 
then $i_l=j_k=l+k-1$ and $\Sigma I+l=\Sigma J-k+1$.
Hence 
\begin{align*}
d_s\CH_n(x)=&\sum_J\sum_{k=1}^{|J|}
(-1)^{\frac{1}{2}(|J|-1)|J|+r|J|-r+\Sigma J+k+1}
\CH_{n+|J|-1}(x_{J-\{j_k\}}|_{X\times D_J})_{\P} \\
&+\sum_I(-1)^{\frac{1}{2}|I|(|I|+1)+r|I|+\Sigma I+r-|I|}
\CH_{n+|I|-1}(\partial x_I)_{\P} \\
=&(-1)^r\sum_J(-1)^{\frac{1}{2}|J|(|J|+1)+r|J|+\Sigma J}\times \\
&\left((-1)^{|J|}\CH_{n+|J|-1}(\partial x_J)_{\P}+
(-1)^{|J|}\sum_{k=1}^{|J|}(-1)^{k-1}
\CH_{n+|J|-1}(x_{J-\{j_k\}}|_{X\times D_J})_{\P}\right).
\end{align*}
On the other hand, the definition of the boundary map of 
the total chain complex of the $\scC$-complex 
$\WT{\Q}\WH{C}_*(X\times \PP{r}; X\times \partial \PP{r})$ 
together with Prop.3.7 implies that 
\begin{align*}
\CH_{n+|J|-1}&\, ((\partial x)_J)_{\P}  \\
=&\, (-1)^{|J|}
\CH_{n+|J|-1}(\partial x_J)_{\P}+(-1)^{|J|}\sum_{k=1}^{|J|}
(-1)^{k-1}\CH_{n+|J|-1}(x_{J-\{j_k\}}|_{X\times D_J})_{\P}.
\end{align*}
Hence $(-1)^rd_s\CH_n(x)=\CH_{n-1}(\partial x)$, 
which completes the proof.
\qed

\vskip 1pc
Let $f:X\to Y$ be a morphism of proper complex algebraic manifolds.
It follows from Prop.3.7 that  that the diagram 
$$
\begin{CD}
\WT{\Q}\WH{C}_*(Y\times \PP{r}; Y\times \partial \PP{r})[r] @>{\CH_*}>> 
\underset{p}{\oplus }\WT{\scD}_{\A,\P}^{2p-*}(Y, p)   \\
@V{f^*}VV  @VV{f^*}V  \\
\WT{\Q}\WH{C}_*(X\times \PP{r}; X\times \partial \PP{r})[r] @>{\CH_*}>> 
\underset{p}{\oplus }\WT{\scD}_{\A,\P}^{2p-*}(X, p)
\end{CD}
$$
is commutitive. 
Hence the diagram 
$$
\begin{CD}
\QQ{K_n(Y\times \PP{r}; Y\times \partial \PP{r})} @>{\CH_n}>> 
\underset{p}{\oplus }H_{\scD}^{2p-n-r}(Y, \R(p))   \\
@V{f^*}VV  @VV{f^*}V  \\
\QQ{K_n(X\times \PP{r}; X\times \partial \PP{r})} @>{\CH_n}>> 
\underset{p}{\oplus }H_{\scD}^{2p-n-r}(X, \R(p))
\end{CD}
$$
is also commutative.

The proposition below can be easily verified.

\vskip 1pc
\begin{prop}
The diagram 
$$
\begin{CD}
\WT{\Q}\WH{C}_*(X) @>{i_X}>> 
\WT{\Q}\WH{C}_*(X\times \PP{r}; X\times \partial \PP{r})[r]  \\
@V{\CH_{*,\P}}VV  @VV{\CH_*}V \\
\underset{p}{\oplus }\WT{\scD}_{\P}^{2p-*}(X, p) @>>> 
\underset{p}{\oplus }\WT{\scD}_{\A,\P}^{2p-*}(X, p)
\end{CD}
$$
is commutative.
In particular, the diagram 
$$
\xymatrix{
 \QQ{K_{n+r}(X)}\ar[rr] \ar[rd]_{\CH_{n+r}^p} & & 
 \QQ{K_n(X\times \PP{r}; X\times \partial \PP{r})} 
 \ar[ld]^{\CH_n^p} \\
 & H^{2p-n-r}_{\scD }(X, \R(p)) & 
 }
$$
is commutative.
\end{prop}

\setcounter{equation}{0}
\vskip 2pc
\section{Chern form of a hermitian vector bundle 
on an iterated double}
\vskip 1pc

\subsection{Hermitian vector bundles on an iterated double}
In this subsection we will introduce a scheme called 
{\it iterated double} and construct a theory of Chern forms of 
hermitian vector bundles on it.
Let $X$ be a scheme and $Y_1, \ldots , Y_r$ closed subschemes 
of $X$.
Denote by $D(X; Y_1, \ldots , Y_r)$ the iterated double 
defined by Levine in \cite{levine}.
As a topological space, it is a union of $2^r$ copies of $X$  
indexed by the set of all subsets of $\{1, \ldots , r\}$.
To be more precise, if we denote by $X_I$ the closed subscheme 
corresponding to $I\subset \{1, \ldots , r\}$, then 
$$
D(X; Y_1, \ldots , Y_r)=
\underset{I\subset \{1, \ldots , r\}}{\bigcup }X_I, 
$$
and for any $j\notin I$, $X_I$ is glued with $X_{I\cup \{j\}}$ 
along $Y_j$ transversally.

\vskip 1pc
\begin{defn}
Suppose that $X$ is a complex algebraic manifold and 
that $Y_1, \ldots , Y_r$ are closed submanifolds of $X$.
Let $\caF$ be a vector bundle on the iterated double 
$D(X; Y_1, \ldots , Y_r)$.
A smooth hermitian metric $h=(h_I)$ on $\caF$ is 
a family of smooth hermitian metrics $h_I$ on the restrictions 
$\caF|_{X_I}$ such that 
$h_I|_{Y_j}=h_{I\cup \{j\}}|_{Y_j}$ for any $j\notin I$.
A smooth hermitian metric $h$ on $\caF$ is said to be smooth 
at infinity if each $h_I$ is a smooth at infinity metric 
on $\caF|_{X_I}$.
\end{defn}

\vskip 1pc
\begin{prop}
If $D(X; Y_1, \ldots , Y_r)$ is quasi-projective, then any vector 
bundle $\caF$ on $D(X; Y_1, \ldots , Y_r)$ admits 
a smooth hermitian metric which is smooth at infinity.
\end{prop}

{\it Proof}: As shown in \cite[\S 3.2]{fulton1}, for any vector 
bundle $\caF$ there is a morphism $\varphi :D(X; Y_1, \ldots , Y_r)\to G$ 
to a projective manifold $G$ and a vector bundle $\caG$ on $G$ 
such that $\varphi^*\caG\simeq \caF$.
If we choose a smooth hermitian metric $h$ on $\caG$, then 
the pull-back metric $\varphi^*h$ on $\caF$ is smooth at infinity.
\qed

\vskip 1pc
Let $X$ be a complex algebraic manifold and assume $X$ to be 
projective. 
Consider the normal crossing divisor 
$$
X\times \partial \square^r=X\times D_1+\cdots +X\times D_r
\subset X\times \square^r
$$
introduced in \S 3.6.
Let $T=D(X\times \square^r; X\times \partial \square^r)$ be 
the associated iterated double.
Then $T$ is quasi-projective, because it is isomorphic to 
$X\times D(\square^r; \partial \square^r)$ and 
$D(\square^r; \partial \square^r)$ is affine and of finite 
type over $\C$.
Hence it follows from Prop.4.2 that any vector bundle on $T$ admits 
a smooth at infinity metric.

Let $(X\times \PP{r})_I\subset T$ be the irreducible 
component corresponding to $I\subset \{1, \ldots , r\}$, and 
$i_I:X\times \PP{r}\to T$ the embedding onto 
$(X\times \PP{r})_I$.
Let $\OV{\caF}$ be a hermitian vector bundle on $T$ with 
a smooth at infinity metric.
Consider the alternating sum of the Chern form of $i_I^*\OV{\caF}$:
$$
\CH_{T,0}(\OV{\caF})=\sum_I(-1)^{|I|}\CH_0(i_I^*\OV{\caF})
\in \underset{p}{\oplus }\WT{\scD}_{\A}^{2p-r}(X, p).
$$
Then the isometry $(\delta_j^i)^*i_I^*\OV{\caF}\simeq 
(\delta_j^i)^*i_{I\cup \{j\}}^*\OV{\caF}$ 
for any $j\notin I$ and for $i=0$ or $\infty $ implies that 
$\delta_{\A}\CH_{T,0}(\OV{\caF})=0$.
Hence $d_s\CH_{T,0}(\OV{\caF})=0$.
The form $\CH_{T,0}(\OV{\caF})$ is called 
the {\it Chern form of $\OV{\caF}$}.

Let 
$$
\OV{\mathscr E}:0\to \OV{\caF}_{-1}\to \OV{\caF}_0\to 
\OV{\caF}_1\to 0
$$
be a short exact sequence of hermitian vector bundles on $T$ 
with smooth at infinity metrics.
Then as shown in \cite[Def.3.8]{BW}, the metric on 
the transgression bundle $\TRNS{1}{i_I^*\OV{\caE}}$ is 
smooth at infinity, therefore we have 
$$
\CH_1(i_I^*\OV{\caE})_{\P}=\CH_0(\TRNS{1}{i_I^*\OV{\caE}})_{\P}\in 
\underset{p}{\oplus }\WT{\scD}_{\A,\P}^{2p-r-1}(X, p).
$$
Let 
$$
\CH_{T,1}(\OV{\caE})=\sum_I(-1)^{|I|}
\CH_1(i_I^*\OV{\caE})_{\P}\in \underset{p}{\oplus }
\WT{\scD}_{\A,\P}^{2p-r-1}(X, p).
$$
Then the isometry $(\delta_j^i)^*i_I^*\OV{\caE}\simeq 
(\delta_j^i)^*i_{I\cup \{j\}}^*\OV{\caE}$ 
for any $j\notin I$ and for $i=0$ or $\infty $ implies that 
$\delta_{\A}\CH_{T,1}(\OV{\caE})=0$.
Hence 
$$
d_s\CH_{T,1}(\OV{\caE})=(-1)^r\delta_{\P}
\CH_{T,1}(\OV{\caE})=(-1)^r\left(\CH_{T,0}(\OV{\caF}_{-1})
+\CH_{T,0}(\OV{\caF}_1)-\CH_{T,0}(\OV{\caF}_0)\right).
$$
The form $\CH_{T,1}(\OV{\caE})$ is called the {\it Bott-Chern form 
of $\OV{\caE}$}.

Summing up the results obtained in this subsection, 
we have the following:

\vskip 1pc
\begin{thm}
The Chern form $\CH_{T,0}(\OV{\caF})$ of a hermitian vector 
bundle $\OV{\caF}$ on $T$ with a smooth at infinity metric 
is $d_s$-closed.
Moreover, for a short exact sequence 
$$
\OV{\mathscr E}:0\to \OV{\caF}_{-1}\to \OV{\caF}_0\to 
\OV{\caF}_1\to 0
$$
of hermitian vector bundles with smooth at infinity metrics on $T$, 
$$
\CH_{T,0}(\OV{\caF}_{-1})+\CH_{T,0}(\OV{\caF}_1)
-\CH_{T,0}(\OV{\caF}_0)=(-1)^rd_s\CH_{T,1}(\OV{\caE}).
$$
In particular, the element of the Deligne cohomology 
represented by $\CH_{T,0}(\OV{\caF})$ is independent of 
the choice of metric, and it induces a map of abelian groups 
$$
\CH_{T,0}=\underset{p}{\oplus }\CH_{T,0}^p:K_0(T)
\to \underset{p}{\oplus }
H_{\scD }^{2p-r}(X, \R(p)).
$$
\end{thm}

\vskip 1pc
\subsection{Relations with the Beilinson's regulator}
Let 
$$
T_j=\underset{j\in I}{\cup }(X\times \PP{r})_I\subset T
$$
for $1\leq j\leq r$, and $\iota_j:T_j\hookrightarrow T$ 
the closed embedding.
Define a morphism $p_j:T\to T_j$ as follows:
For any $I\subset \{1, \ldots , r\}$ with $j\notin I$, 
$p_j|_{(X\times \PP{r})_I}$ and 
$p_j|_{(X\times \PP{r})_{I\cup \{j\}}}$ are given by 
\begin{align*}
p_j|_{(X\times \PP{r})_I}:&\ (X\times \PP{r})_I
\overset{\IId}{\to }(X\times \PP{r})_{I\cup \{j\}}\subset T_j,  \\
p_j|_{(X\times \PP{r})_{I\cup \{j\}}}:&\ 
(X\times \PP{r})_{I\cup \{j\}}\overset{\IId}{\to }
(X\times \PP{r})_{I\cup \{j\}}\subset T_j.
\end{align*}
Then $p_j$ satisfies $p_j\iota_j=\IId_{T_j}$ and 
$p_j(T_i)=T_i\cap T_j$, therefore we have 
maps of $K$-groups 
$$
K_*(T; T_1, \ldots , T_{j-1})\overset{p_j^*}{
\underset{\iota_j^*}{\leftrightarrows }}
K_*(T_j; T_1\cap T_j, \ldots , T_{j-1}\cap T_j)
$$
which satisfy $\iota_j^*p_j^*=\IId $.
Hence $\iota_j^*$ is split surjective, and 
the canonical map 
$$
K_*(T; T_1, \ldots , T_j)\to K_*(T; T_1, \ldots , T_{j-1})
$$
is split injective, and the canonical map  
\begin{equation}
K_*(T, T_1, \ldots , T_r)\hookrightarrow K_*(T) \label{sec4:map1}
\end{equation}
is splitting injective.

Let us now apply the results obtained in \S 2 to the maps of 
$\scC$-complexes 
$$
\WT{\Q}\WALT{*}{T; T_1, \ldots , T_{j-1}}
\overset{p_j^*}{\underset{\iota_j^*}{\leftrightarrows }}
\WT{\Q}\WALT{*}{T_j; T_1\cap T_j, \ldots , T_{j-1}\cap T_j}.
$$
It follows from Prop.2.20 that there is a homotopy from the identity 
to $\iota_j^*p_j^*$, and Cor.2.16 says that 
$\WT{\Q}\WALT{*}{T; T_1, \ldots , T_j}$ is 
isomorphic to the simple complex of $\iota_j^*$.
Hence Prop.2.6 says that there is a left inverse 
map up to homotopy of the canonical map \linebreak 
$\WT{\Q}\WALT{*}{T; T_1, \ldots , T_j}\to 
\WT{\Q}\WALT{*}{T; T_1, \ldots , T_{j-1}}$, 
which we denote by 
$$
t_j:\WT{\Q}\WALT{*}{T; T_1, \ldots , T_{j-1}}\to 
\WT{\Q}\WALT{*}{T; T_1, \ldots , T_j}. 
$$
The description of the map $t_j$ given in Prop.2.6 together with 
Prop.3.7 implies the following:

\vskip 1pc
\begin{prop}
Take an element 
$$
x=(x_I)\in \underset{I\subset \{1, \ldots , j-1\}}
{\underset{|I|=m}{\oplus }}\WT{\Q}\WALT{*+|I|}{T_I}
$$
and consider its image by the map 
$$
t_j^{m,n}:\underset{I\subset \{1, \ldots , j-1\}}
{\underset{|I|=m}{\oplus }}\WT{\Q}\WALT{*+|I|}{T_I}\to 
\underset{I\subset \{1, \ldots , j\}}
{\underset{|J|=n}{\oplus }}\WT{\Q}\WALT{*+|J|}{T_J}.
$$
If $m<n$, then $t_j^{m,n}(x)_J$ is isometrically equivalent 
to a degenerate element.
On the other hand, 
$$
t_j^{m,m}(x)_J=x_J-p_j^*\iota_j^*x_J.
$$
\end{prop}

\vskip 1pc
\begin{cor}
Let 
$$
q:\WT{Q}\WALT{*}{T; T_1, \ldots , T_r}\to \WT{\Q}\WALT{*}{T}
$$
be the canonical map, and 
$$
t=t_rt_{r-1}\cdots t_1:\WT{\Q}\WALT{*}{T}\to 
\WT{\Q}\WALT{*}{T; T_1, \ldots , T_r}.
$$
Then $t$ is a left inverse map of the map $q$ up to homotopy.
For $x\in \WT{\Q}\WALT{*}{T}$, $t^{0,n}(x)$ is isometrically 
equivalent to a degenerate element if $0<n$, and 
$$
t^{0,0}(x)=(1-p_r^*\iota_r^*)\cdots (1-p_1^*\iota_1^*)(x).
$$
\end{cor}

\vskip 1pc
Let 
\begin{align*}
t:&\, \QQ{K_0(T)}\to \QQ{K_0(T; T_1, \ldots , T_r)},  \\
q:&\, \QQ{K_0(T; T_1, \ldots , T_r)}\to \QQ{K_0(T)}, 
\end{align*}
be the maps given by $t$ and $q$ respectively.
Note that $q$ agrees with the map \eqref{sec4:map1} and 
$t$ is a left inverse map of $q$.
Cor.4.5 says that 
$$
qt=(1-p_r^*\iota_r^*)\cdots (1-p_1^*\iota_1^*):\QQ{K_0(T)}\to 
\QQ{K_0(T)}.
$$
Levine have shown in \cite[Thm.1.10]{levine} that the embedding 
$$
i_{\emptyset }:(X\times \PP{r}; X\times \partial \PP{r})
\hookrightarrow (T; T_1, \ldots , T_r)
$$
induces an isomorphism of $K_0$-groups: 
$$
i_{\emptyset }^*:K_0(T; T_1, \ldots , T_r)
\overset{\sim }{\longrightarrow }
K_0(X\times \PP{r}; X\times \partial \PP{r}).
$$

\vskip 1pc
\begin{prop}
The diagram 
$$
\begin{CD}
\QQ{K_0(T)} @>{t}>> \QQ{K_0(T; T_1, \ldots , T_r)} \\ 
@V{\CH_{T,0}^p}VV  @VV{i_{\emptyset }^*}V  \\
H_{\scD}^{2p-r}(X, \R(p)) @<{\CH_0^p}<<
\QQ{K_0(X\times \PP{r}; X\times \partial \PP{r})} 
\end{CD}
$$
is commutative.
\end{prop}

{\it Proof}: 
Let $\OV{\caF}$ be a hermitian vector bundle on $T$ 
with a smooth at infinity metric, and consider 
$$
i_{\emptyset }^*t(\OV{\caF})\in 
\WT{\Q}\WALT{0}{X\times \PP{r}; X\times \partial \PP{r}}=
\underset{I}{\oplus }\, \WT{\Q}\WALT{|I|}{X\times D_I}.
$$
Cor.4.5 implies that $i_{\emptyset }^*t(\OV{\caF})_I$ is 
isometrically equivalent to a degenerate element if 
$I\not= \emptyset $, and that 
$$
i_{\emptyset }^*t(\OV{\caF})_{\emptyset }=i_{\emptyset }^*
(1-p_r^*\iota_r^*)\cdots (1-p_1^*\iota_1^*)\OV{\caF}\in 
\WT{\Q}\WALT{0}{X\times \PP{r}}. 
$$
The commutative diagram of schemes 
\begin{equation}
\xymatrix{
T \ar[r]^{p_j} & T_j\ar[r]^{\iota_j} & T  \\
& X\times \PP{r} \ar[ul]^{i_I} \ar[ur]_{i_{I\cup \{j\}}} &  \\
} \label{sec4:cd1}
\end{equation}
implies that $i_I^*(1-p_j^*\iota_j^*)\OV{\caF}$ is isometric 
to $(i_I^*-i_{I\cup \{j\}}^*)\OV{\caF}$ as virtual hermitian 
vector bundles.
Hence $i_{\emptyset }^*t(\OV{\caF})_{\emptyset }$ is isometric 
to $\underset{I}{\sum }(-1)^{|I|}i_I^*\OV{\caF}$ and 
$$
\CH_0^p(i_{\emptyset }^*t(\OV{\caF}))=
\CH_0^p(i_{\emptyset }^*t(\OV{\caF})_{\emptyset })
=\sum_I(-1)^{|I|}\CH_0^p(i_I^*\OV{\caF})=\CH_{T,0}^p(\OV{\caF}), 
$$
which completes the proof.
\qed

\vskip 1pc
Prop.3.9 and Prop.4.6 imply the following:

\vskip 1pc
\begin{cor}
The diagram 
$$
\xymatrix{
\QQ{K_r(X)}\ar[r]^{\sim \qquad \quad } \ar[dr]_{\CH_r^p} & 
\QQ{K_0(X\times \PP{r}; X\times \partial \PP{r})}& 
\QQ{K_0(T; T_1, \ldots , T_r)}\ar[r] 
\ar[l]^{\quad \sim }_{\quad i_{\emptyset }^*} & 
\QQ{K_0(T)}\ar[dll]^{\CH_{T,0}^p}  \\
& \underset{p}{\oplus }H_{\scD}^{2p-r}(X, \R(p)) 
}
$$
is commutative.
\end{cor}

\vskip 1pc
\begin{prop}
The diagram 
$$
\begin{CD}
\QQ{K_0(T)} @>{t}>>  \QQ{K_0(T; T_1, \ldots , T_r)}  \\
@V{\CH_{T,0}^p}VV @VV{q}V  \\
H_{\scD}^{2p-r}(X, \R(p)) @<{\CH_{T,0}^p}<< \QQ{K_0(T)}
\end{CD}
$$
is commutative.
\end{prop}

{\it Proof}: 
Let $\OV{\caF}$ be a hermitian vector bundle on $T$ with a smooth 
at infinity metric.
Then 
$$
i_I^*qt(\OV{\caF})=i_I^*(1-p_r^*\iota_r^*)\cdots 
(1-p_1^*\iota_1^*)\OV{\caF}.
$$
The commutative diagram \eqref{sec4:cd1} implies that 
$i_I^*qt(\OV{\caF})$ is isometric to zero as virtual hermitian vector 
bundle if $I\not= \emptyset $, and that $i_{\emptyset }^*qt(\OV{\caF})$ 
is isomorphic to $\underset{I}{\sum }(-1)^{|I|}i_I^*\OV{\caF}$.
Hence 
$$
\CH_{T,0}^p(qt(\OV{\caF}))=\sum_I(-1)^{|I|}
\CH_0^p(i_I^*qt(\OV{\caF}))=
\sum_I(-1)^{|I|}\CH_0^p(i_I^*\OV{\caF})=\CH_{T,0}^p(\OV{\caF}), 
$$
which completes the proof.
\qed

\setcounter{equation}{0}
\vskip 2pc
\section{Several arithmetic $K$-groups}

\vskip 1pc
\subsection{Multi-relative arithmetic $K$-theory}
We begin by introducing some terminology and notations which 
are used in Arakelov geometry.
An {\it arithmetic ring} is a triple 
$(A, \Sigma, F_{\infty })$, where $A$ is a Noetherian integral 
domain, $\Sigma $ is a finite set of embeddings 
$A\hookrightarrow \C$, and 
$$
F_{\infty }:\underset{\iota \in \Sigma }{\oplus }\C=\C^{\Sigma }
\to \C^{\Sigma }
$$
is a conjugate-linear involution whose restriction to $A$ 
is the identity.
Here we see $A$ as a subalgebra of $\C^{\Sigma }$ by 
$a\mapsto (\iota (a))_{\iota \in \Sigma }$.
By an {\it arithmetic variety}, we mean a separated regular scheme 
which is flat and of finite type over an arithmetic ring.
For an arithmetic variety $X$ defined over an arithmetic ring $A$, 
denote by $X(\C)$ the complex algebraic manifold associated with 
the smooth algebraic variety $X\underset{A}{\otimes }\C^{\Sigma }$ 
over $\C $, and by $F_{\infty }:X(\C )\to X(\C )$ 
the anti-holomorphic involution induced by 
$F_{\infty }:\C^{\Sigma }\to \C^{\Sigma }$.
A {\it hermitian vector bundle} $\OV{\caF}=(\caF, h)$ on $X$ is 
a vector bundle $\caF $ equipped with an $F_{\infty }$-invariant 
smooth hermitian metric $h$ on $\caF(\C)$.

For a differential form $\eta $ on $X(\C)$, the involution 
on the space $E^n(X(\C))$ given by 
$\eta \mapsto \OV{F}_{\infty }(\eta )=\OV{F_{\infty }^*\eta }$ 
respects the bigrading 
$$
E^n(X(\C))=\underset{p+q=n}{\oplus }E^{p,q}(X(\C)).
$$
Hence it induces involutions on the complexes of 
differential forms and currents introduced in \S 3.
Set 
\begin{align*}
\mmD^*(X, p)&=\mmD^*(X(\C ), p)^{\OV{F}_{\infty }=\IId }, \\
\mmD_D^*(X, p)&=\mmD_D^*(X(\C ), p)^{\OV{F}_{\infty }=\IId },  \\
\WT{\scD}_{\A}^*(X, p)&=
\WT{\scD}_{\A}^*(X(\C ), p)^{\OV{F}_{\infty }=\IId }, \\
\WT{\scD}_{\P}^*(X, p)&=
\WT{\scD}_{\P}^*(X(\C ), p)^{\OV{F}_{\infty }=\IId }, \\
\WT{\scD}_{\A,\P}^*(X, p)&=
\WT{\scD}_{\A,\P}^*(X(\C ), p)^{\OV{F}_{\infty }=\IId }, 
\end{align*}
and 
\begin{align*}
\mmD_*(X)&=\underset{p}{\oplus }\mmD^{2p-*}(X, p), \\
\mmD_{D,*}(X)&=\underset{p}{\oplus }\mmD_D^{2p-*}(X, p),  \\
\WT{\scD}_{\A,*}(X)&=
\underset{p}{\oplus }\WT{\scD}_{\A}^{2p-*}(X, p), \\
\WT{\scD}_{\P,*}(X)&=
\underset{p}{\oplus }\WT{\scD}_{\P}^{2p-*}(X, p), \\
\WT{\scD}_{\A,\P,*}(X)&=
\underset{p}{\oplus }\WT{\scD}_{\A,\P}^{2p-*}(X, p).
\end{align*}
Note that for an exact hermitian $n$-cube $\OV{\caF}$ on $X$, 
we have 
\begin{gather*}
\CH_n(\OV{\caF})_{\P}\in \WT{\scD}_{\P,n}(X),  \\
\CH_n(\OV{\caF})=\kappa_{\P}(\CH_n(\OV{\caF})_{\P})\in \mmD_n(X).
\end{gather*}

We now recall the arithmetic $K_0$-group of $X$ \cite{gilletsoule2}.
Suppose $X$ is proper over an arithmetic ring.
Define the {\it arithmetic $K_0$-group} $\WH{K}_0(X)$ of $X$ 
to be the abelian group generated by pairs $(\OV{\caF}, \WT{\eta })$ 
of a hermitian vector bundle $\OV{\caF}$ on $X$ with 
$\WT{\eta }\in \mmD_1(X)/\IIm d_{\scD}$, subject to the relation 
$$
(\OV{\caF}_{-1}, \WT{\eta }_{-1})+(\OV{\caF}_1, \WT{\eta }_1)-
(\OV{\caF}_0, \WT{\eta }_{-1}+\WT{\eta }_1+\WT{\CH}_1(\OV{\caE}))
$$
for any short exact sequence 
$\OV{\caE}:0\to \OV{\caF}_{-1}\to \OV{\caF}_0\to \OV{\caF}_1\to 0$ 
and for any 
$\WT{\eta }_{-1}, \WT{\eta }_1\in \mmD_1(X)/\IIm d_{\scD}$.
There is a map 
$$
\CH_0:\WH{K}_0(X)\to \mmD_0(X)
$$
which sends a pair $(\OV{\caF}, \WT{\eta })$ to 
$\CH_0(\OV{\caF})+d_{\scD}\eta $.

We next recall the definition of the higher arithmetic $K$-groups.
The original definition given in \cite{takeda} requires 
homotopy theory.
But in this paper we will employ a simpler definition, because 
we only need $K$-groups with rational coefficient.
For $r\geq 1$, define the {\it $r$-th higher rational 
arithmetic $K$-group of $X$}, which we denote by 
$\QQ{\WH{K}_r(X)}$, to be the homology group of 
the simple complex of the higher Bott-Chern form, namely, 
$$
\QQ{\WH{K}_r(X)}=H_r\left(s\left(\CH_*:\WT{\Q}\WH{C}_*(X)
\to \mmD_*(X)\right)\right).
$$
Then there is a long exact sequence 
\begin{align*}
\cdots \to \underset{p}{\oplus }H_{\scD}^{2p-r-1}(X, \R(p))\to 
&\ \QQ{\WH{K}_r(X)}\to \QQ{K_r(X)}\overset{\CH_r}{\to }
\underset{p}{\oplus }H_{\scD}^{2p-r}(X, \R(p))\to \cdots \\
\cdots \to \QQ{K_1(X)}&\ \to \mmD_1(X)/\IIm d_{\scD}\to 
\QQ{\WH{K}_0(X)}\to \QQ{K_0(X)}\to 0.
\end{align*}

For $r\geq 1$, set 
$$
\QQ{\WH{K}_{\P,r}(X)}=H_r\left(s\left(\CH_{*,\P}:
\WT{\Q}\WH{C}_*(X)\to \WT{\scD}_{*,\P}(X)\right)\right).
$$
Since $\kappa_{\P}:\WT{\scD}_{\P,*}(X)\to \mmD_*(X)$ is 
a quasi-isomorphism and $\CH_*=\kappa_{\P}\CH_{*,\P}$, 
there is a natural isomorphism 
\begin{equation}
\QQ{\WH{K}_{\P,r}(X)}\to \QQ{\WH{K}_r(X)}. \label{sec5:iso1}
\end{equation}

Let us now define multi-relative arithmetic $K$-theory.
Let $X$ be a proper arithmetic variety and suppose $r\geq 1$.
Taking the $\OV{F}_{\infty }$-invariant part of the map 
$\CH_*$ defined in Prop.3.8 yields the map of complexes 
$$
\CH_*:\WT{\Q}\WH{C}_*(X\times \PP{r}; X\times \partial \PP{r})[r]
\to \WT{\scD}_{\A,\P,*}(X).
$$
For $n\geq 0$, define {\it multi-relative rational arithmetic 
$K$-theory} of $(X\times \PP{r}; X\times \partial \PP{r})$ 
to be the homology group of the simple complex of this map: 
$$
\QQ{\WH{K}_n(X\times \PP{r}; X\times \partial \PP{r})}=
H_{n+r}(s(\CH_*)). 
$$

\vskip 1pc
\begin{prop}
There is a canonical isomorphism 
$$
\QQ{\WH{K}_{\P,n+r}(X)}\simeq 
\QQ{\WH{K}_n(X\times \PP{r}; X\times \partial \PP{r})}.
$$
\end{prop}

{\it Proof}: 
The commutative diagram of chain complexes in Prop.3.9 
yields the map $s(\CH_{*,\P})\to s(\CH_*)$, 
which leads to the map of arithmetic $K$-groups 
$$
\QQ{\WH{K}_{\P,n+r}(X)}\to 
\QQ{\WH{K}_n(X\times \PP{r}; X\times \partial \PP{r})}
$$
which fits into the commutative diagram 
\begin{multline*}
\begin{CD}
\cdots @>>> \QQ{K_{n+r+1}(X)} @>>> 
\underset{p}{\oplus }H_{\scD}^{2p-n-r-1}(X, \R(p)) \\
@. @VVV @V{\IId }VV  \\
\cdots @>>> \QQ{K_{n+1}(X\times \PP{r}; X\times \partial \PP{r})} @>>> 
\underset{p}{\oplus }H_{\scD}^{2p-n-r-1}(X, \R(p)) 
\end{CD}  \\
\begin{CD}
@>>> \QQ{\WH{K}_{\P,n+r}(X)} @>>> \QQ{K_{n+r}(X)}  @>>> \cdots \\
@. @VVV  @VVV @.  \\
@>>> \QQ{\WH{K}_n(X\times \PP{r}; X\times \partial \PP{r})} 
@>>> \QQ{K_n(X\times \PP{r}; X\times \partial \PP{r})} @>>> \cdots .
\end{CD}
\end{multline*}
The proposition follows from this diagram and the homotopy 
invariant property of $K$-theory.
\qed

\vskip 1pc
\subsection{Arithmetic $K$-group of an iterated double}
In this subsection we will define arithmetic $K_0$-group of 
an iterated double $T=D(X\times \PP{r}; X\times \partial \PP{r})$ 
when $X$ is a projective arithmetic variety over an arithmetic ring.

\vskip 1pc
\begin{defn}
Define arithmetic $K_0$-group $\WH{K}^M_0(T)$ of 
the iterated double $T$ to be the abelian group generated by 
pairs $(\OV{\caF}, \WT{\eta })$, where $\OV{\caF}$ is a vector 
bundle on $T$ with a smooth at infinity metric 
and $\WT{\eta }\in \WT{\scD}_{\A,\P,r+1}(X)/\IIm d_s$, 
subject to the relation 
$$
(\OV{\caF}_{-1}, \WT{\eta }_{-1})+(\OV{\caF}_1, \WT{\eta }_1)-
(\OV{\caF}_0, \WT{\eta }_{-1}+\WT{\eta }_1
+(-1)^r\WT{\CH}_{T,1}(\OV{\caE}))
$$
for any short exact sequence 
$\OV{\caE}:0\to \OV{\caF}_{-1}\to \OV{\caF}_0\to \OV{\caF}_1\to 0$ 
and for any 
$\WT{\eta }_{-1}, \WT{\eta }_1\in \WT{\scD}_{\A,\P,r+1}(X)/\IIm d_s$.
\end{defn}

\vskip 1pc
There is a surjective map 
$$
\zeta :\WH{K}_0^M(T)\to K_0(T)
$$
which sends $[(\OV{\caF}, \WT{\eta })]$ to $[\caF ]$.
It is easy to see that it satisfies the exact sequence 
$$
\WT{\scD}_{\A,\P,r+1}(X)/\IIm d_s\to \WH{K}_0^M(T)
\overset{\zeta }{\to }K_0(T)\to 0.
$$

\vskip 1pc
\begin{defn}
By Thm.4.3 we can define a map  
$$
\CH_{T,0}:\WH{K}^M_0(T)\to \WT{\scD}_{\A,\P,r}(X)
$$
which sends $(\OV{\caF}, \WT{\eta })$ to 
$\CH_{T,0}(\OV{\caF})+d_s\eta$.
Let 
$$
\WH{K}_0(T)=\KER \left(\CH_{T,0}:\WH{K}^M_0(T)\to 
\WT{\scD}_{\A,\P,r}(X)\right).
$$
Moreover, define $\WH{K}^M_0(T; T_1, \ldots , T_r)$ 
to be the cartesian product of the diagram 
$$
\xymatrix{
\WH{K}_0^M(T) \ar[dr]_{\zeta } & & 
K_0(T; T_1, \ldots , T_r) \ar[dl]  \\
 & \quad K_0(T). & \\
}
$$
In other words, $\WH{K}^M_0(T; T_1, \ldots , T_r)$ is a 
subgroup of $\WH{K}^M_0(T)$ given as follows:
$$
\WH{K}_0^M(T; T_1, \ldots , T_r)=\left\{x\in \WH{K}_0^M(T); 
\zeta (x)\in \IIm (K_0(T; T_1, \ldots , T_r)\to K_0(T))\right\}.
$$
Finally, let 
$$
\WH{K}_0(T; T_1, \ldots , T_r)=
\WH{K}_0^M(T; T_1, \ldots , T_r)\cap \WH{K}_0(T)\subset 
\WH{K}_0^M(T).
$$
\end{defn}

\vskip 1pc
Let us recall the maps 
\begin{align*}
t:&\, \WT{\Q}\WALT{*}{T}\to \WT{\Q}\WALT{*}{T; T_1, \ldots , T_r},  \\
q:&\, \WT{\Q}\WALT{*}{T; T_1, \ldots , T_r}\to \WT{\Q}\WALT{*}{T}
\end{align*}
given in Cor.4.5.
In the same way as in the proof of Prop.4.8, we can show that \linebreak 
$\CH_{T,1}(qt(\OV{\caE}))=\CH_{T,1}(\OV{\caE})$ for any short 
exact sequence $\OV{\caE}$ of hermitian vector bundles on $T$.
Hence we obtain a map 
$$
\WH{qt}:\WH{K}_0^M(T)\to \WH{K}_0^M(T)
$$
which sends $[(\OV{\caF}, \WT{\eta })]$ to $[(qt(\OV{\caF}), \WT{\eta })]$.
It is obvious that the image of this map is contained in \linebreak 
$\WH{K}_0^M(T; T_1, \ldots , T_r)$, hence it induces 
\begin{equation}
\WH{t}:\WH{K}_0^M(T)\to \WH{K}_0^M(T; T_1, \ldots , T_r), 
\label{sec5:map0}
\end{equation}
which turns out to be a splitting map of the inclusion 
$\WH{K}_0^M(T; T_1, \ldots , T_r)\subset \WH{K}_0^M(T)$.
Prop.4.8 implies that 
$$
\xymatrix{
\WH{K}_0^M(T) \ar[rr]^{\WH{qt}}\ar[rd]_{\CH_{T,0}} & & 
\WH{K}_0^M(T) \ar[ld]^{\CH_{T,0}}  \\
& \WT{\scD}_{\A,\P,r}(X)
}
$$
is commutative.
Hence the map \eqref{sec5:map0} induces 
$$
\WH{t}:\WH{K}_0(T)\to \WH{K}_0(T; T_1, \ldots , T_r), 
$$
which is also a splitting map of the inclusion 
$\WH{K}_0(T; T_1, \ldots , T_r)\subset \WH{K}_0(T)$.

\vskip 1pc
\begin{prop}
The embedding 
$$
i_{\emptyset }:(X\times \PP{r}; X\times \partial \PP{r})
\hookrightarrow (T; T_1, \ldots , T_r)
$$
induces a surjection 
$$
\WH{i}_{\emptyset }^*:\QQ{\WH{K}_0(T; T_1, \ldots , T_r)}
\twoheadrightarrow 
\QQ{\WH{K}_0(X\times \PP{r}; X\times \partial \PP{r})}
$$
which fits into the commutative diagram up to sign:
\begin{multline*}
\begin{CD}
\underset{p}{\oplus }H_{\scD}^{2p-r-1}(X, \R(p)) @>>> 
\QQ{\WH{K}_0(T; T_1, \ldots , T_r)} \\
@VV{\IId }V  @VV{\WH{i}_{\emptyset }^*}V   \\
\underset{p}{\oplus }H_{\scD}^{2p-r-1}(X, \R(p)) @>>> 
\QQ{\WH{K}_0(X\times \PP{r}; X\times \partial \PP{r})} 
\end{CD} \\
\begin{CD}
@>{\zeta }>> \QQ{K_0(T; T_1, \ldots , T_r)} @>>> 
\underset{p}{\oplus }H_{\scD}^{2p-r}(X, \R(p)) \\
@. @VV{i_{\emptyset }^*}V @VV{\IId }V  \\
@>{\zeta }>> \QQ{K_0(X\times \PP{r}; X\times \partial \PP{r})} 
@>>> \underset{p}{\oplus }H_{\scD}^{2p-r}(X, \R(p)).
\end{CD}
\end{multline*}
\end{prop}

{\it Proof}:
For a short exact sequence of hermitian vector bundles 
$\OV{\caE}:0\to \OV{\caF}_{-1}\to \OV{\caF}_0\to \OV{\caF}_1\to 0$ 
on $T$, consider the element 
$$
(i_{\emptyset }^*t(\OV{\caE}), 0)\in 
\WT{\Q}\WALT{1}{X\times \PP{r}; X\times \partial \PP{r}}\oplus 
\WT{\scD}_{\A,\P,r+2}(X)=s(\CH_*)_{r+1}.
$$
Then we can show in the same way as in the proof of Prop.4.6 that 
$\CH_1(i_{\emptyset }^*t(\OV{\caE}))=\CH_{T,1}(\OV{\caE})$.
Hence 
\begin{align*}
(-1)^r\partial (i_{\emptyset }^*t(\OV{\caE}), 0)=&\, 
(\partial i_{\emptyset }^*t(\OV{\caE}), 
(-1)^r\CH_1(i_{\emptyset }^*t(\OV{\caE})))  \\
=&\, (i_{\emptyset }^*t(\OV{\caF}_{-1})+
i_{\emptyset }^*t(\OV{\caF}_1)-i_{\emptyset }^*t(\OV{\caF}_0), 
(-1)^r\CH_{T,1}(\OV{\caE})).
\end{align*}
Therefore 
$$
(\OV{\caF}, \WT{\eta })\mapsto 
(i_{\emptyset }^*t(\OV{\caF}), -\eta )\in 
\WT{\Q}\WH{C}_0(X\times \PP{r}; X\times \partial \PP{r})\oplus 
\WT{\scD}_{\A,\P,r+1}(X)=s(\CH_*)_0
$$
gives rise to the map 
\begin{equation}
\WH{i_{\emptyset}^*t}:\WH{K}_0^M(T)\to s(\CH_*)_r/\IIm \partial . 
\label{sec5:map1} 
\end{equation}
The equalities 
$$
\partial (i_{\emptyset }^*t(\OV{\caF}), -\eta )=
(0, \CH_0(i_{\emptyset }^*t(\OV{\caF}))+d_s\eta )=
(0, \CH_{T,0}(\OV{\caF})+d_s\eta )
$$
imply that the map \eqref{sec5:map1} induces 
\begin{equation}
\WH{i_{\emptyset}^*t}:\WH{K}_0(T)\to 
\QQ{\WH{K}_0(X\times \PP{r}; X\times \partial \PP{r})}. 
\label{sec5:map2}
\end{equation}
Define the map $\WH{i}_{\emptyset }^*$ to be the composite 
$$
\WH{i}_{\emptyset }^*:\QQ{\WH{K}_0(T; T_1, \ldots , T_r)}\subset 
\QQ{\WH{K}_0(T)}\overset{\WH{i_{\emptyset}^*t}}{\to }
\QQ{\WH{K}_0(X\times \PP{r}; X\times \partial \PP{r})}.
$$
The commutativity of the diagram is easily verified, and 
the surjectivity of $\WH{i}_{\emptyset }^*$ follows 
from the commutative diagram and the bijectivity of 
$i_{\emptyset }^*$ on $K$-groups.
\qed

\vskip 1pc
Summing up the results in this section, we obtain 
the following sequence of maps of arithmetic $K$-groups:
$$
\QQ{\WH{K}_r(X)}\simeq \QQ{\WH{K}_{\P,r}(X)}\simeq 
\QQ{\WH{K}_0(X\times \PP{r}; X\times \partial \PP{r})}
\overset{\WH{i}_{\emptyset }^*}{\twoheadleftarrow }
\QQ{\WH{K}_0(T; T_1, \ldots , T_r)}
\underset{\WH{t}}{\overset{\hookrightarrow }{\leftarrow }}
\QQ{\WH{K}_0(T)}.
$$

\setcounter{equation}{0}
\vskip 2pc
\section{Arithmetic Chern character of a hermitian vector bundle 
on an iterated double}

\vskip 1pc
\subsection{Simple complex of a diagram of complexes}
In this subsection we will introduce simple complex associated 
with a diagram of complexes \cite[\S 1]{BF}.
Let 
$$
\caD =\left(
\text{
\setlength{\unitlength}{1mm}
\begin{picture}(41,12)
 \put(-2,-8){$A_*^1$}
 \put(11,7){$B_*^1$}
 \put(24,-8){$A_*^2$}
 \put(37,7){$B_*^2$}
 \put(3,1){$\scriptstyle{f_1}$}
 \put(22,1){$\scriptstyle{g_1}$}
 \put(29,1){$\scriptstyle{f_2}$}
 \put(3,-3){\vector(1,1){7}}
 \put(24,-3){\vector(-1,1){7}}
 \put(29,-3){\vector(1,1){7}}
\end{picture}
}
\right)
$$
be a diagram of chain complexes.
Consider the map 
$$
\varphi :A_*^1\oplus A_*^2\to B_*^1\oplus B_*^2
$$
defined by 
$\varphi (a_1, a_2)=(f_1(a_1)-g_1(a_2), f_2(a_2))$.
Define the {\it simple complex associated with $\caD$}, 
which we denote by $s(\caD)_*$, to be the simple complex of $\varphi $.
To be more precise, 
$$
s(\caD)_n=A_n^1\oplus A_n^2\oplus B_{n+1}^1\oplus B_{n+1}^2
$$
and if we write an element of $s(\caD)_n$ in the way that 
$$
\left(
\text{
\setlength{\unitlength}{1mm}
\begin{picture}(12,7)
 \put(-2,-4){$a_1$}
 \put(2,3){$b_1$}
 \put(6,-4){$a_2$}
 \put(10,3){$b_2$}
\end{picture}
}
\right), 
$$
then the boundary map is given by 
$$
\partial \left(
\text{
\setlength{\unitlength}{1mm}
\begin{picture}(12,7)
 \put(-2,-4){$a_1$}
 \put(2,3){$b_1$}
 \put(6,-4){$a_2$}
 \put(10,3){$b_2$}
\end{picture}
}
\right)=
\left(
\text{
\setlength{\unitlength}{1mm}
\begin{picture}(67,7)
 \put(-2,-4){$\partial a_1$}
 \put(4,3){$f_1(a_1)-g_1(a_2)-\partial b_1$}
 \put(40,-4){$\partial a_2$}
 \put(47,3){$f_2(a_2)-\partial b_2$}
\end{picture}
}
\right).
$$

\vskip 1pc
\begin{prop}\cite[Cor.1.16]{BF}
If $g_1$ is a quasi-isomorphism in the above diagram $\caD_*$, then 
there is a long exact sequence 
$$
\cdots \to H_n(s(\caD)_*)\to H_n(A_*^1)\to H_{n-1}(B_*^2)\to 
H_{n-1}(s(\caD)_*)\to \cdots .
$$
\end{prop}

\vskip 1pc
\subsection{Cubical higher Chow groups and higher arithmetic 
Chow groups}
In \cite{goncharov1} Goncharov constructed a map from 
the simplicial cycle complex to the complex of currents 
$\mmD_D^*(X, p)$, and defined higher arithmetic Chow groups as 
the homology groups of the simple complex of this map.
Afterward, Burgos, Feliu and the author constructed in 
\cite{BFT} a similar map by using the cubical cycle complex.
In this subsection we will recall their construction.

First we recall the cubical cycle complex.
Let $X$ be an equidimensional variety defined over a field.
A closed subscheme of $X\times \PP{r}$ is said to be 
{\it admissible} if it intersects properly with $X\times D_I$ 
for any $I\subset \{1, \ldots , r\}$.
Denote by $Z^p(X, r)$ the $\Q$-vector space generated by 
all admissible and integral subschemes of $X\times \PP{r}$ 
of codimension $p$.
Then the cocubical scheme structure on $(\PP{*})$ defined in 
\S 3.1 induces the maps 
\begin{align*}
(\delta_j^0)^*, (\delta_j^{\infty })^*:&\, Z^p(X, r)
\to Z^p(X, r-1), \\
\sigma_j^*:&\, Z^p(X, r-1)\to Z^p(X, r)
\end{align*}
for $1\leq j\leq r$ , by which 
$(Z^p(X, *), (\delta_j^0)^*, (\delta_j^{\infty })^*, s_j^*)$ 
is a cubical abelian group.
Define the {\it higher Chow groups} of $X$ to be the homology 
groups of its normalized subcomplex:
$$
CH^p(X, r)=H_r(Z^p(X, *)_0).
$$

Suppose $X$ is a compact complex algebraic manifold.
For an admissible and integral subscheme 
$V\subset X\times \PP{r}$ of codimension $p$, let 
$$
\caP^p(V)={\pi_X}_*(\delta_{\OV{V}}\wedge [\pi_{\P}^*W_r])\in 
\mmD_D^{2p-r}(X, p), 
$$
where $\OV{V}$ is the closure of $V$ in $X\times (\P^1)^r$.
In other words, $\caP^p(V)$ is the current on $X$ defined by 
the integral 
$$
\caP^p(V)(\omega )=\frac{1}{(2\pi i)^{d_X+r-p}}\int_{\WT{V}}
\iota^*\pi_X^*\omega \wedge \iota^*\pi_{\P}^*W_r, 
$$
where $\iota :\WT{V}\to \OV{V}$ is a desingularization.

\vskip 1pc
\begin{thm}\cite[Thm.7.4, Thm.7.8]{BFT}
The above integral is convergent.
The map 
$$
\caP^p:Z^p(X, *)_0\to \mmD_D^{2p-*}(X, p)
$$
given by this integral is a map of complexes, and 
the induced map on homology 
$$
H_r(\caP^p):CH^p(X, r)\to H^{2p-r}(X, \R(p))
$$
agrees with the regulator map of $X$.
\end{thm}

\vskip 1pc
Let $X$ be a smooth proper variety defined over an arithmetic field.
Let $\WHmmD_D^*(X, p)$ be the complex given by 
$$
\WHmmD_D^n(X, p)=\begin{cases}
\mmD_D^n(X, p), &\ n<2p,  \\
\mmD_D^{2p}(X, p)/\mmD^{2p}(X, p), &\ n=2p.
\end{cases}
$$
Then taking the $\OV{F}_{\infty }$-invariant part of $\caP^p$ 
yields a map of complexes 
\begin{equation}
\caP^p:Z^p(X, *)_0\to \WHmmD_D^{2p-*}(X, p). 
\label{sec6:map1}
\end{equation}

\vskip 1pc
\begin{defn}
Let $\WH{Z}^p_D(X, *)_0$ be the simple complex of the map 
\eqref{sec6:map1}, and define the higher arithmetic Chow groups 
of $X$ to be the homology groups of $\WH{Z}^p_D(X, *)_0$:
$$
\WH{CH}^p_D(X, r)=H_r(\WH{Z}^p_D(X, *)_0).
$$
\end{defn}

\vskip 1pc
It is obvious that $\WH{CH}^p_D(X, 0)$ agrees with the arithmetic 
Chow groups defined by Gillet and Soul\'{e} in \cite{gilletsoule1}.
The natural map 
$$
\mmD^{2p-r-1}(X, p)\subset \WHmmD_D^{2p-r-1}(X, p)\to 
\WH{Z}^p_D(X, r)_0
$$
gives 
$$
a:H_{\scD}^{2p-r-1}(X, \R(p))\to \WH{CH}_D^p(X, r)
$$
if $r\geq 1$, and 
$$
a:\mmD^{2p-1}(X, p)/\IIm d_{\scD}\to \WH{CH}_D^p(X, 0)
$$
when $r=0$.
Concerning these maps the following long exact sequence holds:
\begin{align*}
\to &\, H_{\scD}^{2p-r-1}(X, \R(p))\overset{a}{\to }
\WH{CH}^p_D(X, r)\to CH^p(X, r)\to H_{\scD}^{2p-r}(X, \R(p))
\overset{a}{\to }\cdots \\
&\cdots \to CH^p(X, 1)\to \mmD^{2p-1}(X, p)/\IIm d_{\scD}
\overset{a}{\to }\WH{CH}^p_D(X, 0)\to CH^p(X, 0)\to 0.
\end{align*}

\vskip 1pc
\subsection{Another definition of higher arithmetic Chow groups}
In \cite{gilletsoule1}, Gillet and Soul\'{e} defined 
intersection product in the arithmetic Chow groups.
It is quite natural to seek for a similar product structure in 
their higher analogues.
However, since the definition of $\WH{CH}_D^p(X, r)$ involves 
the space of currents, it seems impossible to put a product 
structure on $\WH{CH}_D^p(X, r)$.
Burgos and Feliu gave in \cite{BF} another 
definition of higher arithmetic Chow groups in which 
one can define intersection product.
In this subsection we will recall their definition.
To do this we first introduce several complexes of 
differential forms.

For an equidimensional variety $X$ defined over a field, 
denote by $\scZ^p_{X,r}$ the set of all admissible subschemes 
of $X\times \PP{r}$ of codimension $p$.
We abbreviate $\scZ^p_{X,r}$ to $\scZ_r^p$ if no confusion occurs.
For a complex algebraic manifold $X$, set 
$$
\scD_{\log}^*(X\times \PP{r}-\scZ_r^p, p)=
\underset{Z\in \scZ_r^p}{\varinjlim }
\scD_{\log}^*(X\times \PP{r}-Z, p), 
$$
and take the simple complex 
$$
\scD_{\log,\scZ_r^p}^*(X\times \PP{r}, p)=s\left(
\scD_{\log}^*(X\times \PP{r}, p)\to 
\scD_{\log}^*(X\times \PP{r}-\scZ_r^p, p)\right)
$$
and its truncated subcomplex 
$$
\mmD_{\log,\scZ_r^p}^*(X\times \PP{r}, p)=\tau_{\leq 2p}
\scD_{\log,\scZ_r^p}^*(X\times \PP{r}, p).
$$
Set 
$$
\scH^p(X, r)=H_{\scD,\scZ_r^p}^{2p}(X\times \PP{r}, \R(p)), 
$$
then it holds that 
$$
\scH^p(X, r)\simeq 
H^{2p}(\mmD_{\log,\scZ_r^p}^*(X\times \PP{r}, p)).
$$
Since $\scH^p(X, r)$ and $\mmD_{\log,\scZ_r^p}^*(X\times \PP{r}, p)$ 
have cubical structures with respect to the index $r$, we can obtain 
the normalized subcomplexes, which we denote by $\scH^p(X, *)_0$ 
and $\mmD_{\log,\scZ_r^p}^*(X\times \PP{r}, p)_0$ respectively.
Moreover, the cycle class map in Deligne cohomology gives a map 
$$
\chi_1:Z^p(X, *)_0\to \scH^p(X, *)_0
$$
such that 
$\chi_1\otimes \R:Z^p(X, *)_0\otimes \R\to \scH^p(X, *)_0$ 
is an isomorphism.
Set 
$$
\scD_{\A,\scZ^p}^{s,-r}(X, p)_0=
\mmD_{\log,\scZ_r^p}^s(X\times \PP{r}, p)_0
$$
and denote by $\scD_{\A,\scZ^p}^*(X, p)_0$ the associated single complex.
Then we can write any element of $\scD_{\A,\scZ^p}^{2p-r}(X, p)_0$ 
as $\left((\omega_r, g_r), \ldots , (\omega_0, g_0)\right)$, where 
$$
(\omega_i, g_i)\in 
\mmD_{\log,\scZ_i^p}^{2p-r+i}(X\times \PP{i}, p)_0, 
$$
in other words, 
$$
\omega_i\in \mmD_{\log}^{2p-r+i}(X\times \PP{i}, p)_0, \ 
g_i\in \mmD_{\log}^{2p-r+i-1}(X\times \PP{i}-\scZ_i^p, p)_0
$$
such that $d_{\scD}g_r=\omega_r$.

\vskip 1pc
\begin{prop}\cite[Prop.2.13]{BF}
The map 
$$
\chi_2:\scD_{\A,\scZ^p}^{2p-*}(X, p)_0\to \scH^p(X, *)_0
$$
defined by 
$$
\left((\omega_r, g_r), \ldots , (\omega_0, g_0)\right)
\mapsto [(\omega_r, g_r)]
$$
is a quasi-isomorphism.
\end{prop}

\vskip 1pc
Define a map of complexes 
$$
\scD_{\A,\scZ^p}^*(X, p)_0\to \scD_{\A}^*(X, p)_0
$$
by 
$$
\left((\omega_r, g_r), \ldots , (\omega_0, g_0)\right)
\mapsto (\omega_r, \ldots , \omega_0), 
$$
and denote by $\rho $ the composite below:
$$
\rho :\scD_{\A,\scZ^p}^*(X, p)_0\to \scD_{\A}^*(X, p)_0
\overset{\sim }{\to}\WT{\scD}_{\A}^*(X, p)\to 
\WT{\scD}_{\A,\P}^*(X, p).
$$

Let us now give another definition of higher arithmetic Chow 
groups \cite{BF}.
Let $X$ be a smooth variety defined over an arithmetic field.
Denote by $\scH^p(X, *)_0, \scD_{\A,\scZ^p}^*(X, p)_0$ 
and $\WT{\scD}_{\A,\P}^*(X, p)$ the $\OV{F}_{\infty }$-invariant 
part of the complexes 
$\scH^p(X(\C), *)_0, \scD_{\A,\scZ^p}^*(X(\C), p)_0$ and 
$\WT{\scD}_{\A,\P}^*(X(\C), p)$ respectively.
Let $\WH{\scD}_{\A,\P}^*(X, p)$ be the complex defined by 
$$
\WH{\scD}_{\A,\P}^n(X, p)=\begin{cases}
\WT{\scD}_{\A,\P}^n(X, p), &\ n<2p,  \\
0, &\ n\geq 2p.
\end{cases}
$$

\vskip 1pc
\begin{defn}\cite[Def.4.2]{BF}
Let $\WH{Z}^p(X, *)_0$ be the simple complex 
of the diagram 
$$
\xymatrix{
 & \scH^p(X, *)_0 & & \WH{\scD}_{\A,\P}^{2p-*}(X, p) \\
Z^p(X, *)_0 \ar[ru]^{\chi_1} & & 
\scD_{\A,\scZ^p}^{2p-*}(X, p)_0 \ar[lu]_{\chi_2} 
\ar[ru]^{\rho }. & 
 }
$$
Define the {\it higher arithmetic Chow groups} of $X$ to be 
the homology groups of this complex:
$$
\WH{CH}^p(X, r)=H_r(\WH{Z}^p(X, *)_0).
$$
\end{defn}

\vskip 1pc
\begin{defn}
For $r\geq 1$, define a map 
$$
a:\WT{\scD}_{\A,\P}^{2p-r-1}(X, p)/\IIm d_s\to 
\WH{Z}^p(X, r)_0/\IIm \partial , 
$$
by 
$$
a(\WT{\eta^p})=
\left[\left(
\text{
\setlength{\unitlength}{1mm}
\begin{picture}(15,7)
 \put(-2,-4){$0$}
 \put(2,3){$0$}
 \put(6,-4){$0$}
 \put(9,3){$-\eta ^p$}
\end{picture}
}
\right)\right].
$$
Then it induces a map on (co)homology:
$$
a:H_{\scD}^{2p-r-1}(X, \R(p))\to \WH{CH}^p(X, r).
$$
\end{defn}

\vskip 1pc
If $\eta^p\in \scD_{\A}^{2p-r-1}(X, p)_0$, then 
$$
a(\WT{\eta^p})=\left[\left(
\text{
\setlength{\unitlength}{1mm}
\begin{picture}(21,7)
 \put(-2,-4){$0$}
 \put(2,3){$0$}
 \put(5,-4){$d_s(\eta ^p, 0)$}
 \put(20,3){$0$}
\end{picture}
}
\right)\right].
$$
In particular, if $\delta_{\A}\eta^p=0$, then 
$$
a(\WT{\eta^p})=\left[\left(
\text{
\setlength{\unitlength}{1mm}
\begin{picture}(25,7)
 \put(-2,-4){$0$}
 \put(2,3){$0$}
 \put(5,-4){$(d_{\scD}\eta ^p, \eta ^p)$}
 \put(24,3){$0$}
\end{picture}
}
\right)\right].
$$
In the case that $r=0$, the canonical inclusion 
$$
\mmD^{2p-1}(X, p)/\IIm d_{\scD}\hookrightarrow 
\WT{\scD}_{\A,\P}^{2p-1}(X, p)/\IIm d_s
$$
turns out to be an isomorphism.
Hence it follows from Prop.6.1 and Prop.6.4 that there is 
a long exact sequence 
\begin{align*}
\to H_{\scD}^{2p-r-1}(X, \R(p))&\ \overset{a}{\to }
\WH{CH}^p(X, r)\to CH^p(X, r)\to H_{\scD}^{2p-r}(X, \R(p))
\overset{a}{\to }\cdots  \\
\cdots \to CH^p(X, 1)&\ \to \mmD^{2p-1}(X, p)/\IIm d_{\scD}
\overset{a}{\to }
\WH{CH}^p(X, 0)\to CH^p(X, 0)\to 0.
\end{align*}

Let us compare two chain complexes $\WH{Z}_D^p(X, *)_0$ and 
$\WH{Z}^p(X, *)_0$.
The quasi-isomorphism 
$\kappa =\kappa_{\A,\P}:\WT{\scD}_{\A,\P}^*(X, p)\to \mmD_D^*(X, p)$
defined in \S 3.4 induces a map of complexes 
$$
\kappa :\WH{\scD}_{\A,\P}^*(X, p)\to \WHmmD_D^*(X, p).
$$
Let 
$$
\theta:\scH^p(X, *)_0
\overset{(\chi_1\otimes \R)^{-1}}{\longrightarrow }
Z^p(X, *)_0\otimes \R\overset{\caP^p}{\to }\mmD_D^{2p-*}(X, p), 
$$
which is also a map of complexes.
For $g_i\in \mmD_{\log}^{2p-r+i-1}(X\times \PP{i}-\scZ^p, p)$, 
$g_i\bullet \pi_{\P}^*W_i$ is locally integrable on 
$X\times (\P^1)^i$ such that 
\begin{equation}
d_{\scD}{\pi_X}_*[g_i\bullet \pi_{\P}^*W_i]=
{\pi_X}_*[d_{\scD}g_i\bullet \pi_{\P}^*W_i]
+(-1)^{r+i-1}{\pi_X}_*[\delta_{\A}g_i\bullet \pi_{\P}^*W_{i-1}] 
\label{sec6:eq1}
\end{equation}
if $i<r$, and 
\begin{equation}
d_{\scD}{\pi_X}_*[g_r\bullet \pi_{\P}^*W_r]=
{\pi_X}_*[d_{\scD}g_r\bullet \pi_{\P}^*W_r]
-{\pi_X}_*[\delta_{\A}g_r\bullet \pi_{\P}^*W_{r-1}]
-\caP^p(z) \label{sec6:eq2}
\end{equation}
by \cite[Prop.7.5, Prop.7.6]{BFT}.
In the above, $z\in Z^p(X, r)\otimes \R$ is the unique element 
satisfying $(\chi_1\otimes \R )(z)=[(d_{\scD}g_r, g_r)]$ in 
$\scH^p(X, r)$.

\vskip 1pc
\begin{prop}
Assume $r\geq 1$.
Define a map 
$$
\psi :\scD_{\A,\scZ^p}^{2p-r}(X, p)_0
\to \mmD_D^{2p-r-1}(X, p)
$$
by 
$$
\psi ((\omega_r, g_r), \ldots , (\omega_0, g_0))
=\sum_{i=0}^r{\pi_X}_*[g_i\bullet \pi_{\P}^*W_i].
$$
Then 
$$
\psi d_s+d_{\scD}\psi =\kappa \rho -\theta \chi_2.
$$
\end{prop}

{\it Proof}: 
Let $\alpha =((\omega_r, g_r), \ldots , (\omega_0, g_0))\in 
\scD_{\A,\scZ^p}^{2p-r}(X, p)_0$.
Then by \eqref{sec6:eq1} and \eqref{sec6:eq2}, 
$$
d_{\scD}\psi (\alpha )=\sum_{i=0}^r{\pi_X}_*
[d_{\scD}g_i\bullet \pi_{\P}^*W_i]+\sum_{i=1}^r(-1)^{r+i-1}
{\pi_X}_*[\delta_{\A}g_i\bullet \pi_{\P}^*W_{i-1}]-\caP^p(z), 
$$
where $z\in Z^p(X, r)_0\otimes \R$ such that 
$(\chi_1\otimes \R)(z)=[(\omega_r, g_r)]=\chi_2(\alpha )$.
Moreover, 
\begin{align*}
\psi d_s(\alpha )=&\ 
\sum_{i=0}^{r-1}{\pi_X}_*[(\omega_i-d_{\scD}g_i+
(-1)^{r-i+1}\delta_{\A}g_{i+1})\bullet \pi_{\P}^*W_i], \\
\kappa \rho (\alpha )=&\ \sum_{i=0}^r{\pi_X}_*
[\omega_i\bullet \pi_{\P}^*W_i].
\end{align*}
Then 
$$
d_{\scD}\psi (\alpha )+\psi d_s(\alpha )
-\kappa \rho (\alpha )=-\caP^p(z)
=-\theta (\chi_1\otimes \R)(z)=-\theta \chi_2(\alpha ), 
$$
which completes the proof.
\qed

\vskip 1pc
It follows from Prop.6.7 that the map 
$$
\varGamma :\WH{Z}^p(X, *)_0\to \WH{Z}_D^p(X, *)_0
$$
given by 
$$
\varGamma \left(
\text{
\setlength{\unitlength}{1mm}
\begin{picture}(12,7)
 \put(-2,-4){$y$}
 \put(2,3){$\beta_1$}
 \put(6,-4){$\alpha $}
 \put(10,3){$\beta_2$}
\end{picture}
}
\right)=
(y, -\psi (\alpha )+\theta (\beta_1)+\varphi (\beta_2))
$$
is a map of complexes which fits into the commutative diagram 
\begin{multline*}
\hskip 3pc
\begin{CD}
\cdots @>>> H_{\scD}^{2p-r-1}(X, \R(p)) @>{a}>> \WH{CH}^p(X, r)  \\
@. @V{\IId }VV @V{\varGamma_*}VV  \\
\cdots @>>> H_{\scD}^{2p-r-1}(X, \R(p)) @>{a}>> \WH{CH}_D^p(X, r) 
\end{CD}  \\
\begin{CD}
@>>> CH^p(X, r) @>>> H_{\scD}^{2p-r}(X, \R(p)) @>{a}>> \cdots  \\
@. @V{\IId}VV @V{\IId }VV @. \\
@>>> CH^p(X, r) @>>> H_{\scD}^{2p-r}(X, \R(p)) @>{a}>> \cdots .
\end{CD}
\hskip 3pc
\end{multline*}
This leads to the following theorem:

\vskip 1pc
\begin{thm}
When $r\geq 1$, 
$$
\varGamma_*:\WH{CH}^p(X, r)\to \WH{CH}_D^p(X, r)
$$
is an isomorphism.
\end{thm}

\vskip 1pc
\subsection{Chern character of a vector bundle 
on an iterated double}
Let $X$ be a smooth projective variety defined over a field.
Consider the associated iterated double 
$T=D(X\times \PP{r}; X\times \partial \PP{r})$ and 
a morphism $\varphi :T\to G$ to a smooth projective variety $G$.
For any $I\subset \{1, \ldots , r\}$, denote by 
$\varphi_I:X\times \PP{r}\to G$ the restriction of $\varphi $ to 
the irreducible component corresponding to $I$.
Let $Z^p_{\varphi }(G, n)$ be the subgroup of $Z^p(G, n)$ 
such that $x\in Z_{\varphi }^p(G, n)$ 
if and only if we can take the pull-back cycle of $x$ by the morphism 
$$
X\times D_J\times D_K\hookrightarrow X\times \PP{r}\times \PP{n}
\overset{\varphi_I\times \IId }{\to }G\times \PP{n}
$$
for any $I, J\subset \{1, \ldots , r\}$ and 
$K\subset \{1, \ldots , n\}$.
The moving lemma \cite[Thm.1.3]{hanamura} says that 
the embedding of the normalized subcomplexes 
$$
Z^p_{\varphi }(G, *)_0\hookrightarrow Z^p(G, *)_0
$$
is a quasi-isomorphism.

Let $\caF$ be a vector bundle on $T$.
Then we can obtain a morphism $\varphi :T\to G$ to a smooth projective 
variety $G$ and a vector bundle $\caG$ on $G$ such that 
$\varphi^*\caG \simeq \caF$ \cite[\S 3.2]{fulton1}.
Let $y\in Z^p_{\varphi }(G)$ be a cycle representing the $p$-th 
Chern character $\CH_0^p(\caG)\in CH^p(G)$, and let us denote 
$$
\varphi^*(y)=\sum_I(-1)^{|I|}\varphi_I^*(y)\in Z^p(X\times \PP{r}).
$$
Since $\varphi_I^*(y)|_{\{z_j=i\}}=
\varphi_{I\cup \{j\}}^*(y)|_{\{z_j=i\}}$ 
for any $I\subset \{1, \ldots , r\}$ with $j\notin I$ and 
for $i=0$ or $\infty $, it holds that $(\delta_j^i)^*\varphi^*(y)=0$.
In particular, $\varphi^*(y)\in Z^p(X, r)_0$ such that 
$\partial \varphi^*(y)=0$.

\vskip 1pc
\begin{thm}
The element $[\varphi^*(y)]\in CH^p(X, r)$ depends only on $\caF$, 
and is independent of the choice of morphisms $\varphi $, 
vector bundles $\caG$ and $y\in Z^p_{\varphi }(G)$.
Hence in what follows we denote 
$$
[\varphi^*(y)]=\CH^p_{T,0}(\caF)\in CH^p(X, r).
$$
\end{thm}

{\it Proof}: 
First we show the independence of the choice of 
$y\in Z^p_{\varphi }(G)$.
Let $\PR{y}\in Z^p_{\varphi }(G)$ be another element representing 
$\CH_0^p(\caG)\in CH^p(G)$.
Then there is an element $w\in Z^p_{\varphi }(G, 1)_0$ such that 
$\partial w=y-\PR{y}$.
Let 
$$
\varphi^*(w)=\sum_I(-1)^{|I|}\varphi_I^*(w)
\in Z^p(X\times \PP{r}, 1).
$$
By the identification 
$$
(X\times \PP{r})\times \PP{1}\simeq X\times \PP{r+1}, \quad 
((x, z_1, \ldots , z_r), \PR{z})\mapsto 
(x, z_1, \ldots , z_r, \PR{z}), 
$$
we can see $\varphi^*(w)$ as an element of $Z^p(X\times \PP{r+1})$.
Since $\varphi_I^*(w)|_{\{z_j=i\}}=
\varphi_{I\cup \{j\}}^*(w)|_{\{z_j=i\}}$ for any 
$I\subset \{1, \ldots , r\}$ with $j\notin I$ and 
for $i=0$ or $\infty $, $\varphi^*(w)\in Z^p(X, r+1)_0$ such that 
$$
\partial \varphi^*(w)=(-1)^{r+1}(\delta_{r+1}^0)^*\varphi^*(w)=
(-1)^{r+1}(\varphi^*(y)-\varphi^*(\PR{y})).
$$
Hence $[\varphi^*(y)]=[\varphi^*(\PR{y})]$ in $CH^p(X, r)$.

We next show the independence of the choice of $\varphi:T\to G$ 
and vector bundles $\caG$.
Let $\PR{\varphi }:T\to \PR{G}$ be another morphism to a smooth 
projective variety $\PR{G}$ and $\PR{\caG}$ a vector bundle 
on $\PR{G}$ such that there is an ismorphism 
${\PR{\varphi }}^*\PR{\caG}\simeq \caF$.
Then as shown in \cite[\S 3.2]{fulton1}, there is a commutative 
diagram of morphisms 
$$
\xymatrix{
 & T \ar[ld]_{\PR{\varphi }} \ar[d]^{\PPRR{\varphi }} 
\ar[rd]^{\varphi } &   \\
 \PR{G} & \PPRR{G} \ar[l]^{\PR{\psi }} \ar[r]_{\psi } & G, \\
 }
$$
where $\PPRR{G}$ is a smooth projective variety on which there is 
an isomorphism $\psi^*\caG \simeq {\PR{\psi }}^*\PR{\caG}$ 
which admits the commutative diagram of isomorphisms 
$$
\xymatrix{
\varphi^*\caG \ar@{-}[r]^{\sim } \ar@{-}[rrd]_{\sim } & 
{\PPRR{\varphi }}^*\psi^*\caG \ar@{-}[rr]^{\sim } & & 
{\PPRR{\varphi }}^*{\PR{\psi }}^*\PR{\caG}  
\ar@{-}[r]^{\sim } & {\PR{\varphi }}^*\PR{\caG} 
\ar@{-}[lld]^{\sim } \\
& & \caF. & & 
}
$$
Hence we may assume that $G=\PPRR{G}$ and $\psi =\IId $, 
that is, there is a morphism 
$\psi :G\to \PR{G}$ such that $\psi^*\PR{\caG}\simeq \caG$ 
and $\PR{\varphi }=\psi \varphi $.

Set 
$Z^p_{\psi , \PR{\varphi }}(\PR{G})=Z^p_{\psi }(\PR{G})\cap 
Z^p_{\PR{\varphi }}(\PR{G})$.
Then there is a cycle $\PR{y}\in Z^p_{\psi ,\PR{\varphi }}(\PR{G})$ 
representing $\CH_0^p(\PR{\caG})\in CH^p(\PR{G})$.
Moreover, we can take the pull-back cycle 
$y=\psi^*(\PR{y})\in Z^p_{\varphi }(G)$ which represents 
$\CH_0^p(\caG)\in CH^p(G)$.
Then $\varphi^*(y)={\PR{\varphi }}^*(\psi^*(\PR{y}))=
{\PR{\varphi }}^*(\PR{y})$ in $Z^p(X, r)_0$.
This completes the proof.
\qed

\vskip 1pc
\begin{thm}
The above correspondence $\caF \mapsto [\varphi^*(y)]$ 
gives a map of abelian groups 
$$
\CH_{T,0}^p:K_0(T)\longrightarrow CH^p(X, r).
$$
\end{thm}

{\it Proof}: 
We have only to show that 
$$
\CH_{T,0}^p(\caF_{-1})+\CH_{T,0}^p(\caF_1)=\CH_{T,0}^p(\caF_0)
$$
for any short exact sequence 
$0\to \caF_{-1}\to \caF_0\to \caF_1\to 0$ 
of vector bundles on $T$.
We can show this equality in the same way as the proof 
of Thm.6.9, mimicking the argument in \cite[\S 3.2]{fulton1}.
\qed

\vskip 1pc
\noindent
\begin{defn}
Composing the map $\CH_{T,0}^p$ with the sequence of maps 
$$
K_r(X)\simeq K_0(X\times \PP{r}; X\times \partial \PP{r})
\simeq K_0(T; T_1, \ldots , T_r)\subset K_0(T)
$$
given by Levine in \cite{levine}, we can define a map 
$$
\CH^p_r:K_r(X)\to CH^p(X, r).
$$
\end{defn}

\vskip 1pc
\noindent
{\it Remark}: 
The author do not know if the map defined above agrees with 
the higher Chern character map.
However, we will show later in Thm.7.4 that 
the composite of $\CH^p_r$ with the regulator map 
$$
CH^p(X, r)\to H_{\scD}^{2p-r}(X, \R(p))
$$
agrees with Beilinson's regulator.
This is a strong evidence that these two maps agree.

\vskip 1pc
\subsection{Arithmetic Chern character of a hermitian vector bundle 
on an iterated double}
In \cite{burgos2} Burgos extended the definition of 
arithmetic Chow groups to open varieties, and to this end he gave 
another definiton of Green forms.
We begin by recalling his construction.

Let $X$ be a complex algebraic manifold $X$.
The {\it space of Green forms of codimension $p$} is defined as 
the truncated cohomology groups of Deligne complexes as follows:
$$
GE^p(X)=\WH{H}^p(\scD_{\log}^*(X, p), \scD_{\log}^*(X-\scZ^p, p)).
$$
For any subset $\scZ_*^p\subset \scZ^p$, set 
$$
GE_{\scZ_*^p}^p(X)=\WH{H}^p(\scD_{\log}^*(X, p), 
\scD_{\log}^*(X-\scZ_*^p, p)).
$$
Any element of $GE^p(X)$ (resp. $GE_{\scZ_*^p}^p(X)$) is given 
by a pair $(\omega , \WT{g})$ of $\omega \in \mmD_{\log}^{2p}(X, p)$ 
with $\WT{g}\in \mmD_{\log}^{2p-1}(X-\scZ^p)/\IIm d_{\scD}$ 
(resp. $\WT{g}\in \mmD_{\log}^{2p-1}(X-\scZ_*^p)/\IIm d_{\scD}$) 
such that $\omega =d_{\scD}g$.

For an arithmetic variety $X$ defined over an arithmetic ring, 
let us denote by \linebreak 
$\WH{CH}^p(X, \scD (E_{\log}^*))$ 
the arithmetic Chow group of $X$ by means of $GE^p(X)$ 
\cite[\S 7]{burgos2}.
Note that when $X$ is a proper arithmetic variety, 
$\WH{CH}^p(X, \scD (E_{\log}^*))$ is canonically 
isomorphic to the arithmetic Chow group originally defined 
by Gillet and Soul\'{e} in \cite{gilletsoule1}, since in this case 
the definition of arithmetic Chow groups does not depend 
on the complexes where Green objects lie.

Let $X$ be a smooth projective variety defined over 
an arithmetic field, and 
$T=D(X\times \PP{r}; X\times \partial \PP{r})$ the iterated double.
For a morphism $\varphi :T\to G$ to a smooth projective variety 
and for a differential form $\omega $ on $G$, let us write 
$$
\varphi^*(\omega )=\sum_I(-1)^{|I|}\varphi_I^*(\omega ).
$$
Since 
$\varphi_I^*(\omega )|_{\{z_j=i\}}=
\varphi_{I\cup \{j\}}^*(\omega )|_{\{z_j=i\}}$ 
for any $I\subset \{1, \ldots , r\}$ with $j\notin I$ and 
for $i=0$ or $\infty $, $\varphi^*(\omega )$ is normalized 
such that $\delta_{\A}\varphi^*(\omega )=0$.

Let $\OV{\caF }$ be a vector bundle with a smooth at 
infinity metric on $T$.
Take a morphism $\varphi :T\to G$ to a smooth projective variety 
$G$, a vector bundle $\caG$ such that $\caF \simeq \varphi^*\caG$, 
and a cycle $y\in Z^p_{\varphi }(G)$ representing 
$\CH_0^p(\caG)\in CH^p(G)$.
Given a smooth hermitian metric $h$ on $\caG$, we obtain 
the $p$-th arithmetic Chern character 
$\WH{\CH}_0^p(\OV{\caG})\in \WH{CH}^p(G, \scD (E_{\log}^*))$ of 
$\OV{\caG}=(\caG, h)$.
Take a Green form $(\omega, \WT{g})\in GE_{\scZ_{\varphi }^p}(G)$ 
associated with $y$ in the sense of \cite[\S 5]{burgos2} such 
that the pair $(y, (\omega , \WT{g}))$ represents 
$\WH{\CH}_0^p(\OV{\caG})$.

Let $g\in \mmD_{\log}^{2p-1}(G-\scZ_{\varphi }^p, p)$ be a lift of $\WT{g}$, 
and let us denote 
$$
(\varphi^*(\omega ), \varphi^*(g))=
((\varphi^*(\omega ), \varphi^*(g)), (0, 0), \ldots , (0, 0))\in 
\scD_{\A,\scZ^p}^{2p-r}(X, p)_0.
$$
Then $\delta_{\A}(\varphi^*(\omega ), \varphi^*(g))=0$ and 
$\chi_1(\varphi^*(y))=\chi_2(\varphi^*(\omega ), \varphi^*(g))$.
Moreover, if we take the $(p, p)$-part of the Bott-Chern form 
$$
\CH_{T,1}^p(\OV{\caF}, \varphi^*\OV{\caG})=
\CH_{T,1}^p(\OV{\caF}\overset{\sim }{\to }\varphi^*\OV{\caG}\to 0)
\in \WT{\scD}_{\A,\P}^{2p-r-1}(X, p), 
$$
then 
$$
d_s\CH_{T,1}^p(\OV{\caF}, \varphi^*\OV{\caG})=
(-1)^r\left(\CH_{T,0}^p(\OV{\caF})-\varphi^*(\CH_{T,0}^p(\OV{\caG}))
\right).
$$

\vskip 1pc
\begin{defn}
Let $\WH{Z}^p(X, r)_0^*$ be the subgroup of $\WH{Z}^p(X, r)_0$ defined 
as follows:
$$
\WH{Z}^p(X, r)_0^*=\left\{x\in \WH{Z}^p(X, r)_0; \ 
\partial x=\left(
\text{
\setlength{\unitlength}{1mm}
\begin{picture}(11,7)
 \put(-2,-4){$0$}
 \put(2,3){$0$}
 \put(6,-4){$0$}
 \put(10,3){$*$}
\end{picture}
}
\right)
\right\}.
$$
Then $\IIm \partial \subset \KER \partial \subset \WH{Z}^p(X, r)_0^*$, 
and the quotient group $\WH{Z}^p(X, r)^*_0/\IIm \partial $ can be 
expressed as a homology group of a chain complex.
In fact, if we denote 
$$
(\sigma_{<r}A^*)^n=\begin{cases}A^n, &\ n<r,  \\
0, &\ n\geq r, \end{cases}
$$
then 
$$
\WH{Z}^p(X, r)^*_0/\IIm \partial =
H_r\left(
\text{
\setlength{\unitlength}{1mm}
\begin{picture}(79,12)
 \put(-2,-8){$Z^p(X, *)_0$}
 \put(13,7){$\caH^p(X, *)_0$}
 \put(30,-8){$\scD_{\A,\scZ^p_X}^{2p-*}(X, p)_0$}
 \put(48,7){$\sigma_{<2p-r}\WH{\scD}_{\A,\P}^{2p-*}(X, p)$}
 \put(10,-3){\vector(1,1){7}}
 \put(32,-3){\vector(-1,1){7}}
 \put(44,-3){\vector(1,1){7}}
\end{picture}
}\right).
$$
\end{defn}

\vskip 1pc
Consider the element 
$$
cl(\OV{\caF}, \varphi , \OV{\caG}, (y, (\omega, g))=
\left(
\text{
\setlength{\unitlength}{1mm}
\begin{picture}(72,7)
 \put(-2,-4){$\varphi^*(y)$}
 \put(9,3){$0$}
 \put(13,-4){$(\varphi^*(\omega ), \varphi^*(g))$}
 \put(36,3){$(-1)^{r+1}\CH_{T,1}^p(\OV{\caF}, \varphi^*\OV{\caG})$}
\end{picture}
}
\right)\in\WH{Z}^p(X, r)_0.
$$
Since 
\begin{equation}
\partial \left(cl(\OV{\caF}, \varphi , \OV{\caG}, (y, (\omega, g))\right)
=\left(
\text{
\setlength{\unitlength}{1mm}
\begin{picture}(23,7)
 \put(-2,-4){$0$}
 \put(2,3){$0$}
 \put(6,-4){$0$}
 \put(10,3){$\CH_{T,0}^p(\OV{\caF})$}
\end{picture}
}
\right), \label{sec6:eq3} 
\end{equation}
it holds that $cl(\OV{\caF}, \varphi , \OV{\caG}, (y, (\omega, g))
\in \WH{Z}^p(X, r)_0^*$.

\vskip 1pc
\begin{thm}
In the quotient group, 
$$
[cl]=[cl(\OV{\caF}, \varphi , \OV{\caG}, (y, (\omega, g))]
\in \WH{Z}^p(X, r)^*_0/\IIm \partial 
$$
depends only on $\OV{\caF}$, in other words, $[cl]$ is 
independent of the choice of $\varphi , \OV{\caG}$ 
and $(y, (\omega, g))$.
\end{thm}

{\it Proof}: 
First we show that $[cl]$ is independent of the choice of 
representatives $(y, (\omega , \WT{g}))$ of 
$\WH{\CH}_0^p(\OV{\caG})$ and lifts $g$ of $\WT{g}$.
Fix a morphism $\varphi :T\to G$ and a hermitian vector bundle 
$\OV{\caG}$ such that $\varphi^*\caG \simeq \caF$.
Denote by $\scZ^p_{\varphi ,n}$ the set of admissible 
subschemes of $G\times \PP{n}$ such that 
$Y\in \scZ_{\varphi ,n}^p$ if and only if we can 
take the pull-back cycle of $[Y]$ by the morphism 
$$
X\times D_J\times D_K\hookrightarrow 
X\times \PP{r}\times \PP{n}
\overset{\varphi_I\times \IId }{\longrightarrow }G\times \PP{n}
$$
for any $I, J\subset \{1, \ldots, r\}$ and for any 
$K\subset \{1, \ldots , n\}$.
Let $\scD_{\log,\scZ^p_{\varphi ,n}}^*(G\times \PP{n},p)$ be 
the simple complex of the restriction map 
$$
\scD_{\log}^*(G\times \PP{n}, p)\to 
\scD_{\log}^*(G\times \PP{n}-\scZ^p_{\varphi ,n}, p)
$$
and 
$$
\mmD_{\log,\scZ^p_{\varphi ,n}}^*(G\times \PP{n}, p)=
\tau_{\leq 2p}\scD_{\log,\scZ^p_{\varphi ,n}}^*(G\times \PP{n}, p)
$$
the truncated subcomplex.
Moreover, set 
\begin{gather*}
\scH_{\varphi }^p(G, n)=
H^{2p}_{\scD,\scZ^p_{\varphi ,n}}(G\times \PP{n}, \R(p))=
H^{2p}(\mmD_{\log,\scZ^p_{\varphi ,n}}^*(G\times \PP{n}, p)),  \\
\scD_{\A,\scZ_{\varphi ,n}^p}^{s,-n}(G, p)=
\mmD_{\log,\scZ_{\varphi ,n}^p}^s(G\times \PP{n}, p).
\end{gather*}
These complexes have cubical structures with respect to the index $n$, 
hence we can obtain the normalized subcomplexes, which we denote by 
$\scH_{\varphi }^p(G, *)_0$ and 
$\scD_{\A,\scZ_{\varphi }^p}^{s,-*}(G, p)_0$ respectively.
Let $\scD_{\A,\scZ_{\varphi }^p}^*(G, p)_0$ be the single 
complex associated with $\scD_{\A,\scZ_{\varphi }^p}^{s,-r}(G, p)_0$.
Finally, let $\WH{Z}_{\varphi }^p(G, *)_0$ be the simple 
complex of the diagram 
$$
\xymatrix{
 & \scH_{\varphi }^p(G, *)_0 & & \WH{\scD}_{\A}^{2p-*}(G, p) \\
Z_{\varphi }^p(G, *)_0 \ar[ru]^{\chi_1} & & 
\scD_{\A,\scZ_{\varphi }^p}^{2p-*}(G, p)_0, 
\ar[lu]_{\chi_2} \ar[ru]^{\rho } & 
}
$$
where the maps $\chi_1, \chi_2$ and $\rho $ are defined in a 
similar way to the definition of $\WH{Z}^p(X, *)_0$.
We should note that in the definition of this complex we use 
$\WH{\scD}_{\A}^{2p-*}(G, p)$, not $\WH{\scD}_{\A,\P}^{2p-*}(G, p)$.
Consider the commutative diagram 
$$
\xymatrix{
Z_{\varphi }^p(G, *)_0 \ar[r]^{\chi_1} \ar[d] & 
\scH_{\varphi }^p(G, *)_0 \ar[d] & 
\scD_{\A,\scZ_{\varphi }^p}^{2p-*}(G, p)_0 
\ar[l]_{\chi_2} \ar[d]  \\
Z^p(G, *)_0 \ar[r]^{\chi_1} & \scH^p(G, *)_0 & 
\scD_{\A,\scZ^p}^{2p-*}(G, p)_0. \ar[l]_{\chi_2} 
}
$$
The left vertical arrow in the diagram is 
a quasi-isomorphism by the moving lemma.
Moreover, as shown in \cite[Prop.1.31, Prop.2.13]{BF}, 
the maps $\chi_2$ and $\chi_1\otimes \R$ on 
the both lines are quasi-isomorphisms.
Hence all the vertical maps are quasi-isomorphisms, which 
implies that the natural map 
$$
\WH{Z}_{\varphi }^p(G, *)_0\to \WH{Z}^p(G, *)_0
$$
is also a quasi-isomorphism.

Take another representative 
$(\PR{y}, (\PR{\omega }, \WT{\PR{g}}))$ of $\WH{\CH}_0^p(\OV{\caG})$ 
and a lift $\PR{g}$ of $\WT{\PR{g}}$.
Since the map 
$$
\WH{CH}^p(G, \scD(E_{\log}^*))\to \WH{CH}^p(G, 0)
$$
given by 
$$
[(y, (\omega , \WT{g}))]\mapsto 
\left[\left(
\text{
\setlength{\unitlength}{1mm}
\begin{picture}(17,7)
 \put(-2,-4){$y$}
 \put(2,3){$0$}
 \put(5,-4){$(\omega , g)$}
 \put(16,3){$0$}
\end{picture}
}
\right)\right]
$$
is an isomorphism \cite[Thm.4.8]{BF}, we can obtain an element 
$$
\left(
\text{
\setlength{\unitlength}{1mm}
\begin{picture}(13,7)
 \put(-2,-4){$z$}
 \put(2,3){$\beta_0 $}
 \put(6,-4){$\alpha $}
 \put(10,3){$\beta_1 $}
\end{picture}
}
\right)\in \WH{Z}_{\varphi }^p(G, 1)_0
$$
such that 
$$
\left(
\text{
\setlength{\unitlength}{1mm}
\begin{picture}(42,7)
 \put(-2,-4){$y-\PR{y}$}
 \put(10,3){$0$}
 \put(13,-4){$(\omega -\PR{\omega }, g-\PR{g})$}
 \put(40,3){$0$}
\end{picture}
}
\right)=\partial \left(
\text{
\setlength{\unitlength}{1mm}
\begin{picture}(13,7)
 \put(-2,-4){$z$}
 \put(2,3){$\beta _0$}
 \put(6,-4){$\alpha $}
 \put(10,3){$\beta _1$}
\end{picture}
}
\right).
$$
Consider the morphism 
$$
\varphi_I^T:X\times \PP{r+n}=X\times \PP{n}\times \PP{r}
\overset{\IId \times T}{\longrightarrow }X\times \PP{r}\times \PP{n}
\overset{\varphi_I\times \IId }{\longrightarrow }G\times \PP{n}, 
$$
where $T$ is the map given by 
$$
\PP{n}\times \PP{r}\to \PP{r}\times \PP{n}, \quad 
(a, b)\mapsto (b, a).
$$
Then the maps $(\varphi^T)^*$ which send $x$ to 
$\underset{I}{\sum}(-1)^{|I|}(\varphi_I^T)^*(x)$ 
form the commutative diagram 
$$
\xymatrix{
Z_{\varphi }^p(G, *)_0 \ar[r]^{\chi_1} \ar[d]^{(\varphi^T)^*} & 
\scH_{\varphi }^p(G, *)_0 \ar[d]^{(\varphi^T)^*} & 
\scD_{\A,\scZ_{\varphi ,*}^p}^{2p-*}(G, p)_0 
\ar[l]_{\chi_2}\ar[r]^{\rho } \ar[d]^{(\varphi^T)^*} & 
\WH{\scD}_{\A}^{2p-*}(G, p) \ar[d]^{(\varphi^T)^*}  \\
Z^p(X, *+r)_0 \ar[r]^{\chi_1} & \scH^p(X, *+r)_0 & 
\scD_{\A,\scZ_*^p}^{2p-*-r}(X, p)_0 \ar[l]_{\chi_2} 
\ar[r]^{\rho } & \WH{\scD}_{\A,\P}^{2p-*-r}(X, p).
}
$$
The vertical arrows in this diagram are maps of complexes except 
the right one at $*=1$.
Hence 
\begin{multline*}
\partial \left(
\text{
\setlength{\unitlength}{1mm}
\begin{picture}(58,7)
 \put(-2,-4){$(\varphi^T)^*(z)$}
 \put(12,3){$(\varphi^T)^*(\beta_0)$}
 \put(28,-4){$(\varphi^T)^*(\alpha )$}
 \put(43,3){$(\varphi^T)^*(\beta_1)$}
\end{picture}
}
\right) \\
=\left(
\text{
\setlength{\unitlength}{1mm}
\begin{picture}(83,7)
 \put(-2,-4){$\varphi^*(y)-\varphi^*(\PR{y})$}
 \put(24,3){$0$}
 \put(28,-4){$(\varphi^*(\omega )-\varphi^*(\PR{\omega }), 
\varphi^*(g)-\varphi^*(\PR{g}))$}
 \put(82,3){$0$}
\end{picture}
}
\right), 
\end{multline*}
which shows that $[cl]$ is independent of the choice of 
$(y, (\omega , \WT{g}))$ and lifts of $\WT{g}$.

We next show that the class $[cl]$ is independent of 
the choice of metrics on $\caG $.
To do this we need the following lemma:

\vskip 1pc
\begin{lem}
Let $\OV{\caE}$ be a short exact sequence of hermitian 
vector bundles on $G$, and 
$$
\CH_1^p(\OV{\caE})=\kappa_{\P}(\CH_1^p(\OV{\caE})_{\P})
\in \mmD^{2p-1}(G, p).
$$
Then in $\WT{\scD}_{\A,\P}^*(X, p)$ the difference 
$$
\CH_{T,1}^p(\varphi^*\OV{\caE})-(-1)^r\varphi^*\CH_1^p(\OV{\caE})
$$
is $d_s$-exact.
\end{lem}

{\it Proof}: 
Note that $\CH_1^p(\OV{\caE})_{\P}\in \WT{\scD}_{\P}^{2p,-1}(G, p)$ and 
$\CH_1^p(\OV{\caE})\in \WT{\scD}_{\P}^{2p-1,0}(G, p)$.
Since $\kappa_{\P}$ is a left inverse of the quasi-isomorphism 
$\mmD^*(G, p)\to \WT{\scD}_{\P}^*(G, p)$, 
the difference $\CH_1^p(\OV{\caE})_{\P}-\CH_1^p(\OV{\caE})$ 
is $d_s$-exact in $\WT{\scD}_{\P}^{2p-1}(G, p)$.
Hence there are 
$\alpha_i\in \WT{\scD}_{\P}^{2p+i-2,-i}(G, p)$ for $i=0, 1, 2$ 
such that 
\begin{gather*}
\CH_1^p(\OV{\caE})=\, -d_{\scD}\alpha_0+\delta_{\P}\alpha_1,  \\
\CH_1^p(\OV{\caE})_{\P}=\, d_{\scD}\alpha_1+\delta_{\P}\alpha_2.
\end{gather*}
Consider 
$\varphi^*\alpha_i\in \WT{\scD}_{\A,\P}^{2p+i-2,-r,-i}(X, p)$.
Since $\delta_{\A}\varphi^*\alpha_i=0$, 
\begin{align*}
d_s(\varphi^*\alpha_0+(-1)^r\varphi^*\alpha_1+\varphi^*\alpha_2)
=&\, d_{\scD}\varphi^*\alpha_0-\delta_{\P}\varphi^*\alpha_1+
(-1)^rd_{\scD}\varphi^*\alpha_1+(-1)^r\delta_{\P}\varphi^*\alpha_2  \\
=&\, -\varphi^*\CH_1^p(\OV{\caE})+
(-1)^r\varphi^*\CH_1^p(\OV{\caE})_{\P} \\
=&\, -\varphi^*\CH_1^p(\OV{\caE})
+(-1)^r\CH_{T,1}^p(\varphi^*\OV{\caE}), 
\end{align*}
which completes the proof.
\qed

\vskip 1pc
Let us go back to the proof of Thm.6.12.
Let $\PR{\OV{\caG}}$ be the same vector bundle as $\OV{\caG}$ 
with a different metric.
If we denote by 
$$
\CH_1^p(\PR{\OV{\caG}}, \OV{\caG})=
\CH_1^p(\PR{\OV{\caG}}\overset{\sim }{\to }\OV{\caG}\to 0), 
$$
which is the $(p-1,p-1)$-part of the Bott-Chern form, then 
by \cite[Thm.4.8]{gilletsoule2} we have 
$$
\WH{\CH}_0^p(\PR{\OV{\caG}})-\WH{\CH}_0^p(\OV{\caG})=
a(\WT{\CH}_1^p(\PR{\OV{\caG}}, \OV{\caG}))
$$
in $\WH{CH}^p(G, \scD (E_{\log}^*))$, where 
$$
a:\scD^{2p-1}(G, p)/\IIm d_{\scD}\to 
\WH{CH}^p(G, \scD (E_{\log}^*))
$$
is the map which sends $\WT{h}$ to $[(0, (d_{\scD}h, \WT{h}))]$.
This means that if $\WH{\CH}_0^p(\OV{\caG})$ is represented by 
$(y, (\omega , \WT{g}))$, then $\WH{\CH}_0^p(\PR{\OV{\caG}})$ 
is represented by $(y, (\omega +d_{\scD}h, \WT{g}+\WT{h}))$ 
where $h=\CH_1^p(\PR{\OV{\caG}}, \OV{\caG})$.
Hence 
\begin{align*}
&\left[cl(\OV{\caF}, \varphi , \OV{\caG}, (y, (\omega , g)))
\right]-\left[cl(\OV{\caF}, \varphi , \PR{\OV{\caG}}, 
(y, (\omega +d_{\scD}h, g+h)))\right] \\
&\hskip 3pc =\left[\left(
\text{
\setlength{\unitlength}{1mm}
\begin{picture}(105,7)
 \put(-2,-4){$0$}
 \put(2,3){$0$}
 \put(6,-4){$-(d_{\scD}\varphi^*(h), \varphi^*(h))$}
 \put(37,3){$(-1)^{r+1}
\left(\CH_{T,1}^p(\OV{\caF}, \varphi^*\OV{\caG})
-\CH_{T,1}^p(\OV{\caF}, \varphi^*\PR{\OV{\caG}})\right)$}
\end{picture}
}
\right)\right]  \\
&\hskip 3pc =(-1)^ra\left(\WT{\CH}_{T,1}^p(\OV{\caF}, 
\varphi^*\OV{\caG})-\WT{\CH}_{T,1}^p(\OV{\caF}, \varphi^*\PR{\OV{\caG}})
\right)-a(\varphi^*\WT{\CH}_1^p(\PR{\OV{\caG}}, \OV{\caG})), 
\end{align*}
and by Lem.6.13 it is equal to 
\begin{equation}
(-1)^ra\left(\WT{\CH}_{T,1}^p(\OV{\caF}, 
\varphi^*\OV{\caG})-\WT{\CH}_{T,1}^p(\OV{\caF}, \varphi^*\PR{\OV{\caG}})
-\WT{\CH}_{T,1}^p(\varphi^*\PR{\OV{\caG}}, \varphi^*\OV{\caG})\right).
\label{sec6:ele1}
\end{equation}

Consider the following exact hermitian $2$-cube on $T$: 
$$
\OV{\caC}=\left(
\begin{CD}
\OV{\caF} @>>> \varphi^*\PR{\OV{\caG}} \\
@VV{\IId }V  @VV{\IId }V  \\
\OV{\caF} @>>> \varphi^*\OV{\caG}
\end{CD}
\right), 
$$
and denote by $\CH_{T,2}^p(\OV{\caC})$ the component in 
$\WT{\scD}_{\A,\P}^{2p,r,2}(X)$ of $\CH_{T,2}(\OV{\caC})$.
Then 
\begin{align*}
d_s\CH_{T,2}^p(\OV{\caC})=&\, 
(-1)^r\delta_{\P}\CH_{T,2}^p(\OV{\caC})  \\
=&\, (-1)^r\left(\CH_{T,1}^p(\OV{\caF}, \varphi^*\PR{\OV{\caG}})
-\CH_{T,1}^p(\OV{\caF}, \varphi^*\OV{\caG})
+\CH_{T,1}^p(\varphi^*\PR{\OV{\caG}}, \varphi^*\OV{\caG})\right).
\end{align*}
This means that \eqref{sec6:ele1} is zero.
Hence we conclude that $[cl]$ is independent of the choice of metrics 
on $\caG$.

Finally we show that $[cl]$ is independent of the choice of morphisms 
$\varphi :T\to G$ and hermitian vector bundles $\OV{\caG}$ on $G$.
Take another morphism $\PR{\varphi }:T\to \PR{G}$ and 
a hermitian vector bundle $\OV{\PR{\caG}}$ on $\PR{G}$ such that 
${\PR{\varphi }}^*\PR{\caG}\simeq \caF$.
As shown in the proof of Thm.6.9, we may assume that there is 
a morphism $\psi :G\to \PR{G}$ such that $\psi \varphi =\PR{\varphi }$ 
and $\psi^*\PR{\caG}\simeq \caG$. 
Moreover, since we have shown that $[cl]$ is independent of 
the choice of metrics on $\caG$, we may assume that the isomorphism 
$\psi^*\PR{\caG}\simeq \caG$ preserves the metrics.
If $(\PR{y}, (\PR{\omega }, \WT{\PR{g}}))$ is an representative 
of $\WH{\CH}_0^p(\OV{\PR{\caG}})$ such that 
$\PR{y}\in Z_{\psi, \PR{\varphi }}^p(\PR{G})$ and 
$\PR{g}\in \scD^{2p-1}_{\log}(\PR{G}-\scZ_{\psi, \PR{\varphi }}^p, p)$, 
then $\WH{\CH}_0^p(\OV{\caG})$ is represented by 
$(\psi^*(\PR{y}), (\psi^*(\PR{\omega }), \WT{\psi^*(\PR{g})})$.
Moreover, the isometry $\psi^*\OV{\PR{\caG}}\simeq \OV{\caG}$ 
implies that $\CH_{T,1}^p(\OV{\caF}, {\PR{\varphi }}^*\OV{\PR{\caG}})
=\CH_{T,1}^p(\OV{\caF}, \varphi^*\OV{\caG})$.
Hence it follows that $[cl]$ is independent 
the choice of morphisms $\varphi :T\to G$ and hermitian vector bundles 
$\OV{\caG}$.
\qed

\vskip 1pc
\begin{defn}
We call the element 
$$
\WH{\CH}_{T,0}^p(\OV{\caF})=[cl(\OV{\caF}, \varphi , \OV{\caG}, 
(y, (\omega , g)))]\in \WH{Z}^p(X, r)^*_0/\IIm \partial 
$$
the $p$-th arithmetic Chern character of 
a hermitian vector bundle $\OV{\caF}$ on $T$.
\end{defn}

\setcounter{equation}{0}
\vskip 2pc
\section{Definition of higher arithmetic Chern character}

\vskip 1pc
In this section we construct a map from the higher 
arithmetic $K$-group to the higher arithmetic Chow group.
First we introduce some notations.
For $\eta \in \WT{\scD}_{\A,\P,r}(X)$, write 
$$
\eta =\sum_p\eta^p, \quad \eta^p\in \WT{\scD}_{\A,\P}^{2p-r}(X, p).
$$
Similarly, denote by $\CH_{T,0}^p$ the component in 
$\WT{\scD}_{\A,\P}^{2p-r}(X, p)$ of 
the map $\CH_{T,0}:\WH{K}_0(T)\to \WT{\scD}_{\A,\P,r}(X)$.

\vskip 1pc
\begin{prop}
We can define a map 
$$
\WH{\CH}_{T,0}^p:\WH{K}_0^M(T)\to \WH{Z}^p(X, r)^*_0/\IIm \partial 
$$
by 
$$
\WH{\CH}_{T,0}^p(\OV{\caF}, \WT{\eta })=\WH{\CH}_{T,0}^p(\OV{\caF})
+a(\WT{\eta^p}).
$$
\end{prop}

{\it Proof}: 
We have only to show that 
$$
\WH{\CH}_{T,0}^p(\OV{\caF}_{-1})+\WH{\CH}_{T,0}^p(\OV{\caF}_1)=
\WH{\CH}_{T,0}^p(\OV{\caF}_0)+(-1)^ra(\WT{\CH}_{T,1}^p(\OV{\caE}))
$$
for any short exact sequence of hermitian vector bundles 
$\OV{\caE}:0\to \OV{\caF}_{-1}\to \OV{\caF}_0\to \OV{\caF}_1\to 0$ 
on $T$.
It is shown in \cite[\S 3.2]{fulton1} that for such $\OV{\caE}$ 
there exist a morphism $\varphi :T\to G$ to a smooth projective 
variety $G$ and a short exact sequence of hermitian vector bundles 
$\OV{\PR{\caE}}:0\to \OV{\caG}_{-1}\to \OV{\caG}_0\to 
\OV{\caG}_1\to 0$ on $G$ such that 
$\varphi^*\PR{\caE}\simeq \caE$.
Take representatives $(y, (\omega , \WT{g}))$ of 
$\WH{\CH}_0^p(\OV{\caG}_{-1})$ and 
$(\PR{y}, (\PR{\omega }, \WT{\PR{g}}))$ 
of $\WH{\CH}_0^p(\OV{\caG}_1)$.
Since 
$$
\WH{\CH}_0^p(\OV{\caG}_{-1})+\WH{\CH}_0^p(\OV{\caG}_1)=
\WH{\CH}_0^p(\OV{\caG}_0)+a(\WT{\CH}_1^p(\OV{\PR{\caE}}))
$$
by \cite[Thm.4.8]{gilletsoule2}, 
$\WH{\CH}_0^p(\OV{\caG}_0)$ is represented by 
$$
(y+\PR{y}, (\omega +\PR{\omega }-d_{\scD}\CH_1^p(\OV{\PR{\caE}}), 
\WT{g}+\WT{\PR{g}}-\WT{\CH}_1^p(\OV{\PR{\caE}}))).
$$
Then 
\begin{align*}
\WH{\CH}_{T,0}^p(\OV{\caF}_{-1})&=\left[
\left(
\text{
\setlength{\unitlength}{1mm}
\begin{picture}(71,7)
 \put(-2,-4){$\varphi^*(y)$}
 \put(8,3){$0$}
 \put(12,-4){$(\varphi^*(\omega ), \varphi^*(g))$}
 \put(28,3){$(-1)^{r+1}\CH_{T,1}^p(\OV{\caF}_{-1}, 
\varphi^*\OV{\caG}_{-1})$}
\end{picture}
}
\right)\right],  \\
\WH{\CH}_{T,0}^p(\OV{\caF}_1)&=\left[
\left(
\text{
\setlength{\unitlength}{1mm}
\begin{picture}(67,7)
 \put(-2,-4){$\varphi^*(\PR{y})$}
 \put(8,3){$0$}
 \put(12,-4){$(\varphi^*(\PR{\omega }), \varphi^*(\PR{g}))$}
 \put(29,3){$(-1)^{r+1}\CH_{T,1}^p(\OV{\caF}_1, 
\varphi^*\OV{\caG}_1)$}
\end{picture}
}
\right)\right], \\
\WH{\CH}_{T,0}^p(\OV{\caF}_0)&=\left[
\left(
\text{
\setlength{\unitlength}{1mm}
\begin{picture}(78,7)
 \put(-2,-4){$\varphi^*(y+\PR{y})$}
 \put(15,3){$0$}
 \put(19,-4){$(\varphi^*(\omega +\PR{\omega }), \varphi^*(g+\PR{g}))$}
 \put(40,3){$(-1)^{r+1}\CH_{T,1}^p(\OV{\caF}_0, \varphi^*\OV{\caG}_0)$}
\end{picture}
}
\right)\right]-a(\varphi^*\CH_1^p(\OV{\PR{\caE}})). 
\end{align*}
Lem.6.13 says that $\varphi^*\CH_1^p(\OV{\PR{\caE}})=
(-1)^r\CH_{T,1}^p(\varphi^*\OV{\PR{\caE}})$, therefore 
\begin{align*}
&\WH{\CH}_{T,0}^p(\OV{\caF}_{-1})+\WH{\CH}_{T,0}^p(\OV{\caF}_1)-
\WH{\CH}_{T,0}^p(\OV{\caF}_0)  \\
&\hskip 1pc =(-1)^ra\left(
\WT{\CH}_{T,1}^p(\OV{\caF}_{-1}, \varphi^*\OV{\caG}_{-1})+
\WT{\CH}_{T,1}^p(\OV{\caF}_1, \varphi^*\OV{\caG}_1)-
\WT{\CH}_{T,1}^p(\OV{\caF}_0, \varphi^*\OV{\caG}_0)
+\WT{\CH}_{T,1}^p(\varphi^*\OV{\PR{\caE}})\right).
\end{align*}

Consider the following exact hermitian $2$-cube on $T$:
$$
\OV{\PR{\caC}}=\left(
\begin{CD}
\OV{\caF}_{-1} @>>> \OV{\caF}_0 @>>> \OV{\caF}_1  \\
@VVV @VVV @VVV \\
\varphi^*\OV{\caG}_{-1} @>>> \varphi^*\OV{\caG}_0 @>>> 
\varphi^*\OV{\caG}_1 
\end{CD}
\right).
$$
Then 
\begin{align*}
d_s\CH_{T,2}^p(\OV{\PR{\caC}})=&\, 
(-1)^r\delta_{\P}\CH_{T,2}^p(\OV{\PR{\caC}})  \\
=&\, (-1)^r\left(\CH_{T,1}^p(\OV{\caE})-\CH_{T,1}^p(\varphi^*\OV{\PR{\caE}})
-\CH_{T,1}^p(\OV{\caF}_{-1}, \varphi^*\OV{\caG}_{-1})\right.  \\
&\hskip 8pc \left.+\CH_{T,1}^p(\OV{\caF}_0, \varphi^*\OV{\caG}_0)
-\CH_{T,1}^p(\OV{\caF}_1, \varphi^*\OV{\caG}_1)\right).
\end{align*}
This implies that 
$$
\WH{\CH}_{T,0}^p(\OV{\caF}_{-1})+\WH{\CH}_{T,0}^p(\OV{\caF}_1)
-\WH{\CH}_{T,0}^p(\OV{\caF}_0)=(-1)^ra(\WT{\CH}_{T,1}^p(\OV{\caE})),
$$
which completes the proof.
\qed

\vskip 1pc
\begin{prop}
Suppose $r\geq 1$.
Define a map 
$$
\PR{\partial }:\WH{Z}^p(X, r)_0^*/\IIm \partial \to 
\WT{\scD}_{\A,\P}^{2p-r}(X, p)
$$
so that 
$$
\partial x=\left(
\text{
\setlength{\unitlength}{1mm}
\begin{picture}(13,7)
 \put(-2,-4){$0$}
 \put(2,3){$0$}
 \put(6,-4){$0$}
 \put(9,3){$\PR{\partial }x$}
\end{picture}
}
\right)
$$
for any $x\in \WH{Z}^p(X, r)_0^*$.
Then the diagram 
$$
\xymatrix{
\WH{K}^M_0(T) \ar[rr]^{\WH{\CH}_{T,0}^p} \ar[rd]_{\CH_{T,0}^p} & 
& \WH{Z}^p(X, r)_0^*/\IIm \partial \ar[ld]^{\PR{\partial }}  \\
& \WT{\scD}_{\A,\P}^{2p-r}(X, p) \\
}
$$
is commutative.
Hence taking the kernels of $\CH_{T,0}^p$ and $\PR{\partial }$ 
yields the map 
$$
\WH{\CH}_{T,0}^p:\WH{K}_0(T)\to \WH{CH}^p(X, r).
$$
\end{prop}

{\it Proof}: 
Let $(\OV{\caF}, \WT{\eta })$ be a pair of a hermitian vector 
bundle $\OV{\caF}$ on $T$ with \linebreak 
$\WT{\eta}\in \WT{\scD}_{\A,\P,r+1}(X)/\IIm d_s$.
Then we have seen in \eqref{sec6:eq3} that 
$$
\partial \, \WH{\CH}_{T,0}^p(\OV{\caF})=
\left(\text{
\setlength{\unitlength}{1mm}
\begin{picture}(22,7)
 \put(-2,-4){$0$}
 \put(2,3){$0$}
 \put(6,-4){$0$}
 \put(9,3){$\CH_{T,0}^p(\OV{\caF})$}
\end{picture}
}
\right).
$$
Hence 
$$
\partial \, \WH{\CH}_{T,0}^p([(\OV{\caF}, \WT{\eta })])=
\left(\text{
\setlength{\unitlength}{1mm}
\begin{picture}(34,7)
 \put(-2,-4){$0$}
 \put(2,3){$0$}
 \put(6,-4){$0$}
 \put(9,3){$\CH_{T,0}^p(\OV{\caF})+d_s\eta^p$}
\end{picture}
}
\right), 
$$
which completes the proof.
\qed

\vskip 1pc
\begin{thm}
There is a map 
$$
\WH{\CH}_r^p:\QQ{\WH{K}_r(X)}\to \WH{CH}^p(X, r)
$$
which makes the following diagram commutative:
$$
\xymatrix{
\QQ{\WH{K}_r(X)}\ar[r]^{\sim \qquad \qquad } \ar[dr]_{\WH{\CH}_r^p} & 
\QQ{\WH{K}_0(X\times \PP{r}; X\times \partial \PP{r})} & 
\QQ{\WH{K}_0(T; T_1, \ldots , T_r)} 
\ar@{->>}[l]_{\qquad \WH{i}_{\emptyset }^*} \ar[r] & 
\QQ{\WH{K}_0(T)}\ar[dll]^{\WH{\CH}_{T,0}^p}  \\
 & \WH{CH}^p(X, r). & & \\
}
$$
We call this map the {\it higher arithmetic Chern character} of $X$.
\end{thm}

{\it Proof}: 
We have only to show that the kernel of the surjection 
$\WH{i}_{\emptyset }^*$ goes to zero by the map 
$$
\QQ{\WH{K}_0(T; T_1, \ldots , T_r)}\subset \QQ{\WH{K}_0(T)}
\overset{\WH{\CH}_{T,0}^p}{\to }\WH{CH}^p(X, r).
$$
Recall the commutative diagram in Prop.5.4.
The exactness of the sequence 
$$
\QQ{K_1(X\times \PP{r}; X\times \partial \PP{r})}\to 
\underset{p}{\oplus }H_{\scD}^{2p-r-1}(X, \R(p))\to 
\QQ{\WH{K}_0(X\times \PP{r}; X\times \partial \PP{r})}
$$
with the bijectivity of $i_{\emptyset }^*:\QQ{K_0(T; T_1, \ldots , T_r)}
\to \QQ{K_0(X\times \PP{r}; X\times \partial \PP{r})}$  
implies that $\KER \WH{i}_{\emptyset }^*$ agrees with 
the image of the composite 
$$
\QQ{K_1(X\times \PP{r}; X\times \partial \PP{r})}\to 
\underset{p}{\oplus }H_{\scD}^{2p-r-1}(X, \R(p))\to 
\QQ{\WH{K}_0(T; T_1, \ldots , T_r)}.
$$
The commutative diagram in Cor.3.9 and Thm.3.4 imply that 
any element of $\KER \WH{i}_{\emptyset }^*$ is written as 
$[(0, \WT{\eta })]$ such that the element 
$\WT{\eta }\in \WT{\scD}_{\A,r+1}(X)/\IIm d_s$ 
in contained in the image of Beilinson's regulator 
$$
\QQ{K_{r+1}(X)}\to \underset{p}{\oplus }
H_{\scD}^{2p-r-1}(X, \R(p))\subset \WT{\scD}_{\A,r+1}(X)/\IIm d_s.
$$
On the other hand, recall the exact sequence 
$$
CH^p(X, r+1)\overset{H^{r+1}(\caP^p)}{\to }
H_{\scD}^{2p-r-1}(X, \R(p))\to \WH{CH}^p(X, r), 
$$
which is shown in \cite[Prop.4.4]{BF}.
Since $H^{r+1}(\caP^p)$ agrees with the regulator map 
by \cite[Thm.7.8]{BFT}, we conclude that 
$\WH{\CH}_{T,0}^p([(0, \WT{\eta })])=a(\WT{\eta^p})$ 
is zero in $\WH{CH}^p(X, r)$ if 
$[(0, \WT{\eta })]\in \KER \WH{i}_{\emptyset }^*$, 
which completes the proof.
\qed

\vskip 1pc
\begin{thm}
There is a commutative diagram up to sign:
{\small 
\begin{multline*}
\xymatrix{
\cdots \ar[r] &  \underset{p}{\oplus }H_{\scD}^{2p-r-1}(X, \R(p)) 
\ar[r] \ar[d]^{\PRJ} & \QQ{\WH{K}_r(X)} \ar[d]^{\WH{\CH}_r^p} \\
\cdots \ar[r] & H_{\scD}^{2p-r-1}(X, \R(p)) \ar[r] & \WH{CH}^p(X, r) \\
}   \\
\xymatrix{
\ar[r] & \QQ{K_r(X)} \ar[r]^{\CH_r^p \qquad } \ar[d]^{\CH_r^p} 
\ar@{}[rd] |{(*)} & \underset{p}{\oplus }H_{\scD}^{2p-r}(X, \R(p)) 
\ar[r] \ar[d]^{\PRJ} & \cdots \ar[r] & \QQ{K_1(X)} \ar[d]^{\CH_1^p} \\
\ar[r] & CH^p(X, r) \ar[r]^{H_r(\caP^p) \quad } & H_{\scD}^{2p-r}(X, \R(p)) 
\ar[r] & \cdots \ar[r] & CH^p(X, 1) \\
} \\
\xymatrix{
\ar[r]^{\CH_1^p \qquad \qquad \quad } \ar@{}[rd] |{(**)} & 
\underset{p}{\oplus }\mmD^{2p-1}(X, p)/\IIm d_{\scD} 
\ar[r] \ar[d]^{\PRJ} & \QQ{\WH{K}_0(X)} \ar[r] \ar[d]^{\WH{\CH}_0^p} 
& \QQ{K_0(X)} \ar[r] \ar[d]^{\CH_0^p} & 0 \\
\ar[r]^{H_1(\caP^p) \qquad \qquad } & \mmD^{2p-1}(X, p)/\IIm d_{\scD} 
\ar[r] & \WH{CH}^p(X, 0) \ar[r] & CH^p(X) \ar[r] &  0, 
} 
\end{multline*}
}
where $\CH_r^p:\QQ{K_r(X)}\to CH^p(X, r)$ is the map defined 
in Def.6.11, and $pr$ is the canonical projection.
\end{thm}

{\it Proof}: 
The commutativity of the diagram is straightforward except 
(*) and (**), and by Cor.4.7 the commutativity of them 
follows from the commutativity of the diagram 
$$
\xymatrix{
K_0(T) \ar[rr]^{\CH_{T,0}^p} \ar[rd]_{\CH_{T,0}^p} & &  
CH^p(X, r) \ar[ld]^{H_r(\caP^p)} \\
& H_{\scD}^{2p-r}(X, \R(p)). & 
}
$$
Let $\caF$ be a hermitian vector bundle on $T$.
Take a morphism $\varphi :T\to G$ and a vector bundle $\caG$ 
on $G$ with an isomorphism $\caF\simeq \varphi^*\caG$.
Moreover, let $y\in Z_{\varphi }^p(G)$ be a cycle representing 
$\CH_0^p(\caG)\in CH^p(G)$.
Then $\CH_{T,0}^p(\caF)=[\varphi^*(y)]$ in $CH^p(X, r)$.

Put a hermitian metric on $\caG$, and take the pull-back metric 
on $\caF$ by means of the isomorphism $\caF\simeq \varphi^*\caG$.
Denote by $\OV{\caG}$ and $\OV{\caF}$ the hermitian vector 
bundles obtained in this way.
Let $(\omega , \WT{g})\in GE_{\scZ_{\varphi }^p}^p(G)$ be a Green form 
associated with $y$ such that 
$\WH{\CH}_0^p(\OV{\caG})$ is represented by $(y, (\omega , \WT{g}))$.
Then since $d_{\scD}g=\omega $ and $\delta_{\A}\varphi^*(g)=0$, 
it follows from \eqref{sec6:eq2} that 
$$
d_{\scD}{\pi_X}_*[\varphi^*(g)\bullet \pi_{\P}^*W_r]=
{\pi_X}_*[\varphi^*(\omega )\bullet \pi_{\P}^*W_r]
-\caP^p(\varphi^*(y)).
$$
Since $\omega =\CH_0^p(\OV{\caG})$, it follows that 
$\varphi^*(\omega )=\CH_{T,0}^p(\OV{\caF})$, therefore 
$$
{\pi_X}_*[\varphi^*(\omega )\bullet \pi_{\P}^*W_r]=
\kappa_{\A}(\CH_{T,0}^p(\OV{\caF})).
$$
This means that $\caP^p(\varphi^*(y))$ and 
$\kappa_{\A}(\CH_{T,0}^p(\OV{\caF}))$ give the same cohomology 
class in \linebreak 
$H_{\scD}^{2p-r}(X, \R(p))$.
Since $\kappa_{\A}$ induces the identity on 
cohomology, we conclude that $\caP^p(\varphi^*(y))$ and 
$\CH_{T,0}^p(\OV{\caF})$ give the same cohomology 
class in $H_{\scD}^{2p-r}(X, \R(p))$.
This completes the proof.
\qed

\setcounter{equation}{0}
\vskip 2pc
\section{Compatibility with pull-back maps}

\vskip 1pc
\subsection{Pull-back maps on arithmetic $K$-groups}
Let $f:X\to Y$ be a morphism of smooth projective varieties 
defined over an arithmetic field.
Consider the pull-back map 
$$
f^*:\WT{\Q}\WALT{*}{Y\times \PP{r}; Y\times \partial \PP{r}}
\to \WT{\Q}\WALT{*}{X\times \PP{r}; X\times \partial \PP{r}}
$$
defined in Prop.2.19.
Prop.3.7 says that $(f^*)^{m,n}(x)$ is isometrically equivalent to 
a degenerate element for $m<n$.
Hence the diagram 
$$
\begin{CD}
\WT{\Q}\WALT{*}{Y\times \PP{r}; Y\times \partial \PP{r}}[r] 
@>{\CH_*}>> \WT{\scD}_{\A,\P,*}(Y)  \\
@V{f^*}VV  @VV{f^*}V  \\
\WT{\Q}\WALT{*}{X\times \PP{r}; X\times \partial \PP{r}}[r] 
@>{\CH_*}>> \WT{\scD}_{\A,\P,*}(X)
\end{CD}
$$
is commutative.
This diagram gives the pull-back map 
$$
\WH{f}^*:\QQ{\WH{K}_0(Y\times \PP{r}; Y\times \partial \PP{r})}
\to \QQ{\WH{K}_0(X\times \PP{r}; X\times \partial \PP{r})}.
$$
For two morphisms $f:X\to Y$ and $g:Y\to Z$, we have given 
in Prop.2.19 a homotopy $\varPhi $ from $(gf)^*$ to 
$g^*f^*$.
Prop.3.7 says that $\varPhi ^{m,n}(x)$ is isometrically equivalent 
to a degenerate element for any $m$ and $n$.
This implies that $\WH{g}^*\WH{f}^*=\WH{gf}^*$.
Moreover, Prop.2.20 says that the identity morphism of $X$ induces 
the identity of 
$\WT{\Q}\WALT{*}{X\times \PP{r}; X\times \partial \PP{r}}$, which 
implies that $\WH{\IId}_X^*=\IId $.
Since the alternating part of the commutative diagram in Cor.3.9 
$$
\begin{CD}
\WT{\Q}\WALT{*}{X} @>{i_X}>> 
\WT{\Q}\WALT{*}{X\times \PP{r}; X\times \partial \PP{r}}[r]  \\
@V{\CH_{*,\P}}VV  @VV{\CH_*}V \\
\WT{\scD}_{\P,*}(X) @>>> \WT{\scD}_{\A,\P,*}(X)
\end{CD}
$$
is compatible with the pull back maps $f^*$, 
the diagram 
$$
\begin{CD}
\QQ{\WH{K}_r(Y)} @>{\sim }>> \QQ{\WH{K}_{\P,r}(Y)} @>{\sim }>> 
\QQ{\WH{K}_0(Y\times \PP{r}; Y\times \partial \PP{r})}  \\
@V{\WH{f}^*}VV  @V{\WH{f}^*}VV  @VV{\WH{f}^*}V  \\
\QQ{\WH{K}_r(X)} @>{\sim }>> \QQ{\WH{K}_{\P,r}(X)} @>{\sim }>> 
\QQ{\WH{K}_0(X\times \PP{r}; X\times \partial \PP{r})}
\end{CD}
$$
is commutative.

Consider the iterated doubles 
$T=(X\times \PP{r}; X\times \partial \PP{r})$ and 
$U=D(Y\times \PP{r}; Y\times \partial \PP{r})$.
Then $f$ induces a morphism $f_D:T\to U$.
It is obvious from the definition that the Chern form 
$\CH_{U,0}(\OV{\caF})\in \WT{\scD}_{\A,\P,r}(Y)$ of a hermitian 
vector bundle $\OV{\caF}$ on $U$ satisfies 
$f^*\CH_{U,0}(\OV{\caF})=\CH_{T,0}(f_D^*\OV{\caF})$, 
and that the Bott-Chern form 
$\CH_{U,1}(\OV{\caE})\in \WT{\scD}_{\A,\P,r+1}(Y)$ 
of a short exact sequence $\OV{\caE}$ of hermitian vector 
bundles on $U$ satisfies 
$f^*\CH_{U,1}(\OV{\caE})=\CH_{T,1}(f_D^*\OV{\caE})$.
Hence we can define pull-back map of arithmetic $K$-groups 
$$
\WH{f}_D^*:\WH{K}_0^M(U)\to \WH{K}_0^M(T)
$$
by $[(\OV{\caF}, \WT{\eta })]\mapsto 
[(f_D^*\OV{\caF}, \WT{f^*(\eta )})]$.
Let $T_1, \ldots , T_r\subset T$ and 
$U_1, \ldots , U_r\subset U$ be the closed subschemes 
introduced in \S 4.2, and 
\begin{align*}
\iota_j:&\, U_j\hookrightarrow U, \quad \iota_j: 
T_j\hookrightarrow T,  \\
p_j:&\, U\to U_j, \quad p_j:T\to T_j
\end{align*}
the morphisms defined in \S 4.2.
It is obvious that the pull-back map defined above induces 
\begin{align*}
\WH{f}_D^*:\WH{K}_0^M(U; U_1, \ldots , U_r)&\, \to 
\WH{K}_0^M(T; T_1, \ldots , T_r),  \\
\WH{f}_D^*:\WH{K}_0(U)&\, \to \WH{K}_0(T), \\
\WH{f}_D^*:\WH{K}_0(U; U_1, \ldots , U_r)&\, \to 
\WH{K}_0(T; T_1, \ldots , T_r).
\end{align*}

For $1\leq j\leq r$, consider the diagram 
$$
\begin{CD}
\WT{\Q}\WALT{*}{U; U_1, \cdots , U_{j-1}} @>{\iota_j^*}>> 
\WT{\Q}\WALT{*}{U_j; U_1\cap U_j, \cdots , U_{j-1}\cap U_j} \\
@V{f_D^*}VV  @VV{f_D^*}V  \\
\WT{\Q}\WALT{*}{T; T_1, \cdots , T_{j-1}} @>{\iota_j^*}>> 
\WT{\Q}\WALT{*}{T_j; T_1\cap T_j, \cdots , T_{j-1}\cap T_j}.
\end{CD}
$$
Then Prop.2.19 implies that the family of maps 
$$
\varPhi_{\iota }^{m,n}(x)_J=(-1)^n
\underset{|I|=m}{\sum_{K\coprod I=J}}\sgn \sbinom{K \ I}{J}
\left(\varXi_{K,\iota_j,f_D}^{\alt}(x_I)-
\varXi_{K,f_D,\iota_j}^{\alt}(x_I)\right)
$$
for any $m$ and $n$ gives a homotopy 
$$
\varPhi_{\iota }:\WT{\Q}\WALT{*}{U; U_1, \ldots , U_{j-1}}\to 
\WT{\Q}\WALT{*+1}{T_j; T_1\cap T_j, \ldots , T_{j-1}\cap T_j}
$$
from $f_D^*\iota_j^*$ to $\iota_j^*f_D^*$.
The diagram 
$$
\begin{CD}
\WT{\Q}\WALT{*}{U_j; U_1\cap U_j, \ldots , U_{j-1}\cap U_j} 
@>{p_j^*}>> \WT{\Q}\WALT{*}{U; U_1, \ldots , U_{j-1}}   \\
@V{f_D^*}VV  @VV{f_D^*}V  \\
\WT{\Q}\WALT{*}{T_j; T_1\cap T_j, \ldots , T_{j-1}\cap T_j} 
@>{p_j^*}>> \WT{\Q}\WALT{*}{T; T_1, \ldots , T_{j-1}}, 
\end{CD}
$$
is also commutative up to homotopy, and a homotopy 
$$
\varPhi_p:
\WT{\Q}\WALT{*}{U_j; U_1\cap U_j, \ldots , U_{j-1}\cap U_j}
\to \WT{\Q}\WALT{*+1}{T; T_1, \ldots , T_{j-1}}
$$
from $f_D^*p_j^*$ to $p_j^*f_D^*$ is given by 
$$
\varPhi_p^{m,n}(x)_J=(-1)^n
\underset{|I|=m}{\sum_{K\coprod I=J}}\sgn \sbinom{K \ I}{J}
\left(\varXi_{K,p_j,f_D}^{\alt}(x_I)
-\varXi_{K,f_D,p_j}^{\alt}(x_I)\right).
$$
Moreover, let 
\begin{align*}
\varPsi_U:\WT{\Q}\WALT{*}{U_j; U_1\cap U_j, \ldots , U_{j-1}\cap U_j}
&\, \to \WT{\Q}\WALT{*+1}{U_j; U_1\cap U_j, \ldots , U_{j-1}\cap U_j}, \\
\varPsi_T:\WT{\Q}\WALT{*}{T_j; T_1\cap T_j, \cdots , T_{j-1}\cap T_j}
&\, \to \WT{\Q}\WALT{*+1}{T_j; T_1\cap T_j, \cdots , T_{j-1}\cap T_j}
\end{align*}
be the homotopies from the identity to $\iota_j^*p_j^*$ 
given in Prop.2.20.
Then we have the following:

\vskip 1pc
\begin{prop}
There is a second homotopy 
$$
\varTheta :\WT{\Q}\WALT{*}{U_j; U_1\cap U_j, \cdots , U_{j-1}\cap U_j} 
\to \WT{\Q}\WALT{*+2}{T_j; T_1\cap T_j, \cdots , T_{j-1}\cap T_j}
$$
from $\varPhi_{\iota }p_j^*+\iota_j^*\varPhi_p+f_D^*\varPsi_U$ 
to $\varPsi_Tf_D^*$ in the sense of Def.2.7.
\end{prop}

{\it Proof}: 
We begin by introducing a map of chain complexes of 
exact cubes associated with a sequence of morphisms 
$$
(X_1; Y_{1,1}, \ldots , Y_{1,r})\overset{f_1}{\to }
(X_2; Y_{2,1}, \ldots , Y_{2,r})\overset{f_2}{\to }
(X_3; Y_{3,1}, \ldots , Y_{3,r})\overset{f_3}{\to }
(X_4; Y_{4,1}, \ldots , Y_{4,r}). 
$$
For $1\leq i\leq 4$ and $J\subset \{1, \ldots , r\}$ with 
$k\notin J$, denote by $\iota_k$ the embedding \linebreak 
$Y_{i,J\cup \{k\}}\hookrightarrow Y_{i,J}$.
Define $\varXi_{K,f_1,f_2,f_3}:\WT{\Q}\WH{C}_*(Y_{4,I})\to 
\WT{\Q}\WH{C}_{*+n-m+2}(Y_{1,J})$ by 
\begin{align*}
\varXi_{K,f_1,f_2,f_3}=&\sum_{0\leq p\leq q\leq r\leq n-m}
(-1)^{p+q+r}\times  \\
&\hskip 2pc \sum_{\sigma \in \frS_{n-m}}(\sgn \sigma )
(\ldots , \iota_{k_{\sigma (p)}}, f_1, \ldots , 
\iota_{k_{\sigma (q)}}, f_2, \ldots , \iota_{k_{\sigma (r)}}, 
f_3, \ldots )^*.
\end{align*}
Then we can show in the same way as the proof of Prop.2.12 that 
\begin{align}
\partial \varXi_{K,f_1,f_2,f_3}&(x)+(-1)^{n-m+1}
\varXi_{K,f_1,f_2,f_3}(\partial x) \label{sec8:eq1} \\
=&-\sum_{a=1}^{n-m}\underset{|L|=a}{\sum_{L\coprod \PR{L}=K}}
\sgn \sbinom{L \ \PR{L}}{K}\varXi_L\varXi_{\PR{L},f_1,f_2,f_3}(x) 
\notag \\
&+\sum_{a=0}^{n-m}(-1)^a
\underset{|L|=a}{\sum_{L\coprod \PR{L}=K}}\sgn \sbinom{L \ \PR{L}}{K}
\varXi_{L,f_1}\varXi_{\PR{L},f_2,f_3}(x) \notag \\
&-\sum_{a=0}^{n-m}\underset{|L|=a}{\sum_{L\coprod \PR{L}=K}}
\sgn \sbinom{L \ \PR{L}}{K}\varXi_{L,f_1,f_2}\varXi_{\PR{L},f_3}(x) 
\notag \\
&+\sum_{a=0}^{n-m-1}(-1)^a
\underset{|L|=a}{\sum_{L\coprod \PR{L}=K}}\sgn \sbinom{L \ \PR{L}}{K}
\varXi_{L,f_1,f_2,f_3}\varXi_{\PR{L}}(x) \notag \\
&-\varXi_{K,f_2f_1,f_3}(x)+\varXi_{K,f_1,f_3f_2}(x).  \notag 
\end{align}
Denote
$\varXi_{K,f_1,f_2,f_3}^{\alt}=\Alt_*\varXi_{K,f_1,f_2,f_3}: 
\WALT{*}{Y_{4,I}}\to \WALT{*+n-m+2}{Y_{1,J}}$.

Let us go back to the proof of the proposition.
Define a map of $\scC$-complexes 
$$
\varTheta :\WT{\Q}\WALT{*}{U_j; U_1\cap U_j, \ldots , U_{j-1}\cap U_j}
\to \WT{\Q}\WALT{*}{T_j; T_1\cap T_j, \ldots , T_{j-1}\cap T_j}
$$
by 
$$
\varTheta^{m,n}(x)_J=-\sum_{K\coprod I=J}\sgn \sbinom{K \ I}{J}
\left(\varXi_{K,f_D,\iota_j,p_j}^{\alt}(x_I)
-\varXi_{K,\iota_j,f_D,p_j}^{\alt}(x_I)
+\varXi_{K,\iota_j,p_j,f_D}^{\alt}(x_I)\right).
$$
Note that $\varXi_{f_D, \IId }^{\alt}=0$, which we can 
show in the same way as Prop.2.20.
Using the equality \eqref{sec8:eq1}, we have 
\begin{align*}
(-1)^n\partial &\varTheta^{m,n}(x)_J
-(-1)^m\varTheta^{m,n}(\partial x_J) \\
=&-\sum_{a=1}^{n-m}F^{n-a,n}\varTheta^{m,n-a}(x)_J
+\sum_{a=0}^{n-m-1}\varTheta^{n-a,n}F^{m,n-a}(x)_J  \\
&-\sum_{a=0}^{n-m}(f_D^*)^{n-a,n}\varPsi_U^{m,n-a}(x)_J
-\sum_{a=0}^{n-m}(\iota_j^*)^{n-a,n}\varPhi_p^{m,n-a}(x)_J \\
&-\sum_{a=0}^{n-m}\varPhi_{\iota}^{n-a,n}(p_j^*)^{m,n-a}(x)_J
+\sum_{a=0}^{n-m}\varPsi_T^{n-a,n}(f_D^*)^{m,n-a}(x)_J, 
\end{align*}
which says that $\varTheta $ is a second homotopy from 
$\varPhi_{\iota }p_j^*+\iota_j^*\varPhi_p+f_D^*\varPsi_U$ 
to $\varPsi_Tf_D^*$.
\qed

\vskip 1pc
It follows from Prop.8.1 and Prop.2.8 that the diagram 
\begin{equation}
\begin{CD}
\WT{\Q}\WALT{*}{U; U_1, \cdots , U_{j-1}} @>{t_j}>> 
\WT{\Q}\WALT{*}{U; U_1, \cdots , U_j} \\
@V{f_D^*}VV  @VV{f_D^*}V  \\
\WT{\Q}\WALT{*}{T; T_1, \cdots , T_{j-1}} @>{t_j}>> 
\WT{\Q}\WALT{*}{T; T_1, \cdots , T_j}
\end{CD} \label{sec8:cd1}
\end{equation}
is commutative up to homotopy.
Denote by $\varPi_j$ the homotopy from $f_D^*t_j$ to 
$t_jf_D^*$ given in Prop.2.8.
It is obvious that $\varPhi_{\iota }^{m,n}(x), \varPhi_p^{m,n}(x), 
\varPsi_U^{m,n}(x), \varPsi_T^{m,n}(x)$ and $\varTheta^{m,n}(x)$ 
are isometrically equivalent to degenerate elements for any 
$x\in \underset{|I|=m}{\oplus }\WT{\Q}\WALT{*}{U_I}$.
Hence $\varPi_j^{m,n}(x)$ is also isometrically equivalent 
to a degenerate element.
Connecting the diagram \eqref{sec8:cd1} for all $j$, 
we obtain the following diagram 
$$
\begin{CD}
\WT{\Q}\WALT{*}{U} @>{t}>> 
\WT{\Q}\WALT{*}{U; U_1, \cdots , U_r} \\
@V{f_D^*}VV  @VV{f_D^*}V  \\
\WT{\Q}\WALT{*}{T} @>{t}>> 
\WT{\Q}\WALT{*}{T; T_1, \cdots , T_r}
\end{CD}
$$
which is commutative up to homotopy, and a homotopy from $f_D^*t$ to 
$tf_D^*$ is given by 
$$
\varPi=\sum_{j=1}^nt_n\cdots t_{j+1}\varPi_jt_{j-1}\cdots t_1.
$$
It is obvious that $\varPi^{0,n}(x)$ is isometrically equivalent to 
a degenerate element for $x\in \WT{\Q}\WALT{*}{U}$.

By Prop.2.19 the diagram 
$$
\begin{CD}
\WT{\Q}\WALT{*}{U; U_1, \cdots , U_r} @>{i^*_{\emptyset }}>> 
\WT{\Q}\WALT{*}{Y\times \PP{r}; Y\times \partial \PP{r}} \\
@V{f_D^*}VV  @VV{f^*}V  \\
\WT{\Q}\WALT{*}{T; T_1, \cdots , T_r} @>{i^*_{\emptyset }}>> 
\WT{\Q}\WALT{*}{X\times \PP{r}; X\times \partial \PP{r}} \\
\end{CD}
$$
is commutative up to homotopy, and a homotopy from 
$f^*i_{\emptyset }$ to $i_{\emptyset }^*f_D^*$ is given by 
\linebreak 
$\varPsi =\varPhi_{i_{\emptyset }, f_D}-\varPhi_{f, i_{\emptyset }}$.
Then $\varPsi^{m,n}(x)$ is isometrically equivalent to 
a degenerate element for any 
$x\in \underset{|I|=m}{\oplus }\WT{\Q}\WALT{*}{U_I}$.
Hence the diagram 
$$
\begin{CD}
\WT{\Q}\WALT{*}{U} @>{i^*_{\emptyset }t}>> 
\WT{\Q}\WALT{*}{Y\times \PP{r}; Y\times \partial \PP{r}} \\
@V{f_D^*}VV  @VV{f^*}V  \\
\WT{\Q}\WALT{*}{T} @>{i^*_{\emptyset }t}>> 
\WT{\Q}\WALT{*}{X\times \PP{r}; X\times \partial \PP{r}} \\
\end{CD} 
$$
is also commutative up to homotopy, and a homotopy from 
$f^*i_{\emptyset }^*t$ to $i_{\emptyset }^*tf_D^*$ is given by 
$\PR{\varPi }=\varPsi t+i^*_{\emptyset }\varPi $.
It is obvious that ${\PR{\varPi }}^{0,n}(x)$ is isometrically 
equivalent to a degenerate element for any 
$x\in \WT{\Q}\WALT{*}{U}$.
In particular, $\CH_*({\PR{\varPi }}^{0,n}(x))=0$.

Consider the diagram:
\begin{equation}
\begin{CD}
\QQ{\WH{K}_0(U)} @>{\WH{i_{\emptyset }^*t}}>> 
\QQ{\WH{K}_0(Y\times \PP{r}; Y\times \partial \PP{r})}  \\
@V{\WH{f}_D^*}VV  @VV{\WH{f}^*}V  \\
\QQ{\WH{K}_0(T)} @>{\WH{i_{\emptyset }^*t}}>> 
\QQ{\WH{K}_0(X\times \PP{r}; X\times \partial \PP{r})}, 
\end{CD} \label{sec8:cd2}
\end{equation}
where the maps $\WH{i_{\emptyset }^*t}$ are 
defined in \eqref{sec5:map2}.
Let $\OV{\caF}$ be a virtual hermitian vector bundle on $U$ and 
$\eta \in \WT{\scD}_{\A,\P,r+1}(Y)$ such 
that $\CH_{T,0}(\OV{\caF})+d_s\eta =0$.
Then 
\begin{align*}
\WH{f}^*\WH{i_{\emptyset }^*t}\, [(\OV{\caF}, \WT{\eta })]
=&\, [(f^*i_{\emptyset }^*t(\OV{\caF}), -f^*(\eta ))],  \\
\WH{i_{\emptyset }^*t}\WH{f}_D^*[(\OV{\caF}, \WT{\eta })]
=&\, [(i_{\emptyset }^*tf^*(\OV{\caF}), -f^*(\eta ))].
\end{align*}
On the other hand, since $\CH_*({\PR{\varPi }}^{0,n}(x))=0$ for any 
$x\in \WT{\Q}\WALT{*}{U}$, 
$$
\partial (\PR{\varPi }(\OV{\caF}), 0)=
(f_D^*i_{\emptyset }^*t(\OV{\caF})-
i_{\emptyset }^*tf^*(\OV{\caF}), 0)
$$
in $s(\CH_*)_r$.
This means that the diagram \eqref{sec8:cd2} is commutative.
Restricting \eqref{sec8:cd2} to the relative $K$-theories, 
we have the commutative diagram 
$$
\begin{CD}
\QQ{\WH{K}_0(U; U_1, \ldots , U_r)} @>{\WH{i}_{\emptyset }^*}>> 
\QQ{\WH{K}_0(Y\times \PP{r}; Y\times \partial \PP{r})}  \\
@V{\WH{f}_D^*}VV  @VV{\WH{f}^*}V  \\
\QQ{\WH{K}_0(T; T_1, \ldots , T_r)} @>{\WH{i}_{\emptyset }^*}>> 
\QQ{\WH{K}_0(X\times \PP{r}; X\times \partial \PP{r})}.
\end{CD}
$$

To Sum up, we obtain the following proposition:

\vskip 1pc
\begin{prop}
The diagram 
$$
\begin{CD}
\QQ{\WH{K}_r(Y)} @>{\sim }>>  
\QQ{\WH{K}_0(Y\times \PP{r}; Y\times \partial \PP{r})} 
@<{\WH{i}_{\emptyset }^*}<<  \QQ{\WH{K}_0(U; U_1, \ldots , U_r)} 
@>>> \QQ{\WH{K}_0(U)} \\
@V{\WH{f}^*}VV @V{\WH{f}^*}VV @VV{\WH{f}_D^*}V @VV{\WH{f}_D^*}V  \\
\QQ{\WH{K}_r(X)} @>{\sim }>> 
\QQ{\WH{K}_0(X\times \PP{r}; X\times \partial \PP{r})} 
@<<{\WH{i}_{\emptyset }^*}< \QQ{\WH{K}_0(T; T_1, \ldots , T_r)} 
@>>> \QQ{\WH{K}_0(T)} 
\end{CD}
$$
is commutative.
\end{prop}

\vskip 1pc
\subsection{The main theorem}
We begin by recalling the pull-back map of higher arithmetic 
Chow groups \cite{BF}.
Let $f:X\to Y$ be a morphism of smooth projective varieties 
defined over an arithmetic field.
Let $Z_f^p(Y, *)$ be the subcomplex of $Z^p(Y, *)$ such that 
$y\in Z_f^p(Y, *)$ if and only if one can take the pull-back cycle 
$f^*(y)\in Z^p(X, *)$.
Then the moving lemma says that the inclusion 
$Z_f^p(Y, *)_0\hookrightarrow Z^p(Y, *)_0$ is a quasi-isomorphism.
We can define complexes $\scH_f^p(Y, *)_0$ and 
$\scD_{\A,\scZ_f^p}^*(Y, p)_0$ in the same way as in \S 6.3.
Let $\WH{Z}_f^p(Y, *)_0$ be the simple complex of the diagram 
\begin{equation}
\xymatrix{
 & \scH_f^p(Y, *)_0 & & \WH{\scD}_{\A,\P}^{2p-*}(Y, p)_0 \\
Z_f^p(Y, *)_0 \ar[ru]^{\chi_1} & & 
\scD_{\A,\scZ_f^p}^{2p-*}(Y, p)_0. \ar[lu]_{\chi_2} 
\ar[ru]^{\rho } & 
 } \label{sec8:dia1}
\end{equation}
Then the natural inclusion 
$$
\WH{Z}_f^p(Y, *)_0\hookrightarrow \WH{Z}^p(Y, *)_0
$$
is a quasi-isomorphism, and collecting the pull-back maps on 
the complexes in \eqref{sec8:dia1} yields the map 
$$
f^*:\WH{Z}_f^p(Y, *)_0\to \WH{Z}^p(X, *)_0.
$$
Taking the maps on homology we obtain the pull-back map on 
higher arithmetic Chow groups:
\begin{equation}
f^*:\WH{CH}^p(Y, r)\overset{\sim }{\leftarrow }
H_r(\WH{Z}_f^p(Y, *)_0)\overset{f^*}{\to }
\WH{CH}^p(X, r). \label{sec8:map1}
\end{equation}

\vskip 1pc
\begin{thm}
The diagram 
$$
\begin{CD}
\QQ{\WH{K}_r(Y)} @>{\WH{\CH}_r^p}>> \WH{CH}^p(Y, r)  \\
@V{\WH{f}^*}VV  @VV{\WH{f}^*}V  \\
\QQ{\WH{K}_r(X)} @>{\WH{\CH}_r^p}>> \WH{CH}^p(X, r)
\end{CD}
$$
is commutative.
\end{thm}

{\it Proof}: 
Substituting the complex 
$\sigma_{<2p-r}\WH{\scD}_{\A,\P}^{2p-*}(Y, p)_0$ for 
$\WH{\scD}_{\A,\P}^{2p-*}(Y, p)_0$ in the diagram \eqref{sec8:dia1}, 
we obtain a map 
$$
\WH{f}^*:\WH{Z}^p(Y, r)^*_0/\IIm \partial 
\overset{\sim }{\leftarrow } 
\WH{Z}_f^p(Y, r)^*_0/\IIm \partial \overset{\WH{f}^*}{\to }
\WH{Z}^p(X, *)_0/\IIm \partial , 
$$
which is an extension of \eqref{sec8:map1}.
Then the definition of $\WH{\CH}^p_r$ and Prop.8.2 imply that 
the theorem follows from the commutativity of the diagram 
\begin{equation}
\begin{CD}
\WH{K}^M_0(U) @>{\WH{\CH}_{U,0}^p}>> 
\WH{Z}^p(Y, r)_0^*/\IIm \partial  \\
@V{\WH{f}_D^*}VV  @VV{\WH{f}^*}V  \\
\WH{K}^M_0(T) @>{\WH{\CH}_{T,0}^p}>> \WH{Z}^p(X, r)_0^*/\IIm \partial , 
\end{CD} 	\label{sec8:cd3}
\end{equation}
and it is equivalent to that 
$\WH{f}^*\WH{\CH}_{U,0}(\OV{\caF})=\WH{\CH}_{T,0}(f_D^*\OV{\caF})$ 
for any hermitian vector bundle $\OV{\caF}$ on $U$.

Take a morphism $\varphi :U\to G$ to a smooth projective variety 
and a hermitian vector bundle $\OV{\caG}$ on $G$ 
with an isomorphism $\varphi^*\caG\simeq \caF$.
Let $Z^p_{f,\varphi }(G)$ be the subgroup of $Z^p(G)$ such that 
$z\in Z^p_{f,\varphi }(G)$ if and only if 
one can define the pull-back cycle of $z$ by the morphism 
$$
Y\times D_J\hookrightarrow Y\times \PP{r}
\overset{\varphi_I}{\to }G, 
$$
and also by the morphism 
$$
X\times D_J\hookrightarrow X\times \PP{r}\overset{f\times 1}{\to }
Y\times \PP{r}\overset{\varphi_I}{\to }G
$$
for any $I, J\subset \{1, \ldots , r\}$.
Let $(y, (\omega , \WT{g}))$ be a pair of $y\in Z^p_{f,\varphi }(G)$ 
and a Green form $(\omega , \WT{g})\in GE^p_{\scZ^p_{f,\varphi }}(G)$ 
associated with $y$ representing 
$\WH{\CH}_0^p(\OV{\caG})\in \WH{CH}^p(G, \scD (E_{\log}^*))$.
Then 
$$
\WH{\CH}_{U,0}^p(\OV{\caF})=\left[
\left(
\text{
\setlength{\unitlength}{1mm}
\begin{picture}(69,7)
 \put(-2,-4){$\varphi^*(y)$}
 \put(9,3){$0$}
 \put(13,-4){$(\varphi^*(\omega ), \varphi^*(g))$}
 \put(33,3){$(-1)^{r+1}\CH_{U,1}^p(\OV{\caF}, \varphi^*\OV{\caG})$}
\end{picture}
}
\right)\right], 
$$
which is an element of $\WH{Z}_f^p(Y, r)_0^*/\IIm \partial $.
Since $f_D^*\varphi^*\caG\simeq f_D^*\caF$,  
\begin{align*}
\WH{f}^*\WH{\CH}_{U,0}^p(\OV{\caF})=&\, \left[
\left(
\text{
\setlength{\unitlength}{1mm}
\begin{picture}(83,7)
 \put(-2,-4){$f^*\varphi^*(y)$}
 \put(13,3){$0$}
 \put(16,-4){$(f^*\varphi^*(\omega ), f^*\varphi^*(g))$}
 \put(39,3){$(-1)^{r+1}\CH_{T,1}^p
(f_D^*\OV{\caF}, f_D^*\varphi^*\OV{\caG})$}
\end{picture}
}
\right)\right]  \\
=&\, \WH{\CH}_{T,0}^p(f_D^*\OV{\caF}), 
\end{align*}
which completes the proof.
\qed

\setcounter{equation}{0}
\vskip 2pc
\section{A tensor product structure on the multi-relative complexes 
of exact cubes}

\vskip 1pc
\subsection{An exact cube $\Sps{\caF_0, \ldots, \caF_l}$}
Let $\frA$ be a small exact category.
For a sequence of isomorphisms 
$\caF_0\simeq \caF_1\simeq \cdots \simeq \caF_l$ 
of objects of $\frA$, define an exact $l$-cube 
$\Sps{\caF_0, \ldots , \caF_l}$ of $\frA$ as follows:
If $\alpha_j=1$ for some $j$, then 
$$
\Sps{\caF_0, \ldots , \caF_l}_{\alpha_1, \ldots , \alpha_l}
=0.
$$
On the other hand, if $\alpha_j=-1$ and 
$\alpha_{j+1}=\cdots =\alpha_l=0$, then 
$$
\Sps{\caF_0, \ldots , \caF_l}_{\alpha_1, \ldots , \alpha_l}
=\caF_{l-j}.
$$
The maps in $\Sps{\caF_0, \ldots , \caF_l}$ are the zero maps, 
the identities, or composites of the isomorphisms 
in the sequence.
When $l=0$, $\Sps{\caF_0}$ is supposed to be 
the $0$-cube $\caF_0$.
For instance, $\Sps{\caF_0, \caF_1, \caF_2}$ is described as 
$$
\xymatrix{
\caF_0 \ar[d] \ar[r] & \caF_1 \ar[d] \ar[r] & 0 \ar[d] \\
\caF_0 \ar[d] \ar[r] & \caF_2 \ar[d] \ar[r] & 0 \ar[d] \\
0 \ar[r] & 0 \ar[r] & \ 0.  \\
}
$$

Let us consider the faces of $\Sps{\caF_0, \ldots , \caF_l}$.
Firstly, 
$$
\partial_j^1\Sps{\caF_0, \ldots , \caF_l}=0
$$
for any $1\leq j\leq l$.
Moreover, 
\begin{align*}
\partial_1^{-1}\Sps{\caF_0, \ldots , \caF_l}
=&\Sps{\caF_0, \ldots , \caF_{l-1}},  \\
\partial_1^0\Sps{\caF_0, \ldots , \caF_l}
=&\Sps{\caF_0, \ldots , \caF_{l-2}, \caF_l}, 
\end{align*}
and if $j\geq 2$, then 
$\partial_j^{-1}\Sps{\caF_0, \ldots , \caF_l}$ 
is a degenerate cube and 
$$
\partial_j^0\Sps{\caF_0, \ldots , \caF_l}=\Sps{\caF_0, 
\ldots , \caF_{l-j-1}, \caF_{l-j+1}, \ldots , \caF_l}.
$$
Hence it holds in $\WT{\Q}C_*(\frA)$ that 
$$
\partial \Sps{\caF_0, \ldots , \caF_l}=
\sum_{j=0}^l(-1)^j\Sps{\caF_0, \ldots , 
\caF_{l-j-1}, \caF_{l-j+1}, \ldots , \caF_l}.
$$

We can generalize this construction to a sequence of exact cubes.
If $\caF_0\simeq \cdots \simeq \caF_l$ is a sequence of 
isomorphisms of exact $n$-cubes of $\frA$, then we can obtain 
an exact $(n+l)$-cube $\Sps{\caF_0, \ldots , \caF_l}$ so that 
$$
\partial_{l+1}^{\alpha_1}\cdots \partial_{l+n}^{\alpha_n}
\Sps{\caF_0, \ldots , \caF_l}=\Sps{
(\caF_0)_{\alpha_1, \ldots , \alpha_n}, \ldots , 
(\caF_l)_{\alpha_1, \ldots , \alpha_n}}.
$$
Then it holds in $\WT{\Q}C_*(\frA)$ that 
\begin{align*}
\partial \Sps{\caF_0, \ldots , \caF_l}=&\sum_{j=0}^l
(-1)^j\Sps{\caF_0, \ldots , \caF_{l-j-1}, 
\caF_{l-j+1}, \ldots , \caF_l}  \\
&+\sum_{j=1}^n\sum_{i=-1}^1(-1)^{l+i+j}\Sps{
\partial_j^i\caF_0, \ldots , \partial_j^i\caF_l}.
\end{align*}

\vskip 1pc
\subsection{An exact cube 
$\Sps{\caF; \varphi_1, \varphi_2, \ldots , \varphi_l}(\caG)$}
We begin by recalling the definition of tensor product of exact cubes.
Let $\frA$ be a small exact category, and assume that $\frA$ 
is equipped with tensor product.
Given an exact $n$-cube $\caF$ and an exact $m$-cube $\caG$ 
of $\frA$, define an exact $(n+m)$-cube $\caF\otimes \caG$ by 
$$
(\caF\otimes \caG)_{\alpha_1, \ldots , \alpha_{n+m}}=
\caF_{\alpha_1, \ldots , \alpha_n}\otimes 
\caG_{\alpha_{n+1}, \ldots , \alpha_{n+m}}.
$$
Then it gives a product of chain complexes 
$$
\otimes :\WT{\Q}C_*(\frA)\otimes \WT{\Q}C_*(\frA)\to 
\WT{\Q}C_*(\frA).
$$

Let $\frA_0, \ldots , \frA_l$ and $\frB$ be small exact 
categories, and consider the diagram of functors 
$$
\xymatrix{
\frA_0 & \frA_1\ar[l]_{\varphi_1} & \cdots \ar[l]_{\varphi_2} 
& \frA_l\ar[l]_{\varphi_l}  \\
& \frB, \ar[lu]^{\pi_0}\ar[u]^{\pi_1}\ar[rru]_{\pi_l} \\
}
$$
where each $\pi_p$ is an exact functor and each 
$\varphi_p$ is an exact functor from $\frA_p$ to the category 
of exact $s_p$-cubes of $\frA_{p-1}$.
Then we can extend $\varphi_p$ to an exact functor from 
the category of exact $n$-cubes of $\frA_p$ to that 
of exact $(s_p+n)$-cubes of $\frA_{p-1}$ in the way that 
$$
\partial_{s_p+1}^{\alpha_1}\cdots 
\partial_{s_p+n}^{\alpha_n}\varphi_p(\caG)=
\varphi_p(\caG_{\alpha_1, \ldots , \alpha_n})
$$
for any exact $n$-cube $\caG$ of $\frA_p$.

Suppose that each $\frA_p$ is equipped with tensor product.
Fix an object $\caF$ of $\frB$, and we abbreviate 
$\pi_p(\caF)\otimes \caG$ to $\caF\otimes \caG$ for any 
exact cube $\caG$ of $\frA_p$.
Moreover, assume that there is a natural transformation 
$$
\caF\otimes \varphi_p(\caG)\simeq \varphi_p(\caF\otimes \caG).
$$
Set $s=\sum_{p=1}^{l}s_p$.
Then for any exact $n$-cube $\caG$ of $\frA_l$, we have 
a sequence of isomorphisms of exact $(n+s)$-cubes 
$$
\caF\otimes \varphi_1\varphi_2\cdots \varphi_l(\caG)\simeq 
\varphi_1(\caF\otimes \varphi_2\cdots \varphi_l(\caG))\simeq 
\ldots \simeq 
\varphi_1\varphi_2\cdots \varphi_l(\caF\otimes \caG)
$$
of $\frA_0$, and the associated exact 
$(n+s+l)$-cube of $\frA_0$: 
\begin{align*}
&\Sps{\caF; \varphi_1, \varphi_2, \ldots , \varphi_l}(\caG) \\
&\hskip 2pc =
\Sps{\caF\otimes \varphi_1\varphi_2\cdots \varphi_l(\caG), 
\varphi_1(\caF\otimes \varphi_2\cdots \varphi_l(\caG)), 
\cdots , 
\varphi_1\varphi_2\cdots \varphi_l(\caF\otimes \caG)}.
\end{align*}
When $l=0$, $\Sps{\caF; \ }(\caG)$ is supposed to be 
the tensor product $\caF\otimes \caG$.
Since \linebreak 
$\Sps{\caF; \varphi_1, \varphi_2, \ldots , \varphi_l}(\caG)$ 
is degenerate if so is $\caG$, it induces a map 
$$
\Sps{\caF; \varphi_1, \varphi_2, \ldots , \varphi_l}: 
\WT{\Q}C_*(\frA_l)\to \WT{\Q}C_{*+s+l}(\frA_0).
$$

We can generalize this construction to the case that $\varphi_p$ 
is a linear sum of exact functors from $\frA_p$ to the category 
of exact cubes of $\frA_{p-1}$.
In particular we have 
\begin{align*}
&\Sps{\caF; \varphi_1, \ldots , \partial \varphi_p, 
\ldots , \varphi_l}  \\
&\hskip 2pc =\sum_{j=1}^{s_p}\sum_{i=-1}^1(-1)^{i+j}
\Sps{\caF; \varphi_1, \ldots , \partial_j^i\varphi_p, 
\ldots , \varphi_l}:\WT{\Q}C_*(\frA_l)\to 
\WT{\Q}C_{*+s+l-1}(\frA_0).
\end{align*}

Let us consider the faces of the exact cube 
$\Sps{\caF; \varphi_1, \varphi_2, \ldots , \varphi_l}(\caG)$.
As we have seen in the previous subsection, 
$\partial_j^i\Sps{\caF; \varphi_1, \varphi_2, \ldots , \varphi_l}(\caG)$ 
is degenerate if $1\leq j\leq l$ and $i=1$, or if $2\leq j\leq l$ 
and $i=0$.
On the other hand, 
$$
\partial_1^{-1}\Sps{\caF; \varphi_1, \ldots , 
\varphi_l}(\caG)=\Sps{\caF; \varphi_1, \ldots , 
\varphi_{l-1}}(\varphi_l(\caG)), 
$$
and for $1\leq j\leq l-1$, 
$$
\partial_j^0\Sps{\caF; \varphi_1, \ldots , \varphi_l}
(\caG)=\Sps{\caF; \varphi_1, \ldots , 
\varphi_{l-j}\varphi_{l-j+1}, \ldots , \varphi_l}(\caG).
$$
However, $\partial_l^0\Sps{\caF; \varphi_1, \ldots , 
\varphi_l}(\caG)$ is not equal to 
$\varphi_1(\Sps{\caF; \varphi_2, \ldots , \varphi_l}(\caG))$.  
In fact, 
$$
\partial_l^0\Sps{\caF; \varphi_1, \ldots , \varphi_l}(\caG)
=\sigma \left(\varphi_1(\Sps{\caF; \varphi_2, \ldots , 
\varphi_l}(\caG))\right)
$$
where $\sigma \in \frS_{n+s+l-1}$ is the transposition of 
the sequence $\{1, 2, \ldots , l-1\}$ with the adjacent one 
$\{l, \ldots , l+s_1-1\}$.
Hence we have 
\begin{align*}
\partial \left(\Sps{\caF; \varphi_1, \ldots , \varphi_l}(\caG)
\right)=&\Sps{\caF; \varphi_1, \ldots , \varphi_{l-1}}
(\varphi_l(\caG))  \\
&+\sum_{j=1}^{l-1}(-1)^j\Sps{\caF; \varphi_1, \cdots , 
\varphi_{l-j}\varphi_{l-j+1}, \ldots , \varphi_l}(\caG) \\
&+(-1)^l\sigma \left(\varphi_1(\Sps{\caF; \varphi_2, \ldots , 
\varphi_l}(\caG))\right)   \\
&+(-1)^l\sum_{j=1}^l(-1)^{s_1+\cdots +s_{j-1}}
\Sps{\caF; \varphi_1, \ldots , \partial \varphi_j, 
\ldots , \varphi_l}(\caG)  \\
&+(-1)^{s+l}\Sps{\caF; \varphi_1, \ldots , \varphi_l}
(\partial \caG).
\end{align*}

If we use the chain complex $\WT{\Q}\ALT{*}{\frA}$, then 
we can get rid of the action of $\sigma \in \frS_*$ 
from the above expression.
Let 
$$
\varphi_p^{\alt }=\Alt_*\varphi_p:\WT{\Q}\ALT{*}{\frA_p}\to 
\WT{\Q}\ALT{*+s_p}{\frA_{p-1}}, 
$$
and 
$$
\Sps{\caF; \varphi_1, \varphi_2, \ldots , \varphi_l}^{\alt }=
\Alt_*\Sps{\caF; \varphi_1, \varphi_2, \ldots , \varphi_l}:
\WT{\Q}\ALT{*}{\frA_l}\to \WT{\Q}\ALT{*+s+l}{\frA_0}.
$$
Since the signature of $\sigma $ is $(-1)^{s_1(l-1)}$, 
we have the following:

\vskip 1pc
\begin{prop}
For $x\in \WT{\Q}\ALT{*}{\frA_l}$, 
\begin{align*}
\partial \left(\Sps{\caF; \varphi_1, \ldots , \varphi_l}^{\alt }
(x)\right)=&\Sps{\caF; \varphi_1, \ldots , 
\varphi_{l-1}}^{\alt }(\varphi_l^{\alt }(x))  \\
&+\sum_{j=1}^{l-1}(-1)^j\Sps{\caF; \varphi_1, \cdots , 
\varphi_{l-j}\varphi_{l-j+1}, \ldots , \varphi_l}^{\alt }(x) \\
&+(-1)^{s_1(l-1)+l}\varphi_1^{\alt }(\Sps{\caF; \varphi_2, 
\ldots , \varphi_l}^{\alt }(x))   \\
&+(-1)^l\sum_{j=1}^l(-1)^{s_1+\cdots +s_{j-1}}
\Sps{\caF; \varphi_1, \ldots , \partial \varphi_j, 
\ldots , \varphi_l}^{\alt }(x)  \\
&+(-1)^{s+l}\Sps{\caF; \varphi_1, \ldots , \varphi_l}^{\alt }
(\partial x).
\end{align*}
\end{prop}

\vskip 1pc
\subsection{A tensor product structure on 
$\WT{\Q}\ALT{*}{X; Y_1, \ldots , Y_r}$}
Let $X$ be a scheme and \linebreak 
$Y_1, \cdots , Y_r$ closed subschemes of $X$.
Assume that $X$ is defined over a base scheme $S$.
Fix a vector bundle $\caF$ on $S$.
For any exact cube $\caG$ on $Y_I$, we write $\caF \otimes \caG$ 
for $\pi_I^*\caF\otimes \caG$, where $\pi_I:Y_I\to S$ is 
the structure map.
In this subsection, we will construct a map of $\scC$-complexes 
$$
(\caF \otimes \ ):\WT{\Q}\ALT{*}{X; Y_1, \ldots , Y_r}\to 
\WT{\Q}\ALT{*}{X; Y_1, \ldots , Y_r}
$$
such that 
$$
(\caF \otimes \ )^{m, m}:
\underset{|I|=m}{\oplus }\WT{\Q}\ALT{*}{Y_I}\to 
\underset{|I|=m}{\oplus }\WT{\Q}\ALT{*}{Y_I}
$$
is given by the tensor product with $\caF$.

Let $K\coprod I=J$ be a division of subsets $\{1, \ldots , r\}$ 
with $|I|=m$ and $|J|=n$, and write 
$K=\{k_1, \ldots , k_{n-m}\}$ with $k_1<\cdots <k_{n-m}$.
Assume that $K$ is not empty.
Let us recall the map defined in \S 2.4:
$$
\varXi_K=\sum_{\sigma \in \frS_{n-m}}
(\sgn \sigma )(\iota_{k_{\sigma (1)}}, \ldots , 
\iota_{k_{\sigma (n-m)}})^*:\WT{\Q}C_*(Y_I)\to 
\WT{\Q}C_{*+n-m-1}(Y_J).
$$
If we see this map as a linear sum of exact functors from 
$\frP(Y_I)$ to the category of exact $(n-m-1)$-cubes of 
$\frP(Y_J)$, then Prop.2.12 says that 
\begin{equation}
\partial \varXi_K=\sum_{a=1}^{n-m-1}(-1)^{a+1}
\underset{|L|=a}{\sum_{L\coprod \PR{L}=K}}
\sgn \sbinom{L \ \PR{L}}{K}\varXi_L\varXi_{\PR{L}}. 
\label{sec9:eq1}
\end{equation}

Consider a division 
$K=K_1\coprod K_2\coprod \cdots \coprod K_l$ such that each 
$K_j$ is not empty.
Let $|K_p|=s_p$.
Then we have the following diagram of functors:
$$
\xymatrix{
\frP (Y_J) & \frP (Y_{K_2\cup \cdots \cup K_l\cup I})
\ar[l]_{\varXi_{K_1} \qquad } & \cdots 
\ar[l]_{\qquad \quad \varXi_{K_2}} & 
\frP (Y_I)\ar[l]_{\varXi_{K_l}}  \\
& \frP (S). \ar[lu]^{\pi^*} \ar[u]_{\pi^*} \ar[rru]_{\pi^*} \\
}
$$
Here we see $\varXi_{K_p}$ as a linear sum of exact functors from 
$\frP(Y_{K_{p-1}\cup \cdots \cup K_l\cup I})$ 
to the category of exact $(s_p-1)$-cubes of 
$Y_{K_p\cup \cdots \cup K_l\cup I}$.
With this diagram we can associate a map 
$$
\Sps{\caF; \varXi_{K_1}, \ldots , \varXi_{K_l}}^{\alt }:
\WT{\Q}\ALT{*}{Y_I}\to \WT{\Q}\ALT{*+n-m}{Y_J}.
$$
Then it follows from Prop.9.1 that for 
$x_I\in \WT{\Q}\ALT{*}{Y_I}$, 
\begin{align*}
\partial &\left(\Sps{\caF; \varXi_{K_1}, \ldots , 
\varXi_{K_l}}^{\alt }(x_I)\right)=\Sps{\caF; \varXi_{K_1}, 
\ldots , \varXi_{K_{l-1}}}^{\alt }(\varXi_{K_l}^{\alt }(x_I))  \\
&\hskip 4pc +\sum_{p=1}^{l-1}(-1)^p\Sps{\caF; \varXi_{K_1}, 
\cdots , \varXi_{K_{l-p}}\varXi_{K_{l-p+1}}, \ldots , 
\varXi_{K_l}}^{\alt }(x_I)  \\
&\hskip 4pc +(-1)^{s_1(l-1)+1}\varXi_{K_1}^{\alt }
\Sps{\caF; \varXi_{K_2}, \ldots , \varXi_{K_l}}^{\alt }(x_I)   \\
&\hskip 4pc +(-1)^l\sum_{p=1}^l(-1)^{s_1+\cdots +s_{p-1}-p+1}
\Sps{\caF; \varXi_{K_1}, \ldots , \partial \varXi_{K_p}, 
\ldots , \varXi_{K_l}}^{\alt }(x_I)   \\
&\hskip 4pc +(-1)^{n-m}\Sps{\caF; \varXi_{K_1}, \ldots , 
\varXi_{K_l}}^{\alt }(\partial x_I).
\end{align*}

For any division $K=K_1\coprod K_2\coprod \cdots \coprod K_l$, 
define the signature 
$$
\sgn \binom{K_1 \ \cdots \ K_l \ I}{J}
$$
as follows:
If we write 
\begin{align*}
K_p&=\{k_{p,1}, \ldots , k_{p,s_p}\}, \quad 
k_{p,1}< \ldots <k_{p,s_p},  \\
I&=\{i_1, \ldots , i_m\},  \qquad \ \ i_1< \ldots <i_m,  \\
J&=\{j_1, \ldots , j_n\},  \qquad \ \ j_1< \ldots <j_n, 
\end{align*}
then 
$$
\sgn \binom{K_1 \ \cdots \ K_l \ I}{J}=\sgn \scriptstyle{
\begin{pmatrix}k_{1,1} & \cdots \ k_{1,s_1} \ k_{2,1} \ \cdots \ 
k_{2,s_2} & \cdots \ k_{l,s_l} & i_1 \ \cdots & i_m \\
j_1 & \hdotsfor{3} & j_n\end{pmatrix}}.
$$
Define a map 
$$
\varphi^{m,n}=(\caF\otimes \ )^{m,n}:\underset{|I|=m}{\oplus }
\WT{\Q}\ALT{*}{Y_I}\to \underset{|J|=n}{\oplus }
\WT{\Q}\ALT{*+n-m}{Y_J}
$$
as follows: 
For $x=(x_I)\in \underset{|I|=m}{\oplus }\WT{\Q}\ALT{*}{Y_I}$, 
$$
\varphi^{m,n}(x)_J=
\underset{K_1\coprod \cdots \coprod K_l\coprod I=J}
{\sum }\sgn \sbinom{K_1 \ \cdots \ K_l \ I}{J}
(-1)^{b(s_1, \ldots , s_l)}\Sps{\caF; \varXi_{K_1}, 
\ldots , \varXi_{K_l}}^{\alt }(x_I), 
$$
where $s_j=|K_j|$ and 
$$
b(s_1, \ldots , s_l)=\begin{cases}
s_{l-1}+s_{l-3}+\cdots +s_2, &\ \text{$l$ is odd}, \\
s_{l-1}+s_{l-3}+\cdots +s_1, &\ \text{$l$ is even}.
\end{cases}
$$
The lemma below follows easily from the definition of 
$b(s_1, \ldots , s_l)$.

\vskip 1pc
\begin{lem}
For $l\geq 2$ and $1\leq p\leq l-1$, 
$$
b(s_1, \ldots , s_l)+b(s_1, \ldots , s_p+s_{p+1}, \ldots , s_l)+
\sum_{a=1}^ps_a
$$
is an even number.
Moreover, 
\begin{align*}
b(s_1, \ldots, s_l)+b(s_1, \ldots , s_{l-1})&=m-n-s_l, \\
b(s_1, \ldots, s_l)+b(s_2, \ldots , s_l)&=
\begin{cases}2b(s_2, \ldots , s_l), &\ \text{$l$ is odd},  \\
2b(s_2, \ldots , s_l)+s_1, &\ \text{$l$ is even}.  \end{cases}
\end{align*}
\end{lem}

\vskip 1pc
Using Lem.9.2 and \eqref{sec9:eq1} we can show that 
\begin{multline}
\underset{K_1\coprod \cdots \coprod K_l\coprod I=J}
{\sum }\sgn \sbinom{K_1 \ \cdots \ K_l \ I}{J}
(-1)^{b(s_1, \ldots , s_l)}\times \\ 
\left(\sum_{p=1}^{l-1}(-1)^p\Sps{\caF; \varXi_{K_1}, 
\cdots , \varXi_{K_{l-p}}\varXi_{K_{l-p+1}}, \ldots , 
\varXi_{K_l}}^{\alt }(x_I)\right. \\
+\left.(-1)^l\sum_{p=1}^l(-1)^{s_1+\cdots +s_{p-1}-p+1}
\Sps{\caF; \varXi_{K_1}, \ldots , \partial \varXi_{K_p}, 
\ldots , \varXi_{K_l}}^{\alt }(x_I)\right) \label{sec9:ele1} 
\end{multline}
is equal to zero.
Hence 
\begin{align*}
\partial &\varphi^{m,n}(x)_J+(-1)^{n-m-1}
\varphi^{m,n}\partial (x)_J \\
=&\underset{K_1\coprod \cdots \coprod K_l\coprod I=J}
{\sum }\sgn \sbinom{K_1 \ \cdots \ K_l \ I}{J}
(-1)^{b(s_1, \ldots , s_l)}\Sps{\caF; \varXi_{K_1}, 
\ldots , \varXi_{K_{l-1}}}^{\alt }(\varXi_{K_l}^{\alt }(x_I))  \\
&+\underset{K_1\coprod \cdots \coprod K_l\coprod I=J}
{\sum }\sgn \sbinom{K_1 \ \cdots \ K_l \ I}{J}
(-1)^{b(s_1, \ldots , s_l)+s_1(l-1)+1}
\varXi_{K_1}^{\alt }\Sps{\caF; \varXi_{K_2}, 
\ldots , \varXi_{K_l}}^{\alt }(x_I).   \\
\end{align*}
Applying Lem.9.2 to this equality we have 
\begin{align*}
\partial &\varphi^{m,n}(x)_J+(-1)^{n-m-1}
\varphi^{m,n}\partial (x)_J \\
=&\underset{K_1\coprod \cdots \coprod K_l\coprod I=J}
{\sum }\sgn \sbinom{K_1 \ \cdots \ K_l \ I}{J}
(-1)^{b(s_1, \ldots , s_{l-1})+m+n+s_l}
\Sps{\caF; \varXi_{K_1}, \ldots , \varXi_{K_{l-1}}}^{\alt }
(\varXi_{K_l}^{\alt }(x_I))  \\
&+\underset{K_1\coprod \cdots \coprod K_l\coprod I=J}
{\sum }\sgn \sbinom{K_1 \ \cdots \ K_l \ I}{J}
(-1)^{b(s_2, \ldots , s_l)+1}\varXi_{K_1}^{\alt }
\Sps{\caF; \varXi_{K_2}, \ldots , \varXi_{K_l}}^{\alt }(x_I)   \\
=&(-1)^n\sum_{s_l=1}^{n-m-1}\varphi^{m+s_l,n}(F^{m,m+s_l}(x))_J
+(-1)^{n+1}\sum_{s_1=1}^{n-m-1}F^{n-s_1,n}(\varphi^{m,n-s_1}(x))_J, 
\end{align*}
which leads to the following theorem:

\vskip 1pc
\begin{prop}
$$
(\caF\otimes \ )=\varphi =(\varphi^{m,n}):
\WT{\Q}\ALT{*}{X; Y_1, \ldots , Y_r}
\to \WT{\Q}\ALT{*}{X; Y_1, \ldots , Y_r}
$$
is a map of $\scC$-complexes.
\end{prop}

\subsection{A homotopy from $(\caF \otimes \ )f^*$ to 
$f^*(\caF \otimes \ )$}
Let $T$ be another scheme with closed subschemes 
$D_1, \ldots , D_r$, and 
$f:(X; Y_1, \ldots , Y_r)\to (T; D_1, \ldots , D_r)$ 
a morphism.
Assume that $X$ and $T$ are defined over a base 
scheme $S$, and $f$ is defined over $S$.
Let $\caF$ be a vector bundle on $S$.
We abbreviate $\pi_X^*\caF \otimes \ $ and 
$\pi_T^*\caF \otimes \ $ to $\caF \otimes \ $, where 
$\pi_X:X\to S$ and $\pi_T:T\to S$ are the structure morphisms.
Then we have the diagram 
$$
\begin{CD}
\WT{\Q}\ALT{*}{T; D_1, \ldots , D_r} @>{f^*}>> 
\WT{\Q}\ALT{*}{X; Y_1, \ldots , Y_r}  \\
@V{(\caF \otimes \ )}VV @VV{(\caF \otimes \ )}V  \\
\WT{\Q}\ALT{*}{T; D_1, \ldots , D_r} @>{f^*}>> 
\WT{\Q}\ALT{*}{X; Y_1, \ldots , Y_r}.
\end{CD}
$$
The aim of this subsection is to construct a homotopy 
$\varPhi_f$ from $(\caF \otimes \ )f^*$ to 
$f^*(\caF \otimes \ )$.

Let $K\coprod I=J$ be a division of subsets of $\{1, \ldots , r\}$ 
and $K=K_1\coprod \cdots \coprod K_l$ a division of $K$. 
Consider the following diagram:
$$
\xymatrix{
\frP (Y_J) & \cdots \ar[l]_{\quad \varXi_{K_1}} & 
\frP (Y_{K_p\cup \cdots \cup K_l\cup I}) 
\ar[l]_{\varXi_{K_{p-1}}\qquad \ } & 
\frP (D_{K_{p+1}\cup \cdots \cup K_l\cup I}) \ar[l]_{\varXi_{{K_p},f}} 
& \cdots \ar[l]_{\qquad \quad \varXi_{K_{p+1}}} & 
\frP (D_I)\ar[l]_{\varXi_{K_1}}  \\
& & & \frP (S), \ar[lllu]^{\pi_X^*} \ar[lu]_{\pi_X^*} 
\ar[u]_{\pi_T^*} \ar[rru]_{\pi_T^*} \\
}
$$
where $\varXi_{K_p,f}$ is the linear sum of exact functors defined 
in \S 2.4.
Unlike the previous case, $K_p$ may be empty in this case.
Define a map 
$$
\varPhi_f^{m,n}:\underset{|I|=m}{\oplus }
\WT{\Q}\ALT{*}{D_I}\to \underset{|J|=n}{\oplus }
\WT{\Q}\ALT{*+n-m+1}{Y_J}
$$
by 
\begin{align*}
\varPhi_f^{m,n}&(x)_J=
\underset{K_1\coprod \cdots \coprod K_l\coprod I=J}
{\sum }\sgn \sbinom{K_1 \ \cdots \ K_l \ I}{J}\times \\  
&\sum_{p=1}^l(-1)^{b(s_1, \ldots , s_{p-1}+s_p, \ldots , s_l)
+n+p+l+1}\Sps{\caF; \varXi_{K_1}, \ldots , 
\varXi_{{K_p},f}, \ldots , \varXi_{K_l}}^{\alt }(x_I)
\end{align*}
for $x=(x_I)\in \underset{|I|=m}{\oplus }\WT{\Q}\ALT{*}{D_I}$.
In the above, $|K_j|=s_j$ and $s_0$ is supposed to be zero.

Let us calculate $\partial \varPhi_f^{m,n}(x)$ using Prop.9.1 
and Lem.9.2.
Since a similar cancellation of terms to \eqref{sec9:ele1} 
occurs in this case,  
\begin{align*}
(-1)^n&\partial \varPhi_f^{m,n}(x)_J+(-1)^m
\varPhi_f^{m,n}(\partial x)_J   \\
=&\underset{K_1\coprod \cdots \coprod K_l\coprod I=J}
{\sum }\sgn \sbinom{K_1 \ \cdots \ K_l \ I}{J}\times  \\
&\left\{\sum_{p=1}^{l-1}
(-1)^{b(\ldots , s_{p-1}+s_p, \ldots )+p+l+1}
\Sps{\caF; \ldots , \varXi_{{K_p},f}, \ldots , 
\varXi_{K_{l-1}}}^{\alt }\varXi_{K_l}^{\alt }(x_I)\right.  \\
&\hskip 1pc +(-1)^{b(s_1, \ldots , s_{l-1}+s_l)+1}
\Sps{\caF; \varXi_{K_1}, \ldots , \varXi_{K_{l-1}}}^{\alt }
\varXi_{{K_l},f}^{\alt }(x_I)   \\
&\hskip 1pc +(-1)^{b(s_1, \ldots , s_l)+s_1(l-1)}
\varXi_{{K_1},f}^{\alt }\Sps{\caF; \varXi_{K_2}, \ldots , 
\varXi_{K_l}}^{\alt }(x_I)   \\
&\hskip 1pc \left.+\sum_{p=2}^l
(-1)^{b(\ldots , s_{p-1}+s_p, \ldots )+p+(s_1-1)l+s_1}
\varXi_{K_1}^{\alt }\Sps{\caF; \varXi_{K_2}, \ldots , 
\varXi_{{K_p},f}, \ldots }^{\alt }(x_I)\right\}  \\
=&-\sum_{s_l=1}^{n-m}\varPhi_f^{m+s_l,n}F^{m,m+s_l}(x)_J
-\sum_{s_l=0}^{n-m}(\caF \otimes \ )^{m+s_l,n}
(f^*)^{m,m+s_l}(x)_J   \\
&+\sum_{s_1=0}^{n-m}(f^*)^{n-s_1,n}
(\caF \otimes \ )^{m,n-s_1}(x)_J-\sum_{s_1=1}^{n-m}F^{n-s_1,n}
\varPhi_f^{m,n-s_1}(x)_J.
\end{align*}
Hence we have the following:

\vskip 1pc
\begin{prop}
$$
\varPhi_f:\WT{\Q}\ALT{*}{T; D_1, \ldots , D_r}\to 
\WT{\Q}\ALT{*+1}{X; Y_1, \ldots , Y_r}
$$
is a homotopy from 
$(\caF \otimes \ )f^*$ to $f^*(\caF \otimes \ )$.
\end{prop}

\vskip 1pc
The closed immersion $\iota_r:Y_r\hookrightarrow X$ induces 
a map of $\scC$-complexes 
$$
\iota_r^*:\WT{\Q}\ALT{*}{X; Y_1, \ldots , Y_{r-1}}\to 
\WT{\Q}\ALT{*}{Y_r; Y_1\cap Y_r, \ldots , Y_{r-1}\cap Y_r}
$$
and a homotopy $\varPhi_{\iota_r}$ from 
$(\caF \otimes \ )\iota_r^*$ to $\iota_r^*(\caF \otimes \ )$.
If we identify the simple complex of $\iota_r^*$ with 
the complex $\WT{\Q}\ALT{*}{X; Y_1, \ldots , Y_r}$ by Cor.2.16, 
then Prop.2.5 says that the homotopy $\varPhi_{\iota_r}$ 
gives a map of $\scC$-complexes 
$$
\varphi_s:\WT{\Q}\ALT{*}{X; Y_1, \ldots , Y_r}\to 
\WT{\Q}\ALT{*}{X; Y_1, \ldots , Y_r}.
$$

\vskip 1pc
\begin{prop}
The map $\varphi_s$ agrees with $(\caF \otimes \ )$.
\end{prop}

{\it Proof}: 
Let $x=(x_I)\in \underset{|I|=m}{\oplus }\WT{\Q}
\ALT{*}{Y_I}$.
It follows from the definition of $\varphi_s$ that 
$\varphi_s^{m,n}(x)$ is written as 
$$
\varphi_s^{m,n}(x)_J=
\underset{I\subset J}{\sum }\varphi_s^{I,J}(x_I).
$$
In the case that $r\notin J$ or $r\in I$, $\varphi_s^{I,J}(x_I)$ 
comes from $(\caF \otimes \ )$, that is, 
$$
\varphi_s^{I,J}(x_I)=
\underset{K_1\coprod \cdots \coprod K_l\coprod I=J}{\sum }
\sgn \sbinom{K_1 \ \cdots \ K_l \ I}{J}(-1)^{b(s_1, \ldots , s_l)}
\Sps{\caF; \varXi_{K_1}, \ldots , \varXi_{K_l}}^{\alt }(x_I), 
$$
where $|K_j|=s_j$.
Assume $r\in J$ and $r\notin I$.
In this case, $\varphi_s^{I,J}(x_I)$ comes from the homotopy 
$\varPhi_{\iota_r}$.
To be more precise, if we write $J^r=J-\{r\}$, then 
\begin{align*}
\varphi_s^{I,J}(x_I)=&
\underset{K_1\coprod \cdots \coprod K_l\coprod I=J^r}
{\sum }\sgn \sbinom{K_1 \ \cdots \ K_l \ I}{J^r}\times \\
&\hskip 1pc \sum_{p=1}^l
(-1)^{b(\ldots , s_{p-1}+s_p, \ldots )+(n-1)+p+l+1}
\Sps{\caF; \varXi_{K_1}, \ldots , \varXi_{K_p,\iota_r}, 
\ldots , \varXi_{K_l}}^{\alt }(x_I).
\end{align*}
Write $\PR{K}_p=K_p\cup \{r\}$ and 
$\PR{s}_p=|\PR{K}_p|=s_p+1$.
Then the following equalities hold:
\begin{align*}
\varXi_{K_p,\iota_r}=&\ (-1)^{s_p}\varXi_{\PR{K}_p},  \\
\sgn \sbinom{K_1 \ \cdots \ \PR{K}_p \ \cdots \ K_l \ I}{J}=&\ 
(-1)^{m+s_{p+1}+\cdots +s_l}
\sgn \sbinom{K_1 \ \cdots \ K_p \ \cdots \ K_l \ I}{J^r}.
\end{align*}
Hence using Lem.9.2 we can show that 
\begin{align*}
\varphi_s^{I,J}(x_I)
=&\underset{K_1\coprod \cdots \coprod K_l\coprod I=J^r}
{\sum }\sum_{p=1}^l\sgn \sbinom{K_1 \ \cdots \ \PR{K}_p \ 
\cdots \ K_l \ I}{J}\times \\
&\hskip 1pc (-1)^{b(s_1, \ldots , s_l)+p+l}
\Sps{\caF; \varXi_{K_1}, \ldots , \varXi_{\PR{K}_p}, 
\ldots , \varXi_{K_l}}^{\alt }(x_I).
\end{align*}
It follows from the definition of $b(s_1, \ldots , s_l)$ 
that 
$$
(-1)^{b(s_1, \ldots , s_l)+p+l}=
(-1)^{b(s_1, \ldots , \PR{s}_p, \ldots , s_l)}.
$$
Hence if we change the symbol $\PR{K}_p$ to $K_p$ and $\PR{s}_p$ to $s_p$, 
then 
$$
\varphi_s^{I,J}(x_I)=
\underset{K_1\coprod \cdots \coprod K_l\coprod I=J}{\sum }
\sgn \sbinom{K_1 \ \cdots \ K_l \ I}{J}
(-1)^{b(s_1, \ldots , s_l)}
\Sps{\caF; \varXi_{K_1}, \ldots , \varXi_{K_l}}^{\alt }(x_I), 
$$
which completes the proof.
\qed

\vskip 1pc
\subsection{A second homotopy arising from a section of 
a closed immersion}
Let $f:(X; Y_1, \ldots , Y_r)\to (T; D_1, \ldots , D_r)$ 
be a closed immersion defined over a base scheme $S$, and 
suppose that there is a morphism 
$g:(T; D_1, \ldots , D_r)\to (X; Y_1, \ldots , Y_r)$ 
also defined over $S$ such that $gf=\IId_X$.
Then we have the diagram of $\scC$-complexes 
$$
\begin{CD}
\WT{\Q}\ALT{*}{X; Y_1, \ldots , Y_r} @>{g^*}>>
\WT{\Q}\ALT{*}{T; D_1, \ldots , D_r} @>{f^*}>>
\WT{\Q}\ALT{*}{X; Y_1, \ldots , Y_r} \\
@VV{(\caF \otimes \ )}V @VV{(\caF \otimes \ )}V 
@VV{(\caF \otimes \ )}V \\ 
\WT{\Q}\ALT{*}{X; Y_1, \ldots , Y_r} @>{g^*}>>
\WT{\Q}\ALT{*}{T; D_1, \ldots , D_r} @>{f^*}>>
\WT{\Q}\ALT{*}{X; Y_1, \ldots , Y_r}
\end{CD}
$$
and the homotopies 
\begin{align*}
\varPhi_g:&\ \WT{\Q}\ALT{*}{X; Y_1, \ldots , Y_r}\to 
\WT{\Q}\ALT{*}{T; D_1, \ldots , D_r}, \\
\varPhi_f:&\ \WT{\Q}\ALT{*}{T; D_1, \ldots , D_r}\to 
\WT{\Q}\ALT{*}{X; Y_1, \ldots , Y_r}, \\
\varPsi:&\ \WT{\Q}\ALT{*}{X; Y_1, \ldots , Y_r}\to 
\WT{\Q}\ALT{*}{X; Y_1, \ldots , Y_r}
\end{align*}
from $(\caF \otimes \ )g^*$ to $g^*(\caF \otimes \ )$, 
from $(\caF \otimes \ )f^*$ to $f^*(\caF \otimes \ )$, 
and from the identity to $f^*g^*$ respectively.
Hence we have two homotopies 
$$
\varPhi_fg^*+f^*\varPhi_g+(\caF \otimes \ )\varPsi , \ 
\varPsi (\caF \otimes \ ): 
\WT{\Q}\ALT{*}{X; Y_1, \ldots , Y_r}\to 
\WT{\Q}\ALT{*+1}{X; Y_1, \ldots , Y_r}
$$
from $(\caF \otimes \ )$ to $f^*g^*(\caF \otimes \ )$.
In this subsection we will construct a second homotopy 
between them which admits the condition of Def.2.7.

Define a map 
$$
\varTheta_1^{m,n}:\underset{|I|=m}{\oplus }
\WT{\Q}\ALT{*}{Y_I}\to \underset{|J|=n}{\oplus }
\WT{\Q}\ALT{*+n-m+2}{Y_J}
$$
by 
\begin{align*}
\varTheta_1^{m,n}(x)_J=& 
\underset{K_1\coprod \cdots \coprod K_l\coprod I=J}
{\sum }\sgn \sbinom{K_1 \ \cdots \ K_l \ I}{J}\times \\
&\hskip 1pc \sum_{p=1}^l(-1)^{b(s_1, \ldots , s_l)+1}
\Sps{\caF; \varXi_{K_1}, \ldots , \varXi_{{K_p},f,g}, 
\ldots , \varXi_{K_l}}^{\alt }(x_I)
\end{align*}
for $x=(x_I)\in \underset{|I|=m}{\oplus }\WT{\Q}\ALT{*}{Y_I}$, 
where $|K_j|=s_j$.
In the above, $K_p$ may be empty. 
Then a similar cancellation of terms to \eqref{sec9:ele1} 
occurs, therefore 
\begin{align*}
(-1)^n&\partial \varTheta_1^{m,n}(x)_J-(-1)^m
\varTheta_1^{m,n}(\partial x)_J   \\
=&\sum_{s_l=1}^{n-m}\varTheta_1^{m+s_l,n}F^{m,m+s_l}(x)_J
-\sum_{s_l=0}^{n-m}(\caF \otimes \ )^{m+s_l,n}
\varPsi^{m,m+s_l}(x)_J   \\
&+\sum_{s_1=0}^{n-m}\varPsi^{n-s_1,n}
(\caF \otimes \ )^{m,n-s_1}(x)_J-\sum_{s_1=1}^{n-m}F^{n-s_1,n}
\varTheta_1^{m,n-s_1}(x)_J  \\
&+\underset{K_1\coprod \cdots \coprod K_l\coprod I=J}
{\sum }\sgn \sbinom{K_1 \ \cdots \ K_l \ I}{J}\times \\
&\hskip 2pc \sum_{p=1}^{l-1}
(-1)^{b(s_1, \ldots , s_l)+n+p+l+s_p+1}
\Sps{\caF; \ldots , \varXi_{{K_p},f}\varXi_{K_{p+1},g}, 
\ldots }^{\alt }(x_I).
\end{align*}
Meanwhile, define a map 
$$
\varTheta_2^{m,n}:\underset{|I|=m}{\oplus }
\WT{\Q}\ALT{*}{Y_I}\to \underset{|J|=n}{\oplus }
\WT{\Q}\ALT{*+n-m+2}{Y_J}
$$
by 
\begin{align*}
\varTheta_2^{m,n}(x)_J&=
\underset{K_1\coprod \cdots \coprod K_l\coprod I=J}
{\sum }\sgn \sbinom{K_1 \ \cdots \ K_l \ I}{J}\times \\
&\sum_{1\leq p<q\leq l}(-1)^{c_{p,q}}
\Sps{\caF; \varXi_{K_1}, \ldots , \varXi_{{K_p},f}, 
\ldots , \varXi_{{K_q},g}, \ldots , \varXi_{K_l}}^{\alt }(x_I)
\end{align*}
for $x=(x_I)\in \underset{|I|=m}{\oplus }\WT{\Q}\ALT{*}{Y_I}$, 
where 
$$
c_{p,q}=b(s_1, \ldots , s_l)+\sum_{j=p}^{q-1}s_j+p+q+1
$$
and $|K_j|=s_j$.
Here $K_p$ and $K_q$ may be empty.
Then 
\begin{align*}
(-1)^n&\partial \varTheta_2^{m,n}(x)_J-(-1)^m
\varTheta_2^{m,n}(\partial x)_J   \\
=&\sum_{s_l=1}^{n-m}\varTheta_2^{m+s_l,n}F^{m,m+s_l}(x)_J-
\sum_{s_l=0}^{n-m}\varPhi_f^{m+s_l,n}(g^*)^{m,m+s_l}(x)_J   \\
&-\sum_{s_1=0}^{n-m}(f^*)^{n-s_1,n}\varPhi_g^{m,n-s_1}(x)_J-
\sum_{s_1=1}^{n-m}F^{n-s_1,n}\varTheta_2^{m,n-s_1}(x)_J  \\
&+\underset{K_1\coprod \cdots \coprod K_l\coprod I=J}
{\sum }\sgn \sbinom{K_1 \ \cdots \ K_l \ I}{J}\times \\
&\hskip 2pc \sum_{p=1}^{l-1}
(-1)^{b(s_1, \ldots , s_l)+n+p+l+s_p}
\Sps{\caF; \ldots , \varXi_{{K_p},f}\varXi_{K_{p+1},g}, 
\ldots }^{\alt }(x_I).
\end{align*}
Hence we have the following:

\vskip 1pc
\begin{prop}
If we set 
$$
\varTheta^{m,n}=\varTheta_1^{m,n}+\varTheta_2^{m,n}:
\underset{|I|=m}{\oplus }\WT{\Q}\ALT{*}{Y_I}\to 
\underset{|J|=n}{\oplus }\WT{\Q}\ALT{*+n-m+2}{Y_J}, 
$$
then for $x=(x_I)\in \underset{|I|=m}{\oplus }\WT{\Q}\ALT{*}{Y_I}$, 
\begin{align*}
(-1)^n&\partial \varTheta^{m,n}(x)_J+\sum_{k_1=1}^{n-m}F^{n-k_1,n}
\varTheta^{m,n-k_1}(x)_J  \\
&-(-1)^m\varTheta^{m,n}(\partial x)_J-
\sum_{k_l=1}^{n-m}\varTheta^{m+k_l,n}F^{m,m+k_l}(x)_J \\
=&\sum_{k_1=0}^{n-m}\varPsi^{n-k_1,n}
(\caF \otimes \ )^{m,n-k_1}(x)_J-\sum_{k_l=0}^{n-m}
(\caF \otimes \ )^{m+k_l,n}\varPsi^{m,m+k_l}(x)_J  \\
&-\sum_{k_l=0}^{n-m}\varPhi_f^{m+k_l,n}(g^*)^{m,m+k_l}(x)_J
-\sum_{k_1=0}^{n-m}(f^*)^{n-k_1,n}\varPhi_g^{m,n-k_1}(x)_J.
\end{align*}
In other words, $\Theta =(\Theta^{m,n})$ is a second homotopy from 
$\varPhi_fg^*+f^*\varPhi_g+(\caF \otimes \ )\varPsi $ to 
$\varPsi (\caF \otimes \ )$.
\end{prop}

\setcounter{equation}{0}
\vskip 2pc
\section{$\WH{K}_0(X)$-module structures on 
arithmetic $K$-groups}

\vskip 1pc
\subsection{$\WH{K}_0(X)$-module structures on 
$\WH{K}_r(X)$ and on $\WH{K}_0(T)$}
Let $X$ be a smooth proper variety defined over an arithmetic ring.
Let us first recall the $\WH{K}_0(X)$-module structure on 
$\WH{K}_r(X)$ given in \cite[\S 5]{takeda}.
Define an operation on the Deligne complexes 
$$
\TRI :\mmD^{2p-n}(X, p)\times \mmD^{2q-m}(X, q)
\to \mmD^{2p+2q-n-m-1}(X, p+q)
$$
as follows: 
Let 
$$
a^{n,m}_{i,j}=1-2\sbinom{n+m}{n}^{-1}\sum_{\alpha =0}^{i-1}
\sbinom{n+m-i-j+1}{n-\alpha }\sbinom{i+j-1}{\alpha }, 
$$
and for $\omega \in \mmD^{2p-n}(X, p)$ and 
$\tau \in \mmD^{2q-m}(X, q)$,  
$$
\omega \TRI \tau =
\underset{1\leq j\leq m}{\underset{1\leq i\leq n}{\sum }}
a^{n,m}_{i,j}\omega^{(p-n+i-1,p-i)}\wedge 
\tau^{(q-m+j-1,q-j)}
$$
if $n, m\geq 1$, and $\omega \TRI \tau =0$ if $n=0$ or $m=0$.

\vskip 1pc
\begin{prop}\cite[Thm.5.2]{takeda}
Let $\OV{\caF}$ and $\OV{\caG}$ be an exact hermitian $n$-cube 
and $m$-cube on $X$.
Then 
\begin{align*}
\CH_{n+m}(\OV{\caF}\otimes \OV{\caG})=&\CH_n(\OV{\caF})\bullet 
\CH_m(\OV{\caG})+(-1)^{n+1}
d_{\scD}(\CH_n(\OV{\caF})\TRI \CH_m(\OV{\caG})) \\
&+(-1)^nd_{\scD}\CH_n(\OV{\caF})\TRI \CH_m(\OV{\caG})-
\CH_n(\OV{\caF})\TRI d_{\scD}\CH_m(\OV{\caG}).
\end{align*}
\end{prop}

\vskip 1pc
Suppose $r\geq 1$.
Let $(\OV{\caF}, \WT{\eta })$ be a pair of a hermitian vector bundle 
$\OV{\caF}$ on $X$ with $\WT{\eta }\in \mmD_1(X)/\IIm d_{\scD}$ and 
$$
(y, \tau )\in \WT{\Q}\WH{C}_r(X)\oplus \mmD_{r+1}(X)=s(\CH_*)_r
$$
such that $\partial (y, \tau )=0$.
Note that 
$\CH_0(\OV{\caF})\bullet \tau =\CH_0(\OV{\caF})\wedge \tau $ 
by definition, and in what follows we write 
$\CH_0(\OV{\caF})\bullet \tau $ instead of $\CH_0(\OV{\caF})\wedge \tau $.

\vskip 1pc
\begin{prop}
Define a product of above pairs by 
$$
(\OV{\caF}, \WT{\eta })\cdot (y, \tau )=
(\OV{\caF}\otimes y, \CH_0(\OV{\caF})\bullet \tau )\in s(\CH_*)_r.
$$
Then this product gives rise to a map of arithmetic $K$-groups 
\begin{equation}
\WH{K}_0(X)\times \QQ{\WH{K}_r(X)}\to \QQ{\WH{K}_r(X)}, 
\label{sec10:map1}
\end{equation}
by which $\QQ{\WH{K}_r(X)}$ is a $\WH{K}_0(X)$-module.
\end{prop}

{\it Proof}: 
It is obvious that 
$\partial (\OV{\caF}\otimes y, \CH_0(\OV{\caF})\bullet \tau )=0$ 
and that $\partial (\OV{\caF}\otimes \PR{y}, \CH_0(\OV{\caF})
\bullet \PR{\tau })=(\OV{\caF}\otimes y, \CH_0(\OV{\caF})\bullet \tau )$ 
if $(y, \tau )=\partial (\PR{y}, \PR{\tau })$.
Let $\OV{\caE}:0\to \OV{\caF}_{-1}\to \OV{\caF}_0\to \OV{\caF}_1\to 0$
be a short exact sequence of hermitian vector bundles on $X$.
We see $\OV{\caE}$ as an element of $\WT{\Q}\WH{C}_1(X)$.
Then Prop.10.1 implies that 
$$
\CH_{r+1}(\OV{\caE}\otimes y)=\CH_1(\OV{\caE})\bullet \CH_r(y)
+d_{\scD}(\CH_1(\OV{\caE})\TRI \CH_r(y))-\CH_1(\OV{\caE})\TRI 
d_{\scD}\CH_r(y).
$$
Since $\CH_r(y)=d_{\scD}\tau $, we have 
\begin{align*}
\partial (\OV{\caE}\otimes y, &\, -\CH_1(\OV{\caE})\bullet \tau 
+\CH_1(\OV{\caE})\TRI \CH_n(y))  \\
&=(\partial \OV{\caE}\otimes y, \CH_{r+1}(\OV{\caE}\otimes y)
+d_{\scD}(\CH_1(\OV{\caE})\bullet \tau 
-\CH_1(\OV{\caE})\TRI \CH_n(y)))  \\
&=(\OV{\caF}_{-1}\otimes y+\OV{\caF}_1\otimes y-
\OV{\caF}_0\otimes y, \CH_0(\OV{\caF}_{-1})\bullet \tau 
+\CH_0(\OV{\caF}_1)\bullet \tau -\CH_0(\OV{\caF}_0)\bullet \tau ), 
\end{align*}
which means that the map \eqref{sec10:map1} is well-defined.

We next show that the above product gives a $\WH{K}_0(X)$-module 
structure on $\QQ{\WH{K}_r(X)}$.
Take two elements 
$[(\OV{\caF}_1, \WT{\eta }_1)], [(\OV{\caF}_2, \WT{\eta }_2)]
\in \WH{K}_0(X)$.
Then their product in $\WH{K}_0(X)$ is given by 
$$
[(\OV{\caF}_1\otimes \OV{\caF}_2, \CH_0(\OV{\caF}_1)\wedge \WT{\eta }_2
+\WT{\eta }_1\wedge \CH_0(\OV{\caF}_2)+d_{\scD}\eta_1\wedge \WT{\eta }_2)]
\in \WH{K}_0(X).
$$
For any exact hermitian $r$-cube $\OV{\caG}$ on $X$, consider 
the following exact hermitian $(r+1)$-cube 
\begin{equation}
\Sps{\OV{\caF}_1, \OV{\caF}_2}(\OV{\caG})=\left(
\OV{\caF}_1\otimes (\OV{\caF}_2\otimes \OV{\caG})
\overset{\sim }{\to }(\OV{\caF}_1\otimes \OV{\caF}_2)
\otimes \OV{\caG}\to 0\right). \label{sec10:eq1}
\end{equation}
Then it gives a map 
$$
\Sps{\OV{\caF}_1, \OV{\caF}_2}:\WT{\Q}\WH{C}_r(X)\to 
\WT{\Q}\WH{C}_{r+1}(X)
$$
which satisfies 
$$
\partial \left(\Sps{\OV{\caF}_1, \OV{\caF}_2}(x)\right)=
\OV{\caF}_1\otimes (\OV{\caF}_2\otimes x)
-(\OV{\caF}_1\otimes \OV{\caF}_2)\otimes x
-\Sps{\OV{\caF}_1, \OV{\caF}_2}(\partial x)
$$
for $x\in \WT{\Q}\WH{C}_r(X)$.
Note that $\Sps{\OV{\caF}_1, \OV{\caF}_2}(x)$ is 
isometrically equivalent to a degenerate element since 
the isomorphism in \eqref{sec10:eq1} is an isometry.
Hence we have 
$\CH_{r+1}\left(\Sps{\OV{\caF}_1, \OV{\caF}_2}(x)\right)=0$.

Let $(y, \tau )\in s(\CH_*)_r$ such that 
$\partial (y, \tau )=0$.
Then 
\begin{align*}
&(\OV{\caF}_1, \WT{\eta }_1)\cdot ((\OV{\caF}_2, \WT{\eta }_2)
\cdot (y, \tau ))-\left(\OV{\caF}_1\otimes \OV{\caF}_2, 
\CH_0(\OV{\caF}_1)\wedge \WT{\eta }_2+\WT{\eta }_1\wedge 
\CH_0(\OV{\caF}_2)+d_{\scD}\eta_1\wedge \WT{\eta }_2\right)
\cdot (y, \tau ) \\
&\hskip 2pc =\left(\OV{\caF}_1\otimes (\OV{\caF}_2\otimes y)
-(\OV{\caF}_1\otimes \OV{\caF}_2)\otimes y, 0\right).
\end{align*}
Since $\partial y=0$ and 
$\CH_{r+1}(\Sps{\OV{\caF}_1, \OV{\caF}_2}(y))=0$, 
it is equal to 
$\partial \left(\Sps{\OV{\caF}_1, \OV{\caF}_2}(y), 0\right)$.
This means that the product \eqref{sec10:map1} is associative.
The identity condition 
$[(\OV{\caO}, 0)]\cdot [(y, \tau )]=[(y, \tau )]$, 
where $\OV{\caO}$ is the structure sheaf with the canonical 
metric, can be shown in a similar way by using the canonical 
isometry $\OV{\caO}\otimes \OV{\caG}\simeq \OV{\caG}$.
\qed

\vskip 1pc
We next consider a $\WH{K}_0(X)$-module structure on 
$\QQ{\WH{K}_{\P,r}(X)}$.
Let $(\OV{\caF}, \WT{\eta })$ be a pair of 
a hermitian vector bundle $\OV{\caF}$ on $X$ with 
$\WT{\eta }\in \mmD_1(X)/\IIm d_{\scD}$.
Suppose $r\geq 1$ and let 
$$
(y, \tau )\in \WT{\Q}\WH{C}_r(X)\oplus \WT{\scD}_{\P,r+1}(X)
=s(\CH_{*,\P})_r, 
$$
then define their product by a similar expression:
\begin{equation}
(\OV{\caF}, \WT{\eta })\cdot (y, \tau )=(\OV{\caF}\otimes y, 
\CH_0(\OV{\caF})\PRP \tau ), \label{sec10:eq2}
\end{equation}
where $\PRP $ is the product introduced in \S 3.2.
It is easy to show that the canonical isomorphism 
$\QQ{\WH{K}_{\P,r}(X)}\to \QQ{\WH{K}_r(X)}$, which is 
\eqref{sec5:iso1} given in \S 5.1, is compatible with 
the product with the pair $(\OV{\caF}, \WT{\eta })$, 
hence the product \eqref{sec10:eq2} gives a map 
$$
\WH{K}_0(X)\times \QQ{\WH{K}_{\P,r}(X)}\to \QQ{\WH{K}_{\P,r}(X)}, 
$$
by which $\QQ{\WH{K}_{\P,r}(X)}$ is a $\WH{K}_0(X)$-module.

Here we assume that $X$ is projective over an arithmetic 
field.
Let $T=D(X\times \PP{r}; X\times \partial \PP{r})$ be 
the associated iterated double.
Let $(\OV{\caF}, \WT{\eta })$ be a pair of a hermitian vector 
bundle $\OV{\caF}$ on $X$ with $\WT{\eta }\in \mmD_1(X)/\IIm d_{\scD}$, 
and $(\OV{\caG}, \WT{\tau })$ a pair of a hermitian vector bundle 
$\OV{\caG}$ on $T$ with 
$\WT{\tau }\in \WT{\scD}_{\A,\P,r+1}(X)/\IIm d_s$.
Define a product of such pairs by 
\begin{equation}
(\OV{\caF}, \WT{\eta })\cdot (\OV{\caG}, \WT{\tau })=
(\OV{\caF}\otimes \OV{\caG}, \CH_0(\OV{\caF})\PRAP \WT{\tau }
+\WT{\eta }\PRAP \CH_{T,0}(\OV{\caG})+d_{\scD}\eta \PRAP \WT{\tau }), 
\label{sec10:pro1}
\end{equation}
where $\PRAP $ is the product introduced in the end of \S 3.2.

\vskip 1pc
\begin{prop}
The above product induces a pairing of arithmetic $K$-groups 
$$
\WH{K}_0(X)\times \WH{K}_0^M(T)\to \WH{K}_0^M(T), 
$$
by which $\WH{K}_0^M(T)$ is a $\WH{K}_0(X)$-module.
\end{prop}

{\it Proof}: 
It is easy to see that \eqref{sec10:pro1} is compatible with 
the relation coming from a short exact sequence of hermitain 
vector bundles on $T$.
Let $\OV{\caE}:0\to \OV{\caF}_{-1}\to \OV{\caF}_0\to \OV{\caF}_1\to 0$ 
be a short exact sequence of hermitian vector bundles on $X$.
Then
$$ 
(\OV{\caF}_1+\OV{\caF}_{-1}-\OV{\caF}_0, 
-\WT{\CH}_1(\OV{\caE}))\cdot (\OV{\caG}, \WT{\tau })=
(\OV{\caF}_1\otimes \OV{\caG}+\OV{\caF}_{-1}\otimes \OV{\caG}
-\OV{\caF}_0\otimes \OV{\caG}, 
-\WT{\CH}_1(\OV{\caE})\PRAP \CH_{T,0}(\OV{\caG})).
$$
Note that $\CH_1(\OV{\caE})\PRAP \CH_{T,0}(\OV{\caG})$ is in 
$\WT{\scD}_{\A,r+1}(X)$.
On the other hand, 
\begin{align*}
\CH_{T,1}(\OV{\caE}\otimes \OV{\caG})=&\, 
\sum_{I}(-1)^{|I|}\CH_1(i_I^*(\OV{\caE}\otimes \OV{\caG}))_{\P} \\
=&\, \sum_{I}(-1)^{|I|}\pi_{X,\P}^*\CH_1(\OV{\caE})_{\P}\wedge 
\pi_{X,\A}^*\CH_0(i_I^*\OV{\caG}) \\
=&\, \pi_{X,\P}^*\CH_1(\OV{\caE})_{\P}\wedge 
\pi_{X,\A}^*\CH_{T,0}(\OV{\caG}), 
\end{align*}
where $\pi_{X,\P}:X\times \PP{r}\times \P^1\to X\times \P^1$ and 
$\pi_{X,\A}:X\times \PP{r}\times \P^1\to X\times \PP{r}$ are 
the projections.
Moreover, the difference $\CH_1(\OV{\caE})_{\P}-\CH_1(\OV{\caE})$ is 
$d_s$-exact in $\WT{\scD}_{\P,2}(X)$ by Lem.6.13, and it follows 
that $\CH_1(\OV{\caE})\PRAP \CH_{T,0}(\OV{\caG})=
\pi_X^*\CH_1(\OV{\caE})\wedge \CH_{T,0}(\OV{\caG})$ in 
$\WT{\scD}_{\A,r+1}(X)$.
Then taking the $d_s$-closedness of $\CH_{T,0}(\OV{\caG})=0$ 
into consideration, we can show in the same way as the proof 
of Lem.6.13 that 
\begin{align*}
\CH_{T,1}(\OV{\caE}\otimes \OV{\caG})&\, -
(-1)^r\CH_1(\OV{\caE})\PRAP \CH_{T,0}(\OV{\caG}) \\
=&\, \pi_{X,\P}^*\CH_1(\OV{\caE})_{\P}\wedge 
\pi_{X,\A}^*\CH_{T,0}(\OV{\caG})
-(-1)^r\pi_X^*\CH_1(\OV{\caE})\wedge \CH_{T,0}(\OV{\caG})
\end{align*}
is $d_s$-exact in $\WT{\scD}_{\A,\P,r+1}(X)$.
This means that 
\begin{align*}
(\OV{\caF}_1+&\, \OV{\caF}_{-1}-\OV{\caF}_0, 
-\WT{\CH}_1(\OV{\caE}))\cdot (\OV{\caG}, \WT{\tau }) \\
&=(\OV{\caF}_1\otimes \OV{\caG}+\OV{\caF}_{-1}\otimes \OV{\caG}
-\OV{\caF}_0\otimes \OV{\caG}, 
-(-1)^r\WT{\CH}_{T,1}(\OV{\caE}\otimes \OV{\caG})), 
\end{align*}
therefore the product \eqref{sec10:pro1} is compatible with 
the relation coming from a short exact sequence of hermitian 
vector bundles on $X$.

Finally we show that the product admits the associative law.
We should note that the product $\bullet $ on $\mmD_{\log}^*(X, *)$ 
does not satisfies the associative law in general, whereas if at least 
one of three elements $\omega , \eta , \tau $ is in 
$\mmD_{\log}^{2p}(X\times \PP{r}\times (\P^1)^s, p)$, then 
$\omega \bullet (\eta \bullet \tau )
=(\omega \bullet \eta )\bullet \tau $ holds.
Using this fact we can show the associativity of 
the product in the same way as the proof of the associativity 
of the product in $\WH{K}_0(X)$ in \cite[Thm.7.3.2]{gilletsoule2}.
\qed

\vskip 1pc
It is easy to see that the diagram 
$$
\begin{CD}
\WH{K}_0(X)\times \WH{K}_0^M(T) @>>> \WH{K}_0^M(T)  \\
@V{\CH_0\times \CH_{T,0}}VV  @VV{\CH_{T,0}}V  \\
\mmD_0(X)\times \WT{\scD}_{\A,\P,r}(X)  @>{\PRAP }>> 
\WT{\scD}_{\A,\P,r}(X)
\end{CD}
$$
is commutative.
This implies that $\WH{K}_0(T)\subset \WH{K}_0^M(T)$ is 
a $\WH{K}_0(X)$-submodule.
Moreover, since 
$$
\begin{CD}
\WH{K}_0(X)\times \WH{K}_0^M(T) @>>> \WH{K}_0^M(T)  \\
@V{\zeta \times \zeta }VV  @VV{\zeta }V  \\
K_0(X)\times K_0(T) @>>> K_0(T)
\end{CD}
$$
is commutative, $\WH{K}_0^M(T; T_1, \ldots , T_r)$ 
and $\WH{K}_0(T; T_1, \ldots , T_r)$ 
are also $\WH{K}_0(X)$-submodules of $\WH{K}_0^M(T)$.

\vskip 1pc
\begin{prop}
The splitting maps 
\begin{align*}
\WH{t}:&\ \WH{K}^M_0(T)\to \WH{K}^M_0(T; T_1, \ldots , T_r), \\
\WH{t}:&\ \WH{K}_0(T)\to \WH{K}_0(T; T_1, \ldots , T_r)
\end{align*}
respect $\WH{K}_0(X)$-module structures.
\end{prop}

{\it Proof}: 
Since the second map is the restriction of the first map, 
it suffices to proof the claim for the first one.
We will prove that the composite 
$$
\WH{K}^M_0(T)\overset{\WH{t}}{\to }
\WH{K}^M_0(T; T_1, \ldots , T_r)\overset{\WH{q}}{\to }
\WH{K}^M_0(T)
$$
is a map of $\WH{K}_0(X)$-modules.
Let $\OV{\caF}$ be a hermitian vector bundle on $X$, and 
$\OV{\caG}$ a hermitian vector bundle on $T$.
Since $qt(\OV{\caG})$ is equal to 
$(1-p_r^*\iota_r^*)\cdots (1-p_1^*\iota_1^*)\OV{\caG}$, 
there is a canonical isometry 
$qt(\OV{\caF}\otimes \OV{\caG})\simeq \OV{\caF}\otimes qt(\OV{\caG})$ 
as virtual hermitian vector bundles on $T$.
Hence the proposition follows.
\qed

\vskip 1pc
\subsection{A $\WH{K}_0(X)$-module structure on 
$\QQ{\WH{K}_0(X\times \PP{r}; X\times \partial \PP{r})}$}
Let $X$ be a smooth proper variety defined over an arithmetic 
field, and $\OV{\caF}$ a hermitain vector bundle on $X$.
Applying the construction in \S 9.3 we can obtain 
a tensor product functor 
$$
(\OV{\caF}\otimes \ ):
\WT{\Q}\WALT{*}{X\times \PP{r}; X\times \partial \PP{r}}\to 
\WT{\Q}\WALT{*}{X\times \PP{r}; X\times \partial \PP{r}}.
$$

\vskip 1pc
\begin{prop}
The diagram 
$$
\begin{CD}
\WT{\Q}\WALT{*}{X\times \PP{r}; X\times \partial \PP{r}}[r]
@>{(\OV{\caF}\otimes \ )}>>
\WT{\Q}\WALT{*}{X\times \PP{r}; X\times \partial \PP{r}}[r]  \\
@V{\CH_*}VV  @VV{\CH_*}V  \\
\WT{\scD}_{\A,\P,*}(X) @>{\CH_0(\OV{\caF})\PRAP }>> 
\WT{\scD}_{\A,\P,*}(X)
\end{CD}
$$
is commutative.
\end{prop}

{\it Proof}: 
For any $x_I\in \WT{\Q}\WALT{n}{X\times D_I}$ with 
$I\subset \{1, \ldots , r\}$, it holds that 
$\CH_n(\OV{\caF}\otimes x_I)_{\P}=
\CH_0(\OV{\caF})\PRAP \CH_n(x_I)_{\P}$ in 
$\WT{\scD}_{\A,\P,r+n}(X)$.
Moreover, for $I\subset J$ with $I\not= J$ and for a division 
$K_1\coprod \cdots \coprod K_l\coprod I=J$, 
$\CH_*\left(
\Sps{\OV{\caF}; \varXi_{K_1}, \ldots , \varXi_{K_l}}(x_I)\right)=0$ 
because $\Sps{\OV{\caF}; \varXi_{K_1}, \ldots , \varXi_{K_l}}(x_I)$ 
is isometrically equivalent to a degenerate element.
Hence the proposition follows.
\qed

\vskip 1pc
\begin{cor}
Let $(\OV{\caF}, \WT{\eta })$ be a pair of a hermitian vector bundle 
$\OV{\caF}$ on $X$ with $\WT{\eta }\in \mmD_1(X)/\IIm d_{\scD}$. 
Then for 
$$
(y, \tau )\in 
\WT{\Q}\WALT{0}{X\times \PP{r}; X\times \partial \PP{r}}
\oplus \WT{\scD}_{\A,\P,r+1}(X)=s(\CH_*)_r, 
$$
the product 
\begin{equation}
(\OV{\caF}, \WT{\eta })\cdot (y, \tau )=(\OV{\caF}\otimes y, 
\CH_0(\OV{\caF})\PRAP \tau ). \label{sec10:eq3}
\end{equation}
gives a map of complexes 
$$
(\OV{\caF}, \WT{\eta })\, \cdot \ :\QQ{\WH{K}_0(X\times \PP{r}; X\times 
\partial \PP{r})}\to \QQ{\WH{K}_0(X\times \PP{r}; X\times 
\partial \PP{r})}.
$$
\end{cor}

\vskip 1pc
\begin{prop}
The isomorphism 
$$
\QQ{\WH{K}_{\P,r}(X)}\simeq 
\QQ{\WH{K}_0(X\times \PP{r}; X\times \partial \PP{r})}
$$
given in Prop.5.1 is compatible with the product with the pair 
$(\OV{\caF}, \WT{\eta })$.
Hence the product \eqref{sec10:eq3} gives a 
$\WH{K}_0(X)$-module structure on 
$\QQ{\WH{K}_0(X\times \PP{r}; X\times \partial \PP{r})}$, 
and the above isomorphism respects the 
$\WH{K}_0(X)$-module structures.
\end{prop}

{\it Proof}: 
Let 
$$
i_X:\WT{\Q}\WALT{*}{X}\to 
\WT{\Q}\WALT{*}{X\times \PP{r}; X\times \partial \PP{r}}[r]
$$
be the restriction to the alternating part of the map 
\eqref{sec3:map1} defined in \S 3.6.
Then we can easily verify that $i_X$ commutes with 
the tensor product functor, from which the proposition follows.
\qed

\vskip 1pc
\subsection{Comparison of the $\WH{K}_0(X)$-module structures 
on $\WH{K}_r(X)$ and on $\WH{K}_0(T)$}
In this subsection we will prove the following theorem:

\vskip 1pc
\begin{prop}
The map 
$$
\WH{i}_{\emptyset }^*:\QQ{\WH{K}_0(T; T_1, \ldots , T_r)}\to 
\QQ{\WH{K}_0(X\times \PP{r}; X\times \partial \PP{r})}
$$
respects the $\WH{K}_0(X)$-module structure.
\end{prop}

{\it Proof}: 
It suffices to show that the map \eqref{sec5:map2}
$$
\WH{i_{\emptyset}^*t}:\WH{K}_0(T)\to 
\QQ{\WH{K}_0(X\times \PP{r}; X\times \partial \PP{r})}
$$
respects the $\WH{K}_0(X)$-module structures.
Recall the morphisms 
\begin{align*}
\iota_j:&\, (T_j; T_1\cap T_j, \ldots , T_{j-1}\cap T_j)\to 
(T; T_1, \ldots , T_{j-1}),  \\
p_j:&\, (T; T_1, \ldots , T_{j-1})\to 
(T_j; T_1\cap T_j, \ldots , T_{j-1}\cap T_j)
\end{align*}
introduced in \S 4.2, and the maps of complexes induced by them 
$$
\WT{\Q}\WALT{*}{T; T_1, \ldots , T_{j-1}}
\overset{\iota_j^*}{\underset{p_j^*}{\rightleftarrows }}
\WT{\Q}\WALT{*}{T_j; T_1\cap T_j, \ldots , T_{j-1}\cap T_j}.
$$
Let $\OV{\caF}$ be a hermitian vector bundle on $X$.
Prop.9.4 says that we can obtain a homotopy $\varPhi_{\iota }$ from 
$(\OV{\caF}\otimes \ )\iota_j^*$ to $\iota_j^*(\OV{\caF}\otimes \ )$ 
and a homotopy $\varPhi_p$ from $(\OV{\caF}\otimes \ )p_j^*$ 
to $p_j^*(\OV{\caF}\otimes \ )$.
Let $\varPsi $ be the homotopy from the identity to $\iota_j^*p_j^*$ 
given in Prop.2.20.
Then by Prop.9.6 we have a second homotopy 
$$
\varTheta :\WT{\Q}\WALT{*}{T_j; T_1\cap T_j, \ldots , T_{j-1}\cap T_j}
\to \WT{\Q}\WALT{*+2}{T_j; T_1\cap T_j, \ldots , T_{j-1}\cap T_j}
$$
from 
$\varPhi_{\iota }p_j^*+\iota_j^*\varPhi_p+(\caF \otimes \ )\varPsi $ 
to $\varPsi (\caF \otimes \ )$.
Hence it follows from Prop.2.8 that the diagram 
\begin{equation}
\begin{CD}
\WT{\Q}\WALT{*}{T; T_1, \ldots , T_{j-1}} @>{t_j}>> 
\WT{\Q}\WALT{*}{T; T_1, \ldots , T_j} \\
@V{(\caF \otimes \ )}VV  @VV{(\caF \otimes \ )}V  \\
\WT{\Q}\WALT{*}{T; T_1, \ldots , T_{j-1}} @>{t_j}>> 
\WT{\Q}\WALT{*}{T; T_1, \ldots , T_j}
\end{CD} \label{sec10:cd1}
\end{equation}
is commutative up to homotopy. 
Denote by $\varPi_j$ the homotopy from 
$(\OV{\caF}\otimes \ )t_j$ to 
$t_j(\OV{\caF}\otimes \ )$ given in Prop.2.8.
Since the images of the homotopies 
$\varPhi_{\iota }, \varPhi_p, \varPsi $ and 
the second homotopy $\varTheta $ are isometrically equivalent 
to degenerate elements, $\varPi_j^{m,n}(x)$ is isometrically 
equivalent to a degenerate element for any $m$ and $n$.
Connecting \eqref{sec10:cd1} for all $j$, we can obtain 
the following commutative diagram up to homotopy 
\begin{equation}
\begin{CD}
\WT{\Q}\WALT{*}{T} @>{t}>> 
\WT{\Q}\WALT{*}{T; T_1, \ldots , T_r} \\
@V{(\caF \otimes \ )}VV  @VV{(\caF \otimes \ )}V  \\
\WT{\Q}\WALT{*}{T} @>{t}>> 
\WT{\Q}\WALT{*}{T; T_1, \ldots , T_r},
\end{CD} \label{sec10:cd2}
\end{equation}
and a homotopy $\varPi $ 
from $(\OV{\caF}\otimes \ )t$ to $t(\OV{\caF}\otimes \ )$ is 
given by 
$$
\varPi =\sum_{j=1}^rt_r\cdots t_{j+1}\varPi_jt_{j-1}\cdots t_1.
$$
It is obvious that $\varPi^{0,n}(x)$ is isometrically equivalent 
to a degenerate element for any $n$.
On the other hand, Prop.9.4 says that the diagram 
\begin{equation}
\begin{CD}
\WT{\Q}\WALT{*}{T; T_1, \ldots , T_r} @>{i^*_{\emptyset }}>> 
\WT{\Q}\WALT{*}{X\times \PP{r}; X\times \partial \PP{r}} \\
@V{(\caF \otimes \ )}VV  @VV{(\caF \otimes \ )}V  \\
\WT{\Q}\WALT{*}{T; T_1, \ldots , T_r} @>{i^*_{\emptyset }}>> 
\WT{\Q}\WALT{*}{X\times \PP{r}; X\times \partial \PP{r}} \\
\end{CD} \label{sec10:cd3}
\end{equation}
is also commutative up to homotopy, and if we denote 
by $\varPhi_{i_{\emptyset }}$ the homotopy given 
in Prop.9.4, then $\varPhi_{i_{\emptyset }}^{m,n}(x)$ is 
isometrically equivalent to a degenerate element 
for any $m$ and $n$.
Combining \eqref{sec10:cd2} with \eqref{sec10:cd3} yields 
the diagram 
$$
\begin{CD}
\WT{\Q}\WALT{*}{T} @>{i^*_{\emptyset }t}>> 
\WT{\Q}\WALT{*}{X\times \PP{r}; X\times \partial \PP{r}} \\
@V{(\caF \otimes \ )}VV  @VV{(\caF \otimes \ )}V  \\
\WT{\Q}\WALT{*}{T} @>{i^*_{\emptyset }t}>> 
\WT{\Q}\WALT{*}{X\times \PP{r}; X\times \partial \PP{r}} \\
\end{CD}
$$
which is also commutative up to homotopy, and a homotopy from 
$(\OV{\caF}\otimes \ )i^*_{\emptyset }t$ to 
$i^*_{\emptyset }t(\OV{\caF}\otimes \ )$ is given by 
$\PR{\varPi }=\varPhi_{i_{\emptyset }}t+i^*_{\emptyset }\varPi $.
It is obvious that ${\PR{\varPi }}^{0,n}(x)$ is isometrically 
equivalent to a degenerate element as well.

Any element of $\WH{K}_0(T)$ is represented by a pair 
$(\OV{\caG}, \WT{\tau })$ of a virtual hermitian vector bundle 
$\OV{\caG}$ on $T$ with 
$\WT{\tau }\in \WT{\scD}_{\A,\P,r+1}(X)/\IIm d_s$ 
such that $\CH_{T,0}(\OV{\caG})+d_s\tau =0$.
Since the map \eqref{sec5:map2} 
$$
\WH{i_{\emptyset}^*t}:\WH{K}_0(T)\to 
\QQ{\WH{K}_0(X\times \PP{r}; X\times \partial \PP{r})}
$$
sends $[(\OV{\caG}, \WT{\tau })]$ to 
$[(i_{\emptyset }^*t(\OV{\caG}), -\tau )]$, we have 
$$
[(\OV{\caF}, \WT{\eta })]\cdot \WH{i_{\emptyset}^*t}
\left([(\OV{\caG}, \WT{\tau })]\right)=[(\OV{\caF}\otimes 
i_{\emptyset }^*t(\OV{\caG}), -\CH_0(\OV{\caF})\PRAP \tau )].
$$
On the other hand, since $\CH_{T,0}(\OV{\caG})+d_s\tau =0$, 
we have 
$$
\WH{i_{\emptyset}^*t}\left([(\OV{\caF}, \WT{\eta })]\cdot 
[(\OV{\caG}, \WT{\tau })]\right)=\WH{i_{\emptyset }^*t}
\left([(\OV{\caF}\otimes \OV{\caG}, \CH_0(\OV{\caF})\PRAP \WT{\tau })]
\right)=[(i_{\emptyset }^*t(\OV{\caF}\otimes \OV{\caG}), 
-\CH_0(\OV{\caF})\PRAP \tau )].
$$
Since $\CH_*({\PR{\varPi }}^{0,n}(\OV{\caG}))=0$, we have 
$$
(i_{\emptyset }^*t(\OV{\caF}\otimes \OV{\caG})
-\OV{\caF}\otimes i_{\emptyset }^*t(\OV{\caG}), 0)
=\partial (\PR{\varPi }(\OV{\caG}), 0), 
$$
which means that 
$[(\OV{\caF}, \WT{\eta })]\cdot \WH{i_{\emptyset}^*t}
[(\OV{\caG}, \WT{\tau })]=\WH{i_{\emptyset}^*t}
\left([(\OV{\caF}, \WT{\eta })]\cdot 
[(\OV{\caG}, \WT{\tau })]\right)$.
This completes the proof.
\qed

\setcounter{equation}{0}
\vskip 2pc
\section{Compatibility with product}

\vskip 1pc
\subsection{Product of $\WH{CH}^q(X, r)$ with $\WH{CH}^p(X)$}
In \cite{BF}, Burgos and Feliu constructed a product in 
higher arithmetic Chow groups.
In this subsection, we will recall their construction 
in the case of the product with $\WH{CH}^p(X)$.
We begin by fixing some notations.

Throughout this section, fix $r\geq 1$.
For a variety $X$ defined over a field, 
denote by $\scZ_X^p$ the set of all closed subschemes of $X$ of 
codimension $p$, and by $\scZ^q_{X,r}$ the set of all 
admissible subschemes of $X\times \PP{r}$ of codimension $q$.
For another variety $Y$, write 
$\scZ^{p,q}_{X,Y,r}=\scZ_X^p\times \scZ_{Y,r}^q$.
We regard $\scZ^{p,q}_{X,Y,r}$ as a subset of 
$\scZ^{p+q}_{X\times Y,r}$ in the way that 
$(Z, W)\mapsto Z\times W$.
In what follows, we identify $\scZ_X^p$ with 
$\scZ_{X,Y,r}^{p,0}$ by $Z\mapsto Z\times Y\times \PP{r}$ 
and $\scZ_{Y,r}^q$ with $\scZ_{X,Y,r}^{0,q}$ by $W\mapsto X\times W$.

Suppose $X$ and $Y$ are compact complex algebraic manifolds.
Write $\PP{r}_{X,Y}=X\times Y\times \PP{r}$ for short.
Then it is shown in \cite[Lem.5.5]{BF} that 
\begin{multline*}
0\to \scD_{\log}^*(\PP{r}_{X,Y}-\scZ_{X,Y,r}^{p,q}, p+q)\to 
\scD_{\log}^*(\PP{r}_{X,Y}-\scZ_X^p, p+q)\oplus 
\scD_{\log}^*(\PP{r}_{X,Y}-\scZ_{Y,r}^q, p+q) \\
\overset{j_{X,Y,r}^{p,q}}{\longrightarrow }
\scD_{\log}^*(\PP{r}_{X,Y}-\scZ_X^p\cup \scZ_{Y,r}^q, p+q)\to 0
\end{multline*}
is a short exact sequence, where the map $j_{X,Y,r}^{p,q}$ 
is given by $(\omega , \tau )\mapsto -\omega +\tau $.
This implies that the natural map 
\begin{equation}
\scD_{\log}^*(\PP{r}_{X,Y}-\scZ_{X,Y,r}^{p,q}, p+q)
\to s(-j_{X,Y,r}^{p,q})^* \label{sec11:iso1}
\end{equation}
is a quasi-isomorphism.
Let 
$$
i_{X,Y,r}^{p,q}:\scD_{\log}^*(\PP{r}_{X,Y}, p+q)\to 
s(-j_{X,Y,r}^{p,q})^*
$$
be the composite of the restriction to 
$\scD_{\log}^*(\PP{r}_{X,Y}-\scZ_{X,Y,r}^{p,q}, p+q)$ with 
\eqref{sec11:iso1}, and take the simple complex 
$s(i^{p,q}_{X,Y,r})^*$.
To be more precise, 
\begin{align*}
s(i_{X,Y,r}^{p,q})^n=&\, \scD_{\log}^n(\PP{r}_{X,Y}, p+q)\oplus 
\scD_{\log}^{n-1}(\PP{r}_{X,Y}-\scZ_X^p, p+q)  \\
&\oplus \scD_{\log}^{n-1}(\PP{r}_{X,Y}-\scZ_{Y,r}^q, p+q)\oplus 
\scD_{\log}^{n-2}(\PP{r}_{X,Y}-\scZ_X^p\cup \scZ_{Y,r}^q, p+q) 
\end{align*}
with the differential given by 
$$
d_s(\omega_0, (\omega_1, \omega_2), \omega_3)=
(d_{\scD}\omega_0, (\omega_0-d_{\scD}\omega_1, 
\omega_0-d_{\scD}\omega_2), -\omega_1+\omega_2+d_{\scD}\omega_3).
$$
If we set 
$$
\scD_{\log,\scZ_{X,Y,r}^{p,q}}^*(\PP{r}_{X,Y}, p+q)=
s(\scD_{\log}^*(\PP{r}_{X,Y}, p+q)\to 
\scD_{\log}^*(\PP{r}_{X,Y}-\scZ_{X,Y,r}^{p,q}, p+q)), 
$$
then the natural map 
\begin{equation}
\scD_{\log,\scZ_{X,Y,r}^{p,q}}^*(\PP{r}_{X,Y}, p+q)\to 
s(i_{X,Y,r}^{p,q})^* \label{sec11:iso2}
\end{equation}
is a quasi-isomorphism.
Take the truncated subcomplex 
$$
\scD_{\A,\scZ_{X,Y}^{p,q}}^{*,-r}(X\times Y, p+q)=
\tau_{\leq 2(p+q)}
\scD_{\log,\scZ_{X,Y,-r}^{p,q}}^*(\PP{-r}_{X,Y}, p+q).
$$
Let $\scD_{\A,\scZ_{X,Y}^{p,q}}^*(X\times Y, p+q)_0$ be 
the single complex of the normalized subcomplex of \linebreak 
$\scD_{\A,\scZ_{X,Y}^{p,q}}^{s,-r}(X\times Y, p+q)$ with 
respect to the cubical structure on the index $r$.
In the same way, we can obtain the single complex of the normalized 
subcomplex of $\tau_{\leq 2(p+q)}s(i_{X,Y,-r}^{p,q})^*$ with 
respect to the cubical structure on the index $r$, which 
we denote by $s_{\A}(i_{X,Y}^{p,q})_0^*$.
Then the map \eqref{sec11:iso2} induces a quasi-isomorphism 
$$
\scD_{\A,\scZ_{X,Y}^{p,q}}^*(X\times Y, p+q)_0\to 
s_{\A}(i_{X,Y}^{p,q})_0^*.
$$

Let 
\begin{align*}
\pi_X&:\PP{r}_{X,Y}=X\times Y\times \PP{r}\to X, \\
\pi_Y&:\PP{r}_{X,Y}=X\times Y\times \PP{r}\to Y\times \PP{r}
\end{align*}
be the projections, and define an exterior product of 
$\omega \in \mmD^*(X, p)$ with 
$\PR{\omega }\in \mmD_{\log}^*(Y\times \PP{r}, q)$ as follows:
$$
\omega \EPRA \PR{\omega }=\pi_X^*\omega \bullet 
\pi_Y^*\PR{\omega }\in \mmD_{\log}^*(\PP{r}_{X,Y}, p+q).
$$
Then it induces a map of complexes 
$$
\EPRA :\mmD^*(X, p)\times \mmD_{\log}^*(Y\times \PP{r}, q)\to 
\mmD_{\log}^*(\PP{r}_{X,Y}, p+q).
$$
Similarly, for $g\in \mmD^*_{\log}(X-\scZ_X^p, p)$ and 
$\PR{g}\in \mmD_{\log}^*(Y\times \PP{r}-\scZ_{Y,r}^q, q)$, 
define 
\begin{align*}
g\EPRA \PR{g}&=\pi_X^*g\bullet \pi_Y^*\PR{g}\in 
\mmD_{\log}^*(\PP{r}_{X,Y}-\scZ_X^p\cup \scZ_{Y,r}^q, p+q),  \\
\omega \EPRA \PR{g}&=\pi_X^*\omega \bullet \pi_Y^*\PR{g}
\in \mmD_{\log}^*(\PP{r}_{X,Y}-\scZ_{Y,r}^q, p+q),  \\
g\EPRA \PR{\omega }&=\pi_X^*g \bullet \pi_Y^*\PR{\omega }\in 
\mmD_{\log}^*(\PP{r}_{X,Y}-\scZ_X^q, p+q).
\end{align*}
Using these products, we can define 
\begin{equation}
\diamond :\mmD_{\scZ^p}^n(X, p)\times 
\mmD_{\scZ_{Y,r}^q}^m(Y\times \PP{r}, q)\to 
\tau_{\leq 2p+2q}s(i_{X,Y,r}^{p,q})^{n+m} \label{sec11:map1}
\end{equation}
by 
$$
(\omega , g)\diamond (\PR{\omega }, \PR{g})\mapsto 
(\omega \EPRA \PR{\omega }, (g\EPRA \PR{\omega }, 
(-1)^n\omega \EPRA \PR{g}), (-1)^{n-1}g\EPRA \PR{g}), 
$$
which turns out to be a map of complexes.
Taking the cohomology yields 
$$
H_{\scD,\scZ^p}^{2p}(X, \R(p))\times 
H_{\scD,\scZ_{Y,r}^q}^{2q}(Y\times \PP{r}, \R(q))\to 
H^{2p+2q}(s(i_{X,Y,r}^{p,q})^*).
$$
Combining this map with 
$$
H^{2p+2q}(s(i_{X,Y,r}^{p,q})^*)\simeq 
H_{\scD,\scZ^{p,q}_{X,Y,r}}^{2p+2q}(\PP{r}_{X,Y}, \R(p+q))\to 
H_{\scD,\scZ^{p+q}_{X\times Y,r}}^{2p+2q}(\PP{r}_{X,Y}, \R(p+q))
$$
induced by \eqref{sec11:iso2} and the inclusion 
$\scZ^{p,q}_{X,Y,r}\subset \scZ^{p+q}_{X\times Y,r}$, and taking 
the normalized subcomplexes, we obtain a map 
$$
\times :\caH^p(X, 0)\times \caH^q(Y, r)_0\to 
\caH^{p+q}(X\times Y, r)_0.
$$
Prop.5.13 in \cite{BF} says that 
$$
\begin{CD}
Z^p(X, 0)\times Z^q(Y, r)_0 @>{\times }>> 
Z^{p+q}(X\times Y, r)_0 \\
@VV{\chi_1\times \chi_1}V  @VV{\chi_1}V  \\
\caH^p(X, 0)\times \caH^q(Y, r)_0 @>{\times }>> 
\caH^{p+q}(X\times Y, r)_0
\end{CD}
$$
is commutative, where the upper horizontal arrow is the exterior 
product of cycles.
On the other hand, taking the single complexes of the normalized 
subcomplexes on the both sides of \eqref{sec11:map1}, we have 
$$
\diamond :\mmD_{\scZ^p}^*(X, p)\times 
\scD_{\A,\scZ^q}^*(Y, q)_0\to s_{\A}(i_{X,Y}^{p,q})_0^*.
$$
An element $a\in s_{\A}(i^{p,q}_{X,Y})_0^{2p+2q-r}$ is written as 
$$
a=\sum_{j=0}^ra_j, \qquad 
a_j\in \tau_{\leq 2p+2q}s(i^{p,q}_{X,Y,j})_0^{2p+2q-r+j}.
$$
Then by the quasi-isomorphism \eqref{sec11:iso2} 
we have the cohomology class 
$$
[a_r]\in 
H_{\scD,\scZ^{p,q}_{X,Y,r}}^{2(p+q)}(\PP{r}_{X\times Y}, p+q)_0.
$$
Since $\scZ^{p,q}_{X,Y,r}\subset \scZ^{p+q}_{X\times Y,r}$, 
the correspondence $a\mapsto [a_r]$ gives 
$$
\chi_2:s_{\A}(i^{p,q}_{X,Y})_0^{2p+2q-*}\to 
\caH^{p+q}(X\times Y, *)_0, 
$$
which turns out be a map of complexes.
Cor.5.14 in \cite{BF} says that 
$$
\begin{CD}
\mmD_{\scZ_p}^{2p}(X, p)\times \scD_{\A,\scZ^q}^{2q-*}(Y, q)_0 
@>{\diamond }>> s_{\A}(i^{p,q}_{X,Y})_0^{2p+2q-*}  \\
@VV{\chi_2\times \chi_2}V  @VV{\chi_2}V  \\
\caH^p(X, 0)\times \caH^q(Y, *)_0 @>{\times }>> 
\caH^{p+q}(X\times Y, *)_0
\end{CD}
$$
is commutative.

Denote the following canonical maps by the same symbol $\rho $:
\begin{gather*}
\rho :\mmD_{\scZ^p}^*(X, p)\to \mmD^*(X, p),  \\
\rho :\scD_{\A,\scZ^q}^*(Y, q)_0\to \scD_{\A}^*(Y, q)_0
\simeq \WT{\scD}_{\A}^*(Y, q)\to \WT{\scD}_{\A,\P}^*(Y, q), \\
\rho :s_{\A}(i_{X,Y}^{p,q})_0^*\to 
\scD_{\A}^*(X\times Y, p+q)_0\simeq 
\WT{\scD}_{\A}^*(X\times Y, p+q)\to \WT{\scD}_{\A,\P}^*(X\times Y, p+q).
\end{gather*}
Moreover, define a map of complexes 
\begin{equation}
\EPRAP :\mmD^*(X, p)\times \WT{\scD}_{\A,\P}^*(Y, q)\to 
\WT{\scD}_{\A,\P}^*(X\times Y, p+q) \label{sec11:pro0}
\end{equation}
by 
$\omega \EPRAP \PR{\omega }=
\pi_X^*\omega \bullet \pi_Y^*\PR{\omega }$, where 
\begin{align*}
\pi_X:&\, X\times Y\times \PP{r}\times (\P^1)^s\to X,  \\
\pi_Y:&\, X\times Y\times \PP{r}\times (\P^1)^s\to 
Y\times \PP{r}\times (\P^1)^s
\end{align*}
are the projections.
Then the diagram 
$$
\begin{CD}
\mmD_{\scZ^p}^*(X, p)\times \scD_{\A,\scZ^q}^*(Y, q)_0 
@>{\diamond }>> s_{\A}(i_{X,Y}^{p,q})_0^*   \\
@V{\rho \times \rho }VV   @VV{\rho }V  \\
\mmD^*(X, p)\times \WT{\scD}^*_{\A,\P}(X\times Y, p+q) 
@>{\EPRAP }>> \WT{\scD}^*_{\A,\P}(X\times Y, p+q)
\end{CD}
$$
is commutative.

Under the above preparation, we can define an exterior product 
of $\WH{CH}^p(X)$ with $\WH{CH}^q(Y, r)$.
Let $(y, (\omega , \WT{g}))$ be a pair of $y\in Z^p(X)$ with 
a Green form $(\omega , \WT{g})$ associated with $y$, and 
take a lift $g\in \mmD_{\log}^{2p-1}(X-\scZ^p, p)$.
Let 
$$
\left(
\text{
\setlength{\unitlength}{1mm}
\begin{picture}(12,7)
 \put(-2,-4){$z$}
 \put(2,3){$\beta_1$}
 \put(6,-4){$\alpha $}
 \put(10,3){$\beta_2$}
\end{picture}
}
\right)\in \WH{Z}^q(Y, r)_0.
$$
Define their product as follows:
\begin{equation}
(y, (\omega , g))\cdot 
\left(
\text{
\setlength{\unitlength}{1mm}
\begin{picture}(12,7)
 \put(-2,-4){$z$}
 \put(2,3){$\beta_1$}
 \put(6,-4){$\alpha $}
 \put(10,3){$\beta_2$}
\end{picture}
}
\right)=\left(
\text{
\setlength{\unitlength}{1mm}
\begin{picture}(53,7)
 \put(-2,-4){$y\times z$}
 \put(7,3){$\chi_1(y)\times \beta_1$}
 \put(23,-4){$(\omega , g)\diamond \alpha $}
 \put(39,3){$\omega \EPRAP \beta_2$}
\end{picture}
}
\right).  \label{sec11:pro1}
\end{equation}
Then as shown in \cite[\S 5]{BF}, this product gives a map 
\begin{align}
\WH{CH}^p(X)&\times \WH{CH}^q(Y, r) \notag \\
&\longrightarrow H_r\left(
\text{
\setlength{\unitlength}{1mm}
\begin{picture}(95,12)
 \put(-2,-8){$Z^{p+q}(X\times Y, *)_0$}
 \put(17,7){$\caH^{p+q}(X\times Y, *)_0$}
 \put(40,-8){$s_{\A}(i_{X,Y}^{p,q})_0^{2p+2q-*}$}
 \put(57,7){$\WH{\scD}_{\A,\P}^{2p+2q-*}(X\times Y, p+q)$}
 \put(10,0){$\scriptstyle{\chi_1}$}
 \put(43,0){$\scriptstyle{\chi_2}$}
 \put(62,0){$\scriptstyle{\rho }$}
 \put(12,-3){\vector(1,1){7}}
 \put(45,-3){\vector(-1,1){7}}
 \put(62,-3){\vector(1,1){7}}
\end{picture}
}
\right), \label{sec11:map3}
\end{align}
where $\WH{\scD}_{\A,\P}^*(X\times Y, p+q)$ is 
the complex defined in \S 6.3.
On the other hand, it is obvious that the diagram 
\begin{equation}
\text{
\setlength{\unitlength}{1mm}
\begin{picture}(120,13)
 \put(0,8){$s_{\A}(i^{p,q}_{X,Y})_0^*\overset{\sim }{\leftarrow }
\scD_{\A,\scZ^{p,q}_{X,Y}}^*(X\times Y, p+q)_0\to 
\scD_{\A,\scZ^{p+q}_{X\times Y}}^*(X\times Y, p+q)_0$}
 \put(30,-8){$\WH{\scD}_{\A,\P}^*(X\times Y, p+q)$}
 \put(15,5){\vector(2,-1){17}}
 \put(80,3){\vector(-2,-1){14}}
 \put(25,1){$\scriptstyle{\rho }$}
 \put(71,0){$\scriptstyle{\rho }$}
\end{picture}
} \label{sec11:map4}
\end{equation}
\vskip 2pc
\noindent
is commutative, hence it gives 
\begin{align}
H_r&\left(
\text{
\setlength{\unitlength}{1mm}
\begin{picture}(93,12)
 \put(-2,-8){$Z^{p+q}(X\times Y, *)_0$}
 \put(18,7){$\caH^{p+q}(X\times Y, *)_0$}
 \put(40,-8){$s_{\A}(i_{X,Y}^{p,q})_0^{2p+2q-*}$}
 \put(55,7){$\WH{\scD}_{\A,\P}^{2p+2q-*}(X\times Y, p+q)$}
 \put(14,-3){\vector(1,1){7}}
 \put(46,-3){\vector(-1,1){7}}
 \put(60,-3){\vector(1,1){7}}
\end{picture}
}
\right) \notag \\
\simeq &H_r\left(
\text{
\setlength{\unitlength}{1mm}
\begin{picture}(93,12)
 \put(-2,-8){$Z^{p+q}(X\times Y, *)_0$}
 \put(18,7){$\caH^{p+q}(X\times Y, *)_0$}
 \put(40,-8){
$\scD_{\A,\scZ^{p,q}_{X,Y}}^{2p+2q-*}(X\times Y, p+q)_0$}
 \put(55,7){$\WH{\scD}_{\A,\P}^{2p+2q-*}(X\times Y, p+q)$}
 \put(14,-3){\vector(1,1){7}}
 \put(46,-3){\vector(-1,1){7}}
 \put(60,-3){\vector(1,1){7}}
\end{picture}
}
\right) \label{sec11:map5} \\
&\longrightarrow H_r\left(
\text{
\setlength{\unitlength}{1mm}
\begin{picture}(93,12)
 \put(-2,-8){$Z^{p+q}(X\times Y, *)_0$}
 \put(18,7){$\caH^{p+q}(X\times Y, *)_0$}
 \put(40,-8){
$\scD_{\A,\scZ^{p+q}_{X\times Y}}^{2p+2q-*}(X\times Y, p+q)_0$}
 \put(55,7){$\WH{\scD}_{\A,\P}^{2p+2q-*}(X\times Y, p+q)$}
 \put(14,-3){\vector(1,1){7}}
 \put(46,-3){\vector(-1,1){7}}
 \put(60,-3){\vector(1,1){7}}
\end{picture}
}
\right).\notag
\end{align}
\vskip 1pc
\noindent
Define an exterior product of higher arithmetic Chow groups 
$$
\cup :\WH{CH}^p(X)\times \WH{CH}^q(Y, r)
\to \WH{CH}^{p+q}(X\times Y, r)
$$
to be the composite of \eqref{sec11:map3} with \eqref{sec11:map5}.
Substituting $\sigma_{<2q-r}\WH{\scD}_{\A,\P}^*(Y, q)$ for 
$\WH{\scD}_{\A,\P}^*(Y, q)$ and 
$\sigma_{<2p+2q-r}\WH{\scD}_{\A,\P}^*(X\times Y, q)$ for 
$\WH{\scD}_{\A,\P}^*(X\times Y, q)$, 
we can also define an extended exterior product 
$$
\cup :\WH{CH}^p(X)\times \WH{Z}^q(Y, r)_0^*/\IIm \partial 
\to \WH{Z}^{p+q}(X\times Y, r)_0^*/\IIm \partial .
$$

Suppose that $X$ is projective over an arithmetic field.
As we have seen in \S 8.2, we can define the pull-back map 
$$
\WH{\Delta}^*:\WH{CH}^{p+q}(X\times X, r)\to \WH{CH}^{p+q}(X, r)
$$
by the diagonal map.
Define intersection product of higher arithmetic Chow groups 
$$
\WH{CH}^p(X)\times \WH{CH}^q(X, r)\to \WH{CH}^{p+q}(X, r) 
$$
to be the composite of the exterior product in the case that 
$X=Y$ with the pull-back map $\WH{\Delta}^*$.

\vskip 1pc
\begin{prop}
Let $\eta^p\in \mmD^{2p-1}(X, p)$.
If $r\geq 1$, then for any $x\in \WH{CH}^q(Y, r)$ the exterior 
product $a(\WT{\eta^p})\cup x$ is equal to zero.
\end{prop}

{\it Proof}: 
Suppose that $x$ is represented by 
\begin{equation}
\left(
\text{
\setlength{\unitlength}{1mm}
\begin{picture}(12,7)
 \put(-2,-4){$y$}
 \put(2,3){$\beta_1$}
 \put(6,-4){$\alpha $}
 \put(10,3){$\beta_2$}
\end{picture}
}
\right)\in \KER \partial \subset \WH{Z}^p(Y, r)_0. 
\label{sec11:ele1}
\end{equation}
Then in 
$$
H_r\left(\text{
\setlength{\unitlength}{1mm}
\begin{picture}(96,13)
 \put(-2,-7){$Z^{p+q}(X\times Y, *)_0$}
 \put(17,8){$\caH^{p+q}(X\times Y, *)_0$}
 \put(40,-7){$s_{\A}(i_{X,Y}^{p,q})_0^{2p+2q-*}$}
 \put(57,8){$\WH{\scD}_{\A,\P}^{2p+2q-*}(X\times Y, p+q)$}
 \put(10,1){$\scriptstyle{\chi_1}$}
 \put(43,1){$\scriptstyle{\chi_2}$}
 \put(62,1){$\scriptstyle{\rho }$}
 \put(12,-2){\vector(1,1){7}}
 \put(45,-2){\vector(-1,1){7}}
 \put(62,-2){\vector(1,1){7}}
\end{picture}
}
\right), 
$$
the product of $(0, (d_{\scD}\eta^p, \eta^p))$ with 
the element \eqref{sec11:ele1} is given by 
\begin{equation}
\left(
\text{
\setlength{\unitlength}{1mm}
\begin{picture}(44,7)
 \put(-2,-4){$0$}
 \put(2,3){$0$}
 \put(6,-4){$(d_{\scD}\eta^p, \eta^p)\diamond \alpha $}
 \put(25,3){$d_{\scD}\eta^p\EPRAP \beta_2$}
\end{picture}
}
\right). \label{sec11:ele2}
\end{equation}
If we write $\alpha =(\alpha_1, \alpha_2)\in 
\scD^{2q-r}_{\A,\scZ^q}(Y, q)_0$ 
with $\alpha_1\in \scD^{2p-r}_{\A}(Y, q)_0$ and 
$\alpha_2\in \scD^{2p-r-1}_{\A}(Y-\scZ_Y^q, q)_0$, then it holds in 
$s(i_{X,Y,r}^{p,q})^{2p+2q-r}_0$ that 
\begin{align*}
(d_{\scD}\eta^p, \eta^p)\diamond \alpha =&\, 
(d_{\scD}\eta^p, \eta^p)\diamond (\alpha_1, \alpha_2) \\
=&\, (d_{\scD}\eta^p\EPRA \alpha_1, (\eta^p\EPRA \alpha_1, 
d_{\scD}\eta^p\EPRA \alpha_2), -\eta^p\EPRA \alpha_2).
\end{align*}
Since $\eta^p\in \mmD^{2p-1}(X, p)$ and $d_s\alpha =0$, 
it is equal to 
$$
d_s(\eta^p\EPRA \alpha_1, (0, -\eta^p\EPRA \alpha_2), 0).
$$
It follows from the definition of $\chi_2$ that 
$\chi_2(\eta^p\EPRA \alpha_1, (0, -\eta^p\EPRA \alpha_2), 0)=0$.
Hence \eqref{sec11:ele2} is homologous to 
\begin{equation}
\left(\text{
\setlength{\unitlength}{1mm}
\begin{picture}(49,7)
 \put(-2,-4){$0$}
 \put(2,3){$0$}
 \put(6,-4){$0$}
 \put(10,3){$d_{\scD}\eta^p\EPRAP \beta_2-\eta^p\EPRAP \alpha_1$}
\end{picture}
}
\right).  \label{sec11:ele3}
\end{equation}
The $\partial $-closedness of \eqref{sec11:ele1} implies 
that $d_s\beta_2=\rho (\alpha )=\alpha_1$, therefore 
$$
d_{\scD}\eta^p\EPRAP \beta_2-\eta^p\EPRAP \alpha_1=
d_s(\eta^p\EPRAP \beta_2), 
$$
which means that \eqref{sec11:ele3} is $\partial $-exact.
Hence the proposition follows.
\qed

\vskip 1pc
\subsection{The main theorem}
In this subsection we prove the following theorem:

\vskip 1pc
\begin{thm}
Let $X$ be a smooth projective variety defined over an arithmetic 
field.
Then 
$$
\begin{CD}
\WH{K}_0(X)\times \QQ{\WH{K}_r(X)} @>>> 
\QQ{\WH{K}_r(X)}  \\
@V{\underset{p+q=n}{\oplus }\WH{\CH}_0^p\times \WH{\CH}^q_r}VV  
@VV{\WH{\CH}^n_r}V  \\
\underset{p+q=n}{\oplus }\WH{CH}^p(X)\times \WH{CH}^q(X, r)
@>>> \WH{CH}^n(X, r)
\end{CD}
$$
is commutative.
\end{thm}

{\it Proof}:
We have shown in \S 10 that the sequence of maps 
$$
\QQ{\WH{K}_r(X)}\simeq 
\QQ{\WH{K}_0(X\times \PP{r}; X\times \partial \PP{r})}
\twoheadleftarrow 
\QQ{\WH{K}_0(T; T_1, \ldots , T_r)}\subset \QQ{\WH{K}_0(T)}
$$
respects the $\WH{K}_0(X)$-module structures.
Hence we have only to show that the diagram 
\begin{equation}
\begin{CD}
\WH{K}_0(X)\times \WH{K}_0(T) @>>> \WH{K}_0(T)  \\
@V{\underset{p+q=n}{\oplus }\WH{\CH}_0^p\times \WH{\CH}^q_{T,0}}VV  
@VV{\WH{\CH}^n_{T,0}}V  \\
\underset{p+q=n}{\oplus }\WH{CH}^p(X)\times \WH{CH}^q(X, r)
@>>> \WH{CH}^n(X, r)
\end{CD} \label{sec11:cd1}
\end{equation}
is commutative.

First we consider the product with 
$x_1=[(0, \WT{\eta })]\in \WH{K}_0(X)$ where $\eta \in \mmD_1(X)$.
In this case it follows from the definition of the product of 
arithmetic $K$-groups described as \eqref{sec10:pro1} that 
$x_1\cdot x_2=0$ for any $x_2\in \WH{K}_0(T)$, and Prop.11.1 says 
that $a(\WT{\eta^p})\cdot \WH{\CH}_{T,0}^q(x_2)=0$.
Hence the theorem follows in this case.

Let us now introduce exterior product of arithmetic $K_0$-groups.
Let $Y$ be a smooth projective variety, and 
$T=D(Y\times \PP{r}; Y\times \partial \PP{r})$ the associated iterated double.
We identify $X\times T$ with 
$\PR{T}=T(X\times Y\times \PP{r}; X\times Y\times \partial \PP{r})$.
Define an exterior product 
$$
\cup :\WH{K}_0(X)\times \WH{K}^M_0(T)\to \WH{K}^M_0(\PR{T})
$$
by 
$$
[(\OV{\caF}, \WT{\eta })]\cup [(\OV{\PR{\caF}}, \WT{\tau })]=
[(\OV{\caF}\boxtimes \OV{\PR{\caF}}, \CH_0(\OV{\caF})\EPRAP \WT{\tau }
+\WT{\eta }\EPRAP \CH_{T,0}(\OV{\PR{\caF}})
+d_{\scD}\eta \EPRAP \WT{\tau })].
$$
In the above, $\OV{\caF}\boxtimes \OV{\PR{\caF}}=
\pi_X^*\OV{\caF}\otimes \pi_T^*\OV{\PR{\caF}}$, where 
$\pi_X:\PR{T}=X\times T\to X$ and $\pi_T:|PR{T}=X\times T\to T$ 
are the projections, and the product $\EPRAP $ is defined 
in \eqref{sec11:pro0}.
We can show in the same way as the proof of Prop.10.3 that 
this product is well-defined, and it gives a map 
$$
\cup :\WH{K}_0(X)\times \WH{K}_0(T)\to \WH{K}_0(\PR{T}).
$$

Here we assume that $X=Y$.
Let $\Delta :X\to X\times X$ be the diagonal morphism, 
and $\Delta_D:T\to \PR{T}$ the morphism between the iterated 
doubles induced by $\Delta $. 
Then the diagram \eqref{sec11:cd1} is decomposed as follows:
$$
\begin{CD}
\WH{K}_0(X)\times \WH{K}_0(T) @>{\cup }>> \WH{K}_0(\PR{T}) 
@>{\Delta_D^*}>> \WH{K}_0(T)  \\
@V{\underset{p+q=n}{\oplus }\WH{\CH}_0^p\times \WH{\CH}^q_{T,0}}VV  
@VV{\WH{\CH}^n_{\PR{T},0}}V  @VV{\WH{\CH}^n_{T,0}}V  \\
\underset{p+q=n}{\oplus }\WH{CH}^p(X)\times \WH{CH}^q(X, r)
@>{\cup }>> \WH{CH}^n(X\times X, r) @>{\Delta^*}>> \WH{CH}^n(X, r).
\end{CD}
$$
The right diagram is commutative by \eqref{sec8:cd3}.
Hence the theorem follows from the lemma below.
\qed

\vskip 1pc
\begin{lem}
For a hermitian vector bundle $\OV{\caF}$ on $X$, denote by 
$$
\OV{\caF}\, \boxtimes \ :\WH{K}_0(T)\to \WH{K}_0(\PR{T})
$$
the exterior product with $[(\OV{\caF}, 0)]\in \WH{K}_0(X)$.
Then the diagram 
$$
\begin{CD}
\WH{K}_0(T) @>{\OV{\caF}\boxtimes }>> \WH{K}_0(\PR{T})  \\
@V{\underset{p}{\oplus }\WH{\CH}^{n-p}_{T,0}}VV  
@VV{\WH{\CH}^n_{\PR{T},0}}V  \\
\underset{p}{\oplus }\WH{CH}^{n-p}(X, r)
@>>{\underset{p}{\sum }\WH{\CH}_0^p(\OV{\caF})\cup }> 
\WH{CH}^n(X\times Y, r)
\end{CD}
$$
is commutative.
\end{lem}

{\it Proof}: 
In the previous subsection we have shown that the exterior product 
is extended to 
$$
\WH{CH}^p(X)\times \WH{Z}^q(Y, r)_0^*/\IIm \partial \to 
\WH{Z}^{p+q}(X\times Y, r)_0^*/\IIm \partial .
$$
Hence it suffices to prove the commutativity of the extended diagram 
$$
\begin{CD}
\WH{K}_0^M(T) @>{\OV{\caF}\boxtimes }>> \WH{K}^M_0(\PR{T})  \\
@V{\underset{p}{\oplus }\WH{\CH}^{n-p}_{T,0}}VV  
@VV{\WH{\CH}^n_{\PR{T},0}}V  \\
\underset{p}{\oplus }\WH{Z}^{n-p}(Y, r)_0^*/\IIm \partial 
@>>{\underset{p}{\sum }\WH{\CH}_0^p(\OV{\caF})\cup }> 
\WH{Z}^n(X\times Y, r)_0^*/\IIm \partial .
\end{CD}
$$

First we consider the case that $x_2=[(0, \WT{\tau })]\in K_0^M(T)$ 
with $\tau \in \WT{\scD}_{\A,\P,r+1}(X)$.
In this case, 
$$
\WH{\CH}_{\PR{T},0}^n(\OV{\caF}\boxtimes x_2)=
\WH{\CH}_{\PR{T},0}^n([(0, \CH_0(\OV{\caF})\EPRAP \WT{\tau })])
=\underset{p+q=n}{\sum }a(\CH_0^p(\OV{\caF})\EPRAP \WT{\tau^q}).
$$
On the other hand, if $\WH{\CH}_0^p(\OV{\caF})$ is represented 
by $(y^p, (\omega^p, \WT{g^p}))$ such that $y^p\in Z^p(X)$ and 
$(\omega^p, \WT{g^p})$ is a Green form associated with $y^p$, 
then it follows from the definition of the product 
\eqref{sec11:pro1} that 
$$
\WH{\CH}_0^p(\OV{\caF})\cup a(\WT{\tau^q})
=\left[\left(
\text{
\setlength{\unitlength}{1mm}
\begin{picture}(28,7)
 \put(-2,-4){$0$}
 \put(2,3){$0$}
 \put(6,-4){$0$}
 \put(10,3){$-\omega^p\EPRAP \tau^q$}
\end{picture}
}
\right)\right]=
a(\CH_0^p(\OV{\caF})\EPRAP \WT{\tau^q})
$$
since $\CH_0^p(\OV{\caF})=\omega^p$.
Hence the lemma follows in this case.

Finally let us assume that $x_2=[(\OV{\PR{\caF}}, 0)]$.
In this case what we need to show is the equality 
\begin{equation}
\WH{\CH}_{\PR{T},0}^n(\OV{\caF}\boxtimes \OV{\PR{\caF}})=
\sum_{p+q=n}\WH{\CH}_0^p(\OV{\caF})\cup 
\WH{\CH}_{T,0}^q(\OV{\PR{\caF}})  \label{sec11:eq1}
\end{equation}
in $\WH{Z}^n(X\times Y, r)_0^*/\IIm \partial $.
Take a morphism $\varphi :T\to G$ to a smooth projective variety 
$G$ and a hermitian vector bundle $\OV{\PR{\caG}}$ on $G$ such 
that $\varphi^*\PR{\caG}\simeq \PR{\caF}$.
Let $(y^p, (\omega^p, \WT{g^p}))$ and 
$(z^q, (\tau^q, \WT{h^q}))$ be representatives of 
$\WH{\CH}_0^p(\OV{\caF})$ and 
$\WH{\CH}_0^q(\OV{\PR{\caG}})$ respectively.
Assume that $z^q\in Z^q_{\varphi }(G)$.
Then the exterior product of $\WH{\CH}_0^p(\OV{\caF})$ with 
$\WH{\CH}_{T,0}^q(\OV{\PR{\caF}})$ in 
{\small
$$
H_r\left(
\text{
\setlength{\unitlength}{1mm}
\begin{picture}(101,12)
 \put(-2,-8){$Z^{p+q}(X\times Y, *)_0$}
 \put(16,7){$\caH^{p+q}(X\times Y, *)_0$}
 \put(38,-8){$s_{\A}(i_{X,Y}^{p,q})_0^{2p+2q-*}$}
 \put(50,7){$\sigma_{<2p+2q-r}
\WH{\scD}_{\A}^{2p+2q-*}(X\times Y, p+q)$}
 \put(14,-3){\vector(1,1){7}}
 \put(46,-3){\vector(-1,1){7}}
 \put(60,-3){\vector(1,1){7}}
\end{picture}
}
\right)
$$}
is equal to 
\begin{equation}
\left[\left(\text{
\setlength{\unitlength}{1mm}
\begin{picture}(101,7)
 \put(-2,-4){$y^p\times \varphi^*(z^q)$}
 \put(18,3){$0$}
 \put(21,-4){$(\omega^p, g^p)\diamond (\varphi^*(\tau^q), \varphi^*(h^q))$}
 \put(44,3){$(-1)^{r+1}\CH_0^p(\OV{\caF})\EPRAP 
\CH_{T,1}^q(\OV{\PR{\caF}}, \varphi^*\OV{\PR{\caG}})$}
\end{picture}
}
\right)\right]. \label{sec11:ele4}
\end{equation}
On the other hand, since there is an isomorphism 
$(\IId_X\times \varphi )^*(\caF\boxtimes \PR{\caG})\simeq 
\caF\boxtimes \PR{\caF}$, 
\begin{multline*}
\WH{\CH}^n(\OV{\caF}\boxtimes \OV{\PR{\caF}})=  \\
\left[\left(
\text{
\setlength{\unitlength}{1mm}
\begin{picture}(125,7)
 \put(-2,-4){$\underset{p+q=n}{\ssum }(1\times \varphi )^*
(y^p\times z^q)$}
 \put(38,3){$0$}
 \put(40,-4){$\underset{p+q=n}{\ssum }(1\times \varphi )^*
\left(\pi_X^*(\omega^p, g^p)\ast \pi_G^*(\tau^q, h^q)\right)$}
 \put(71,3){$(-1)^{r+1}\CH_{\PR{T},1}^n(\OV{\caF}\boxtimes 
\PR{\OV{\caF}}, \OV{\caF}\boxtimes \varphi^*\OV{\PR{\caG}})$}
\end{picture}
}
\right)\right], 
\end{multline*}
\vskip 1pc
\noindent
where $\pi_X^*(\omega^p, g^p)\ast \pi_G^*(\tau^q, h^q)\in 
\scD^{2n}_{\A,\scZ^n_{X\times G,1\times \varphi }}
(X\times G, n)_0$ is a representative of the star product 
of the Green forms $\pi_X^*(\omega^p, \WT{g^p})$ and 
$\pi_G^*(\tau^q, \WT{h^q})$.
It is easy to see that 
$$
\underset{p+q=n}{\sum }\CH_0^p(\OV{\caF})\EPRAP 
\CH_{T,1}^q(\OV{\PR{\caF}}, \varphi^*\OV{\PR{\caG}})=
\CH_{\PR{T},1}^n(\OV{\caF}\boxtimes 
\PR{\OV{\caF}}, \OV{\caF}\boxtimes \varphi^*\OV{\PR{\caG}}).
$$

Let $\scZ_{X,G,\varphi }^{p,q}=\scZ_X^p\times \scZ_{G,\varphi }^q
\subset \scZ_{X\times G}^n$, and we identify 
$\scZ_X^p=\scZ_{X,G,\varphi }^{p,0}$ and 
$\scZ_{G,\varphi }^q=\scZ_{X,G,\varphi }^{0,q}$ in 
the same way as in \S 11.1.
Let 
$$
j_{X,G,\varphi }^{p,q}:\scD_{\log}^*(X\times G-\scZ_X^p, n)\oplus 
\scD_{\log}^*(X\times G-\scZ_{G,\varphi }^q, n)\to 
\scD_{\log}^*(X\times G-\scZ_X^p\cup \scZ_{G,\varphi }^q, n)
$$
be the map given by $(\omega , \tau )\mapsto -\omega +\tau $.
Then the natural map 
\begin{equation}
\scD_{\log}^*(X\times G-\scZ_{X,G,\varphi }^{p,q}, n)\to 
s(-j_{X,G,\varphi }^{p,q})^*  \label{sec11:map7}
\end{equation}
is a quasi-isomorphism \cite[Lem.5.5]{BF}.
Let 
$$
i_{X,G,\varphi }^{p,q}:\scD_{\log}^*(X\times G, n)\to 
s(-j_{X,G,\varphi }^{p,q})^* 
$$
be the composite of the restriction map to 
$\scD_{\log}^*(X\times G-\scZ_{X,G,\varphi }^{p,q}, n)$ with 
\eqref{sec11:map7}, and take the truncated subcomplex of 
the simple complex of this map, which we denote by 
$\tau_{\leq 2n}s(i_{X,G,\varphi }^{p,q})^*$.
Then in the same way as the definition of the map \eqref{sec11:map1} 
we can define a map of complexes 
$$
\diamond :\mmD_{\scZ^p}^*(X, p)\times 
\mmD_{\scZ_{G,\varphi }^q}^*(G, q)\to 
\tau_{\leq 2n}s(i_{X,G,\varphi }^{p,q})^*
$$
which induces an exterior product of Green forms 
$$
\diamond :GE^p(X)\times GE^q_{\scZ^q_{\varphi }}(G)\to 
\WH{H}^{2n}(\scD^*(X\times G, n), s(-j_{X,G,\varphi }^{p,q})^*).
$$
Moreover, we have a commutative diagram 
$$
{\small 
\begin{CD}
\WH{H}^{2n}(\scD^*(X\times G, n), s(-j_{X,G,\varphi }^{p,q})^*)  
@<{\sim }<< GE_{\scZ_{X,G,\varphi }^n}^n(X\times G) @>>> 
GE_{\scZ_{1\times \varphi }^n}^n(X\times G)   \\
@AAA  @AAA @AAA \\
\tau_{\leq 2n}s(i^{p,q}_{X,G,\varphi })^{2n} @<<< 
\mmD^{2n}_{\scZ_{X,G,\varphi }^{p,q}}(X\times G, n) @>>> 
\mmD^{2n}_{\scZ_{X\times G,1\times \varphi }^n}(X\times G, n) \\
@V{1\times \varphi^*}VV @VV{1\times \varphi^*}V 
@VV{1\times \varphi^*}V \\
s_{\A}(i^{p,q}_{X,Y})_0^{2n-r} @<<< 
\scD_{\A,\scZ^{p,q}_{X,Y}}^{2n-r}(X\times Y, n)_0 @>>> 
\scD_{\A,\scZ^n_{X\times Y}}^{2n-r}(X\times Y, n)_0, 
\end{CD}
}
$$
where the upper vertical maps are surjective, 
and Thm.5.12 in \cite{burgos2} says that the upper horizontal map 
sends the exterior product 
$(\omega^p, \WT{g^p})\diamond (\tau^q, \WT{h^q})$ 
to the star product of 
$\pi_X^*(\omega^p, \WT{g^p})$ with $\pi_G^*(\tau^q, \WT{h^q})$.
Since 
$$
(\omega^p, g^p)\diamond (\varphi^*(\tau^q), \varphi^*(h^q))=
(1\times \varphi )^*((\omega^p, g^p)\diamond (\tau^q, h^q)), 
$$
the sequence of maps \eqref{sec11:map5} sends the element 
\eqref{sec11:ele4} to 
$$
\left[\left(\text{
\setlength{\unitlength}{1mm}
\begin{picture}(101,7)
 \put(-2,-4){$y^p\times \varphi^*(z^q)$}
 \put(18,3){$0$}
 \put(21,-4){$\pi_X^*(\omega^p, g^p)\ast \pi_G^*(\tau^q, h^q)$}
 \put(44,3){$(-1)^{r+1}\CH_0^p(\OV{\caF})\EPRAP 
\CH_{T,1}^q(\OV{\PR{\caF}}, \varphi^*\OV{\PR{\caG}})$}
\end{picture}
}
\right)\right], 
$$
which means that 
$$
\sum_{p+q=n}\WH{\CH}_0^p(\OV{\caF})\cup 
\WH{\CH}_{T,0}^q(\OV{\PR{\caF}})
=\WH{\CH}_{\PR{T},0}^n(\OV{\caF}\boxtimes \OV{\PR{\caF}}).
$$
This completes the proof of \eqref{sec11:eq1}.
\qed


\end{document}